\documentclass[opre,sglanonrev]{informs4}

\RequirePackage{tgtermes}
\RequirePackage{newtxtext}
\RequirePackage{newtxmath}
\RequirePackage{endnotes}

\OneAndAHalfSpacedXII
\usepackage{anyfontsize}
\usepackage{mathtools}
\usepackage{enumitem}
\usepackage{dsfont}
\usepackage{booktabs}
\usepackage{caption}
\usepackage[labelfont=sf]{subcaption}
\usepackage{placeins}
\usepackage{stackengine,scalerel}
\usepackage{xcolor}
\usepackage{algorithm}
\usepackage{algpseudocode}
\usepackage{tikz}
\usetikzlibrary{calc,positioning,arrows.meta,quotes,decorations.pathreplacing,intersections,patterns,shapes.geometric}
\usepackage{pgfplots}
\pgfplotsset{compat=1.17}
\usepgfplotslibrary{groupplots}

\colorlet{Daf}{blue!72!black}
\colorlet{Dag}{orange!88!black}
\colorlet{Dial}{green!55!black}
\colorlet{MZ}{violet!72!black}
\colorlet{SBA}{red!72!black}
\pgfplotsset{
  objectiveaxis/.style={
    width=0.44\textwidth,
    height=0.26\textwidth,
    xmode=log,
    ymode=log,
    xmin=0.0001,
    xmax=3600,
    ymin=0.000001,
    ymax=2,
    xtick={0.0001,0.01,1,100,3600},
    xticklabels={$10^{-4}$,$10^{-2}$,$1$,$10^2$,$3600$},
    ytick={0.000001,0.0001,0.01,1},
    yticklabels={$10^{-6}$,$10^{-4}$,$10^{-2}$,$1$},
    tick label style={font=\scriptsize},
    title style={font=\small,align=center},
    axis line style={line width=0.35pt},
    tick style={line width=0.35pt},
    grid=both,
    major grid style={gray!28,line width=0.3pt},
    minor grid style={gray!12,line width=0.2pt},
    every axis plot/.append style={mark=none,line cap=round},
    clip mode=individual
  },
  daf/.style={Daf,solid,line width=0.8pt,mark=none},
  dag/.style={Dag,dashed,line width=0.8pt,mark=none},
  dial/.style={Dial,dash dot,line width=0.8pt,mark=none},
  mz/.style={MZ,dotted,line width=1.0pt,mark=none},
  sba/.style={SBA,densely dashed,line width=0.9pt,mark=none},
  bta/.style={black,solid,line width=1.25pt,mark=none}
}

\newcommand{\ToleranceLine}[1][3600]{%
  \addplot[gray!70,densely dashed,line width=0.45pt,forget plot]
  coordinates {(1e-9,0.00001) (#1,0.00001)};
  \node[anchor=south west,font=\tiny,text=gray!70]
  at ({rel axis cs:0,0}|-{axis cs:1,0.00001}) {$10^{-5}$};
}

\usepackage{natbib}
\bibpunct[, ]{(}{)}{,}{a}{}{,}%
\def\bibfont{\small}%

\newcommand{\defi}{\vcentcolon=}
\newcommand{\dom}{\mathrm{dom}\,}
\newcommand{\inner}[2]{\langle #1,#2\rangle}
\newcommand{\ignore}[1]{}
\newcommand{\ledot}{\mathrel{\ThisStyle{\ensurestackMath{\stackengine{-.5\LMex}{\SavedStyle\lessdot}{\SavedStyle-}{U}{c}{F}{F}{S}}}}}
\newcommand{\jacob}[1]{\textcolor{blue}{Jacob: #1}}

\def\diff{\mathop{}\!\mathrm{d}}

\makeatletter
\renewcommand{\ALG@beginalgorithmic}{\footnotesize}

\makeatother

\newtheorem{manuallemmainner}{Lemma}
\newenvironment{manuallemma}[1]{%
  \renewcommand\themanuallemmainner{#1}%
  \manuallemmainner
}{\endmanuallemmainner}

\colorlet{ExampleFrontier}{green!45!black}
\colorlet{ExampleInactive}{red!70!black}
\colorlet{ExamplePruned}{red!75!black}
\colorlet{ExampleWarm}{orange!80!black}
\colorlet{ExampleMid}{green!60!black}
\colorlet{ExampleChoiceOne}{blue!18}
\colorlet{ExampleChoiceTwo}{green!18}
\colorlet{ExampleChoiceThree}{orange!20}
\colorlet{ExampleGuide}{teal!55!black}
\pgfplotsset{
  exampleaxis/.style={
    grid=none,
    axis background/.style={fill=white},
    axis line style={draw=black, line width=0.35pt},
    tick style={draw=black, line width=0.35pt},
    tick align=outside,
    tick label style={font=\scriptsize},
    label style={font=\small},
    title style={font=\small},
    every axis plot/.append style={line cap=round, line join=round},
    clip mode=individual
  },
  frontierline/.style={very thick, color=ExampleFrontier, densely dotted},
  frontierguide/.style={densely dashed, line width=0.55pt, color=gray!58},
  frontierpoint/.style={only marks, mark=*, mark size=2.4pt, draw=black, fill=black},
  inactivepoint/.style={only marks, mark=o, mark size=2.1pt, draw=ExampleInactive, fill=white, line width=0.7pt},
  auxiliarypoint/.style={only marks, mark=o, mark size=2.1pt, draw=gray!70, fill=white, line width=0.7pt},
  prunedpoint/.style={only marks, mark=*, mark size=2.6pt, color=ExamplePruned}
}
\tikzset{
  thresholdline/.style={densely dashed, gray!72, line width=0.65pt}
}

\graphicspath{{figures/}}
\allowdisplaybreaks

\EquationsNumberedThrough
\TheoremsNumberedThrough
\ECRepeatTheorems
\MANUSCRIPTNO{}

\makeatletter
\def\theARTICLETOP{}%
\RRHSecondLine{}%
\LRHSecondLine{}%
\def\theARTICLEABSTRACT{%
  \HOOKb
  \vspace*{18pt}%
  \noindent\begin{minipage}[t]{\textwidth}\parindent1em
    \ABSfont
    \noindent\theABSTRACT\endgraf
    \vskip5pt
    \theFUNDING
    \theKEYWORDS
    \theSUBJECTCLASS
    \theAREAOFREVIEW
    \theMSCCLASS
    \theORMSCLASS
    \theHISTORY
  \noindent\hrulefill
  \end{minipage}%
  \vspace*{0pt}}
\makeatother

\begin{document}

\RUNAUTHOR{Aguirre, Kleywegt, and Monteiro}
\RUNTITLE{Bicriteria Traffic Assignment}
\TITLE{Bicriteria Traffic Assignment}

\ARTICLEAUTHORS{%
\AUTHOR{Jacob M. Aguirre, Anton J. Kleywegt, and Renato D. C. Monteiro}
\AFF{School of Industrial and Systems Engineering, Georgia Institute of Technology, Atlanta, GA 30332-0205, USA\\
\EMAIL{aguirre@gatech.edu}, \EMAIL{anton@isye.gatech.edu}, \EMAIL{monteiro@isye.gatech.edu}}
}

\FUNDING{Jacob M. Aguirre is supported by the National Science Foundation Graduate Research Fellowship under Grant No.~\mbox{DGE-2039655}.
Renato D.C. Monteiro is supported by AFOSR Grant \mbox{FA9550-25-1-0131}.}

\ABSTRACT{%
The classical traffic assignment problem determines the equilibrium flows on a congested network with many infinitesimal users who choose paths to minimize their own travel times.
Travelers care about cost (including tolls) as well as time, and different travelers put different weights on these two criteria.
We study a bicriteria traffic assignment problem in which travelers differ in their tradeoff between a flow-independent cost, such as tolls, and a flow-dependent cost, such as travel time.
Existing formulations require either finitely many user classes or positive preference densities, including formulations for positive preference densities that are either infinite-dimensional or nonconvex.
We formulate a finite-dimensional convex optimization problem for arbitrary preference distributions with bounded support, including distributions with a mixture of atoms, densities, and gaps in their support,
and we show how to convert an optimal solution to an equilibrium of the bicriteria traffic assignment problem.
We propose an algorithm that alternates between solving a convex problem for a (typically small) number of paths, and solving a parametric shortest-path problem.
Computational experiments on benchmark networks illustrate the performance of the algorithm relative to existing methods.
}%

\maketitle

\section{Bicriteria Traffic Assignment Problem}

Individual users, who want to travel from their origins to their destinations in a network, choose their preferred paths (we use the term \emph{user} instead of traveler/decision maker/player/agent).
An important problem is the determination of Wardrop equilibrium flows in the network, in which no user wants to change their chosen path.
We consider a problem with many infinitesimal (nonatomic) users, so that an individual user's path choice does not affect the desirability of the paths, that is, each individual user can regard the desirability of the paths as independent of the user's own path choice.
In addition, different users have different path preferences.
For example, one type of user may prefer a path from origin to destination with the least travel time, and another type of user may prefer a path from the same origin to the same destination with the least monetary travel cost (such as transit fees or tolls).
We consider the bicriteria traffic assignment problem (BiTAP) in which users choose paths from their origins to their destinations that minimize a weighted combination of two criteria, one called (travel) time and the other called cost or toll.
The travel time along a path depends on the flows on each of the arcs of the path, typically due to traffic congestion.
As pointed out in \citet{marc:96b} and \citet{marc:97}, a necessary condition for the equilibrium conditions to coincide with the first-order optimality conditions of an optimization problem, is that one of the criteria does not depend on flow rates.
Therefore, in the standard BiTAP, the cost/toll along a path does not depend on the flows on the path.
Different users use different weights for the two criteria, and users are partitioned into types according to their relative weights.
If the number of user types is finite, that is, the user type distribution has finite support, then a modification, with additional user type indices, of the various formulations for the single criterion traffic assignment problem can be used \citep[see, e.g.,][]{dafe:72}.
If the number of user types is infinite, then alternative approaches are needed.

Next we set notation that will be used throughout, and then we present various formulations of the BiTAP.
Let $(\mathcal{N},\mathcal{A})$ denote the finite directed network.
For each origin-destination (OD) pair $(u,v) \in \mathcal{N}^2$, let $D^{(u,v)} \ge 0$ denote the demand rate (flow rate per unit time) for flow from $u$ to $v$, and let $\Pi \defi \{(u,v) \in \mathcal{N}^2 \, : \, D^{(u,v)} > 0\}$ denote the OD pairs with positive demand.
We use $(u,v)$ to denote an OD pair when it is necessary to be specific about the origin or destination, otherwise we use $\pi$.
Unless noted otherwise, it is assumed that $D^{\pi}$ is input; this is called the fixed demand case.
For each OD pair $\pi \in \Pi$, let $\mathcal{R}^{\pi}$ denote the set of paths in the network for OD pair~$\pi$, for each arc~$a \in \mathcal{A}$, let $\mathcal{R}^{\pi}_{a}$ denote the set of paths for OD pair~$\pi$ that use arc~$a$, and for each $r \in \mathcal{R}^{\pi}$, let $\mathcal{A}^{\pi}_{r}$ denote the set of arcs on path~$r$.
For each arc~$a \in \mathcal{A}$, let $s_{a} \in \mathbb{R}_{+}$ denote the cost per unit flow on arc~$a$ that does not depend on the flow rates, and let $t_{a} : \mathbb{R}_{+} \mapsto \mathbb{R}_{+}$ denote the travel time on arc~$a$ as a function of the total flow rate on arc~$a$.
We assume that $t_{a}$ is nondecreasing.
A user with weight~$w \in \mathbb{R}_{+}$ for the cost relative to the travel time has preference consistent with bicriteria cost $c_{a}(w,x_{a}) \defi t_{a}(x_{a}) + w s_{a}$ when $x_{a}$ is the total flow rate on arc~$a$.
Let $x \defi (x_{a}, a \in \mathcal{A})$.
It is assumed that the bicriteria cost $C^{\pi}_{r}(w,x)$ per unit flow for a path~$r \in \mathcal{R}^{\pi}$ is the sum of the arc costs per unit flow, $C^{\pi}_{r}(w,x) = \sum_{a \in \mathcal{A}^{\pi}_{r}} c_{a}(w,x_{a})$.
For each OD pair $\pi \in \Pi$, and each path~$r \in \mathcal{R}^{\pi}$, let $S^{\pi}_{r} \defi \sum_{a \in \mathcal{A}^{\pi}_{r}} s_{a}$ be the sum of the arc costs per unit flow, and let $T^{\pi}_{r}(x) \defi \sum_{a \in \mathcal{A}^{\pi}_{r}} t_{a}(x_{a})$ be the total travel time as a function of arc flows~$x$.
For each OD pair $\pi \in \Pi$, and each path~$r \in \mathcal{R}^{\pi}$, let $f^{\pi}_{r}$ denote the flow on path~$r$, and let $f \defi (f^{\pi}_{r}, \pi \in \Pi, r \in \mathcal{R}^{\pi})$.
Given path flows $f$, the flow~$x_{a}$ on each arc~$a$ is uniquely given by $x_{a} = \sum_{\pi \in \Pi} \sum_{r \in \mathcal{R}^{\pi}_{a}} f^{\pi}_{r}$.
However, given arc flows~$x$, the path flows~$f$ are not unique in general.

As mentioned above, different users have different weights; let $G^{\pi}$ denote the cumulative distribution function of weight~$w$ among travelers for OD pair $\pi \in \Pi$, and let $[\underline{w}^{\pi},\overline{w}^{\pi}] \subset \mathbb{R}_{+}$ denote the smallest closed interval that contains the support of $G^{\pi}$.

For each OD pair $\pi \in \Pi$, we choose a total order $\prec$ on the set of paths $\mathcal{R}^{\pi}$ such that for any $r,r' \in \mathcal{R}^{\pi}$ it holds that if $r \prec r'$ (implying $r \neq r'$) then $S^{\pi}_{r} \ge S^{\pi}_{r'}$.
(Unlike previous papers on the BiTAP, we do not require the order to be unique, that is, we do not require that if $r \prec r'$ then $S^{\pi}_{r} > S^{\pi}_{r'}$.)
The notation $\prec$ does not show the dependence of the order on $\pi$, because it will be clear from the context.

Given arc flows~$x$, a path $r \in \mathcal{R}^{\pi}$ is called efficient if there is an interval $[\underline{w}^{\pi}_{r},\overline{w}^{\pi}_{r}] \subset [\underline{w}^{\pi},\overline{w}^{\pi}]$ with $\underline{w}^{\pi}_{r} \le \overline{w}^{\pi}_{r}$ such that $C^{\pi}_{r}(w,x) \le C^{\pi}_{r'}(w,x)$ for all $r' \in \mathcal{R}^{\pi}$ and $w \in [\underline{w}^{\pi}_{r},\overline{w}^{\pi}_{r}]$, and a path $r \in \mathcal{R}^{\pi}$ is called strictly efficient if there is an interval $[\underline{w}^{\pi}_{r},\overline{w}^{\pi}_{r}] \subset [\underline{w}^{\pi},\overline{w}^{\pi}]$ with $\underline{w}^{\pi}_{r} < \overline{w}^{\pi}_{r}$ such that $C^{\pi}_{r}(w,x) < C^{\pi}_{r'}(w,x)$ for all $r' \in \mathcal{R}^{\pi} \setminus \{r\}$ and $w \in (\underline{w}^{\pi}_{r},\overline{w}^{\pi}_{r})$.
Next we point out some basic properties that have been useful in previous research on the BiTAP.
Consider two efficient paths $r,r' \in \mathcal{R}^{\pi}$, and suppose that $r \prec r'$ (thus $S^{\pi}_{r} \ge S^{\pi}_{r'}$).
Consider any $w \in [\underline{w}^{\pi}_{r},\overline{w}^{\pi}_{r}]$ and $w' \in [\underline{w}^{\pi}_{r'},\overline{w}^{\pi}_{r'}]$.
Then $T^{\pi}_{r}(x) + w S^{\pi}_{r} \le T^{\pi}_{r'}(x) + w S^{\pi}_{r'}$ and thus $T^{\pi}_{r'}(x) - T^{\pi}_{r}(x) \ge w \left(S^{\pi}_{r} - S^{\pi}_{r'}\right) \ge 0$, hence $T^{\pi}_{r'}(x) \ge T^{\pi}_{r}(x)$.
Similarly, $T^{\pi}_{r'}(x) - T^{\pi}_{r}(x) \le w' \left(S^{\pi}_{r} - S^{\pi}_{r'}\right)$.
Previous papers assumed that different paths have different costs so that $S^{\pi}_{r} > S^{\pi}_{r'}$, in which case it can be concluded that $w \le w'$.
The assumption that $S^{\pi}_{r} > S^{\pi}_{r'}$ is awkward when these costs represent tolls (as in the canonical application), because in many applications there are many paths with no tolls.
Although the path costs could be perturbed for the purpose of numerical work, we show that this is unnecessary, and more importantly, that the focus on weight intervals $[\underline{w}^{\pi}_{r},\overline{w}^{\pi}_{r}]$ is unnecessarily restrictive without contributing to efficient computation.
For the purpose of presentation and comparison with previous work, we will also identify weight intervals $[\underline{w}^{\pi}_{r},\overline{w}^{\pi}_{r}]$, and we will show that in any case (even if $S^{\pi}_{r} = S^{\pi}_{r'}$), the intervals can be chosen such that $\overline{w}^{\pi}_{r} \le \underline{w}^{\pi}_{r'}$.
The choices of intervals $[\underline{w}^{\pi}_{r},\overline{w}^{\pi}_{r}]$ and $[\underline{w}^{\pi}_{r'},\overline{w}^{\pi}_{r'}]$ may not be unique; for example, if $r,r' \in \mathcal{R}^{\pi}$ are both efficient and $S^{\pi}_{r} = S^{\pi}_{r'}$ (and thus $T^{\pi}_{r}(x) = T^{\pi}_{r'}(x)$), then there is a common interval $[\underline{w}^{\pi}_{r,r'},\overline{w}^{\pi}_{r,r'}] \subset [\underline{w}^{\pi},\overline{w}^{\pi}]$ with $\underline{w}^{\pi}_{r,r'} \le \overline{w}^{\pi}_{r,r'}$ such that $C^{\pi}_{r}(w,x) = C^{\pi}_{r'}(w,x) \le C^{\pi}_{r''}(w,x)$ for all $r'' \in \mathcal{R}^{\pi}$ and $w \in [\underline{w}^{\pi}_{r,r'},\overline{w}^{\pi}_{r,r'}]$.
In such a case, if $r \prec r'$, then we will choose intervals $[\underline{w}^{\pi}_{r},\overline{w}^{\pi}_{r}], [\underline{w}^{\pi}_{r'},\overline{w}^{\pi}_{r'}] \subset [\underline{w}^{\pi}_{r,r'},\overline{w}^{\pi}_{r,r'}]$ such that $\overline{w}^{\pi}_{r} \le \underline{w}^{\pi}_{r'}$.
Also, if paths $r,r' \in \mathcal{R}^{\pi}$ are both efficient, and $r \prec r'$, and there is no strictly efficient path $r'' \in \mathcal{R}^{\pi}$ such that $r \prec r'' \prec r'$, then the intervals $[\underline{w}^{\pi}_{r},\overline{w}^{\pi}_{r}]$ and $[\underline{w}^{\pi}_{r'},\overline{w}^{\pi}_{r'}]$ can be chosen such that $\overline{w}^{\pi}_{r} = \underline{w}^{\pi}_{r'}$.
This allows the possibility that there is an efficient path $r'' \in \mathcal{R}^{\pi}$ such that $r \prec r'' \prec r'$ and $\overline{w}^{\pi}_{r} = \underline{w}^{\pi}_{r''} = \overline{w}^{\pi}_{r''} = \underline{w}^{\pi}_{r'}$.
Thus, given arc flows~$x$, there is a nonempty set $\{r^{1},\ldots,r^{k}\} \subset \mathcal{R}^{\pi}$ of efficient paths and a corresponding sequence of thresholds $\underline{w}^{\pi} = w^{\pi}_{r^{0}} \le w^{\pi}_{r^{1}} \le \cdots \le w^{\pi}_{r^{k-1}} \le w^{\pi}_{r^{k}} = \overline{w}^{\pi}$ such that for all $i=1,\ldots,k$ it holds that $C^{\pi}_{r^{i}}(w,x) \le C^{\pi}_{r'}(w,x)$ for all $r' \in \mathcal{R}^{\pi}$ and $w \in [w^{\pi}_{r^{i-1}},w^{\pi}_{r^{i}}]$.
As mentioned above, the set $\{r^{1},\ldots,r^{k}\} \subset \mathcal{R}^{\pi}$ and the thresholds $w^{\pi}_{r^{1}} \le \cdots \le w^{\pi}_{r^{k-1}}$ do not have to be unique, as illustrated in Example~\ref{ex:efficient frontier example}.

\begin{figure}[tbp]
  \centering
  \begin{tikzpicture}
    \begin{axis}[
      exampleaxis,
      width=0.58\textwidth,
      height=0.42\textwidth,
      axis lines=left,
      axis line style={-{Stealth[length=2mm]}, draw=black, line width=0.35pt},
      xmin=0.75, xmax=5.25,
      ymin=0.03, ymax=0.34,
      xtick={1.0,2.0,3.0,4.0,5.0},
      ytick={0.05,0.10,0.15,0.20,0.25,0.30},
      xlabel={Cost/toll $S$},
      ylabel={Time $T$},
      legend style={
        font=\scriptsize,
        draw=none,
        fill=white,
        fill opacity=0.88,
        text opacity=1,
        at={(0.98,0.98)},
        anchor=north east
      },
      legend cell align=left,
      clip=false
    ]
      \addplot[
        frontierline
      ] coordinates {
        (5.00,0.040)
        (4.00,0.050)
        (3.00,0.120)
        (2.00,0.190)
        (1.00,0.310)
      };
      \addlegendentry{efficient frontier}
      \addplot[
        only marks,
        mark=diamond*,
        mark size=3.0pt,
        draw=black,
        fill=ExampleFrontier
      ] coordinates {
        (5.00,0.040)
        (4.00,0.050)
        (3.00,0.120)
        (2.00,0.190)
        (1.00,0.310)
      };
      \addlegendentry{efficient path}

      \addplot[
        only marks,
        mark=none,
        mark size=2.4pt,
        draw=blue!70!black,
        fill=blue!70!black
      ] coordinates {
        (3.70,0.135)
        (3.35,0.160)
        (2.60,0.180)
        (2.25,0.235)
      };
      \addlegendentry{inefficient path}

      \addplot[
        only marks,
        mark=o,
        mark size=2.7pt,
        draw=ExamplePruned,
        fill=white,
        line width=0.75pt,
        forget plot
      ] coordinates {
        (3.00,0.134)
      };

      \node[font=\scriptsize, anchor=west] at (axis cs:5.00,0.040) {$r_1$};
      \node[font=\scriptsize, anchor=north] at (axis cs:4.04,0.05) {$r_2,r_3$};
      \node[font=\scriptsize, anchor=north east] at (axis cs:3.00,0.120) {$r_6$};
      \node[font=\scriptsize, anchor=south east] at (axis cs:2.00,0.174) {$r_9$};
      \node[font=\scriptsize, anchor=south] at (axis cs:1.00,0.310) {$r_{10}$};

      \node[font=\scriptsize, anchor=south west] at (axis cs:3.70,0.135) {$r_4$};
      \node[font=\scriptsize, anchor=south west] at (axis cs:3.35,0.160) {$r_5$};
      \node[font=\scriptsize, anchor=south west] at (axis cs:2.60,0.180) {$r_7$};
      \node[font=\scriptsize, anchor=south west] at (axis cs:2.25,0.235) {$r_8$};
      \draw[
        densely dashed,
        line width=0.55pt,
        color=gray!58,
        -{Stealth[length=1.5mm]}
      ]
        (axis cs:4.27,0.08) -- (axis cs:4.04,0.054);
      \node[
        font=\scriptsize,
        anchor=west,
        align=left,
        fill=white,
        inner sep=1pt
      ] at (axis cs:4.28,0.082) {exact tie};

      \draw[
        ->,
        color=ExamplePruned,
        line width=0.55pt
      ]
        (axis cs:3.00,0.123) -- (axis cs:3.00,0.132);

      \draw[draw=none]
        (axis cs:4.00,0.050) -- (axis cs:5.00,0.04)
        node[pos=0.65, sloped, above, font=\scriptsize] {slope $-0.01$};
      \draw[draw=none]
        (axis cs:3.00,0.120) -- (axis cs:4.00,0.050)
        node[midway, sloped, below, font=\scriptsize] {slope $-0.07$};
      \draw[draw=none]
        (axis cs:2.00,0.190) -- (axis cs:3.00,0.120)
        node[midway, sloped, below, font=\scriptsize] {slope $-0.07$};
      \draw[draw=none]
        (axis cs:1.00,0.310) -- (axis cs:2.00,0.190)
        node[midway, sloped, below, font=\scriptsize] {slope $-0.12$};

    \end{axis}
  \end{tikzpicture}
  \caption{Bicriteria cost pairs $(S^{\pi}_{r_{i}},T^{\pi}_{r_{i}}(x))$ for different paths~$r_{i}$ for a given~$x$, and the corresponding efficient frontier, as in Example~\ref{ex:efficient frontier example}.
}
  \label{fig:efficient-frontier-tie-interior}
\end{figure}
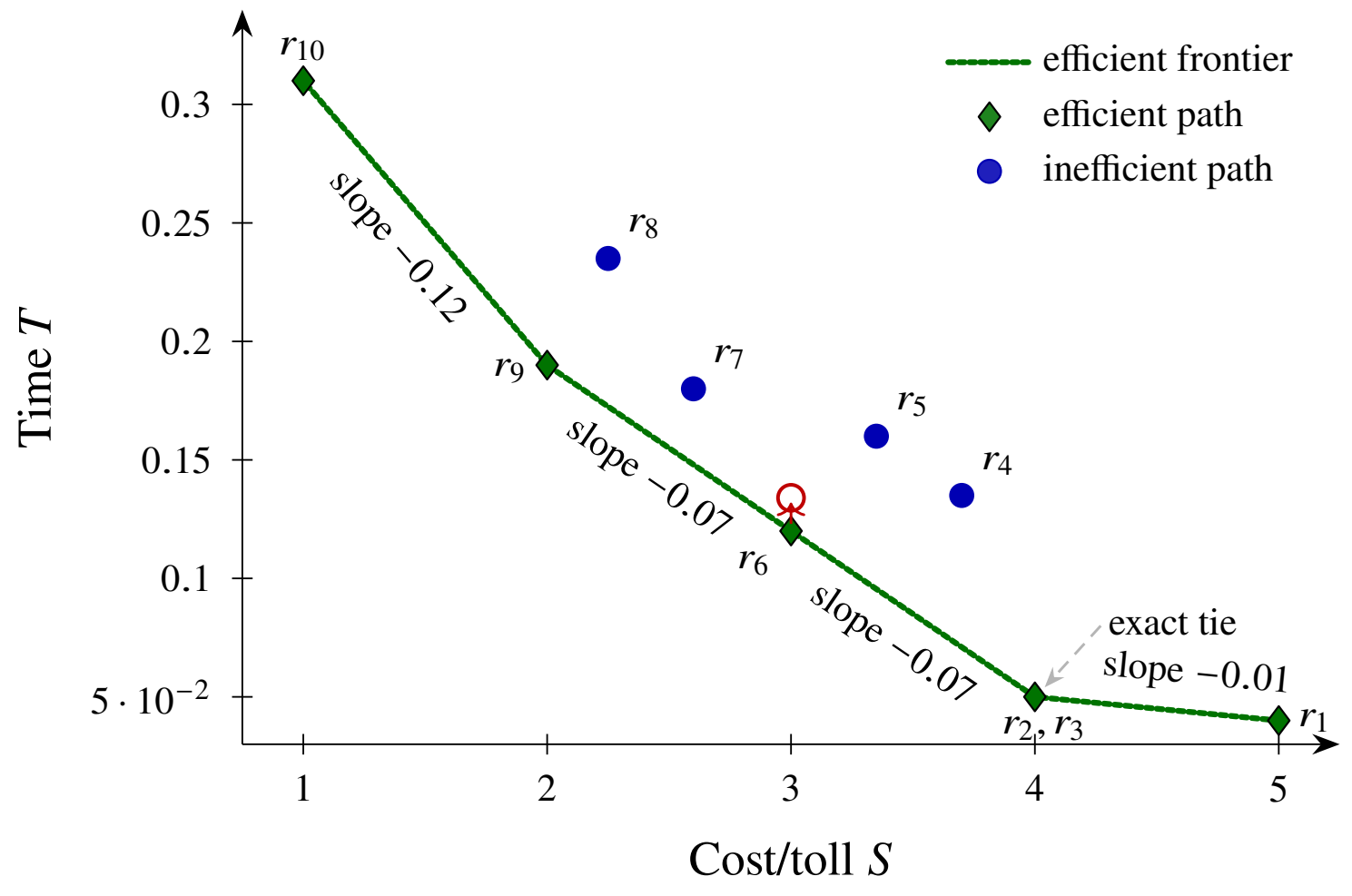

\begin{example}
\label{ex:efficient frontier example}
Fix an OD pair $\pi \in \Pi$ with support interval $[\underline{w}^{\pi},\overline{w}^{\pi}] = [0,0.14]$, arc flows~$x$, and a path set $\mathcal{R}^{\pi} = \{r_{1},r_{2},\ldots,r_{10}\}$ indexed so that $r_{1} \prec r_{2} \prec \cdots \prec r_{10}$.
Let $(S^{\pi}_{r_{1}},S^{\pi}_{r_{2}},\ldots,S^{\pi}_{r_{10}}) = (5,4,4,3.7,3.35,3,2.6,2.25,2,1)$ and $(T^{\pi}_{r_{1}}(x),T^{\pi}_{r_{2}}(x),\ldots,T^{\pi}_{r_{10}}(x)) = (0.04,0.05,0.05,0.135,0.16,0.12,0.18,0.235,0.19,0.31)$.
Figure~\ref{fig:efficient-frontier-tie-interior} shows the points $(S^{\pi}_{r_{i}},T^{\pi}_{r_{i}}(x))$ for $i = 1,\ldots,10$.
As $x$ changes, $T^{\pi}(x)$ changes, and thus the points move up or down.
Note that, at the given~$x$, the cost pairs $(S^{\pi}_{r_{i}},T^{\pi}_{r_{i}}(x))$ for paths $r_{1}$, $r_{2}$, $r_{3}$, $r_{6}$, $r_{9}$, and $r_{10}$ are on the efficient frontier, whereas the cost pairs $(S^{\pi}_{r_{i}},T^{\pi}_{r_{i}}(x))$ for paths $r_{4}$, $r_{5}$, $r_{7}$, and $r_{8}$ are not on the efficient frontier.

Next we identify thresholds $\underline{w}^{\pi} = w^{\pi}_{0} \le w^{\pi}_{1} \le \cdots \le w^{\pi}_{k-1} \le w^{\pi}_{k} = \overline{w}^{\pi}$ such that for all efficient paths $r_{i}$ it holds that $C^{\pi}_{r_{i}}(w,x) \le C^{\pi}_{r'}(w,x)$ for all $r' \in \mathcal{R}^{\pi}$ and $w \in [w^{\pi}_{i-1},w^{\pi}_{i}]$.
Note that $C^{\pi}_{r_{1}}(w,x) \le C^{\pi}_{r'}(w,x)$ for all $r' \in \mathcal{R}^{\pi}$ and $w \in [w^{\pi}_{0},w^{\pi}_{1}] = [0,0.01]$, and that the slope between $(S^{\pi}_{r_{1}},T^{\pi}_{r_{1}}(x))$ and $(S^{\pi}_{r_{2}},T^{\pi}_{r_{2}}(x))$ is (the negative of) $0.01$.
Next, $(S^{\pi}_{r_{2}},T^{\pi}_{r_{2}}(x)) = (S^{\pi}_{r_{3}},T^{\pi}_{r_{3}}(x))$, that is, paths $r_{2}$ and $r_{3}$ are exactly tied.
Thus, $C^{\pi}_{r_{2}}(w,x) = C^{\pi}_{r_{3}}(w,x) \le C^{\pi}_{r'}(w,x)$ for all $r' \in \mathcal{R}^{\pi}$ and $w \in [0.01,0.07]$, where the slope between $(S^{\pi}_{r_{2}},T^{\pi}_{r_{2}}(x))$ and $(S^{\pi}_{r_{6}},T^{\pi}_{r_{6}}(x))$ is (the negative of) $0.07$.
Because $r_{2} \prec r_{3}$, their intervals may be chosen as $[w^{\pi}_{1},w^{\pi}_{2}] = [0.01,\eta]$ and $[w^{\pi}_{2},w^{\pi}_{3}] = [\eta,0.07]$ for any $\eta \in [0.01,0.07]$.
Next, because paths $r_{4}$ and $r_{5}$ are inefficient, $w^{\pi}_{4}$ and $w^{\pi}_{5}$ have to be chosen as $w^{\pi}_{4} = w^{\pi}_{5} = w^{\pi}_{3} = 0.07$.
However, we have to keep in mind that paths $r_{4}$ and $r_{5}$ are inefficient, and thus it does not hold that $C^{\pi}_{r_{4}}(w,x) \le C^{\pi}_{r'}(w,x)$ for all $r' \in \mathcal{R}^{\pi}$ and $w \in [w^{\pi}_{3},w^{\pi}_{4}] = [0.07,0.07]$, and similarly for path~$r_{5}$.
Next, note that point $(S^{\pi}_{r_{6}},T^{\pi}_{r_{6}}(x))$ lies on the line segment between points $(S^{\pi}_{r_{3}},T^{\pi}_{r_{3}}(x))$ and $(S^{\pi}_{r_{9}},T^{\pi}_{r_{9}}(x))$.
Thus, $w^{\pi}_{6}$ also has to be chosen as $w^{\pi}_{6} = w^{\pi}_{5} = 0.07$.
However, in this case path $r_{6}$ is efficient, and thus it holds that $C^{\pi}_{r_{6}}(w,x) \le C^{\pi}_{r'}(w,x)$ for all $r' \in \mathcal{R}^{\pi}$ and $w \in [w^{\pi}_{5},w^{\pi}_{6}] = [0.07,0.07]$.
Thus, path $r_{6}$ is not strictly efficient.
Nevertheless, if the distribution $G^{\pi}$ has a point mass at $w^{\pi}_{6}$, then path~$r_{6}$ may be preferred by a positive amount of demand, even though $[w^{\pi}_{5},w^{\pi}_{6}]$ is a singleton.
That does not mean that all the demand corresponding to the point mass at $w^{\pi}_{6}$ has to flow on path~$r_{6}$ --- the demand corresponding to the point mass at $w^{\pi}_{6}$ is partitioned among paths $r_{2}$, $r_{3}$, $r_{6}$, and $r_{9}$.
In general, the demand corresponding to any point mass of $G^{\pi}$ at $w$ is partitioned among all the paths $r_{i}$ such that $C^{\pi}_{r_{i}}(w,x) \le C^{\pi}_{r'}(w,x)$ for all $r' \in \mathcal{R}^{\pi}$.
Also note that the status of path~$r_{6}$ is sensitive to small perturbations of $T^{\pi}_{r_{6}}(x)$.
For example, if $x$ changes by a small amount to $x'$ and as a result $T^{\pi}_{r_{6}}(x)$ increases to $T^{\pi}_{r_{6}}(x')$ that is slightly above the line segment between points $(S^{\pi}_{r_{3}},T^{\pi}_{r_{3}}(x'))$ and $(S^{\pi}_{r_{9}},T^{\pi}_{r_{9}}(x'))$, as shown by the red circle in Figure~\ref{fig:efficient-frontier-tie-interior}, then path~$r_{6}$ ceases to be efficient and the flow on path~$r_{6}$ can jump from a positive amount to zero.
Continuing in the same way, paths $r_{7}$ and $r_{8}$ are inefficient, and thus $w^{\pi}_{7}$ and $w^{\pi}_{8}$ have to be chosen as $w^{\pi}_{7} = w^{\pi}_{8} = w^{\pi}_{6} = 0.07$.
Next, note that $C^{\pi}_{r_{9}}(w,x) \le C^{\pi}_{r'}(w,x)$ for all $r' \in \mathcal{R}^{\pi}$ and $w \in [w^{\pi}_{8},w^{\pi}_{9}] = [0.07,0.12]$, and that the slope between $(S^{\pi}_{r_{9}},T^{\pi}_{r_{9}}(x))$ and $(S^{\pi}_{r_{10}},T^{\pi}_{r_{10}}(x))$ is (the negative of) $0.12$.
Finally, $C^{\pi}_{r_{10}}(w,x) \le C^{\pi}_{r'}(w,x)$ for all $r' \in \mathcal{R}^{\pi}$ and $w \in [0.12,\infty)$, and thus $w^{\pi}_{10}$ is chosen as $w^{\pi}_{10} = \overline{w}^{\pi} = 0.14$.
\hfill\Halmos
\end{example}

\section{Related Literature}
\label{sec:literature}

\citet{dafe:72} considered a BiTAP with a finite number of user types.
It was pointed out that it is a special case of the traffic assignment problem in which $t_{a}$ depends on $x$ and not just on $x_{a}$, that is, the travel cost on an arc~$a$ depends on the flow on other arcs and not just the flow on arc~$a$.
The results of \citet{dafe:71} were used to obtain the following sufficient condition for the user-equilibrium to be given by a solution of an optimization problem:
Let $w = 1,\ldots,W$ index the user types, let $x^{(u,v)}_{a}(w)$ denote the flow of type~$w$ users from origin~$u$ to destination~$v$ on arc~$a$, and let $\mathcal{X}$ denote the set of arc flows $(x^{(u,v)}_{a}(w), (u,v) \in \Pi, a \in \mathcal{A}, w = 1,\ldots,W) \ge 0$ such that
\begin{equation}
\label{eqn:arc flow demand satisfaction finite types}
\sum_{\{j \in \mathcal{N} \, : \, (i,j) \in \mathcal{A}\}} x^{(u,v)}_{(i,j)}(w) - \sum_{\{j \in \mathcal{N} \, : \, (j,i) \in \mathcal{A}\}} x^{(u,v)}_{(j,i)}(w) \ = \ \left\{
\begin{array}{ll}
D^{(u,v)}(w) & \mbox{if } i = u \\
- D^{(u,v)}(w) & \mbox{if } i = v \\
0 & \mbox{otherwise}
\end{array}\right.
\hspace{1mm} \forall \ i \in \mathcal{N}, \ \ \forall \ (u,v) \in \Pi, \ \ \forall \ w = 1,\ldots,W
\end{equation}
Let $x_{a}(w) \defi \sum_{(u,v) \in \Pi} x^{(u,v)}_{a}(w)$, and let $x_{a} \defi (x_{a}(1), \ldots, x_{a}(W))$.
Let
\[
\overline{\mathcal{X}} \ \ \defi \ \ \left\{(x_{a}, a \in \mathcal{A}) \; : \; x_{a}(w) = \sum_{(u,v) \in \Pi} x^{(u,v)}_{a}(w), \; \left(x^{(u,v)}_{a}(w), (u,v) \in \Pi, a \in \mathcal{A}, w = 1,\ldots,W\right) \in \mathcal{X}\right\}
\]
denote the set of feasible arc flows for each user type.
Let $c_{a}(w,x_{a})$ denote the unit travel cost of type~$w$ flow on arc~$a$, and let $c(x) \defi (c_{a}(w,x_{a}), a \in \mathcal{A}, w = 1,\ldots,W)$.
Then user equilibrium can be specified with the variational inequality $x^* \in \overline{\mathcal{X}}$, $\langle c(x^*), y - x^*\rangle \ge 0$ for all $y \in \overline{\mathcal{X}}$, that is, $-c(x^*) \in \mathcal{N}_{\overline{\mathcal{X}}}(x^*)$, where $\mathcal{N}_{\overline{\mathcal{X}}}(x^*)$ denotes the normal cone to $\overline{\mathcal{X}}$ at $x^*$.

If, for each arc~$a \in \mathcal{A}$ and all flow vectors $x_{a}$, the matrix with entries $\partial c_{a}(w,x_{a}) / \partial x_{a}(w')$, $w, w' = 1,\ldots,W$, is symmetric and positive definite, then the user-equilibrium can be computed by solving
\[
\min_{x \ge 0} \left\{\sum_{a \in \mathcal{A}} \sum_{w=1}^{W} \int_{0}^{x_{a}(w)} c_{a}\big(w,(z,x_{-a}(-w))\big) \diff z \; : \; \eqref{eqn:arc flow demand satisfaction finite types}\right\}
\]
where $x_{-a}(-w) \defi \big(x_{a'}(w'), a' \in \mathcal{A}, w' = 1,\ldots,W, (a',w') \neq (a,w)\big)$ denotes the vector of all $x_{a'}(w')$ except for $(a',w') = (a,w)$.
Note that in the case of the BiTAP, $c_{a}(w,x_{a}) = t_{a}(x_{a}(1) + \cdots + x_{a}(W)) + w s_{a}$, and thus $\partial c_{a}(w,x_{a}) / \partial x_{a}(w') = \diff t_{a}(x_{a}(1) + \cdots + x_{a}(W)) / \diff x$.
Hence, all entries $\partial c_{a}(w,x_{a}) / \partial x_{a}(w')$ of the matrix are the same positive number $\diff t_{a}(x_{a}(1) + \cdots + x_{a}(W)) / \diff x$, and therefore the matrix for the BiTAP is symmetric and positive semi-definite.
\citet{dafe:72} proposed an algorithm that solves shortest path problems similar to the Frank-Wolfe algorithm for the traffic assignment problem, but instead of the new point being a convex combination of the previous point and the shortest path, the new point is determined by removing flow from the most expensive path with positive flow for an OD pair, and adding the flow to the shortest path for the OD pair.

\citet{dial:79} considered the special case of the BiTAP in which both $s_{a}$ and $t_{a}$ do not depend on the flow rates~$x$.
An algorithm was proposed that, for a given OD pair $\pi$, finds a set $\{r_{1},\ldots,r_{k}\} \subset \mathcal{R}^{\pi}$ of efficient paths, and corresponding thresholds $\underline{w}^{\pi} = w^{\pi}_{0} < w^{\pi}_{1} < \cdots < w^{\pi}_{k-1} < w^{\pi}_{k} = \overline{w}^{\pi}$, such that for all $i=1,\ldots,k$ it holds that $C^{\pi}_{r_{i}}(w,x) \le C^{\pi}_{r'}(w,x)$ for all $r' \in \mathcal{R}^{\pi}$ and $w \in [w^{\pi}_{i-1},w^{\pi}_{i}]$.

\citet{daga:83} considered a traffic assignment problem with a finite number of user types, similar to \citet{dafe:72}.
The model allows more general arc costs of the form $c_{a}(w,x) \defi t_{a}(x) + w s_{a}$ that may depend on the flows~$x$ on all arcs.
It is assumed that the mapping $t : \mathbb{R}^{\mathcal{A}} \mapsto \mathbb{R}^{\mathcal{A}}$ is increasing.
In addition, it is assumed that for each user type $w = 1,\ldots,W$, there is a function $x^{w} : \mathbb{R}^{\mathcal{A}} \mapsto \mathbb{R}^{\mathcal{A}}$ that maps arc costs to arc flows that is continuously differentiable with a negative semidefinite Jacobian.
The assumption that $x^{w}$ is continuous means that the Wardrop equilibrium conditions (only OD paths with least travel cost have positive flow) do not apply.
Existence of a fixed point is established with Brouwer's fixed point theorem and Kakutani's fixed point theorem.
A stochastic approximation algorithm is proposed to compute fixed points.

\citet{leur:93} considered the following BiTAP.
(\citet{leur:93} groups paths $r \in \mathcal{R}^{\pi}$ with the same value of $S^{\pi}_{r}$ into ``monetary expense classes'', but to simplify notation we omit the classification.)
For each OD pair $\pi \in \Pi$, and each path~$r \in \mathcal{R}^{\pi}$, let $F^{\pi}_{r} \defi \sum_{r'=1}^{r} f^{\pi}_{r'}$ denote the cumulative flow on paths~$r' = 1,\ldots,r$, let $F^{\pi}_{0} \defi 0$, and let $F^{\pi} \defi \sum_{r \in \mathcal{R}^{\pi}} f^{\pi}_{r} = F^{\pi}_{R^{\pi}}$ denote the total flow for OD pair $\pi$.
It is assumed that $G^{\pi}$ has a positive density so that $\left(G^{\pi}\right)^{-1}(q)$ is a singleton for all $q \in (0,1)$.
The OD demand is allowed to depend on the mean total cost over all OD paths and all weights, as follows:
\[
F^{\pi} \ \ = \ \ D^{\pi} \left(\sum_{r \in \mathcal{R}^{\pi}} \left[\frac{f^{\pi}_{r}}{F^{\pi}} T^{\pi}_{r}(x) + S^{\pi}_{r} \int_{F^{\pi}_{r-1} / F^{\pi}}^{F^{\pi}_{r} / F^{\pi}} \left(G^{\pi}\right)^{-1}(q) \diff q\right]\right)
\]
for $F^{\pi} > 0$, where $D^{\pi} : \mathbb{R}_{+} \mapsto \mathbb{R}_{+}$ is assumed to be continuous and decreasing.
It is shown that the equilibrium conditions coincide with the first-order optimality conditions of the following optimization problem:
\begin{equation}
\label{eqn:path cumulative flow formulation}
\min_{f \ge 0} \left\{\sum_{a \in \mathcal{A}} \int_{0}^{x_{a}} t_{a}(z) \diff z + \sum_{\pi \in \Pi} F^{\pi} \sum_{r \in \mathcal{R}^{\pi}} S^{\pi}_{r} \int_{F^{\pi}_{r-1} / F^{\pi}}^{F^{\pi}_{r} / F^{\pi}} \left(G^{\pi}\right)^{-1}(q) \diff q - \sum_{\pi \in \Pi} \int_{0}^{F^{\pi}} \left(D^{\pi}\right)^{-1}(q) \diff q\right\}
\end{equation}

It is shown that the objective function is convex.
Similar to \citet{daga:83}, a stochastic approximation algorithm is proposed to solve the optimization problem.
The algorithm enumerates paths.

\citet{leur:96} considered the same problem as in \citet{leur:93}, with the possible addition of arc capacity constraints of the form $x_{a} \le C_{a}$ (here, $C_{a}$ denotes the capacity of arc~$a$).
Three algorithms are compared for the case without capacity constraints.
For the case with capacity constraints, an augmented Lagrangian algorithm is proposed.

\citet{marc:96a} considered the following formulation of the BiTAP.
It is assumed that $G^{\pi}$ has a density, but as before we give a more general formulation here to match our notation.
Let $x^{\pi}_{a}(w)$ denote the (unscaled) flow on arc~$a \in \mathcal{A}$ of travelers for OD pair~$\pi \in \Pi$ with preference weight~$w \in [0,w_{\max}]$.
Also, for each arc~$a \in \mathcal{A}$, let $X_{a} \defi \sum_{\pi \in \Pi} \int_{0}^{w_{\max}} x^{\pi}_{a}(w) \diff G^{\pi}(w)$ denote the total flow on arc~$a$.
Let $\mathcal{L}^{\pi}_{2}$ denote the set of functions $f : [0,w_{\max}] \mapsto \mathbb{R}$ such that $\int_{0}^{w_{\max}} f^2(w) \diff G^{\pi}(w) < \infty$; for any $f^{\pi},g^{\pi} \in \mathcal{L}^{\pi}_{2}$, let $\langle f^{\pi}, g^{\pi} \rangle^{\pi}_{2} \defi \int_{0}^{w_{\max}} f^{\pi}(w) g^{\pi}(w) \diff G^{\pi}(w)$.
Then, for any $f = \left(f^{\pi}_{a} \in \mathcal{L}^{\pi}_{2}, \pi \in \Pi, a \in \mathcal{A}\right) \in \mathcal{L}_{2} \defi \prod_{\pi \in \Pi} \left(\mathcal{L}^{\pi}_{2}\right)^{\mathcal{A}}$, $g = \left(g^{\pi}_{a} \in \mathcal{L}^{\pi}_{2}, \pi \in \Pi, a \in \mathcal{A}\right) \in \mathcal{L}_{2}$, let $\langle f, g \rangle_{2} \defi \sum_{\pi \in \Pi} \sum_{a \in \mathcal{A}} \langle f^{\pi}_{a}, g^{\pi}_{a} \rangle^{\pi}_{2}$.
For any $x \in \mathcal{X}$, let $c^{\pi}_{a}(w,x) \defi t_{a}(X_{a}) + w s_{a}$ (it has the same value for all $\pi \in \Pi$), and let $c(x) \defi (c^{\pi}_{a}(\cdot,x) \in \mathcal{L}^{\pi}_{2}, \pi \in \Pi, a \in \mathcal{A}) \in \mathcal{L}_{2}$.
Then the equilibrium conditions can be expressed as the infinite dimensional variational inequality $x^* \in \mathcal{X}$, $\langle c(x^*), y - x^* \rangle_{2} \ge 0$ for all $y \in \mathcal{X}$.
It is shown that the equilibrium conditions coincide with the first-order optimality conditions of the infinite dimensional optimization problem
\begin{equation}
\min_{x \in \mathcal{X}} \ \ \left\{\sum_{a \in \mathcal{A}} \int_{0}^{X_{a}} t_{a}(z) \diff z + \sum_{a \in \mathcal{A}} s_{a} \sum_{\pi \in \Pi} \int_{\underline{w}}^{\overline{w}} w x^{\pi}_{a}(w) \diff G^{\pi}(w)\right\}
\label{eqn:inf dimensional arc flow formulation}
\end{equation}
It is shown that the objective function is convex.
The Frank-Wolfe subproblem at $\bar{x} \in \mathcal{X}$ is
\begin{equation}
\min_{x \in \mathcal{X}} \ \ \sum_{a \in \mathcal{A}} \sum_{\pi \in \Pi} \int_{\underline{w}}^{\overline{w}} \left[t_{a}(\bar{X}_{a}) + w s_{a}\right] x^{\pi}_{a}(w) \diff G^{\pi}(w)
\end{equation}
The Frank-Wolfe subproblem separates into the following shortest path problem for each $(u,v) \in \Pi$ and each $w \in [\underline{w},\overline{w}]$:
\begin{align}
\label{eqn:parametric shortest path problem}
\min_{x^{(u,v)}(w) \ge 0} & \ \ \left\{\sum_{a \in \mathcal{A}} \left[t_{a}(\bar{X}_{a}) + w s_{a}\right] x^{(u,v)}_{a}(w) \; : \; \eqref{eqn:arc flow demand satisfaction finite types}\right\}
\end{align}
It is mentioned that, for each $(u,v) \in \Pi$, problems~\eqref{eqn:parametric shortest path problem} for all $w \in [0,w_{\max}]$ can be solved by solving one parametric shortest path problem with parameter~$w$.
Thereafter a line search problem is solved to determine the minimizing convex combination between the previous flows~$x$ and the shortest path flows obtained from the Frank-Wolfe subproblem.
\citet{marc:96b} provided more details about the proposed algorithm,
as well as numerical results.

\citet{dial:96} considered the BiTAP in which both arc costs $s_{a}$ and $t_{a}$ are allowed to be functions of the total arc flow~$x_{a}$.
Two algorithms were proposed.
One algorithm is a heuristic based on the efficient paths algorithm of \citet{dial:79} and Frank-Wolfe steps.
Another algorithm is based on the same optimization problem as in \citet{marc:96a} and \citet{marc:96b}, with $s_{a}$ not varying with the arc flow~$x_{a}$.
\citet{dial:97} considered the same problem as in \citet{dial:96}.
An algorithm is proposed to solve the same variational inequality, $x^* \in \mathcal{X}$, $\langle c(x^*), y - x^* \rangle_{2} \ge 0$ for all $y \in \mathcal{X}$, as in \citet{marc:96a}.
Instead of solving the optimization problem~\eqref{eqn:inf dimensional arc flow formulation}, the proposed algorithm alternates between solving the parametric shortest path problem~\eqref{eqn:parametric shortest path problem} at current flow $\bar{x} \in \mathcal{X}$ that produces flow $x^{m} \in \mathcal{X}$ in iteration~$m$, and a linear programming approximation that chooses the next flow as the convex combination $x(\lambda) \defi \sum_{m'=1}^{m} \lambda_{m'} x^{m'}$, $\sum_{m'=1}^{m} \lambda_{m'} = 1$, $\lambda \ge 0$, attempting to minimize $\sum_{m'=1}^{m} \max\left\{0, \left\langle c\big(x(\lambda)\big), x^{m'} - x(\lambda) \right\rangle_{2}\right\}$.
\citet{marc:97} considered the problem in \citet{marc:96a} and \citet{marc:96b} with more general linear constraints and with arc costs given by $c^{\pi}_{a}(x)(w) = t_{a}(X) + w s_{a}(X)$.
\citet{marc:97} showed that, if $s$ and $t$ are Lipschitz continuous and continuously differentiable, then $c$ is a monotone mapping if and only if $t : \mathbb{R}^{\mathcal{A}} \mapsto \mathbb{R}^{\mathcal{A}}$ is monotone and $s : \mathbb{R}^{\mathcal{A}} \mapsto \mathbb{R}^{\mathcal{A}}$ is constant.

\begin{figure}[t]
  \centering
  \begin{subfigure}[t]{0.30\textwidth}
    \centering
    \begin{tikzpicture}
      \begin{axis}[
        exampleaxis,
        width=\linewidth,
        height=5.8cm,
        axis lines=left,
        axis line style={-{Stealth[length=2mm]}, draw=black, line width=0.35pt},
        xmin=-0.02, xmax=1.05,
        ymin=-0.02, ymax=1.05,
        xtick={0,0.5,1},
        ytick={0,0.5,1},
        xlabel={Path-$1$ share},
        ylabel={Path-$2$ share},
        ylabel style={xshift=-2mm},
        clip=false
      ]
        \addplot[
          only marks,
          mark=diamond*,
          mark size=3.0pt,
          draw=black,
          fill=ExampleFrontier
        ] coordinates {
          (1,0)
          (0,1)
        };

        \addplot[
          only marks,
          mark=*,
          mark size=2.6pt,
          draw=ExamplePruned,
          fill=ExamplePruned
        ] coordinates {
          (0.5,0.5)
        };

        \draw[densely dashed, line width=0.55pt, color=gray!58]
          (axis cs:0,1) -- (axis cs:1,0);

        \node[font=\scriptsize, anchor=south east] at (axis cs:1.2,0.01) {$\hat{f}^{\pi},\hat{\bar x}$};
        \node[font=\scriptsize, anchor=south west] at (axis cs:0,1) {$\check f^{\pi},\check{\bar x}$};
        \node[font=\scriptsize, anchor=south west, color=ExamplePruned] at (axis cs:0.5,0.5) {$f^{\pi},\bar x$};

        \node[
          font=\tiny,
          fill=white,
          inner sep=1pt
        ] at (axis cs:0.72,0.72) {$\bar x_{1}+\bar x_{2}=1$};
      \end{axis}
    \end{tikzpicture}
    \subcaption{Solutions and convex combination for formulations of \citet{leur:93} and \citet{marc:99}.}
  \end{subfigure}
  \hfill
  \begin{subfigure}[t]{0.30\textwidth}
    \centering
    \begin{tikzpicture}
      \begin{axis}[
        exampleaxis,
        width=\linewidth,
        height=5.8cm,
        axis lines=left,
        axis line style={-{Stealth[length=2mm]}, draw=black, line width=0.35pt},
        xmin=0, xmax=1.05,
        ymin=-0.02, ymax=1.05,
        xtick={0,0.5,1},
        xticklabels={$0$,$1/2$,$1$},
        ytick={0,0.5,1},
        xlabel={Weight $w$},
        ylabel={Path shares},
        ylabel style={xshift=-1mm},
        clip=false
      ]
        \addplot[
          draw=none,
          fill=ExampleChoiceOne,
          fill opacity=0.55
        ] coordinates {
          (0,0) (0.5,0) (0.5,1) (0,1)
        } \closedcycle;

        \addplot[
          draw=none,
          fill=ExampleChoiceTwo,
          fill opacity=0.55
        ] coordinates {
          (0.5,0) (1,0) (1,1) (0.5,1)
        } \closedcycle;

        \addplot[
          very thick,
          color=ExampleFrontier
        ] coordinates {
          (0,1) (0.5,1) (0.5,0) (1,0)
        };

        \addplot[
          dashed,
          very thick,
          color=ExampleWarm
        ] coordinates {
          (0,0) (0.5,0) (0.5,1) (1,1)
        };

        \addplot[
          thresholdline
        ] coordinates {
          (0.5,0) (0.5,1.02)
        };

        \node[font=\scriptsize, anchor=west, color=ExampleFrontier] at (axis cs:1.01,0.0) {$\tilde{x}_{1}(w)$};
        \node[font=\scriptsize, anchor=west, color=ExampleWarm] at (axis cs:1.01,1.0) {$\tilde{x}_{2}(w)$};

        \node[
          font=\tiny,
          anchor=south,
          fill=white,
          inner sep=1pt
        ] at (axis cs:0.25,0.07) {mass on path $1$};

        \node[
          font=\tiny,
          anchor=south,
          fill=white,
          inner sep=1pt
        ] at (axis cs:0.75,0.07) {mass on path $2$};
      \end{axis}
    \end{tikzpicture}
    \subcaption{Weight-dependent flow solution implied by convex combination for formulations of \citet{leur:93} and \citet{marc:99}.}
  \end{subfigure}
  \hfill
  \begin{subfigure}[t]{0.30\textwidth}
    \centering
    \begin{tikzpicture}
      \begin{axis}[
        exampleaxis,
        width=0.9\linewidth,
        height=5.8cm,
        axis lines=left,
        axis line style={-{Stealth[length=2mm]}, draw=black, line width=0.35pt},
        xmin=0, xmax=1.05,
        ymin=-0.02, ymax=1.05,
        xtick={0,0.5,1},
        ytick={0,0.5,1},
        xlabel={Weight $w$},
        ylabel={Path-$1$ share},
        ylabel style={xshift=-2mm},
        clip=false
      ]
        \addplot[
          very thick,
          color=ExampleFrontier
        ] coordinates {
          (0,1) (1,1)
        };

        \addplot[
          dashed,
          very thick,
          color=blue!65!black
        ] coordinates {
          (0,0) (1,0)
        };

        \addplot[
          dashdotted,
          very thick,
          color=ExamplePruned
        ] coordinates {
          (0,0.5) (1,0.5)
        };

        \node[font=\scriptsize, anchor=west, color=ExampleFrontier] at (axis cs:1.01,1.0) {$\hat{x}_{1}(w) = 1$};
        \node[font=\scriptsize, anchor=west, color=blue!65!black] at (axis cs:1.01,0.0) {$\check{x}_{1}(w) = 0$};
        \node[font=\scriptsize, anchor=west, color=ExamplePruned] at (axis cs:1.01,0.5) {$x_{1}(w) = \frac{1}{2}$};

        \node[
          font=\tiny,
          fill=white,
          inner sep=1pt
        ] at (axis cs:0.43,0.16) {$x_{2}(w) = 1 - x_{1}(w)$};
      \end{axis}
    \end{tikzpicture}
    \subcaption{Solutions and convex combination for formulation of \citet{marc:96a} and \citet{marc:96b}.}
  \end{subfigure}
\caption{Difference between convex combinations for aggregate-flow formulations and convex combinations for function-valued formulations for Example~\ref{ex:different convex combinations}.
}
\label{fig:aggregate-vs-functional-convex-combination}
\end{figure}
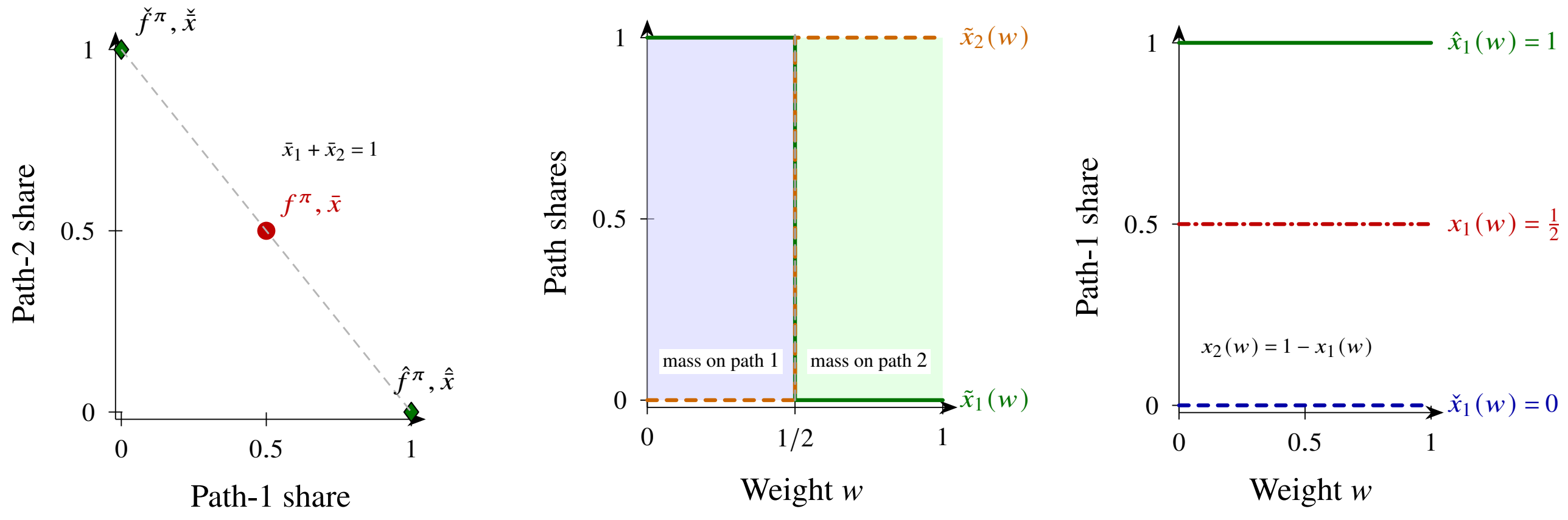

\citet{marc:99} compared a number of formulations for the BiTAP.
It was assumed that $G^{\pi}$ has a positive density $g$ on $[\underline{w},\overline{w}]$ for all $\pi$.
\citet{marc:99} showed that the formulations in \citet{leur:93} (with fixed demand $D^{\pi}$ and with the ratios $f^{\pi}_{r} / F^{\pi}$ replaced by $\bar{x}_{r}$) and in \citet{marc:96a} and \citet{marc:96b} are equivalent. 
However, it is important to note that the Frank-Wolfe framework applied to the formulations does not give equivalent algorithms.
The reason is that convex combinations of equivalent solutions for different formulations do not give equivalent solutions.
To illustrate that, consider the following simple example.
\begin{example}
\label{ex:different convex combinations}
The network has a single OD pair~$\pi$ with fixed demand $D^{\pi} = 1$, two paths indexed~$1$ and~$2$ from the origin to the destination with $1 \prec 2$, and weights~$w$ are uniformly distributed on $[\underline{w}^{\pi},\overline{w}^{\pi}] = [0,1]$ for both paths.
Solution~$1$ puts all flow on path~1.
Solution~$1$ for the formulations of \citet{leur:93} and \citet{marc:99} is depicted by the bottom right diamond in Figure~\ref{fig:aggregate-vs-functional-convex-combination}(a), and solution~$1$ for the formulation of \citet{marc:96a} and \citet{marc:96b} is depicted by the top horizontal line in Figure~\ref{fig:aggregate-vs-functional-convex-combination}(c).
Solution~$2$ puts all flow on path~2.
Solution~$2$ for the formulations of \citet{leur:93} and \citet{marc:99} is depicted by the top left diamond in Figure~\ref{fig:aggregate-vs-functional-convex-combination}(a), and solution~$2$ for the formulation of \citet{marc:96a} and \citet{marc:96b} is depicted by the bottom horizontal line in Figure~\ref{fig:aggregate-vs-functional-convex-combination}(c).

Next, consider the convex combination halfway between the two solutions.
The convex combination for the formulations of \citet{leur:93} and \citet{marc:99} is depicted by the middle dot in Figure~\ref{fig:aggregate-vs-functional-convex-combination}(a), and the convex combination for the formulation of \citet{marc:96a} and \citet{marc:96b} is depicted by the middle horizontal line in Figure~\ref{fig:aggregate-vs-functional-convex-combination}(c).
Also, recall from Example~\ref{ex:efficient frontier example} and the preceding explanation that, given solution $f^{\pi}$ or $\bar{x}$ and taking into account that $1 \prec 2$, there are thresholds $\underline{w}^{\pi} = w^{\pi}_{0} \le w^{\pi}_{1} \le w^{\pi}_{2} = \overline{w}^{\pi}$, such that users with weight $w \in (w^{\pi}_{i-1}, w^{\pi}_{i})$ travel on path~$i$ and thus contribute to flow $f^{\pi}_{i}$ or $\bar{x}_{i}$.
Thus, solution $f^{\pi}$ or $\bar{x}$ can be extended to weight-dependent flow functions $\tilde{x}_{i}(w)$.
(This statement is based on the distribution $G^{\pi}$ having a density, and we will show how this construction can be generalized for arbitrary distributions.)
Figure~\ref{fig:aggregate-vs-functional-convex-combination}(b) shows the resulting weight-dependent flow functions $\tilde{x}_{1}(w) = 1$ and $\tilde{x}_{2}(w) = 0$ for all $w \in (0,1/2)$ and $\tilde{x}_{1}(w) = 0$ and $\tilde{x}_{2}(w) = 1$ for all $w \in (1/2,1)$.
Note that the weight-dependent flow functions $\tilde{x}$ resulting from the convex combination $(1 - \alpha) \hat{f}^{\pi} + \alpha \check{f}^{\pi}$ is not equivalent to the convex combination $(1 - \alpha) \hat{x} + \alpha \check{x}$ depicted by the middle horizontal line in Figure~\ref{fig:aggregate-vs-functional-convex-combination}(c).
Also note that the weight-dependent flow functions $\tilde{x}$ has the interval structure illustrated in Example~\ref{ex:efficient frontier example} and the preceding explanation, while no strict convex combination $(1 - \alpha) \hat{x} + \alpha \check{x}$ has the interval structure.
This difference has important computational consequences.
\hfill\Halmos
\end{example}

\citet{marc:99} also considered the following formulation for the BiTAP based on thresholds $\underline{w}^{\pi} = w^{\pi}_{0} \le w^{\pi}_{1} \le \cdots \le w^{\pi}_{R^{\pi}-1} \le w^{\pi}_{R^{\pi}} = \overline{w}^{\pi}$.
Given~$w = (w^{\pi}_{r}, \pi \in \Pi, r \in \mathcal{R}^{\pi})$, the flow on each path $r \in \mathcal{R}^{\pi}$ is equal to $f^{\pi}_{r}(w^{\pi}) = D^{\pi} \int_{w^{\pi}_{r-1}}^{w^{\pi}_{r}} g^{\pi}(z) \diff z$, and the flow on each arc $a \in \mathcal{A}$ is given by $x_{a} = \sum_{\pi \in \Pi} \sum_{r \in \mathcal{R}^{\pi}_{a}} f^{\pi}_{r}(w^{\pi})$.
Assuming that costs $S^{\pi}_{r}$ do not depend on the flows, and with thresholds as decision variables, the equilibrium flows are given by the optimal solutions of the optimization problem
\begin{align}
\label{eqn:threshold formulation}
\min_{w} \ \ & \left\{\sum_{a \in \mathcal{A}} \int_{0}^{x_{a}} t_{a}(z) \diff z + \sum_{\pi \in \Pi} D^{\pi} \sum_{r \in \mathcal{R}^{\pi}} S^{\pi}_{r} \int_{w^{\pi}_{r-1}}^{w^{\pi}_{r}} z g^{\pi}(z) \diff z\right\} \\
\mbox{s.t.} \ \ & w^{\pi}_{r} - w^{\pi}_{r-1} \ \ \ge \ \ 0 \qquad \forall \ \pi \in \Pi, \ r \in \mathcal{R}^{\pi}
\label{eqn:monotone thresholds} \\
& w^{\pi}_{0} \ \ = \ \ \underline{w}, \qquad w^{\pi}_{R^{\pi}} \ \ = \ \ \overline{w} \qquad \forall \ \pi \in \Pi
\label{eqn:threshold boundaries}
\end{align}
As pointed out by \citet{marc:99}, one issue with this formulation is that the objective is not convex in~$w$ (the variational inequality is not monotone in~$w$).
On the positive side, it is shown that a point that satisfies the first-order optimality conditions of \eqref{eqn:threshold formulation}--\eqref{eqn:threshold boundaries} is an optimal solution of \eqref{eqn:threshold formulation}--\eqref{eqn:threshold boundaries}.

\citet{xie:26} considered the formulation \eqref{eqn:threshold formulation}--\eqref{eqn:threshold boundaries},
and confirmed that the first-order optimality conditions
for \eqref{eqn:threshold formulation}--\eqref{eqn:threshold boundaries} coincide with the equilibrium conditions for the BiTAP.
Similar to \citet{marc:96a} and \citet{marc:96b}, \citet{xie:26} used a parametric shortest path algorithm to generate efficient paths.
Unlike previously proposed algorithms based on Frank-Wolfe, instead of performing a line search between the current solution and an optimal solution of the Frank-Wolfe subproblem, the Single-Boundary Adjustment (SBA) algorithm proposed by \citet{xie:26} adjusted one threshold~$w^{\pi}_{r}$ at a time using an Armijo line search algorithm.

\section{Alternative Formulation of the Bicriteria Traffic Assignment Problem}
\label{sec:alternative formulation}

Recall that the formulation~\eqref{eqn:inf dimensional arc flow formulation} considered in \citet{marc:96a} and \citet{marc:96b} is a convex optimization problem, but that it has infinite dimensional decision variables (weight-dependent flow functions), and that even if two solutions have the interval structure possessed by optimal solutions, a strict convex combination of the two solutions loses the interval structure.
These properties are detrimental to efficient computation.
In contrast, the formulation~\eqref{eqn:threshold formulation}--\eqref{eqn:threshold boundaries} considered in \citet{marc:99} and \citet{xie:26} has finite dimensional decision variables that maintain the interval structure by construction, but the objective function is not convex.
In this section, we propose an alternative formulation.
Similar to the formulation~\eqref{eqn:path cumulative flow formulation} considered in \citet{leur:93} and \citet{leur:96}, our formulation is a convex optimization problem with finite dimensional decision variables.
All the previous formulations either require $G^{\pi}$ to have finite support or require $G^{\pi}$ to have a positive density on $[\underline{w}^{\pi},\overline{w}^{\pi}]$ for each $\pi \in \Pi$.
In contrast, our formulation does not require anything of $G^{\pi}$.
To simplify the presentation and algorithm, we also assume that $G^{\pi}$ has bounded support; this assumption can be removed at the cost of complicating the algorithm.

\subsection{Notation}

For a convex set $S \subset \mathbb{R}^{n}$, the relative interior of $S$ is denoted by $\text{ri}(S)$, and the normal cone to $S$ at $x$ is denoted by $\mathcal{N}_{S}(x)$.
Consider a convex function $\psi : \mathbb{R}^{n} \mapsto \mathbb{R} \cup \{\infty\}$.
The effective domain of $\psi$ is denoted by $\text{dom}(\psi) \defi \{x \in \mathbb{R}^{n} : \psi(x) < \infty\}$.
Let
\begin{equation}
\label{eqn:subdiff}
\partial \psi(z) \ \ \defi \ \ \{u \in \mathbb{R}^{n} \; : \; \psi(z') \ge \psi(z) + \inner{u}{z' - z} \ \forall \; z' \in \mathbb{R}^{n}\}
\end{equation}
denote the subdifferential of $\psi$ at a point $z \in \text{dom}(\psi)$.
For a finite set $\mathcal{I}$, let $\mathbb{R}^{\mathcal{I}}$ denote the vector space indexed by $\mathcal{I}$.
Let
\begin{equation*}
\label{eqn:simplex}
\Delta_{\mathcal{I}} \ \ \defi \ \ \left\{x \in \mathbb{R}^{\mathcal{I}}_+ \; : \; \sum_{i \in \mathcal{I}} x_{i} = 1\right\}
\end{equation*}
denote the standard simplex over $\mathcal{I}$.
Let
\begin{equation*}
\label{eqn:simplex_ri}
\Delta_{\mathcal{I}}^{\circ} \ \ \defi \ \ \left\{x \in \Delta_{\mathcal{I}} \; : \; x_{i} > 0 \ \forall \; i \in \mathcal{I}\right\}
\end{equation*}
denote the relative interior of $\Delta_{\mathcal{I}}$.
For any Borel set $S \subset \mathbb{R}$, let $\mathcal{P}(S)$ denote the set of cumulative distribution functions on $S$.
The quantile map for any $G \in \mathcal{P}([\underline{w}, \overline{w}])$ is defined as the set-valued map $G^{-1} : [0,1] \mapsto 2^{\mathbb{R}}$ given by
\begin{equation*}
\label{eqn:quantile}
G^{-1}(p) \ \ \defi \ \ \{x \in \mathbb{R} \; : \; G(y) \le p \ \forall \; y < x \text{ and } G(y) \ge p \ \forall \; y > x\}.
\end{equation*}
Note that $G^{-1}(p)$ is a nonempty closed interval for all $p \in (0,1)$, and $G^{-1}(0) = (-\infty,w]$ and $G^{-1}(1) = [w',\infty)$ for $\underline{w} \le w \le w' \le \overline{w}$.
The quantile map is non-decreasing in the sense that for any $0 \le a < b \le 1$, $x_{a} \in G^{-1}(a)$, and $x_{b} \in G^{-1}(b)$, it holds that $x_{a} \le x_{b}$.

For any $0 \leq a \leq b \leq 1$, let the truncation of $G \in \mathcal{P}([\underline{w}, \overline{w}])$ be given by
\begin{equation*}
\label{eqn:truncation_def}
G_{a,b}(x) \ \ \defi \ \ \max\{a, \min\{b, G(x)\}\}.
\end{equation*}
Note that $G_{a,b}$ is nondecreasing and right-continuous, and it induces a finite Lebesgue-Stieltjes measure on $\mathbb{R}$, which we also denote by $G_{a,b}$.
Its total mass satisfies
\begin{equation*}
\label{eqn:measure_value}
\int_{\mathbb{R}} \diff G_{a,b} \ \ = \ \ \lim_{x \to \infty} G_{a,b}(x) - \lim_{x \to -\infty} G_{a,b}(x) \ \ = \ \ b - a.
\end{equation*}
Note that for any integrable Borel measurable function $f : \mathbb{R} \mapsto \mathbb{R}$ and any $0 \leq a \leq b \leq c \leq 1$, the additivity property
\begin{equation}
\label{eqn:additivity}
\int_{\mathbb{R}} f(x) \, \diff G_{a,c}(x) \ \ = \ \ \int_{\mathbb{R}} f(x) \, \diff G_{a,b}(x) + \int_{\mathbb{R}} f(x) \, \diff G_{b,c}(x)
\end{equation}
holds.

\subsection{Formulation}

The variables $q^{\pi}_{r}$ denote the fraction of demand $D^{\pi}$ for OD pair $\pi \in \Pi$ that flows on path $r \in \mathcal{R}^{\pi}$.
Let $q \defi \left(q^{\pi}_{r}, \pi \in \Pi, r \in \mathcal{R}^{\pi}\right)$.
As before, the weight thresholds are denoted with $w = (w^{\pi}_{r}, \pi \in \Pi, r \in \mathcal{R}^{\pi})$.
Recall that for each OD pair $\pi \in \Pi$, we choose a total order $\prec$ on the set of paths $\mathcal{R}^{\pi}$ such that for any $r,r' \in \mathcal{R}^{\pi}$ it holds that if $r \prec r'$ (implying $r \neq r'$) then $S^{\pi}_{r} \ge S^{\pi}_{r'}$.
To make expressions more concise, for every $\pi \in \Pi$, and $r \in \mathcal{R}^{\pi}$, let
\[
p^{\pi}_{r} \ \ \defi \ \ \sum_{\{r' \in \mathcal{R}^{\pi} \, : \, r' \preceq r\}} q^{\pi}_{r'}
\]
denote the cumulative fraction of demand $D^{\pi}$ for OD pair $\pi \in \Pi$ that flows on paths $r' \in \mathcal{R}^{\pi}$ with $r' \preceq r$.
Also, let the linear operator $L^{\pi} : \mathbb{R}^{\mathcal{R}^{\pi}} \mapsto \mathbb{R}^{\mathcal{R}^{\pi}}$ be given by
\begin{equation*}
\label{eqn:linear map q to p}
L^{\pi}_{r}(q) \ \ \defi \ \ \sum_{\{r' \in \mathcal{R}^{\pi} \, : \, r' \preceq r\}} q_{r'}
\ \ = \ \ \left(\mathds{1}\{r' \preceq r\}, \; r' \in \mathcal{R}^{\pi}\right)^{\top} \left(q_{r'}, \; r' \in \mathcal{R}^{\pi}\right)
\end{equation*}
so that $p^{\pi} = L^{\pi}(q^{\pi})$.
Also, for each $r \in \mathcal{R}^{\pi}$, $r-1$ denotes the path that is the immediate predecessor of~$r$, so that $p^{\pi}_{r-1} \defi \sum_{\{r' \in \mathcal{R}^{\pi} \, : \, r' \prec r\}} q^{\pi}_{r'}$, with the convention that if $r$ has no predecessor in $\mathcal{R}^{\pi}$, then $p^{\pi}_{r-1} = 0$ and $w^{\pi}_{r-1} = \underline{w}^{\pi}$.
Given~$q$, for each arc~$a \in \mathcal{A}$, let
\begin{equation}
x_{a}(q) \ \ \defi \ \ \sum_{\pi \in \Pi} D^{\pi} \sum_{r \in \mathcal{R}^{\pi}_{a}} q^{\pi}_{r} \hspace{10mm} \forall \ a \in \mathcal{A}
\label{eqn:q reformulated BiTAP arc flows}
\end{equation}
denote the total flow on arc~$a$, and let $x(q) \defi (x_{a}(q), a \in \mathcal{A})$.

\subsubsection{Equilibrium Conditions}

The pair $(q,w)$ satisfies the equilibrium conditions if the following hold for all $\pi \in \Pi$:
\begin{enumerate}[nosep]
\item
\begin{equation}
\label{eqn:q probabilities}
q^{\pi} \ \in \ \Delta_{\mathcal{R}^{\pi}}, \quad \mbox{i.e.} \ \ \sum_{r \in \mathcal{R}^{\pi}} q^{\pi}_{r} \ = \ 1 \ \ \mbox{and} \ \ q^{\pi}_{r} \ \ge \ 0 \ \forall \; r \in \mathcal{R}^{\pi}
\end{equation}
that is, each $q^{\pi}$ is a probability distribution with support on $\mathcal{R}^{\pi}$.
\item
\begin{equation}
\label{eqn:w ordered}
\underline{w}^{\pi} \le \cdots \le w^{\pi}_{r-1} \le w^{\pi}_{r} \le \cdots \le \overline{w}^{\pi}
\end{equation}
that is, the thresholds $w^{\pi}_{r}$ are properly ordered between $\underline{w}^{\pi}$ and $\overline{w}^{\pi}$.
We can also require that the last threshold be equal to $\overline{w}^{\pi}$, but given the other conditions, it is not necessary.
\item
For each $r \in \mathcal{R}^{\pi}$ such that $q^{\pi}_{r} > 0$, it holds that
\begin{equation}
\label{eqn:each used path optimal in its interval}
C^{\pi}_{r}(w,x(q)) \ \le \ C^{\pi}_{r'}(w,x(q)) \quad \forall \; r' \in \mathcal{R}^{\pi}, \quad \forall \; w \in [w^{\pi}_{r-1},w^{\pi}_{r}]
\end{equation}
that is, each used path is user-optimal for all weights~$w$ in its interval.
\item
\begin{equation}
\label{eqn:probabilities and thresholds consistent}
w^{\pi}_{r} \ \in \ \left(G^{\pi}\right)^{-1}\big(L^{\pi}_{r}(q^{\pi})\big) \quad \forall \; r \in \mathcal{R}^{\pi}
\end{equation}
that is, the thresholds and probabilities are consistent with each other.
\end{enumerate}

\begin{example}
\label{ex:equilibrium conditions example}
Consider Example~\ref{ex:efficient frontier example}, and suppose that OD pair~$\pi$ with demand $D^{\pi} = 1$ is the only OD pair with positive demand.
Let
\begin{equation*}
G^{\pi}(w) \ \ = \ \ \left\{
\begin{array}{ll}
0 & \mbox{if } w < 0 \\
0.05 & \mbox{if } 0 \le w < 0.01 \\
0.3 + 5 \left(w - 0.01\right) & \mbox{if } 0.01 \le w < 0.07 \\
0.8 & \mbox{if } 0.07 \le w < 0.14 \\
1 & \mbox{if } w \ge 0.14
\end{array}\right.
\end{equation*}
as shown in Figure~\ref{fig:equilibrium-example-cdf}.
\begin{figure}[tbp]
  \centering
  \begin{tikzpicture}
    \begin{axis}[
      exampleaxis,
      width=0.72\textwidth,
      height=0.5\textwidth,
      axis lines=left,
      axis line style={-{Stealth[length=2mm]}, draw=black, line width=0.35pt},
      xmin=-0.006, xmax=0.155,
      ymin=0, ymax=1.05,
      xtick={0,0.01,0.03,0.07,0.12,0.14},
      xticklabels={$0$,$0.01$,$0.03$,$0.07$,$0.12$,$0.14$},
      ytick={0,0.1,0.3,0.4,0.6,0.7,0.8,1},
      yticklabels={$0$,$0.1$,$0.3$,$0.4$,$0.6$,$0.7$,$0.8$,$1$},
      clip=false
    ]
      \draw[thresholdline] (axis cs:-0.006,0.1) -- (axis cs:0.01,0.1) -- (axis cs:0.01,0);
      \draw[thresholdline] (axis cs:-0.006,0.4) -- (axis cs:0.03,0.4) -- (axis cs:0.03,0);
      \draw[thresholdline] (axis cs:-0.006,0.7) -- (axis cs:0.07,0.7) -- (axis cs:0.07,0);
      \draw[thresholdline] (axis cs:-0.006,0.8) -- (axis cs:0.07,0.8);
      \draw[thresholdline] (axis cs:0.12,0.8) -- (axis cs:0.12,0);
      \draw[thresholdline] (axis cs:-0.006,1.0) -- (axis cs:0.14,1.0);
      \draw[thresholdline] (axis cs:0.14,0.8) -- (axis cs:0.14,0);

      \node[font=\scriptsize, anchor=south west] at (axis cs:0.0112,0.012) {$w^{\pi}_{1}$};
      \node[font=\scriptsize, anchor=south west] at (axis cs:0.0312,0.012) {$w^{\pi}_{2}$};
      \node[font=\scriptsize, anchor=south west] at (axis cs:0.0712,0.012) {$w^{\pi}_{3}{=}\cdots{=}w^{\pi}_{8}$};
      \node[font=\scriptsize, anchor=south east] at (axis cs:0.1188,0.012) {$w^{\pi}_{9}$};
      \node[font=\scriptsize, anchor=south west] at (axis cs:0.1412,0.012) {$w^{\pi}_{10}$};

      \draw[very thick, color=blue!70!black] (axis cs:-0.006,0) -- (axis cs:0,0);
      \draw[color=blue!70!black, line width=0.6pt] (axis cs:0,0) -- (axis cs:0,0.05);
      \draw[very thick, color=blue!70!black] (axis cs:0,0.05) -- (axis cs:0.01,0.05);
      \draw[color=blue!70!black, line width=0.6pt] (axis cs:0.01,0.05) -- (axis cs:0.01,0.3);
      \draw[very thick, color=blue!70!black] (axis cs:0.01,0.3) -- (axis cs:0.07,0.6);
      \draw[color=blue!70!black, line width=0.6pt] (axis cs:0.07,0.6) -- (axis cs:0.07,0.8);
      \draw[very thick, color=blue!70!black] (axis cs:0.07,0.8) -- (axis cs:0.14,0.8);
      \draw[color=blue!70!black, line width=0.6pt] (axis cs:0.14,0.8) -- (axis cs:0.14,1);
      \draw[very thick, color=blue!70!black] (axis cs:0.14,1) -- (axis cs:0.152,1);
      \addplot[only marks, mark=none, mark size=1.7pt, color=blue!70!black] coordinates {(0,0.05) (0.01,0.3) (0.07,0.8) (0.14,1)};
      \addplot[only marks, mark=none, mark size=1.7pt, mark options={draw=blue!70!black, fill=white, line width=0.6pt}] coordinates {(0,0) (0.01,0.05) (0.07,0.6) (0.14,0.8)};

      \node[font=\small, anchor=south] at (axis cs:-0.006,1.07) {$G^{\pi}(w)$};
      \node[font=\small, anchor=west] at (axis cs:0.157,0) {$w$};

      \draw[ExampleFrontier, line width=1pt] (axis cs:-0.0175,0) -- (axis cs:-0.0175,0.1)
        (axis cs:-0.0175,0) -- (axis cs:-0.0163,0) (axis cs:-0.0175,0.1) -- (axis cs:-0.0163,0.1);
      \draw[ExampleWarm, line width=1pt] (axis cs:-0.0175,0.1) -- (axis cs:-0.0175,0.4)
        (axis cs:-0.0175,0.4) -- (axis cs:-0.0163,0.4);
      \draw[ExampleFrontier, line width=1pt] (axis cs:-0.0175,0.4) -- (axis cs:-0.0175,0.7)
        (axis cs:-0.0175,0.7) -- (axis cs:-0.0163,0.7);
      \draw[ExampleWarm, line width=1pt] (axis cs:-0.0175,0.7) -- (axis cs:-0.0175,0.8)
        (axis cs:-0.0175,0.8) -- (axis cs:-0.0163,0.8);
      \draw[ExampleFrontier, line width=1pt] (axis cs:-0.0175,0.8) -- (axis cs:-0.0175,1)
        (axis cs:-0.0175,1) -- (axis cs:-0.0163,1);
      \node[font=\scriptsize, anchor=east, color=ExampleFrontier] at (axis cs:-0.0195,0.05) {$q^{\pi}_{r_{1}} = 0.1$};
      \node[font=\scriptsize, anchor=east, color=ExampleWarm] at (axis cs:-0.0195,0.25) {$q^{\pi}_{r_{2}} = 0.3$};
      \node[font=\scriptsize, anchor=east, color=ExampleFrontier] at (axis cs:-0.0195,0.55) {$q^{\pi}_{r_{3}} = 0.3$};
      \node[font=\scriptsize, anchor=east, color=ExampleWarm] at (axis cs:-0.0195,0.75) {$q^{\pi}_{r_{6}} = 0.1$};
      \node[font=\scriptsize, anchor=east, color=ExampleFrontier] at (axis cs:-0.0195,0.9) {$q^{\pi}_{r_{10}} = 0.2$};
      \node[font=\scriptsize, anchor=east, color=gray!65] at (axis cs:-0.0195,0.69) {$q^{\pi}_{r_{4}}{=}q^{\pi}_{r_{5}}{=}0$};
      \node[font=\scriptsize, anchor=east, color=gray!65] at (axis cs:-0.0195,0.81) {$q^{\pi}_{r_{7}}{=}q^{\pi}_{r_{8}}{=}q^{\pi}_{r_{9}}{=}0$};

      \draw[ExampleFrontier, line width=1pt] (axis cs:0,-0.10) -- (axis cs:0.01,-0.10)
        (axis cs:0,-0.10) -- (axis cs:0,-0.082) (axis cs:0.01,-0.10) -- (axis cs:0.01,-0.082);
      \draw[ExampleWarm, line width=1pt] (axis cs:0.01,-0.10) -- (axis cs:0.03,-0.10)
        (axis cs:0.03,-0.10) -- (axis cs:0.03,-0.082);
      \draw[ExampleFrontier, line width=1pt] (axis cs:0.03,-0.10) -- (axis cs:0.07,-0.10)
        (axis cs:0.07,-0.10) -- (axis cs:0.07,-0.082);
      \draw[gray!65, line width=1pt] (axis cs:0.07,-0.10) -- (axis cs:0.12,-0.10)
        (axis cs:0.12,-0.10) -- (axis cs:0.12,-0.082);
      \draw[ExampleFrontier, line width=1pt] (axis cs:0.12,-0.10) -- (axis cs:0.14,-0.10)
        (axis cs:0.14,-0.10) -- (axis cs:0.14,-0.082);
      \fill[ExampleWarm] (axis cs:0.07,-0.10) circle (1.9pt);
      \node[font=\scriptsize, anchor=north, color=ExampleFrontier] at (axis cs:0.005,-0.125) {$r_{1}$};
      \node[font=\scriptsize, anchor=north, color=ExampleWarm] at (axis cs:0.02,-0.125) {$r_{2}$};
      \node[font=\scriptsize, anchor=north, color=ExampleFrontier] at (axis cs:0.05,-0.125) {$r_{3}$};
      \node[font=\scriptsize, anchor=north west, color=ExampleWarm] at (axis cs:0.0715,-0.125) {$r_{6}$};
      \node[font=\scriptsize, anchor=north, color=gray!65] at (axis cs:0.095,-0.125) {$r_{9}$};
      \node[font=\scriptsize, anchor=north, color=ExampleFrontier] at (axis cs:0.13,-0.125) {$r_{10}$};
    \end{axis}
  \end{tikzpicture}
\caption{The distribution $G^{\pi}$ of Example~\ref{ex:equilibrium conditions example}, the fractions $q^{\pi}_{r}$ as brackets to the left of the vertical axis, the cumulative fractions $p^{\pi}_{r}$ on the vertical axis, the thresholds $w^{\pi}_{r}$ on the horizontal axis, and the intervals $[w^{\pi}_{r-1},w^{\pi}_{r}]$ as brackets below the horizontal axis.
Dashed lines connect each cumulative fraction $p^{\pi}_{r}$ to the thresholds $w^{\pi}_{r} \in (G^{\pi})^{-1}(p^{\pi}_{r})$ of condition~\eqref{eqn:probabilities and thresholds consistent}.
}
\label{fig:equilibrium-example-cdf}
\end{figure}
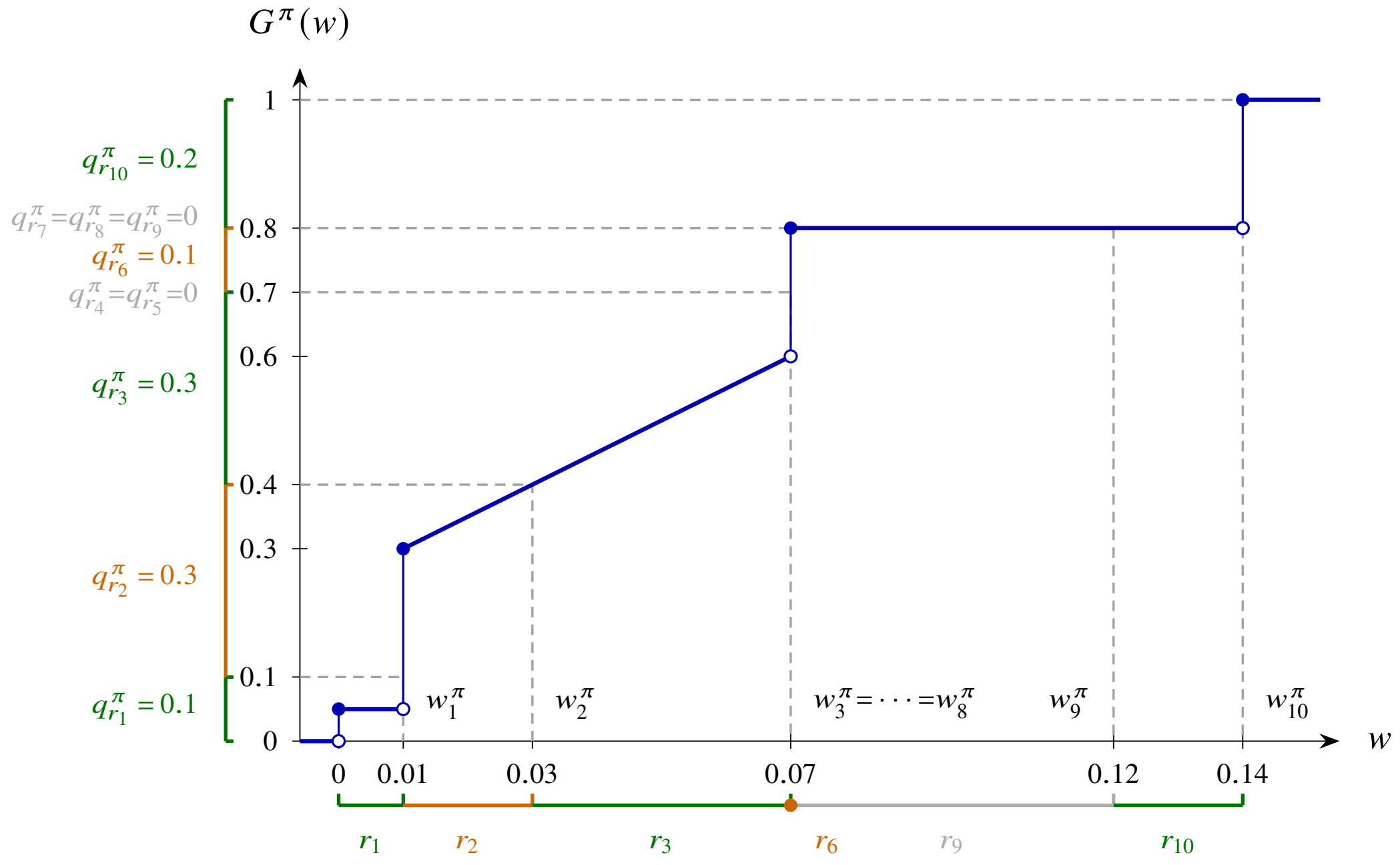
Let $\left(q^{\pi}_{r_{1}},\ldots,q^{\pi}_{r_{10}}\right) = \left(0.1, 0.3, 0.3, 0, 0, 0.1, 0, 0, 0, 0.2\right)$, so that condition~\eqref{eqn:q probabilities} holds and $\left(p^{\pi}_{r_{1}},\ldots,p^{\pi}_{r_{10}}\right) = L^{\pi}(q^{\pi}) = (0.1, 0.4, 0.7, 0.7, 0.7, 0.8, 0.8, 0.8, 0.8, 1)$.
Figure~\ref{fig:equilibrium-example-cdf} shows how $q^{\pi}$ partitions the demand among the used paths.
Consider $\left(w^{\pi}_{1},\ldots,w^{\pi}_{10}\right) = \left(0.01, 0.03, 0.07, 0.07, 0.07, 0.07, 0.07, 0.07, 0.12, 0.14\right)$, which are the thresholds of Example~\ref{ex:efficient frontier example} with the choice $\eta = 0.03$, so that conditions~\eqref{eqn:w ordered} and~\eqref{eqn:each used path optimal in its interval} hold.
Figure~\ref{fig:equilibrium-example-cdf} shows that condition~\eqref{eqn:probabilities and thresholds consistent} also holds, and thus $(q,w)$ is an equilibrium.
Recall that paths $r_{2}$ and $r_{3}$ are tied.
It is clear how their order could be swapped and thus $(q,w)$ is not a unique equilibrium.
Figure~\ref{fig:equilibrium-example-cdf} exhibits some features excluded by formulations that require a positive density and distinct costs: paths $r_{1}$ and $r_{2}$ split the jump at $w^{\pi}_{1}$, paths $r_{3}$ and $r_{6}$ split the jump at $w^{\pi}_{3} = \cdots = w^{\pi}_{8}$, path~$r_{6}$ carries positive flow although its interval is a singleton, and path~$r_{9}$ carries zero flow although its interval $[0.07,0.12]$ is not a singleton.
The example does not show how $T^{\pi}_{r}(x)$ specified in Example~\ref{ex:efficient frontier example} is obtained from the travel time functions $t_{a}$ and the flow partition $q$ --- that will be demonstrated in Example~\ref{ex:bta-psp-equilibrium-example}.
\hfill\Halmos
\end{example}


\subsubsection{Optimization Problem}

In this section we formulate a convex optimization problem that we use to compute an equilibrium.
In a typical solution of the traffic assignment problem, only a small fraction of the paths $r \in \mathcal{R}^{\pi}$ for an OD pair~$\pi \in \Pi$ is used (has $q^{\pi}_{r} > 0$).
For algorithmic efficiency, it is important to maintain sparse solutions, and our proposed algorithm is designed to consider small subsets of paths $\hat{\mathcal{R}}^{\pi} \subset \mathcal{R}^{\pi}$, $\pi \in \Pi$, and find the best solution for the subsets in each iteration.
Therefore, we extend the notation to make provision for the smaller as well as the full optimization problems.
For each $a \in \mathcal{A}$, let $\hat{\mathcal{R}}^{\pi}_{a} \defi \hat{\mathcal{R}}^{\pi} \cap \mathcal{R}^{\pi}_{a}$ denote the paths in $\hat{\mathcal{R}}^{\pi}$ that use arc~$a$.
Note that, for each $\pi \in \Pi$, the same order $\prec$ applies to every subset $\hat{\mathcal{R}}^{\pi} \subset \mathcal{R}^{\pi}$.
When a nonempty subset $\hat{\mathcal{R}}^{\pi} \subset \mathcal{R}^{\pi}$ of paths is considered, we will let $q \defi \left(q^{\pi}_{r}, \pi \in \Pi, r \in \hat{\mathcal{R}}^{\pi}\right)$, we will require that $\sum_{r \in \hat{\mathcal{R}}^{\pi}} q^{\pi}_{r} = 1$ and
\[
p^{\pi}_{r} \ \ \defi \ \ \sum_{\{r' \in \hat{\mathcal{R}}^{\pi} \, : \, r' \preceq r\}} q^{\pi}_{r'},
\]
and thus the linear operator $L^{\pi} : \mathbb{R}^{\hat{\mathcal{R}}^{\pi}} \mapsto \mathbb{R}^{\hat{\mathcal{R}}^{\pi}}$ will be given by
\[
L^{\pi}_{r}(q) \ \ \defi \ \ \sum_{\{r' \in \hat{\mathcal{R}}^{\pi} \, : \, r' \preceq r\}} q_{r'}
\ \ = \ \ \left(\mathds{1}\{r' \preceq r\}, \; r' \in \hat{\mathcal{R}}^{\pi}\right)^{\top} \left(q_{r'}, \; r' \in \hat{\mathcal{R}}^{\pi}\right)
\]
so that $p^{\pi} = L^{\pi}(q^{\pi})$.
Also, for each $r \in \hat{\mathcal{R}}^{\pi}$, $r-1$ will denote the path that is the immediate predecessor in $\hat{\mathcal{R}}^{\pi}$ of~$r$, so that $p^{\pi}_{r-1} \defi \sum_{\{r' \in \hat{\mathcal{R}}^{\pi} \, : \, r' \prec r\}} q^{\pi}_{r'}$.
(The subset $\hat{\mathcal{R}}^{\pi}$ will be clear from the context, and therefore is not displayed in the notation.)
Then, let $\phi : \prod_{\pi \in \Pi} \Delta_{\hat{\mathcal{R}}^{\pi}} \mapsto \mathbb{R}$ be given by
\begin{equation}
\phi(q) \ \ \defi \ \ \sum_{\pi \in \Pi} D^{\pi} \sum_{r \in \hat{\mathcal{R}}^{\pi}} S^{\pi}_{r} \int_{\mathbb{R}} y \, \diff G^{\pi}_{p^{\pi}_{r-1}, p^{\pi}_{r}}(y)
\label{eqn:toll objective}
\end{equation}
where, in this formulation, variables $p^{\pi}_{r}$ are not regarded as separate decision variables, but rather as linear functions of~$q$.
Also, let $\psi : \prod_{\pi \in \Pi} \Delta_{\hat{\mathcal{R}}^{\pi}} \mapsto \mathbb{R}$ be given by
\begin{equation}
\psi(q) \ \ \defi \ \ \sum_{a \in \mathcal{A}} \int_{0}^{x_{a}(q)} t_{a}(y) \, \diff y
\label{eqn:travel time objective}
\end{equation}
where
\begin{equation*}
x_{a}(q) \ \ \defi \ \ \sum_{\pi \in \Pi} D^{\pi} \sum_{r \in \hat{\mathcal{R}}^{\pi}_{a}} q^{\pi}_{r} \hspace{10mm} \forall \ a \in \mathcal{A}
\end{equation*}
denotes the total flow on arc~$a$.
For any subsets $\hat{\mathcal{R}}^{\pi} \subset \mathcal{R}^{\pi}$, $\pi \in \Pi$, the optimization problem is as follows:
\begin{align}
\min_{q} \ & \left\{\Upsilon(q) \ \defi \ \psi(q) + \phi(q) \ : \ q \in \prod_{\pi \in \Pi} \Delta_{\hat{\mathcal{R}}^{\pi}}\right\}
\label{eqn:q reformulated BiTAP objective}
\end{align}
A similar formulation can be written with $p^{\pi}_{r}$ as decision variables, with no apparent advantage of one formulation over the other.
Note that, whereas the equilibrium conditions are stated for a pair $(q,w)$, a solution of the optimization problem is $q$ only.
Thus, unlike previous approaches, the first-order optimality conditions of the optimization problem do not coincide with the equilibrium conditions.
Instead, we will show how, given a solution~$q$ of the optimization problem for $\hat{\mathcal{R}}^{\pi} = \mathcal{R}^{\pi}$, $\pi \in \Pi$, thresholds~$w$ can be determined in run time that is linear in the number of efficient paths, such that $(q,w)$ is an equilibrium.

\section{The Optimization Problem and Equilibria}
\label{sec:algorithm}

\subsection{Preliminary Results}
\label{sec:preliminary results}

First we confirm that the value of $\phi(q)$ is not affected by the choice of total order on $\mathcal{R}^{\pi}$.
\begin{lemma}
\label{lem:ordering_invariance}
For any $\pi \in \Pi$, consider two total orders $\prec$ and $\lessdot$ on $\mathcal{R}^{\pi}$ such that for any $r,r' \in \mathcal{R}^{\pi}$ it holds that if $r \prec r'$ or $r \lessdot r'$ then $S^{\pi}_{r} \ge S^{\pi}_{r'}$.
Consider any $\hat{\mathcal{R}}^{\pi} \subset \mathcal{R}^{\pi}$ and any $q^{\pi} \in \Delta_{\hat{\mathcal{R}}^{\pi}}$.
For each $r \in \hat{\mathcal{R}}^{\pi}$, let
\[
p^{\pi}_{r} \ \ \defi \ \ \sum_{\{r' \in \hat{\mathcal{R}}^{\pi} \, : \, r' \preceq r\}} q^{\pi}_{r'}
\qquad \mbox{and} \qquad
\dot{p}^{\pi}_{r} \ \ \defi \ \ \sum_{\{r' \in \hat{\mathcal{R}}^{\pi} \, : \, r' \ledot r\}} q^{\pi}_{r'}
\]
Then
\begin{equation*}
\sum_{r \in \hat{\mathcal{R}}^{\pi}} S^{\pi}_{r} \int_{\mathbb{R}} y \, \diff G^{\pi}_{p^{\pi}_{r-1}, p^{\pi}_{r}}(y)
\ \ = \ \ \sum_{r \in \hat{\mathcal{R}}^{\pi}} S^{\pi}_{r} \int_{\mathbb{R}} y \, \diff G^{\pi}_{\dot{p}^{\pi}_{r-1}, \dot{p}^{\pi}_{r}}(y)
\end{equation*}
and thus $\phi(q)$ is invariant with respect to the chosen order.
\end{lemma}

Next we establish some pertinent properties of the function $\phi$.
For that purpose, consider any $[\underline{w},\overline{w}] \subset \mathbb{R}_{+}$, any $G \in \mathcal{P}([\underline{w},\overline{w}])$, and the function $\varphi : [0,1] \to \mathbb{R}$ given by
\begin{equation*}
\label{def:varphi}
\varphi(p) \ \ \defi \ \ \int_{\mathbb{R}} w \, \diff G_{0,p}(w)
\end{equation*}
It follows from $[\underline{w},\overline{w}] \subset \mathbb{R}_{+}$ that $\varphi$ is nondecreasing.
Note that $\varphi(0) = 0$ and $\varphi(1) \in [\underline{w},\overline{w}]$.
We will also use the result in Lemma~\ref{lem:varphi technical properties}.

\begin{lemma}
\label{lem:varphi technical properties}
The function $\varphi : [0,1] \to \mathbb{R}$ is convex and Lipschitz continuous, and $\partial \varphi(p) = G^{-1}(p)$ for every $p \in [0,1]$.
\end{lemma}

For every $\pi \in \Pi$, let $\varphi^{\pi} : [0,1] \to \mathbb{R}$ be given by
\begin{equation*}
\label{def:varphi-pi}
\varphi^{\pi}(p) \ \ \defi \ \ \int_{\mathbb{R}} w \, \diff G^{\pi}_{0,p}(w)
\end{equation*}
and given any subset of paths $\hat{\mathcal{R}}^{\pi} \subset \mathcal{R}^{\pi}$ and any $r \in \hat{\mathcal{R}}^{\pi}$, let $\phi^{\pi}_{r} : \Delta_{\hat{\mathcal{R}}^{\pi}} \to \mathbb{R}$ be given by
\begin{equation*}
\label{def:phi-pi}
\phi^{\pi}_{r}(q) \ \ \defi \ \ \varphi^{\pi}\big(L^{\pi}_{r}(q)\big)
\end{equation*}
Then it follows from Lemma~\ref{lem:varphi technical properties} that $\phi^{\pi}_{r} : \Delta_{\hat{\mathcal{R}}^{\pi}} \to \mathbb{R}$ is convex and Lipschitz continuous.
Also, $\text{range}(L^{\pi}_{r}) \cap \text{ri}\big(\text{dom}(\varphi^{\pi})\big) = (0,1) \neq \varnothing$, and thus it follows from \citet[Theorem 23.9]{rockafellar1997convex} that
\[
\partial \phi^{\pi}_{r}(q) \ \ = \ \ \left(\mathds{1}\left\{r' \preceq r\right\}, \; r' \in \hat{\mathcal{R}}^{\pi}\right) \left(G^{\pi}\right)^{-1}\big(L^{\pi}_{r}(q)\big) + \mathcal{N}_{\Delta_{\hat{\mathcal{R}}^{\pi}}}(q)
\]
Let $S^{\pi}_{\min} \defi \min\{S^{\pi}_{r} \, : \, r \in \hat{\mathcal{R}}^{\pi}\}$.
For each $r \in \hat{\mathcal{R}}^{\pi}$ that has a successor in $\hat{\mathcal{R}}^{\pi}$ in the order $\prec$, let $r+1 \in \hat{\mathcal{R}}^{\pi}$ denote the immediate successor in $\hat{\mathcal{R}}^{\pi}$ of $r$.
Note that, if $r$ does not have an immediate successor in $\hat{\mathcal{R}}^{\pi}$, then $S^{\pi}_{r} = S^{\pi}_{\min}$, and $p^{\pi}_{r} = L^{\pi}_{r}(q^{\pi}) = 1$ for any $q^{\pi} \in \Delta_{\hat{\mathcal{R}}^{\pi}}$, and thus $S^{\pi}_{r} \varphi^{\pi}(p^{\pi}_{r}) = S^{\pi}_{\min} \varphi^{\pi}(1)$ is a constant independent of $q^{\pi}$.
Hence
\begin{align*}
\phi(q) \ \ & = \ \ \sum_{\pi \in \Pi} D^{\pi} \sum_{r \in \hat{\mathcal{R}}^{\pi}} S^{\pi}_{r} \left[\varphi^{\pi}(p^{\pi}_{r}) - \varphi^{\pi}(p^{\pi}_{r-1})\right] \\
& = \ \ \sum_{\pi \in \Pi} D^{\pi} \sum_{r \in \hat{\mathcal{R}}^{\pi}} \left[\mathds{1}\{r+1 \in \hat{\mathcal{R}}^{\pi}\} \left(S^{\pi}_{r} - S^{\pi}_{r+1}\right) \phi^{\pi}_{r}(q^{\pi}) + \mathds{1}\{r+1 \notin \hat{\mathcal{R}}^{\pi}\} S^{\pi}_{\min} \varphi^{\pi}(1)\right]
\end{align*}

\begin{lemma}
\label{lem:phi technical properties}
The function $\phi : \prod_{\pi \in \Pi} \Delta_{\hat{\mathcal{R}}^{\pi}} \to \mathbb{R}$ is convex and continuous, and
\[
\partial \phi(q) \ \ = \ \ \prod_{\pi \in \Pi} \left(D^{\pi} \sum_{\{r \in \hat{\mathcal{R}}^{\pi} \, : \, r+1 \in \hat{\mathcal{R}}^{\pi}\}} \left(S^{\pi}_{r} - S^{\pi}_{r+1}\right) \left(\mathds{1}\left\{r' \preceq r\right\}, \; r' \in \hat{\mathcal{R}}^{\pi}\right) \left(G^{\pi}\right)^{-1}\big(L^{\pi}_{r}(q^{\pi})\big) + \mathcal{N}_{\Delta_{\hat{\mathcal{R}}^{\pi}}}(q^{\pi})\right)
\]
\end{lemma}

Recall that $t_{a} : \mathbb{R}_{+} \mapsto \mathbb{R}_{+}$ is assumed to be nondecreasing.
Next we show that it is sufficient for $\psi$ to be convex.
For any $x > 0$, let $t_{a}^{-}(x) \defi \lim_{y \uparrow x} t_{a}(y)$ denote the left limit of $t_{a}$ at~$x$, and for any $x \ge 0$, let $t_{a}^{+}(x) \defi \lim_{y \downarrow x} t_{a}(y)$ denote the right limit of $t_{a}$ at~$x$.
It follows from $t_{a}$ being nondecreasing that $t_{a}^{-}(x) \le t_{a}^{+}(x)$ for all $x > 0$.

\begin{lemma}
\label{lem:psi technical properties}
The function $\psi$ is convex and continuous on $\prod_{\pi \in \Pi} \Delta_{\hat{\mathcal{R}}^{\pi}}$, and
\[
\partial \psi(q) \ \ = \ \ \left\{g = \left(g^{\pi}_{r}, \pi \in \Pi, r \in \hat{\mathcal{R}}^{\pi}\right) \; : \; g^{\pi}_{r} = D^{\pi} \sum_{a \in \mathcal{A}^{\pi}_{r}} \tau_{a}, \; \tau_{a} \in \left[t_{a}^{-}(x_{a}(q)), t_{a}^{+}(x_{a}(q))\right] \ \forall \; a \in \mathcal{A}\right\} + \prod_{\pi \in \Pi} \mathcal{N}_{\Delta_{\hat{\mathcal{R}}^{\pi}}}(q^{\pi})
\]
with $\left[t_{a}^{-}(0), t_{a}^{+}(0)\right] \defi \left(-\infty, t_{a}^{+}(0)\right]$.
\end{lemma}

To convert an optimal solution~$q$ of problem~\eqref{eqn:q reformulated BiTAP objective} to a pair~$(q,w)$ that satisfies the equilibrium conditions~\eqref{eqn:q probabilities}--\eqref{eqn:probabilities and thresholds consistent}, we will assume that $t_{a} : \mathbb{R}_{+} \mapsto \mathbb{R}_{+}$ is continuous, in addition to being nondecreasing.
Then Lemma~\ref{lem:psi technical properties} simplifies to Lemma~\ref{lem:psi continuous technical properties}.

\begin{lemma}
\label{lem:psi continuous technical properties}
If $t_{a} : \mathbb{R}_{+} \mapsto \mathbb{R}_{+}$ is continuous, then
\[
\partial \psi(q) \ \ = \ \ \left\{g = \left(g^{\pi}_{r}, \pi \in \Pi, r \in \hat{\mathcal{R}}^{\pi}\right) \; : \; g^{\pi}_{r} = D^{\pi} \sum_{a \in \mathcal{A}^{\pi}_{r}} t_{a}(x_{a}(q)) = D^{\pi} T^{\pi}_{r}(x(q))\right\} + \prod_{\pi \in \Pi} \mathcal{N}_{\Delta_{\hat{\mathcal{R}}^{\pi}}}(q^{\pi})
\]
\end{lemma}

\subsection{Optimality Conditions}
\label{sec:optimality conditions}

Next we establish necessary and sufficient conditions for $q^{*} \in \prod_{\pi \in \Pi} \Delta_{\hat{\mathcal{R}}^{\pi}}$ to be optimal for problem~\eqref{eqn:q reformulated BiTAP objective}.
For each $\pi \in \Pi$, let $\lambda^{\pi}$ denote the dual variable associated with the simplex upper-bound constraint $\sum_{r \in \hat{\mathcal{R}}^{\pi}} q^{\pi}_{r} = 1$, and let $\lambda \defi (\lambda^{\pi}, \pi \in \Pi)$.

\begin{lemma}
\label{lem:KKT equivalence}
A vector $q^{*} \in \prod_{\pi \in \Pi} \Delta_{\hat{\mathcal{R}}^{\pi}}$ is optimal for problem~\eqref{eqn:q reformulated BiTAP objective} if and only if there exist $\tau^* = (\tau^*_{a}, a \in \mathcal{A})$ with $\tau^*_{a} \in \left[t_{a}^{-}(x_{a}(q^*)), t_{a}^{+}(x_{a}(q^*))\right]$, $w^{*} = \left(w^{\pi*}_{r}, \pi \in \Pi, r \in \hat{\mathcal{R}}^{\pi}\right)$ with $w^{\pi*}_{r} \in \left(G^{\pi}\right)^{-1}\left(L^{\pi}_{r}(q^{\pi*})\right)$, and $\lambda^{*} = (\lambda^{\pi*}, \pi \in \Pi)$, such that the following holds:
\begin{align}
& q^{\pi*}_{r} \left\{D^{\pi} \left[\sum_{a \in \mathcal{A}^{\pi}_{r}} \tau^*_{a} + \sum_{\{r' \in \hat{\mathcal{R}}^{\pi} \, : \, r \preceq r', r'+1 \in \hat{\mathcal{R}}^{\pi}\}} \left(S^{\pi}_{r'} - S^{\pi}_{r'+1}\right) w^{\pi*}_{r'}\right] + \lambda^{\pi*}\right\} \ \ = \ \ 0 \qquad \forall \ \pi \in \Pi, \ \forall \ r \in \hat{\mathcal{R}}^{\pi}
\label{eqn:q reformulated BiTAP optimality1} \\
& D^{\pi} \left[\sum_{a \in \mathcal{A}^{\pi}_{r}} \tau^*_{a} + \sum_{\{r' \in \hat{\mathcal{R}}^{\pi} \, : \, r \preceq r', r'+1 \in \hat{\mathcal{R}}^{\pi}\}} \left(S^{\pi}_{r'} - S^{\pi}_{r'+1}\right) w^{\pi*}_{r'}\right] + \lambda^{\pi*} \ \ \ge \ \ 0 \qquad \forall \ \pi \in \Pi, \ \forall \ r \in \hat{\mathcal{R}}^{\pi}
\label{eqn:q reformulated BiTAP optimality3}
\end{align}
\end{lemma}

If $t_{a} : \mathbb{R}_{+} \mapsto \mathbb{R}_{+}$ is continuous, in addition to being nondecreasing, then Lemma~\ref{lem:KKT equivalence} simplifies to Lemma~\ref{lem:KKT continuous equivalence}.

\begin{lemma}
\label{lem:KKT continuous equivalence}
Suppose that $t_{a} : \mathbb{R}_{+} \mapsto \mathbb{R}_{+}$ is continuous.
Then, a vector $q^{*} \in \prod_{\pi \in \Pi} \Delta_{\hat{\mathcal{R}}^{\pi}}$ is optimal for problem~\eqref{eqn:q reformulated BiTAP objective} if and only if there exist $w^{*} = \left(w^{\pi*}_{r}, \pi \in \Pi, r \in \hat{\mathcal{R}}^{\pi}\right)$ with $w^{\pi*}_{r} \in \left(G^{\pi}\right)^{-1}\left(L^{\pi}_{r}(q^{\pi*})\right)$ and $\lambda^{*} = (\lambda^{\pi*}, \pi \in \Pi)$, such that the following holds:
\begin{align}
& q^{\pi*}_{r} \left\{D^{\pi} \left[T^{\pi}_{r}(x(q^*)) + \sum_{\{r' \in \hat{\mathcal{R}}^{\pi} \, : \, r \preceq r', r'+1 \in \hat{\mathcal{R}}^{\pi}\}} \left(S^{\pi}_{r'} - S^{\pi}_{r'+1}\right) w^{\pi*}_{r'}\right] + \lambda^{\pi*}\right\} \ \ = \ \ 0 \qquad \forall \ \pi \in \Pi, \ \forall \ r \in \hat{\mathcal{R}}^{\pi}
\label{eqn:q continuous BiTAP optimality1} \\
& D^{\pi} \left[T^{\pi}_{r}(x(q^*)) + \sum_{\{r' \in \hat{\mathcal{R}}^{\pi} \, : \, r \preceq r', r'+1 \in \hat{\mathcal{R}}^{\pi}\}} \left(S^{\pi}_{r'} - S^{\pi}_{r'+1}\right) w^{\pi*}_{r'}\right] + \lambda^{\pi*} \ \ \ge \ \ 0 \qquad \forall \ \pi \in \Pi, \ \forall \ r \in \hat{\mathcal{R}}^{\pi}
\label{eqn:q continuous BiTAP optimality3}
\end{align}
\end{lemma}



\subsection{From Optimality to Equilibrium}
\label{sec:optimality to equilibrium}

Recall that optimization problem~\eqref{eqn:q reformulated BiTAP objective} has solution~$q$, whereas the equilibrium conditions~\eqref{eqn:q probabilities}--\eqref{eqn:probabilities and thresholds consistent} are stated for a pair~$(q,w)$.
Next we present Algorithm~\ref{algo:threshold}, that takes any subsets~$\hat{\mathcal{R}}^{\pi}$, $\pi \in \Pi$, of paths, and optimal solution~$q$ for problem~\eqref{eqn:q reformulated BiTAP objective}, as input, and produces the subset $\mathcal{R}^{\pi}_{\min} \subset \hat{\mathcal{R}}^{\pi}$ of paths on the efficient frontier, with thresholds $w = \left(w^{\pi}_{r}, \pi \in \Pi, r \in \{0\} \cup \hat{\mathcal{R}}^{\pi}\right)$, as output.
Algorithm~\ref{algo:threshold} makes provision for multiple tied paths, that is, paths~$r$ and $r'$ with $S^{\pi}_{r} = S^{\pi}_{r'}$ and $T^{\pi}_{r}(x(q)) = T^{\pi}_{r'}(x(q))$, with zero or positive flow.

\begin{algorithm}
\SingleSpacedXI
\caption{Threshold Selection Algorithm $\mathrm{Thresholds}\big((\mathcal{N},\mathcal{A}),\Pi,\big((G^{\pi},[\underline{w}^{\pi},\overline{w}^{\pi}]), \pi \in \Pi\big),\big((t_{a},s_{a}), a \in \mathcal{A}\big),\big(\hat{\mathcal{R}}^{\pi}, \pi \in \Pi\big),q\big)$}
\label{algo:threshold}
\begin{algorithmic}[1]
\Require Network $(\mathcal{N},\mathcal{A})$, set $\Pi$ of OD pairs with cumulative distribution function $G^{\pi}$ with support interval $[\underline{w}^{\pi},\overline{w}^{\pi}]$ for each $\pi \in \Pi$, cost functions $t_{a} : \mathbb{R}_{+} \mapsto \mathbb{R}$ and costs $s_{a} \in \mathbb{R}_{+}$ for each $a \in \mathcal{A}$, any nonempty path sets $\hat{\mathcal{R}}^{\pi}$ indexed with $r = 1,\ldots,\hat{R}^{\pi}$ such that $S^{\pi}_{1} \ge \cdots \ge S^{\pi}_{\hat{R}^{\pi}}$ for each $\pi \in \Pi$, optimal solution $q \defi (q^{\pi}_{r}, \pi \in \Pi, r = 1,\ldots,\hat{R}^{\pi}) \in \prod_{\pi \in \Pi} \Delta_{\hat{\mathcal{R}}^{\pi}}$ for problem~\eqref{eqn:q reformulated BiTAP objective}
\For{$\pi \in \Pi$}
\State
\vspace{-2mm}
\begin{equation}
\hspace{-85mm} r_{1} \ \ \leftarrow \ \ \min\left\{\arg\min\left\{T^{\pi}_{r}(x(q)) + \underline{w}^{\pi} S^{\pi}_{r} \; : \; r \in \hat{\mathcal{R}}^{\pi}\right\}\right\}
\label{eqn:efficient path 1}
\end{equation}
\For{$r = 0,\ldots,r_{1}-1$}
\State $w^{\pi}_{r} \leftarrow \underline{w}^{\pi}$
\EndFor
\State
$\mathcal{R}^{\pi}_{\min} \leftarrow \{r_{1}\}$,
\ $\ell \leftarrow 1$
\While{$r_{\ell} < \hat{R}^{\pi}$}
\If{$\left(S^{\pi}_{r_{\ell}} > S^{\pi}_{r_{\ell}+1}\right)$}
\State
\vspace{-4mm}
\begin{align}
\hspace{-5mm} & w^{\min}_{\ell} \ \ \leftarrow \ \ \min\left\{\frac{T^{\pi}_{r}(x(q)) - T^{\pi}_{r_{\ell}}(x(q))}{S^{\pi}_{r_{\ell}} - S^{\pi}_{r}} \; : \; r \in \{r_{\ell}+1,\ldots,\hat{R}^{\pi}\}\right\}, \qquad w^{\pi}_{r_{\ell}} \ \ \leftarrow \ \ \min\left\{w^{\min}_{\ell}, \, \overline{w}^{\pi}\right\}
\label{eqn:threshold >} \\
\hspace{-5mm} & r_{\ell+1} \ \ \leftarrow \ \ \min\left\{\arg\min\left\{\frac{T^{\pi}_{r}(x(q)) - T^{\pi}_{r_{\ell}}(x(q))}{S^{\pi}_{r_{\ell}} - S^{\pi}_{r}} \; : \; r \in \{r_{\ell}+1,\ldots,\hat{R}^{\pi}\}\right\}\right\}
\label{eqn:efficient path >}
\end{align}
\Else \ $\left(S^{\pi}_{r_{\ell}} = S^{\pi}_{r_{\ell}+1}\right)$
\State
\[
\hspace{-55mm} T^{\min}_{\ell} \ \ \leftarrow \ \ \min\left\{T^{\pi}_{r}(x(q)) \; : \; S^{\pi}_{r} = S^{\pi}_{r_{\ell}}, \; r \in \{r_{\ell}+1,\ldots,\hat{R}^{\pi}\}\right\}
\]
\If{$\left(T^{\min}_{\ell} = T^{\pi}_{r_{\ell}}(x(q))\right)$}
\State
\vspace{-4mm}
\begin{align}
\hspace{-12mm} & w^{\min}_{\ell} \ \ \leftarrow \ \ w^{\pi}_{r_{\ell}} \ \ \leftarrow \ \ \max\left\{w^{\pi}_{r_{\ell}-1}, \; \min\left\{w \in \left(G^{\pi}\right)^{-1}\left(L^{\pi}_{r_{\ell}}(q^{\pi})\right)\right\}\right\}
\label{eqn:threshold ==} \\
\hspace{-12mm} & r_{\ell+1} \ \ \leftarrow \ \ \min\left\{\arg\min\left\{T^{\pi}_{r}(x(q)) \; : \; S^{\pi}_{r} = S^{\pi}_{r_{\ell}}, \; r \in \{r_{\ell}+1,\ldots,\hat{R}^{\pi}\}\right\}\right\}
\label{eqn:efficient path ==}
\end{align}
\Else \ $\left(T^{\min}_{\ell} > T^{\pi}_{r_{\ell}}(x(q))\right)$
\If{$\left(\left\{r \in \hat{\mathcal{R}}^{\pi} \; : \; S^{\pi}_{r_{\ell}} > S^{\pi}_{r}\right\} \neq \varnothing\right)$}
\State
\vspace{-4mm}
\begin{align}
\hspace{14mm} & w^{\min}_{\ell} \ \ \leftarrow \ \ \min\left\{\frac{T^{\pi}_{r}(x(q)) - T^{\pi}_{r_{\ell}}(x(q))}{S^{\pi}_{r_{\ell}} - S^{\pi}_{r}} \; : \; S^{\pi}_{r_{\ell}} > S^{\pi}_{r}, \; r \in \hat{\mathcal{R}}^{\pi}\right\}, \qquad w^{\pi}_{r_{\ell}} \ \ \leftarrow \ \ \min\left\{w^{\min}_{\ell}, \, \overline{w}^{\pi}\right\}
\label{eqn:threshold =>} \\
& r_{\ell+1} \ \ \leftarrow \ \ \min\left\{\arg\min\left\{\frac{T^{\pi}_{r}(x(q)) - T^{\pi}_{r_{\ell}}(x(q))}{S^{\pi}_{r_{\ell}} - S^{\pi}_{r}} \; : \; S^{\pi}_{r_{\ell}} > S^{\pi}_{r}, \; r \in \hat{\mathcal{R}}^{\pi}\right\}\right\}
\label{eqn:efficient path =>}
\end{align}
\Else \ $\left(\left\{r \in \hat{\mathcal{R}}^{\pi} \; : \; S^{\pi}_{r_{\ell}} > S^{\pi}_{r}\right\} = \varnothing\right)$
\State
\vspace{-2mm}
\begin{equation}
\hspace{-60mm} w^{\min}_{\ell} \ \ \leftarrow \ \ \overline{w}^{\pi} + 1, \qquad w^{\pi}_{r_{\ell}} \ \ \leftarrow \ \ \overline{w}^{\pi}
\label{eqn:threshold =0}
\end{equation}
\EndIf
\EndIf
\EndIf
\If{$w^{\min}_{\ell} \le \overline{w}^{\pi}$}
\State
$\mathcal{R}^{\pi}_{\min} \ \ \leftarrow \ \ \mathcal{R}^{\pi}_{\min} \cup \{r_{\ell+1}\}$
\Else
\State
$r_{\ell+1} \leftarrow \hat{R}^{\pi}$
\EndIf
\For{$r = r_{\ell}+1,\ldots,r_{\ell+1}-1$}
\State $w^{\pi}_{r} \ \ \leftarrow \ \ w^{\pi}_{r_{\ell}}$
\EndFor
\State $\ell \ \ \leftarrow \ \ \ell+1$
\EndWhile
\State
\vspace{-2mm}
\begin{equation}
\label{eqn:last threshold}
\hspace{-130mm} w^{\pi}_{\hat{R}^{\pi}} \ \ \leftarrow \ \ \overline{w}^{\pi}
\end{equation}
\EndFor
\State Output $w = \left(w^{\pi}_{r}, \pi \in \Pi, r \in \{0\} \cup \hat{\mathcal{R}}^{\pi}\right)$ and $\mathcal{R}^{\pi}_{\min}$, $\pi \in \Pi$
\end{algorithmic}
\end{algorithm}

A brief summary of the steps of Algorithm~\ref{algo:threshold} is given next.
First, for each OD pair $\pi$, the algorithm selects the path with the least cost $C^{\pi}_{r}(\underline{w}^{\pi},x(q))$, initializes $w^{\pi}_{r} \leftarrow \underline{w}^{\pi}$ for lower indexed (inefficient) paths, and starts the list $\mathcal{R}^{\pi}_{\min}$ of efficient paths.
Second, each iteration applies one of the cases \eqref{eqn:threshold >}, \eqref{eqn:threshold ==}, \eqref{eqn:threshold =>}, or \eqref{eqn:threshold =0} to determine the next threshold value and the next efficient path.
Third, whenever the next threshold value is feasible, the algorithm appends the next efficient path to the list $\mathcal{R}^{\pi}_{\min}$, and fills in the threshold values for inefficient paths with intermediate indices, while preserving monotonicity of the threshold sequence.
Finally, it assigns terminal inefficient paths threshold values of $\overline{w}^{\pi}$, and outputs both the thresholds~$w$ and the lists $\mathcal{R}^{\pi}_{\min}$ of efficient paths.

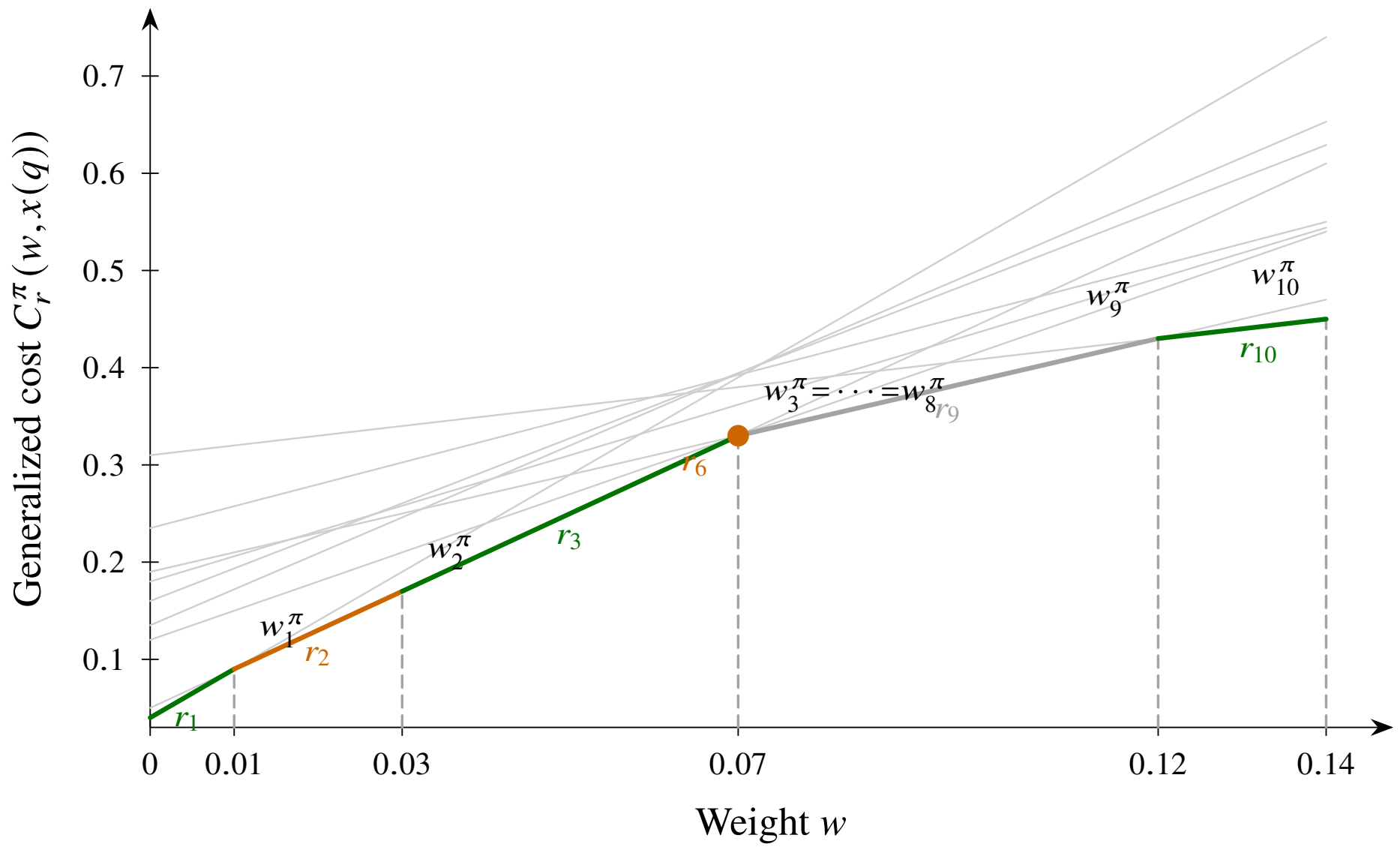
\begin{figure}[tbp]
  \centering
  \begin{tikzpicture}
    \begin{axis}[
      exampleaxis,
      width=0.76\textwidth,
      height=0.48\textwidth,
      axis lines=left,
      axis line style={-{Stealth[length=2mm]}, draw=black, line width=0.35pt},
      xmin=0, xmax=0.148,
      ymin=0.03, ymax=0.77,
      xtick={0,0.01,0.03,0.07,0.12,0.14},
      xticklabels={$0$,$0.01$,$0.03$,$0.07$,$0.12$,$0.14$},
      ytick={0.1,0.2,0.3,0.4,0.5,0.6,0.7},
      xlabel={Weight $w$},
      ylabel={Generalized cost $C^{\pi}_{r}(w,x(q))$},
      clip=false
    ]
      \addplot[domain=0:0.14, samples=2, color=gray!38, line width=0.45pt] {0.04+5*x};
      \addplot[domain=0:0.14, samples=2, color=gray!38, line width=0.45pt] {0.05+4*x};
      \addplot[domain=0:0.14, samples=2, color=gray!38, line width=0.45pt] {0.135+3.7*x};
      \addplot[domain=0:0.14, samples=2, color=gray!38, line width=0.45pt] {0.16+3.35*x};
      \addplot[domain=0:0.14, samples=2, color=gray!38, line width=0.45pt] {0.12+3*x};
      \addplot[domain=0:0.14, samples=2, color=gray!38, line width=0.45pt] {0.18+2.6*x};
      \addplot[domain=0:0.14, samples=2, color=gray!38, line width=0.45pt] {0.235+2.25*x};
      \addplot[domain=0:0.14, samples=2, color=gray!38, line width=0.45pt] {0.19+2*x};
      \addplot[domain=0:0.14, samples=2, color=gray!38, line width=0.45pt] {0.31+x};

      \addplot[thresholdline] coordinates {(0.01,0.03) (0.01,0.09)};
      \addplot[thresholdline] coordinates {(0.03,0.03) (0.03,0.17)};
      \addplot[thresholdline] coordinates {(0.07,0.03) (0.07,0.33)};
      \addplot[thresholdline] coordinates {(0.12,0.03) (0.12,0.43)};
      \addplot[thresholdline] coordinates {(0.14,0.03) (0.14,0.45)};

      \addplot[domain=0:0.01, samples=2, very thick, color=ExampleFrontier] {0.04+5*x};
      \addplot[domain=0.01:0.03, samples=2, very thick, color=ExampleWarm] {0.05+4*x};
      \addplot[domain=0.03:0.07, samples=2, very thick, color=ExampleFrontier] {0.05+4*x};
      \addplot[domain=0.07:0.12, samples=2, very thick, color=gray!72] {0.19+2*x};
      \addplot[domain=0.12:0.14, samples=2, very thick, color=ExampleFrontier] {0.31+x};

      \addplot[only marks, mark=*, mark size=2.5pt, draw=ExampleWarm, fill=ExampleWarm]
        coordinates {(0.07,0.33)};

      \node[font=\scriptsize, anchor=north, color=ExampleFrontier] at (axis cs:0.0045,0.063) {$r_{1}$};
      \node[font=\scriptsize, anchor=north, color=ExampleWarm] at (axis cs:0.020,0.130) {$r_{2}$};
      \node[font=\scriptsize, anchor=north, color=ExampleFrontier] at (axis cs:0.050,0.250) {$r_{3}$};
      \node[font=\scriptsize, anchor=north east, color=ExampleWarm] at (axis cs:0.0685,0.326) {$r_{6}$};
      \node[font=\scriptsize, anchor=north, color=gray!72] at (axis cs:0.095,0.380) {$r_{9}$};
      \node[font=\scriptsize, anchor=north, color=ExampleFrontier] at (axis cs:0.132,0.442) {$r_{10}$};

      \node[font=\scriptsize, anchor=south west] at (axis cs:0.011,0.095) {$w^{\pi}_{1}$};
      \node[font=\scriptsize, anchor=south west] at (axis cs:0.031,0.175) {$w^{\pi}_{2}$};
      \node[font=\scriptsize, anchor=south west] at (axis cs:0.071,0.337) {$w^{\pi}_{3}{=}\cdots{=}w^{\pi}_{8}$};
      \node[font=\scriptsize, anchor=south east] at (axis cs:0.119,0.435) {$w^{\pi}_{9}$};
      \node[font=\scriptsize, anchor=south east] at (axis cs:0.139,0.457) {$w^{\pi}_{10}$};
    \end{axis}
  \end{tikzpicture}
  \caption{Application of Threshold Algorithm~\ref{algo:threshold} to the data of Examples~\ref{ex:efficient frontier example} and~\ref{ex:equilibrium conditions example}.
  The thin gray lines are the distinct bicriteria path-cost functions $C^{\pi}_{r}(w,x)$, $r \in \hat{\mathcal{R}}^{\pi}$ (note that the functions for $r_{2}$ and $r_{3}$ coincide), and the thick segments show the intervals assigned along the lower envelope.
  The values $0.01$, $0.07$, and $0.12$ are intersection thresholds where the slope of the lower envelope changes.
  The value $w^{\pi}_{2} = 0.03$ is not a kink, rather it is the quantile threshold that separates the intervals assigned to the tied paths $r_{2}$ and $r_{3}$.
  The orange point at $w = 0.07$ is the singleton interval assigned to $r_{6}$, and the gray segment is the positive-length optimality interval of the zero-flow path $r_{9}$.}
  \label{fig:threshold-selection-equilibrium-example}
\end{figure}

\begin{example}
\label{ex:threshold-selection-equilibrium-example}
In this example, we apply Algorithm~\ref{algo:threshold} to the path data of Example~\ref{ex:efficient frontier example} and the path fractions~$q$ of Example~\ref{ex:equilibrium conditions example}.
Example~\ref{ex:bta-psp-equilibrium-example} will verify that these fractions~$q$ solve problem~\eqref{eqn:q reformulated BiTAP objective}.

At $\underline{w}^{\pi} = 0$, relation~\eqref{eqn:efficient path 1} selects path~$r_{1}$.
The minimum in~\eqref{eqn:threshold >} is $0.01$, attained by both paths $r_{2}$ and $r_{3}$, so the minimum-index rule in~\eqref{eqn:efficient path >} selects $r_{2}$ and sets $w^{\pi}_{1} = \min\{0.01,0.14\}$.
Because $S^{\pi}_{r_{2}} = S^{\pi}_{r_{3}}$, $T^{\pi}_{r_{2}}(x(q)) = T^{\pi}_{r_{3}}(x(q))$, and $L^{\pi}_{r_{2}}(q^{\pi}) = 0.4$, case~\eqref{eqn:threshold ==} gives $w^{\pi}_{2} = \max\left\{0.01,\min\left\{w \in (G^{\pi})^{-1}(0.4)\right\}\right\} = 0.03$, and relation~\eqref{eqn:efficient path ==} selects $r_{3}$.

Starting from $r_{3}$, the minimum in~\eqref{eqn:threshold >} is $0.07$, attained by $r_{6}$ and $r_{9}$.
The minimum-index rule selects $r_{6}$, sets $w^{\pi}_{3} = \min\{0.07,0.14\}$, and the fill step sets $w^{\pi}_{4} = w^{\pi}_{5} = 0.07$.
Next, starting from $r_{6}$, the minimum in~\eqref{eqn:threshold >} is also $0.07$, attained by $r_{9}$, setting $w^{\pi}_{6} = \min\{0.07,0.14\}$, and the fill step sets $w^{\pi}_{7} = w^{\pi}_{8} = 0.07$.
Starting from $r_{9}$, the minimum in~\eqref{eqn:threshold >} is $0.12$, attained by $r_{10}$, setting $w^{\pi}_{9} = \min\{0.12,0.14\}$.
Finally, the terminal step sets $w^{\pi}_{10} = \overline{w}^{\pi} = 0.14$.
Consequently,
\begin{equation}
\label{eqn:threshold-selection-example-output}
\mathcal{R}^{\pi}_{\min} \ \ = \ \ \{r_{1},r_{2},r_{3},r_{6},r_{9},r_{10}\},
\qquad
(w^{\pi}_{0},\ldots,w^{\pi}_{10}) \ \ = \ \ (0,.01,.03,.07,.07,.07,.07,.07,.07,.12,.14)
\end{equation}
which is the threshold vector of Example~\ref{ex:equilibrium conditions example}.
Figure~\ref{fig:equilibrium-example-cdf} verifies that~\eqref{eqn:probabilities and thresholds consistent} holds, and Figure~\ref{fig:threshold-selection-equilibrium-example} shows the corresponding optimality intervals.
\hfill\Halmos
\end{example}

In the remainder of the paper we assume that $t_{a} : \mathbb{R}_{+} \mapsto \mathbb{R}_{+}$ is continuous.
Example~\ref{ex:discontinuous no equilibrium} shows that continuity of $t_{a}$ is needed to guarantee existence of an equilibrium.

\begin{example}
\label{ex:discontinuous no equilibrium}
Consider a network with a single OD pair~$\pi$ with demand $D^{\pi} = 1$, and two disjoint paths, indexed~$1$ and~$2$.
The weight~$w$ is $0$ with probability~$1$, that is, it is a classical single-criterion traffic assignment problem.
The travel time on path~$1$ is $t_{1}(x_{1}) = 1 + 2 \mathds{1}\{x_{1} \ge 1/2\}$, and the travel time on path~$2$ is $t_{2}(x_{2}) = 2$.
Note that $x_{1}(q) = q_{1}$ and $x_{2}(q) = q_{2}$.
Any solution with $q_{1} < 1/2$ satisfies $q_{2} > 1/2$ but $C^{\pi}_{1}(w,x(q)) < C^{\pi}_{2}(w,x(q))$, and any solution with $q_{1} \ge 1/2$ satisfies $C^{\pi}_{1}(w,x(q)) > C^{\pi}_{2}(w,x(q))$.
Thus there is no solution that satisfies equilibrium condition~\eqref{eqn:each used path optimal in its interval}.
\hfill\Halmos
\end{example}

\begin{lemma}
\label{lem:optimal thresholds}
Assume that $t_{a} : \mathbb{R}_{+} \mapsto \mathbb{R}_{+}$ is continuous for all $a \in \mathcal{A}$.
Consider any optimal solution $q$ for problem~\eqref{eqn:q reformulated BiTAP objective}, and the resulting output $w$ and $\mathcal{R}^{\pi}_{\min}$, $\pi \in \Pi$, of Algorithm~\ref{algo:threshold}.
Then the following hold:
\begin{enumerate}[label=(\arabic*), ref=(\arabic*)]
\item
\label{item:thresholds ordered}
The thresholds $w$ satisfy $\underline{w}^{\pi} = w^{\pi}_{0} \le w^{\pi}_{1} \le \cdots \le w^{\pi}_{\hat{R}^{\pi}} = \overline{w}^{\pi}$ for all $\pi \in \Pi$.
\item
\label{item:positive flow paths are efficient}
Every path with positive flow is an efficient path, that is, $\left\{r \in \hat{\mathcal{R}}^{\pi} \, : \, q^{\pi}_{r} > 0\right\} \subset \mathcal{R}^{\pi}_{\min}$.
\item
\label{item:thresholds quantiles}
The thresholds $w$ satisfy $w^{\pi}_{r} \in \left(G^{\pi}\right)^{-1}\left(L^{\pi}_{r}(q^{\pi})\right)$ for all $\pi \in \Pi$, $r \in \hat{\mathcal{R}}^{\pi}$.
\item
\label{item:thresholds optimal}
For every $\pi \in \Pi$ and every $r_{\ell} \in \mathcal{R}^{\pi}_{\min}$ it holds that
\begin{equation}
\label{eqn:optimal intervals}
T^{\pi}_{r_{\ell}}(x(q)) + w S^{\pi}_{r_{\ell}} \ \ \le \ \ T^{\pi}_{r}(x(q)) + w S^{\pi}_{r}
\qquad \forall \ r \in \hat{\mathcal{R}}^{\pi}, \ \forall \ w \in [w^{\pi}_{r_{\ell-1}},w^{\pi}_{r_{\ell}}]
\end{equation}
where $w^{\pi}_{r_{0}} = w^{\pi}_{0} = \underline{w}^{\pi}$.
\end{enumerate}
\end{lemma}

\begin{theorem}
\label{thm:optimal solution to equilibrium}
Assume that $t_{a} : \mathbb{R}_{+} \mapsto \mathbb{R}_{+}$ is continuous for all $a \in \mathcal{A}$.
Consider any optimal solution $q$ for problem~\eqref{eqn:q reformulated BiTAP objective} with $\hat{\mathcal{R}}^{\pi} = \mathcal{R}^{\pi}$ for all $\pi \in \Pi$, and the resulting output $w$ of Algorithm~\ref{algo:threshold}.
Then $(q,w)$ satisfies the equilibrium conditions~\eqref{eqn:q probabilities}--\eqref{eqn:probabilities and thresholds consistent}.
\end{theorem}

\subsection{Solving the Optimization Problem}
\label{sec:solve optimization}

Section~\ref{sec:optimality to equilibrium} showed how an optimal solution~$q$ of optimization problem~\eqref{eqn:q reformulated BiTAP objective} can be converted to an equilibrium~$(q,w)$ that satisfies conditions~\eqref{eqn:q probabilities}--\eqref{eqn:probabilities and thresholds consistent}.
In this section we propose an algorithm to solve optimization problem~\eqref{eqn:q reformulated BiTAP objective} for the full set $\hat{\mathcal{R}}^{\pi} = \mathcal{R}^{\pi}$, $\pi \in \Pi$, of OD paths.
Lemma~\ref{lem:thresholds optimality conditions} shows that the optimality conditions~\eqref{eqn:q continuous BiTAP optimality1}--\eqref{eqn:q continuous BiTAP optimality3} are satisfied by $(q,\lambda,w)$ for any optimal solution~$q$ of problem~\eqref{eqn:q reformulated BiTAP objective}, $\lambda$ given by~\eqref{eqn:choice of dual}, and with $w$ given by Algorithm~\ref{algo:threshold}.

\ignore{
\begin{lemma}
\label{lem:thresholds optimality conditions}
Consider any optimal solution $q^*$ for problem~\eqref{eqn:q reformulated BiTAP objective}, and the resulting output \linebreak $w^* = \left(w^{(u,v)*}_{r}, (u,v) \in \Pi, r \in \{0\} \cup \mathcal{R}^{(u,v)}\right)$ and $\mathcal{R}^{(u,v)*}_{\min} = \{r_{1},\ldots,r_{L}\}$ of Algorithm~\ref{algo:threshold}, and arc flows $x^*$ given by $q^*$ and~\eqref{eqn:q reformulated BiTAP arc flows}.
Let
\begin{equation}
\lambda^{(u,v)*} \ \ \defi \ \ \left\{
\begin{array}{lcl}
0 & \mbox{if} & q^{(u,v)*}_{R^{(u,v)}} > 0 \\
D^{(u,v)} \left[T^{(u,v)}_{R^{(u,v)}}(x^*) - T^{(u,v)}_{r_{L}}(x^*) + \left(S^{(u,v)}_{R^{(u,v)}} - S^{(u,v)}_{r_{L}}\right) w^{(u,v)*}_{r_{L}}\right] & \mbox{if} & q^{(u,v)*}_{R^{(u,v)}} = 0
\end{array}\right.
\label{eqn:choice of dual}
\end{equation}
Then it holds that $(q^*,\lambda^*,w^*)$ satisfies the optimality conditions~\eqref{eqn:q reformulated BiTAP optimality1}--\eqref{eqn:q reformulated BiTAP complementarity}.
\end{lemma}

\begin{proof}
Proof of Lemma~\ref{lem:thresholds optimality conditions} is deferred to Appendix~\ref{app:threshold-proofs}.
\end{proof}
}

\begin{lemma}
\label{lem:thresholds optimality conditions}
Assume that $t_{a} : \mathbb{R}_{+} \mapsto \mathbb{R}_{+}$ is continuous for all $a \in \mathcal{A}$.
Consider any optimal solution $q$ for problem~\eqref{eqn:q reformulated BiTAP objective}, and the resulting output $w$ and $\mathcal{R}^{\pi}_{\min} = \{r_{1},\ldots,r_{L}\}$ of Algorithm~\ref{algo:threshold}.
Let $\hat{\mathcal{R}}^{\pi}$ be indexed with $r = 1,\ldots,\hat{R}^{\pi}$, and let
\begin{equation}
\lambda^{\pi} \ \ \defi \ \ -D^{\pi} \left[T^{\pi}_{r_{L}}(x(q)) + \overline{w}^{\pi} \left(S^{\pi}_{r_{L}} - S^{\pi}_{\hat{R}^{\pi}}\right)\right]
\qquad \forall \ \pi \in \Pi
\label{eqn:choice of dual}
\end{equation}
Then $(q,\lambda,w)$ satisfies the optimality conditions~\eqref{eqn:q continuous BiTAP optimality1}--\eqref{eqn:q continuous BiTAP optimality3}.
\end{lemma}

Problem~\eqref{eqn:q reformulated BiTAP objective} is a convex optimization problem, and therefore can be solved by general convex optimization algorithms.
However, typical instances of the traffic assignment problem have many paths for most OD pairs, only a small fraction of which have positive flow in an equilibrium.
Thus, we propose Algorithm~BTA (Algorithm~\ref{algo:BTA}) that exploits both convexity as well as the sparsity of typical solutions.
More specifically, Algorithm~BTA maintains a (typically small) subset $\hat{\mathcal{R}}^{\pi} \subset \mathcal{R}^{\pi}$ of paths for each OD pair $\pi \in \Pi$, solves the convex problem~\eqref{eqn:q reformulated BiTAP objective} for the subsets $\hat{\mathcal{R}}^{\pi}$, $\pi \in \Pi$, to produce flow fractions $\hat{q} = \left(\hat{q}^{\pi}_{r}, \pi \in \Pi, r \in \hat{\mathcal{R}}^{\pi}\right)$, update arc flows $x_{a}(\hat{q})$ and travel times $t_{a}(x_{a}(\hat{q}))$, $a \in \mathcal{A}$, and then solves a parametric shortest path problem for each origin $u \in \mathcal{N}$ such that $\mathcal{V}_{u} \defi \{v \in \mathcal{N} \, : \, (u,v) \in \Pi\} \neq \varnothing$.
Algorithms for the parametric shortest path problem have been proposed in \citet{karp1981parametric}.
Algorithm~\ref{algo:BPSP0} solves the parametric shortest path problem for a given origin $u \in \mathcal{N}$.
If the parametric shortest path problems produce paths not yet in the subsets $\hat{\mathcal{R}}^{\pi}$, then the new paths are added to $\hat{\mathcal{R}}^{\pi}$, and Algorithm~BTA returns to solving problem~\eqref{eqn:q reformulated BiTAP objective} for the updated subsets $\hat{\mathcal{R}}^{\pi}$, $\pi \in \Pi$; otherwise Algorithm~BTA returns $\hat{q} = \left(\hat{q}^{\pi}_{r}, \pi \in \Pi, r \in \hat{\mathcal{R}}^{\pi}\right)$.
Note that $\hat{q}$ contains path flow fractions $\hat{q}^{\pi}_{r}$ for paths $r \in \hat{\mathcal{R}}^{\pi}$ only.
Lemma~\ref{lem:optimality conditions} shows that the expanded solution~$q^*$ with $q^{\pi*}_{r} = \hat{q}^{\pi}_{r}$ for $r \in \hat{\mathcal{R}}^{\pi}$ and $q^{\pi*}_{r} = 0$ for $r \in \mathcal{R}^{\pi} \setminus \hat{\mathcal{R}}^{\pi}$ is optimal for problem~\eqref{eqn:q reformulated BiTAP objective} with the complete set of paths.
As shown in Section~\ref{sec:optimality to equilibrium}, Algorithm~\ref{algo:threshold} can be used to convert $q^*$ to an equilibrium~$(q^*,w^*)$ with the complete set of paths.


\begin{algorithm}
\SingleSpacedXI
\caption{Bicriteria Traffic Assignment Algorithm $\mathrm{BTA}\big((\mathcal{N},\mathcal{A}),\Pi,\big((D^{\pi},G^{\pi},[\underline{w}^{\pi},\overline{w}^{\pi}]), \pi \in \Pi\big),\big((t_{a},s_{a}), a \in \mathcal{A}\big)\big)$}
\label{algo:BTA}
\begin{algorithmic}[1]
\Require Network $(\mathcal{N},\mathcal{A})$, set $\Pi$ of OD pairs with demand $D^{\pi}$ and cumulative distribution function $G^{\pi}$ with support $[\underline{w}^{\pi},\overline{w}^{\pi}]$ for each $\pi \in \Pi$, cost functions $t_{a} : \mathbb{R}_{+} \mapsto \mathbb{R}_{+}$ and costs $s_{a} \in \mathbb{R}_{+}$ for each $a \in \mathcal{A}$
\For{$u \in \mathcal{N}$}
\State $\mathcal{V}_{u} \leftarrow \left\{v \in \mathcal{N} \, : \, (u,v) \in \Pi\right\}$
\If{$\mathcal{V}_{u} \neq \varnothing$}
\State $\big((\hat{\mathcal{R}}^{(u,v)},w^{(u,v)}), v \in \mathcal{V}_{u}\big) \leftarrow \mathrm{BPSP}\big((\mathcal{N},\mathcal{A}),u,\mathcal{V}_{u},((t_{a}(0),s_{a}),a \in \mathcal{A}),([\underline{w}^{(u,v)},\overline{w}^{(u,v)}], v \in \mathcal{V}_{u})\big)$
\EndIf
\EndFor
\State $\eta \leftarrow 1$
\While{$\eta = 1$}
\State $\eta \leftarrow 0$
\State $\hat{q} \leftarrow \arg\min\left\{\Upsilon(q) \ : \ q \in \prod_{\pi \in \Pi} \Delta_{\hat{\mathcal{R}}^{\pi}}\right\}$
\For{$u \in \mathcal{N}$ such that $\mathcal{V}_{u} \neq \varnothing$}
\State $\big((\mathcal{R}^{(u,v)+},w^{(u,v)+}), v \in \mathcal{V}_{u}\big) \leftarrow \mathrm{BPSP}\big((\mathcal{N},\mathcal{A}),u,\mathcal{V}_{u},((t_{a}(x_{a}(\hat{q})),s_{a}), a \in \mathcal{A}),([\underline{w}^{(u,v)},\overline{w}^{(u,v)}], v \in \mathcal{V}_{u})\big)$
\For{$v \in \mathcal{V}_{u}$}
\If{$\mathcal{R}^{(u,v)+} \not\subset \hat{\mathcal{R}}^{(u,v)}$}
\State $\hat{\mathcal{R}}^{(u,v)} \leftarrow \hat{\mathcal{R}}^{(u,v)} \cup \mathcal{R}^{(u,v)+}$
\State $\eta \leftarrow 1$
\EndIf
\EndFor
\EndFor
\EndWhile
\State Output $\hat{q}$ and $\hat{\mathcal{R}}^{\pi}$ for all $\pi \in \Pi$
\end{algorithmic}
\end{algorithm}

\begin{algorithm}
\SingleSpacedXI
\caption{Bicriteria Parametric Shortest Path Algorithm from One Origin to All Destinations $\mathrm{BPSP}\big((\mathcal{N},\mathcal{A}),u,\mathcal{V}_{u},((t_{a},s_{a}),a \in \mathcal{A}),([\underline{w}^{(u,v)},\overline{w}^{(u,v)}], v \in \mathcal{V}_{u})\big)$}
\label{algo:BPSP0}
\begin{algorithmic}[1]
\Require Network $(\mathcal{N},\mathcal{A})$, origin~$u$, set $\mathcal{V}_{u}$ of destinations, costs $t_{(i,j)} \in \mathbb{R}_{+}$ and $s_{(i,j)} \in \mathbb{R}_{+}$ for each $(i,j) \in \mathcal{A}$, support intervals $[\underline{w}^{(u,v)},\overline{w}^{(u,v)}]$ for all $v \in \mathcal{V}_{u}$
\State $\underline{w}_{u} \leftarrow \min\left\{\underline{w}^{(u,v)} \, : \, v \in \mathcal{V}_{u}\right\}$ and $\overline{w}_{u} \leftarrow \max\left\{\overline{w}^{(u,v)} \, : \, v \in \mathcal{V}_{u}\right\}$
\State $w \leftarrow w' \leftarrow \underline{w}_{u}$
\State Compute a shortest path tree $(\mathcal{N},\mathcal{A}_{u})$ from origin~$u$ to destinations $v \in \mathcal{V}_{u}$ on network $(\mathcal{N},\mathcal{A})$ with arc costs $c_{(i,j)} = t_{(i,j)} + w s_{(i,j)}$
\For{$i \in \mathcal{N}$}
\State $T_{i} \leftarrow \mbox{the $t$-cost from $u$ to $i$ in the shortest path tree } (\mathcal{N},\mathcal{A}_{u})$
\State $S_{i} \leftarrow \mbox{the $s$-cost from $u$ to $i$ in the shortest path tree } (\mathcal{N},\mathcal{A}_{u})$
\EndFor
\For{$v \in \mathcal{V}_{u}$}
\If{$\underline{w}^{(u,v)} \le w$}
\State $\mathcal{R}^{(u,v)+} \leftarrow \left\{\mbox{path from $u$ to $v$ in the shortest path tree } (\mathcal{N},\mathcal{A}_{u})\right\}$
\State $R^{(u,v)+} \leftarrow 1$
\Else
\State $\mathcal{R}^{(u,v)+} \leftarrow \varnothing$
\State $R^{(u,v)+} \leftarrow 0$
\EndIf
\State $w^{(u,v)+}_{0} \leftarrow \underline{w}^{(u,v)}$
\EndFor
\While{$w' \le \overline{w}_{u}$ and $\left\{(i,j) \in \mathcal{A} \setminus \mathcal{A}_{u} \; : \; S_{i} + s_{(i,j)} < S_{j}\right\} \neq \varnothing$}
\State
\begin{align*}
w' \ \ & \leftarrow \ \ \min\left\{\frac{T_{i} + t_{(i,j)} - T_{j}}{S_{j} - S_{i} - s_{(i,j)}} \; : \; (i,j) \in \mathcal{A} \setminus \mathcal{A}_{u}, \; S_{i} + s_{(i,j)} < S_{j}\right\} \\
\mathcal{A}' \ \ & \leftarrow \ \ \arg\min\left\{\frac{T_{i} + t_{(i,j)} - T_{j}}{S_{j} - S_{i} - s_{(i,j)}} \; : \; (i,j) \in \mathcal{A} \setminus \mathcal{A}_{u}, \; S_{i} + s_{(i,j)} < S_{j}\right\}
\end{align*}
\If{$w' \le \overline{w}_{u}$}
\State Select $(i,j) \in \mathcal{A}'$
\State Find $(i',j) \in \mathcal{A}_{u}$
\State $\mathcal{A}'_{u} \leftarrow \mathcal{A}_{u} \cup \{(i,j)\} \setminus \{(i',j)\}$
\State $\mathcal{N}' \leftarrow \left\{j\right\}$
\State $\Delta T \leftarrow T_{i} + t_{(i,j)} - T_{j}$
\State $\Delta S \leftarrow S_{i} + s_{(i,j)} - S_{j}$
\While{$\mathcal{N}' \neq \varnothing$}
\State Remove $i'' \in \mathcal{N}'$
\State Update $T_{i''} \leftarrow T_{i''} + \Delta T$ and $S_{i''} \leftarrow S_{i''} + \Delta S$
\If{$i'' \in \mathcal{V}_{u}$}
\State Update $(\mathcal{R}^{(u,i'')+}, w^{(u,i'')+}) \leftarrow \mathrm{UpdatePSP}\big((\mathcal{N},\mathcal{A}_{u}),(\mathcal{N},\mathcal{A}'_{u}),u,i'',[\underline{w}^{(u,i'')},\overline{w}^{(u,i'')}],w,w',\mathcal{R}^{(u,i'')+}, w^{(u,i'')+}\big)$
\EndIf
\For{$(i'',j') \in \mathcal{A}_{u}$}
\State $\mathcal{N}' \leftarrow \mathcal{N}' \cup \{j'\}$
\EndFor
\EndWhile
\State $w \leftarrow w'$
\State $\mathcal{A}_{u} \leftarrow \mathcal{A}'_{u}$
\EndIf
\EndWhile
\For{$v \in \mathcal{V}_{u}$}
\If{$w < \underline{w}^{(u,v)}$}
\State $\mathcal{R}^{(u,v)+} \leftarrow \left\{\mbox{path from $u$ to $v$ in the tree } (\mathcal{N},\mathcal{A}_{u})\right\}$
\State $w^{(u,v)+}_{1} \leftarrow \overline{w}^{(u,v)}$
\State $R^{(u,v)+} \leftarrow 1$
\Else \ {($\underline{w}^{(u,v)} \le w$)}
\State $w^{(u,v)+}_{R^{(u,v)+}} \leftarrow \overline{w}^{(u,v)}$
\EndIf
\EndFor
\State Output $\mathcal{R}^{(u,v)+}$ and $w^{(u,v)+}_{0},\ldots,w^{(u,v)+}_{R^{(u,v)+}}$ for all $v \in \mathcal{V}_{u}$
\end{algorithmic}
\end{algorithm}

\begin{algorithm}
\SingleSpacedXI
\caption{Update Parametric Shortest Paths and Thresholds $\mathrm{UpdatePSP}\big((\mathcal{N},\mathcal{A}_{u}),(\mathcal{N},\mathcal{A}'_{u}),u,v,[\underline{w}^{(u,v)},\overline{w}^{(u,v)}],w,w',\mathcal{R}^{(u,v)+},w^{(u,v)+}\big)$}
\label{algo:path update}
\begin{algorithmic}[1]
\Require Trees $(\mathcal{N},\mathcal{A}_{u})$ and $(\mathcal{N},\mathcal{A}'_{u})$, origin~$u$, destination~$v$, support interval $[\underline{w}^{(u,v)},\overline{w}^{(u,v)}]$, current thresholds $w \le w'$, path list $\mathcal{R}^{(u,v)+}$ with $|\mathcal{R}^{(u,v)+}| = R^{(u,v)+}$, threshold list $w^{(u,v)+}$
\If{$w < \underline{w}^{(u,v)} \le \overline{w}^{(u,v)} < w'$}
\State $\mathcal{R}^{(u,v)+} \leftarrow \left\{\mbox{path from $u$ to $v$ in the tree } (\mathcal{N},\mathcal{A}_{u})\right\}$
\State $w^{(u,v)+}_{1} \leftarrow \overline{w}^{(u,v)}$
\State $R^{(u,v)+} \leftarrow 1$
\ElsIf{$w < \underline{w}^{(u,v)} \le w' \le \overline{w}^{(u,v)}$}
\State $\mathcal{R}^{(u,v)+} \leftarrow \left\{\mbox{path from $u$ to $v$ in the tree } (\mathcal{N},\mathcal{A}_{u})\right\}$
\State $w^{(u,v)+}_{1} \leftarrow w'$
\State $\mathcal{R}^{(u,v)+} \leftarrow \mathcal{R}^{(u,v)+} \cup \left\{\mbox{path from $u$ to $v$ in the tree } (\mathcal{N},\mathcal{A}'_{u})\right\}$
\State $R^{(u,v)+} \leftarrow 2$
\ElsIf{$\underline{w}^{(u,v)} \le w \le \overline{w}^{(u,v)} < w'$}
\State $w^{(u,v)+}_{R^{(u,v)+}} \leftarrow \overline{w}^{(u,v)}$ 
\ElsIf{$\underline{w}^{(u,v)} \le w \le w' \le \overline{w}^{(u,v)}$}
\State $w^{(u,v)+}_{R^{(u,v)+}} \leftarrow w'$ 
\State $\mathcal{R}^{(u,v)+} \leftarrow \mathcal{R}^{(u,v)+} \cup \left\{\mbox{path from $u$ to $v$ in the tree } (\mathcal{N},\mathcal{A}'_{u})\right\}$
\State $R^{(u,v)+} \leftarrow R^{(u,v)+} + 1$
\EndIf
\State Output path list $\mathcal{R}^{(u,v)+}$, threshold list $w^{(u,v)+}$
\end{algorithmic}
\end{algorithm}

\begin{remark}
\label{rmk:shortest path intervals not complete}
Algorithm~\ref{algo:BPSP0} returns both shortest paths $\mathcal{R}^{(u,v)+}$ as well as thresholds $w^{(u,v)+}_{0},\ldots,w^{(u,v)+}_{R^{(u,v)+}}$.
The thresholds are not used in Algorithm~BTA, but facilitate the following point.
Lemma~\ref{lem:optimal thresholds} established that, given any optimal solution $q$ for problem~\eqref{eqn:q reformulated BiTAP objective}, the thresholds produced by  Algorithm~\ref{algo:threshold} satisfy~\eqref{eqn:w ordered}, \eqref{eqn:each used path optimal in its interval}, and~\eqref{eqn:probabilities and thresholds consistent}.
In contrast, as shown next, the thresholds produced by Algorithm~\ref{algo:BPSP0} satisfy~\eqref{eqn:w ordered} and~\eqref{eqn:each used path optimal in its interval}, but not~\eqref{eqn:probabilities and thresholds consistent}.
Consider the data of Example~\ref{ex:equilibrium conditions example}, and suppose that each path~$r_{i}$ is associated with a separate arc also denoted with~$r_{i}$ and that the input data with $t_{(i,j)} = t_{(i,j)}(x_{(i,j)}(q))$ is provided to Algorithm~\ref{algo:BPSP0}.
Then $w^{\pi}_{0} = \underline{w} = 0$, and path~$r_{1}$ is the unique shortest path with arc costs $c_{(i,j)} = t_{(i,j)} + \underline{w} s_{(i,j)}$.
Then $w' = 0.01$, $\mathcal{A}' = \{r_{2},r_{3}\}$, and $w^{\pi}_{1} = w' = 0.01$.
Algorithm~\ref{algo:BPSP0} does not specify which arc in $\mathcal{A}'$ is selected to construct the next tree.
Whichever arc in $\mathcal{A}'$ is selected, in the next iteration $w' = 0.07$, $\mathcal{A}' = \{r_{6},r_{9}\}$, and $w^{\pi}_{2} = w' = 0.07$.
Recall from Example~\ref{ex:threshold-selection-equilibrium-example} that $w^{\pi}_{2} = (G^{\pi})^{-1}(L^{\pi}_{r_{2}}(q^{\pi})) = 0.03$, which cannot be obtained by solving the parametric shortest path problem.
The point is that although a parametric shortest path algorithm computes a sequence of thresholds, when there are ties on the efficient frontier, these thresholds may not satisfy~\eqref{eqn:probabilities and thresholds consistent}, because the distribution $G^{\pi}$ is immaterial to the parametric shortest path problem.
Algorithm~\ref{algo:threshold} is designed to produce thresholds that satisfy both the quantile condition~\eqref{eqn:probabilities and thresholds consistent}, as well as, similar to a parametric shortest path algorithm, the minimal cost condition~\eqref{eqn:each used path optimal in its interval}.
\end{remark}

\begin{lemma}
\label{lem:optimality conditions}
Assume that $t_{a} : \mathbb{R}_{+} \mapsto \mathbb{R}_{+}$ is continuous for all $a \in \mathcal{A}$.
For each $\pi \in \Pi$, let $\mathcal{R}^{\pi}$ denote the complete path set and let $\hat{\mathcal{R}}^{\pi} \subset \mathcal{R}^{\pi}$.
Consider any optimal solution $\hat{q}$ for problem~\eqref{eqn:q reformulated BiTAP objective} with path sets $\hat{\mathcal{R}}^{\pi}$, $\pi \in \Pi$.
Let $\hat{w} = \left(\hat{w}^{\pi}_{r}, \pi \in \Pi, r \in \{0\} \cup \hat{\mathcal{R}}^{\pi}\right)$ and $\hat{\mathcal{R}}^{\pi}_{\min} = \{\hat{r}_{1},\ldots,\hat{r}_{\hat{L}}\}$ denote the output of Algorithm~\ref{algo:threshold} with input $\hat{\mathcal{R}}^{\pi}$ and $\hat{q}$.

Let $q^* = \left(q^{\pi*}_{r}, \pi \in \Pi, r \in \mathcal{R}^{\pi}\right)$ be as follows:
For each $r \in \hat{\mathcal{R}}^{\pi}$ let $q^{\pi*}_{r} = \hat{q}^{\pi}_{r}$, and for each $r \in \mathcal{R}^{\pi} \setminus \hat{\mathcal{R}}^{\pi}$ let $q^{\pi*}_{r} = 0$.
Then $x(q^*) = x(\hat{q})$.

Let $w^* = \left(w^{\pi*}_{r}, \pi \in \Pi, r \in \{0\} \cup \mathcal{R}^{\pi}\right)$ and $\mathcal{R}^{\pi*}_{\min} = \{r^*_{1},\ldots,r^*_{L^*}\}$ denote the output of Algorithm~\ref{algo:threshold} with input $\mathcal{R}^{\pi}$ and $q^*$.
Let $\mathcal{R}^{\pi}$ be indexed with $r = 1,\ldots,R^{\pi}$, and let
\begin{equation}
\lambda^{\pi*} \ \ \defi \ \ -D^{\pi} \left[T^{\pi}_{r^*_{L^*}}(x(q^*)) + \overline{w}^{\pi} \left(S^{\pi}_{r^*_{L^*}} - S^{\pi}_{R^{\pi}}\right)\right]
\qquad \forall \ \pi \in \Pi
\label{eqn:choice of dual 2}
\end{equation}
Let $\lambda^{*} \defi (\lambda^{\pi*}, \pi \in \Pi)$.
For each $\pi \in \Pi$, let $\mathcal{R}^{\pi+}$ denote an optimal set of paths for the parametric shortest path problem (over the complete path set $\mathcal{R}^{\pi}$) with arc costs $t_{a}(x_{a}(\hat{q}))$ and $s_{a}$.
That is, $\mathcal{R}^{\pi+}$ can be indexed such that $\mathcal{R}^{\pi+} = \{r^{1},\ldots,r^{k}\}$ and there are thresholds $\underline{w}^{\pi} = w^{\pi}_{r^{0}} \le w^{\pi}_{r^{1}} \le \cdots \le w^{\pi}_{r^{k-1}} \le w^{\pi}_{r^{k}} = \overline{w}^{\pi}$ such that for all $i=1,\ldots,k$ it holds that $C^{\pi}_{r^{i}}(w,x(\hat{q})) \le C^{\pi}_{r'}(w,x(\hat{q}))$ for all $r' \in \mathcal{R}^{\pi}$ and $w \in [w^{\pi}_{r^{i-1}},w^{\pi}_{r^{i}}]$.
Suppose that $\mathcal{R}^{\pi+} \subset \hat{\mathcal{R}}^{\pi}$ for each $\pi \in \Pi$.
Then it holds that $(q^*,\lambda^*,w^*)$ satisfies the optimality conditions for the complete path set:
\begin{align}
& q^{\pi*}_{r} \left\{D^{\pi} \left[T^{\pi}_{r}(x(q^*)) + \sum_{r' = r}^{R^{\pi}-1} \left(S^{\pi}_{r'} - S^{\pi}_{r'+1}\right) w^{\pi*}_{r'}\right] + \lambda^{\pi*}\right\} \ \ = \ \ 0 \qquad \forall \ \pi \in \Pi, \ \forall \ r \in \mathcal{R}^{\pi}
\label{eqn:q complete BiTAP optimality1} \\
& D^{\pi} \left[T^{\pi}_{r}(x(q^*)) + \sum_{r' = r}^{R^{\pi}-1} \left(S^{\pi}_{r'} - S^{\pi}_{r'+1}\right) w^{\pi*}_{r'}\right] + \lambda^{\pi*} \ \ \ge \ \ 0 \qquad \forall \ \pi \in \Pi, \ \forall \ r \in \mathcal{R}^{\pi}
\label{eqn:q complete BiTAP optimality3}
\end{align}
\end{lemma}

\begin{theorem}
\label{thm:BTA algorithm}
Algorithm~BTA returns a solution~$\hat{q} = (\hat{q}^{\pi}_{r}, \pi \in \Pi, r \in \hat{\mathcal{R}}^{\pi})$ in $O\left(\sum_{\pi \in \Pi} \big|\hat{\mathcal{R}}^{\pi}\big|\right)$ steps.
Then $q^* = (q^{\pi*}_{r}, \pi \in \Pi, r \in \mathcal{R}^{\pi})$, where $q^{\pi*}_{r} = \hat{q}^{\pi}_{r}$ for all $\pi \in \Pi$, $r \in \hat{\mathcal{R}}^{\pi}$ and $q^{\pi*}_{r} = 0$ for all $\pi \in \Pi$, $r \in \mathcal{R}^{\pi} \setminus \hat{\mathcal{R}}^{\pi}$, is optimal for problem~\eqref{eqn:q reformulated BiTAP objective} with the complete path set.
Given $q^*$ and the resulting output $w^*$ of Algorithm~\ref{algo:threshold}, $(q^*,w^*)$ is an equilibrium of the BiTAP.
\end{theorem}

\begin{figure}[tbp]
  \centering
  \begin{subfigure}[t]{0.48\textwidth}
    \centering
    \begin{tikzpicture}
      \begin{axis}[
        exampleaxis,
        width=0.98\linewidth,
        height=5.4cm,
        axis lines=left,
        axis line style={-{Stealth[length=1.7mm]}, draw=black, line width=0.35pt},
        xmin=0, xmax=0.142,
        ymin=0.02, ymax=0.46,
        xtick={0,0.01,0.07,0.08,0.12,0.14},
        xticklabels={$0$,$.01$,$.07$,$.08$,$.12$,$.14$},
        ytick={0.1,0.2,0.3,0.4},
        xlabel={$w$},
        ylabel={$C^{\pi}_{r}(w,x)$},
        tick label style={font=\tiny},
        clip=true
      ]
        \addplot[domain=0:0.14,samples=2,color=gray!40,line width=0.4pt] {0.03+5*x};
        \addplot[domain=0:0.14,samples=2,color=gray!40,line width=0.4pt] {0.05+4*x};
        \addplot[domain=0:0.14,samples=2,color=gray!40,line width=0.4pt] {0.04+4*x};
        \addplot[domain=0:0.14,samples=2,color=gray!40,line width=0.4pt] {0.135+3.7*x};
        \addplot[domain=0:0.14,samples=2,color=gray!40,line width=0.4pt] {0.16+3.35*x};
        \addplot[domain=0:0.14,samples=2,color=gray!40,line width=0.4pt] {0.11+3*x};
        \addplot[domain=0:0.14,samples=2,color=gray!40,line width=0.4pt] {0.18+2.6*x};
        \addplot[domain=0:0.14,samples=2,color=gray!40,line width=0.4pt] {0.235+2.25*x};
        \addplot[domain=0:0.14,samples=2,color=gray!40,line width=0.4pt] {0.19+2*x};
        \addplot[domain=0:0.14,samples=2,color=gray!40,line width=0.4pt] {0.31+x};
        \addplot[very thick,color=ExampleFrontier] coordinates {(0,0.03) (0.01,0.08)};
        \addplot[very thick,color=ExampleWarm] coordinates {(0.01,0.08) (0.07,0.32)};
        \addplot[very thick,color=violet!75!black] coordinates {(0.07,0.32) (0.08,0.35)};
        \addplot[very thick,color=ExampleGuide] coordinates {(0.08,0.35) (0.12,0.43)};
        \addplot[very thick,color=black] coordinates {(0.12,0.43) (0.14,0.45)};
        \addplot[thresholdline] coordinates {(0.01,0.02) (0.01,0.08)};
        \addplot[thresholdline] coordinates {(0.07,0.02) (0.07,0.32)};
        \addplot[thresholdline] coordinates {(0.08,0.02) (0.08,0.35)};
        \addplot[thresholdline] coordinates {(0.12,0.02) (0.12,0.43)};
        \node[font=\tiny,color=ExampleFrontier,fill=white,inner sep=0.7pt] at (axis cs:0.004,0.06) {$r_1$};
        \node[font=\tiny,color=ExampleWarm,fill=white,inner sep=0.7pt] at (axis cs:0.038,0.18) {$r_3$};
        \node[font=\tiny,color=violet!75!black,fill=white,inner sep=0.7pt] at (axis cs:0.075,0.335) {$r_6$};
        \node[font=\tiny,color=ExampleGuide,fill=white,inner sep=0.7pt] at (axis cs:0.097,0.385) {$r_9$};
        \node[font=\tiny,fill=white,inner sep=0.7pt] at (axis cs:0.132,0.448) {$r_{10}$};
      \end{axis}
    \end{tikzpicture}
    \subcaption{Parametric shortest path solution with zero initial flow.}
  \end{subfigure}
  \hfill
  \begin{subfigure}[t]{0.48\textwidth}
    \centering
    \begin{tikzpicture}
      \begin{axis}[
        exampleaxis,
        width=0.98\linewidth,
        height=5.4cm,
        axis lines=left,
        axis line style={-{Stealth[length=1.7mm]}, draw=black, line width=0.35pt},
        xmin=0, xmax=0.142,
        ymin=0.02, ymax=0.46,
        xtick={0,0.0055556,0.07,0.12,0.14},
        xticklabels={$0$,$1/180$,$.07$,$.12$,$.14$},
        ytick={0.1,0.2,0.3,0.4},
        xlabel={$w$},
        tick label style={font=\tiny},
        clip=true
      ]
        \addplot[domain=0:0.14,samples=2,color=gray!40,line width=0.4pt] {2/45+5*x};
        \addplot[domain=0:0.14,samples=2,color=gray!40,line width=0.4pt] {0.05+4*x};
        \addplot[domain=0:0.14,samples=2,color=gray!40,line width=0.4pt] {49/900+4*x};
        \addplot[domain=0:0.14,samples=2,color=gray!40,line width=0.4pt] {0.135+3.7*x};
        \addplot[domain=0:0.14,samples=2,color=gray!40,line width=0.4pt] {0.16+3.35*x};
        \addplot[domain=0:0.14,samples=2,color=gray!40,line width=0.4pt] {0.12+3*x};
        \addplot[domain=0:0.14,samples=2,color=gray!40,line width=0.4pt] {0.18+2.6*x};
        \addplot[domain=0:0.14,samples=2,color=gray!40,line width=0.4pt] {0.235+2.25*x};
        \addplot[domain=0:0.14,samples=2,color=gray!40,line width=0.4pt] {0.19+2*x};
        \addplot[domain=0:0.14,samples=2,color=gray!40,line width=0.4pt] {0.31+x};
        \addplot[very thick,color=ExampleFrontier] coordinates {(0,0.0444444) (0.0055556,0.0722222)};
        \addplot[very thick,color=ExamplePruned] coordinates {(0.0055556,0.0722222) (0.07,0.33)};
        \addplot[only marks,mark=*,mark size=2.3pt,color=violet!75!black] coordinates {(0.07,0.33)};
        \addplot[very thick,color=ExampleGuide] coordinates {(0.07,0.33) (0.12,0.43)};
        \addplot[very thick,color=black] coordinates {(0.12,0.43) (0.14,0.45)};
        \addplot[thresholdline] coordinates {(0.0055556,0.02) (0.0055556,0.0722222)};
        \addplot[thresholdline] coordinates {(0.07,0.02) (0.07,0.33)};
        \addplot[thresholdline] coordinates {(0.12,0.02) (0.12,0.43)};
        \node[font=\tiny,color=ExampleFrontier,fill=white,inner sep=0.7pt] at (axis cs:0.002,0.058) {$r_1$};
        \node[font=\tiny,color=ExamplePruned,fill=white,inner sep=0.7pt] at (axis cs:0.034,0.19) {new $r_2$};
        \node[font=\tiny,color=violet!75!black,anchor=south west,fill=white,inner sep=0.7pt] at (axis cs:0.071,0.331) {$r_6$};
        \node[font=\tiny,color=ExampleGuide,fill=white,inner sep=0.7pt] at (axis cs:0.094,0.38) {$r_9$};
        \node[font=\tiny,fill=white,inner sep=0.7pt] at (axis cs:0.132,0.448) {$r_{10}$};
      \end{axis}
    \end{tikzpicture}
    \subcaption{Parametric shortest path solution after the first flow update.}
  \end{subfigure}

  \medskip

  \begin{subfigure}[t]{0.54\textwidth}
    \centering
    \begin{tikzpicture}
      \begin{axis}[
        exampleaxis,
        width=0.98\linewidth,
        height=5.4cm,
        axis lines=left,
        axis line style={-{Stealth[length=1.7mm]}, draw=black, line width=0.35pt},
        xmin=0, xmax=0.142,
        ymin=0.02, ymax=0.46,
        xtick={0,0.01,0.03,0.07,0.12,0.14},
        xticklabels={$0$,$.01$,$.03$,$.07$,$.12$,$.14$},
        ytick={0.1,0.2,0.3,0.4},
        xlabel={$w$},
        tick label style={font=\tiny},
        clip=true
      ]
        \addplot[domain=0:0.14,samples=2,color=gray!40,line width=0.4pt] {0.04+5*x};
        \addplot[domain=0:0.14,samples=2,color=gray!40,line width=0.4pt] {0.05+4*x};
        \addplot[domain=0:0.14,samples=2,color=gray!40,line width=0.4pt] {0.05+4*x};
        \addplot[domain=0:0.14,samples=2,color=gray!40,line width=0.4pt] {0.135+3.7*x};
        \addplot[domain=0:0.14,samples=2,color=gray!40,line width=0.4pt] {0.16+3.35*x};
        \addplot[domain=0:0.14,samples=2,color=gray!40,line width=0.4pt] {0.12+3*x};
        \addplot[domain=0:0.14,samples=2,color=gray!40,line width=0.4pt] {0.18+2.6*x};
        \addplot[domain=0:0.14,samples=2,color=gray!40,line width=0.4pt] {0.235+2.25*x};
        \addplot[domain=0:0.14,samples=2,color=gray!40,line width=0.4pt] {0.19+2*x};
        \addplot[domain=0:0.14,samples=2,color=gray!40,line width=0.4pt] {0.31+x};
        \addplot[very thick,color=ExampleFrontier] coordinates {(0,0.04) (0.01,0.09)};
        \addplot[very thick,color=ExampleWarm] coordinates {(0.01,0.09) (0.07,0.33)};
        \addplot[only marks,mark=*,mark size=2.3pt,color=violet!75!black] coordinates {(0.07,0.33)};
        \addplot[very thick,color=ExampleGuide] coordinates {(0.07,0.33) (0.12,0.43)};
        \addplot[very thick,color=black] coordinates {(0.12,0.43) (0.14,0.45)};
        \addplot[thresholdline] coordinates {(0.01,0.02) (0.01,0.09)};
        \addplot[thresholdline] coordinates {(0.07,0.02) (0.07,0.33)};
        \addplot[thresholdline] coordinates {(0.12,0.02) (0.12,0.43)};
        \addplot[color=ExampleWarm,densely dotted,line width=0.9pt] coordinates {(0.03,0.02) (0.03,0.17)};
        \node[font=\tiny,color=ExampleFrontier,fill=white,inner sep=0.7pt] at (axis cs:0.004,0.066) {$r_1$};
        \node[font=\tiny,color=ExampleWarm,fill=white,inner sep=0.7pt] at (axis cs:0.047,0.23) {$r_2,r_3$};
        \node[font=\tiny,color=violet!75!black,anchor=south west,fill=white,inner sep=0.7pt] at (axis cs:0.071,0.331) {$r_6$};
        \node[font=\tiny,color=ExampleGuide,fill=white,inner sep=0.7pt] at (axis cs:0.094,0.38) {$r_9$};
        \node[font=\tiny,fill=white,inner sep=0.7pt] at (axis cs:0.132,0.448) {$r_{10}$};
        \node[font=\tiny,color=ExampleWarm,anchor=west,fill=white,inner sep=0.7pt] at (axis cs:0.031,0.145) {$0.03$: Threshold only};
      \end{axis}
    \end{tikzpicture}
    \subcaption{Parametric shortest path solution for flow $x(q^*)$.}
  \end{subfigure}
  \caption{The three parametric shortest path solutions produced by Algorithm~\ref{algo:BPSP0} in Example~\ref{ex:bta-psp-equilibrium-example}.
  Thin gray lines are the bicriteria costs of all ten paths and thick colored pieces their lower envelope.
  Panel~(a) shows the first path set $\hat{\mathcal{R}}^{\pi}$.
  Panel~(b) shows the parametric shortest path solution for the updated flow $x(\hat{q})$; path~$r_{6}$ is optimal only at the marked singleton $w = 0.07$.
  Panel~(c) shows the final parametric shortest path solution for the final flow $x(q^*)$, when every parametric shortest path is already in $\hat{\mathcal{R}}^{\pi}$.
}
  \label{fig:bta-psp-equilibrium-example}
\end{figure}
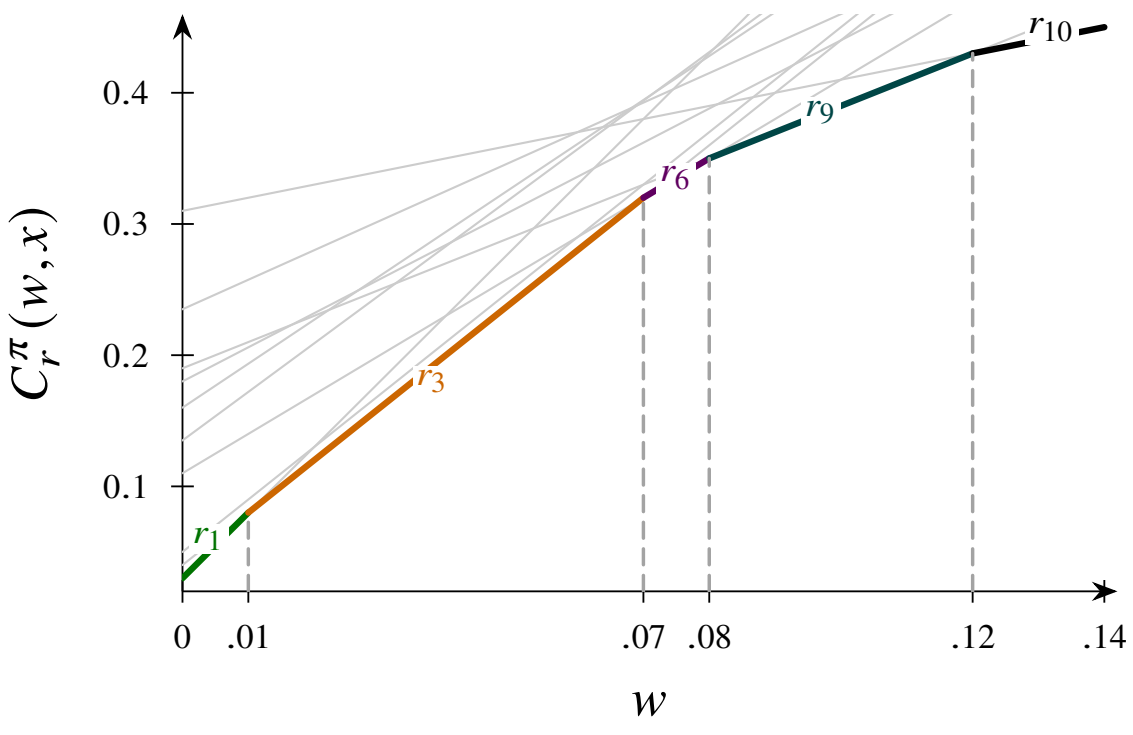
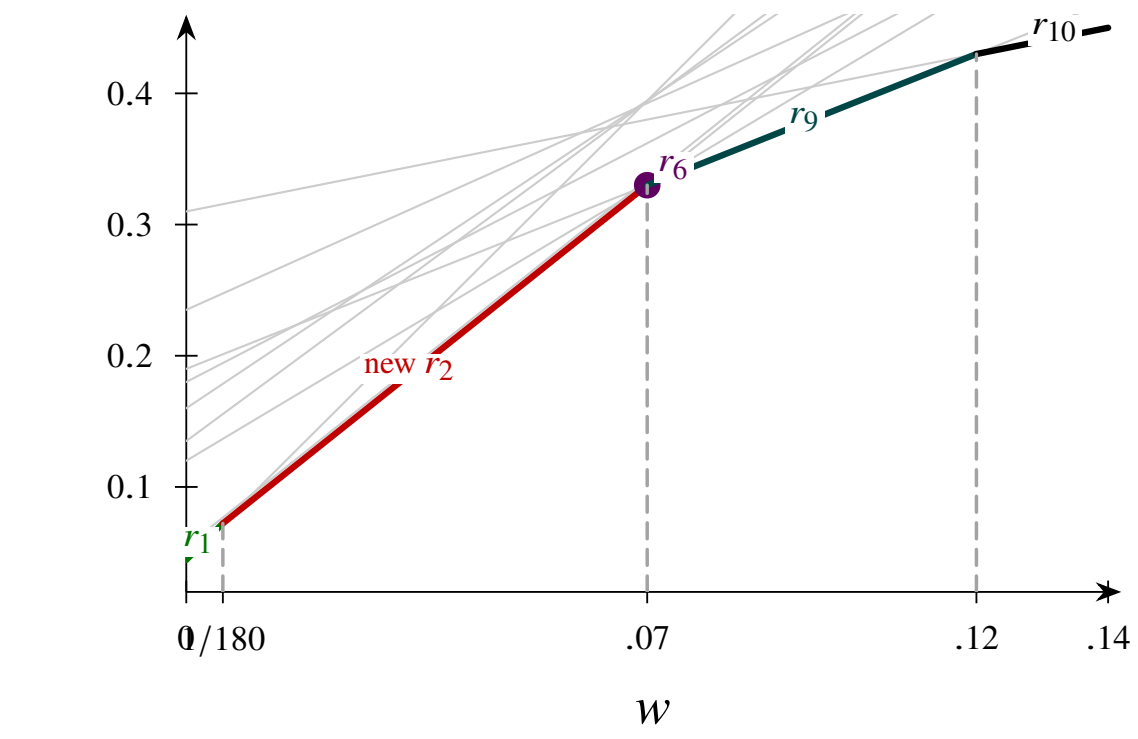
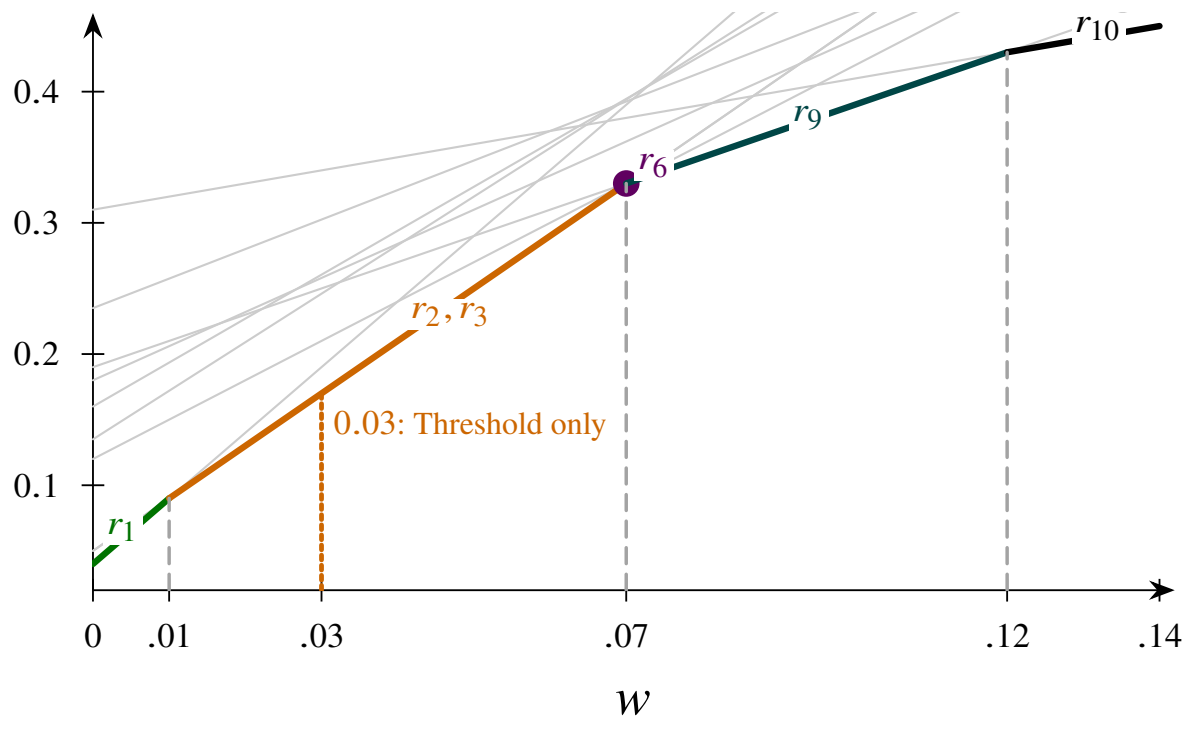

\begin{example}
\label{ex:bta-psp-equilibrium-example}
We present a network for which $(q,w)$ given in Examples~\ref{ex:efficient frontier example} and~\ref{ex:equilibrium conditions example} is an equilibrium.
We also present the steps of Algorithm~BTA.
Let
\[
\mathcal{N} \ \ = \ \ \{u,v,z_{1},\ldots,z_{10}\},
\qquad
\mathcal{A} \ \ = \ \ \{(u,z_i),(z_i,v) \; : \; i=1,\ldots,10\},
\]
let $\Pi = \{\pi\}$ with $\pi = (u,v)$, let $D^{\pi} = 1$, and let path $r_{i} = ((u,z_{i}),(z_{i},v))$.
For every $i$, let $t_{(u,z_{i})}(y) = 0$, $s_{(u,z_{i})} = 0$, and $s_{(z_{i},v)} = S^{\pi}_{r_{i}}$, with $S^{\pi}$ as in Example~\ref{ex:efficient frontier example}, and let
\begin{equation}
\label{eqn:bta-psp-example-travel-times}
\begin{array}{rclcrcl}
t_{(z_{1},v)}(y) & = & 0.03 + 0.1y,
& & t_{(z_{2},v)}(y) & = & 0.05,\\
t_{(z_{3},v)}(y) & = & 0.04 + y/30,
& & t_{(z_{4},v)}(y) & = & 0.135,\\
t_{(z_{5},v)}(y) & = & 0.16,
& & t_{(z_{6},v)}(y) & = & 0.11 + 0.1y,\\
t_{(z_{7},v)}(y) & = & 0.18,
& & t_{(z_{8},v)}(y) & = & 0.235,\\
t_{(z_{9},v)}(y) & = & 0.19,
& & t_{(z_{10},v)}(y) & = & 0.31.
\end{array}
\end{equation}
The paths are arc-disjoint, so~\eqref{eqn:q reformulated BiTAP arc flows} and~\eqref{eqn:bta-psp-example-travel-times} give
\begin{equation}
\label{eqn:bta-psp-example-final-q}
q^* \ \ = \ \ (.1,.3,.3,0,0,.1,0,0,0,.2)
\quad\Rightarrow\quad
T^{\pi}(x(q^*)) \ \ = \ \ (.04,.05,.05,.135,.16,.12,.18,.235,.19,.31),
\end{equation}
as required.
Starting with arc travel times $t_{a}(0)$, Algorithm~\ref{algo:BPSP0} gives initial path set $\hat{\mathcal{R}}^{\pi} = \mathcal{R}^{\pi+} = \{r_{1},r_{3},r_{6},r_{9},r_{10}\}$.
Path~$r_{2}$ is omitted because $S^{\pi}_{r_{2}} = S^{\pi}_{r_{3}}$ and $t_{(z_{2},v)}(0) = .05 > .04 = t_{(z_{3},v)}(0)$; Figure~\ref{fig:bta-psp-equilibrium-example}(a) shows that the other omitted paths also lie above the efficient frontier.

The optimal solution $\hat{q}$ for problem~\eqref{eqn:q reformulated BiTAP objective} with initial path set $\hat{\mathcal{R}}^{\pi}$ is $(\hat{q}^{\pi}_{r_{1}},\hat{q}^{\pi}_{r_{3}},\hat{q}^{\pi}_{r_{6}},\hat{q}^{\pi}_{r_{9}},\hat{q}^{\pi}_{r_{10}}) = \left(\frac{13}{90},\frac{13}{30},\frac1{10},\frac{11}{90},\frac15\right)$, with thresholds $(\hat{w}^{\pi}_{r_{1}},\hat{w}^{\pi}_{r_{3}},\hat{w}^{\pi}_{r_{6}},\hat{w}^{\pi}_{r_{9}},\hat{w}^{\pi}_{r_{10}}) = \left(\frac{1}{100},\frac{59}{900},\frac{7}{100},\frac{3}{25},\frac{7}{50}\right)$.
Optimality conditions~\eqref{eqn:q continuous BiTAP optimality1}--\eqref{eqn:q continuous BiTAP optimality3} hold with $\lambda^{\pi} = -31/100$.

At this solution, $T^{\pi}_{r_{2}}(x(\hat{q})) = 1/20 < 49/900 = T^{\pi}_{r_{3}}(x(\hat{q}))$.
Algorithm~\ref{algo:BPSP0} with arc travel times $t_{a}(x_{a}(\hat{q}))$ gives path set $\mathcal{R}^{\pi+} = \{r_{1},r_{2},r_{6},r_{9},r_{10}\}$, and thus $\hat{\mathcal{R}}^{\pi} \leftarrow \hat{\mathcal{R}}^{\pi} \cup \{r_2\}$.
The optimal solution for problem~\eqref{eqn:q reformulated BiTAP objective} with updated path set $\hat{\mathcal{R}}^{\pi}$ gives $q^*$ in~\eqref{eqn:bta-psp-example-final-q}.
Algorithm~\ref{algo:BPSP0} with arc travel times $t_{a}(x_{a}(q^*))$ gives path set $\mathcal{R}^{\pi+} = \{r_{1},r_{2},r_{6},r_{9},r_{10}\} \subset \hat{\mathcal{R}}^{\pi}$.
If Algorithm~\ref{algo:BPSP0} selects $r_{3}$ instead of $r_{2}$ at $w = 0.01$, or skips path $r_{6}$ at $w = 0.07$, it still holds that $\mathcal{R}^{\pi+} \subset \hat{\mathcal{R}}^{\pi}$.
Algorithm~BTA terminates, and Lemma~\ref{lem:optimality conditions} verifies optimality for the full path set.
Algorithm~\ref{algo:threshold} then gives $w^*$ in~\eqref{eqn:threshold-selection-example-output}.
\hfill\Halmos
\end{example}

\section{Numerical Results}
\label{sec:numerical_experiments}

In this section we compare Algorithm~BTA with the methods of \citet{dafe:72}, \citet{daga:83}, \citet{dial:96}, \citet{marc:96b}, and \citet{xie:26}.
Specifically, from \citet{dafe:72} we use the path-correction algorithm for the multiclass problem, in its TNTP surrogate cost mode.
From \citet{daga:83} we use the stochastic approximation algorithm, with the feasible Ces\`aro average.
From \citet{dial:96} we use T2-ETA, which calls T2-MPA as its minimum-path assignment subproblem.
From \citet{marc:96b} we use the parametric shortest path descent with the paper's exact line search (not its Armijo variants).
From \citet{xie:26} we use SBA rather than the Frank-Wolfe or T2-$\theta$ algorithms that accompany it.

Table~\ref{tab:exp_networks} gives summary network data such as the numbers of nodes, arcs, zones, OD pairs, and total demand, for each benchmark network used in the numerical tests.
Link travel times are given by the BPR function $t_{a}(x_{a}) = t^{0}_{a} \bigl(1 + \alpha_{a} (x_{a} / c_{a})^{\beta_{a}}\bigr)$, with $\alpha_{a}$ and $\beta_{a}$ given in Table~\ref{tab:exp_master_params}.
The arc flow monetary costs are given by $s_{a} = \alpha_{0} x_{a}^{\mathrm{SO}} t'_{a}(x_{a}^{\mathrm{SO}})$, called {\em marginal-cost tolls}, rounded to the nearest $0.01$ for every $a \in \mathcal A$, with $\alpha_{0} = 1$ (except for the results discussed in Section~\ref{sec:exp_toll_impact}) and $x^{\mathrm{SO}}_{a}$ being the system optimal arc flows for the same instance.
Tests were performed with two types of instances: smooth instances, meaning that $\Upsilon$ is smooth (a sufficient condition is that $G^{\pi}$ has a positive density on $[\underline{w}^{\pi},\overline{w}^{\pi}]$, as the methods of \citealt{marc:96b} and \citealt{xie:26} require), and nonsmooth instances.
For the smooth instances, all methods used $w \sim \mathcal{U}[6,30]$, and for the nonsmooth instances $w$ followed a ``staircase'' distribution with equal point masses at $w = 0, 2.5, 5, 7.5, 10$.

All methods, except the SBA method, are implemented in C11 with POSIX threads, and compiled via Clang.
The SBA code, obtained from the authors of \citet{xie:26}, was written in C++14, built from GCC 12.3.0.
The reported computation times exclude the computation of $x^{\mathrm{SO}}$ to determine input arc flow monetary costs $s_{a}$.
We ran BTA with two solvers for solving problem~\eqref{eqn:q reformulated BiTAP objective} with subsets of paths~$\hat{\mathcal{R}}^{\pi}$ in Line~10 of Algorithm~\ref{algo:BTA}:
The FISTA method of \citet{sujanani2024efficient} was used for the smooth instances, and the two-cut adaptive proximal bundle method (Ad-GPB) of \citet{monteiro2024parameter} was used for the nonsmooth instances.
Table~\ref{tab:exp_master_params} gives the parameters used for these solvers.

\subsection{Benchmarking Metrics}
\label{sec:exp_metrics}

For each instance~$i$ and method~$m$, let $q_{i,m}(t')$ denote the incumbent solution at CPU time~$t'$, and let
\begin{equation}
V_{i,m}(t) \ \ \defi \ \ \min_{t' \le t} \ \Upsilon_{i}\!\left(q_{i,m}(t')\right)
\label{eqn:exp_best_value}
\end{equation}
denote the best objective value attained by time~$t$.
For each instance~$i$ we also compute a lower bound on the optimal objective value $\Upsilon_{i}^{\ast}$, as follows.
It follows from convexity of $\Upsilon_{i}$ that
\begin{equation}
\Upsilon_{i}^{\ast} \ \ \ge \ \ \Upsilon_{i}(q) \, + \, \min\left\{\inner{g}{y - q} \ : \ y \in \prod_{\pi \in \Pi} \Delta_{\mathcal{R}^{\pi}}\right\}
\label{eqn:exp_lower_bound}
\end{equation}
for any feasible~$q$ and any subgradient $g \in \partial \Upsilon_{i}(q)$ at~$q$.
For instance~$i$, let $L_{i}$ denote the largest lower bound of~\eqref{eqn:exp_lower_bound} over all iterates~$q$ produced by all methods~$m$.
Thus, $L_{i}$ is the best lower bound on $\Upsilon_{i}^{\ast}$ computed during the test.
Let $q_{i}^{0}$ denote an initial feasible solution for instance~$i$ obtained after Line~4 of Algorithm~\ref{algo:BTA}.
Then the relative gap for any method~$m$ on instance~$i$ at time~$t$ is
\begin{equation}
\gamma_{i,m}(t) \ \ \defi \ \ \frac{V_{i,m}(t) - L_{i}}{\Upsilon_{i}(q_{i}^{0}) - L_{i}}.
\label{eqn:exp_relative_excess}
\end{equation}
Thus, $\gamma_{i,\mathrm{BTA}}(0) = 1$.
Note that the denominator does not depend on the method~$m$, and is used only for scaling.
For any tolerance $\varepsilon_{i} > 0$ and maximum compute time~$T_{i}$ for instance~$i$, let
\[
t_{i,m} \ \ \defi \ \ \begin{cases}
\inf\{t \in [0,T_{i}] \, : \, \gamma_{i,m}(t) \le \varepsilon_{i}\} & \mbox{ if } \{t \in [0,T_{i}] \, : \, \gamma_{i,m}(t) \le \varepsilon_{i}\} \neq \varnothing \\
T_{i} & \mbox{ otherwise}
\end{cases}
\]
denote the time at which method~$m$ first reaches the tolerance for instance~$i$, capped at~$T_{i}$.
We used $\varepsilon_{i} = 10^{-5}$ for smooth instances and $\varepsilon_{i} = 10^{-3}$ for nonsmooth instances.
We used $T_{i} = 3600$s for all instances, except for the smooth Goldcoast, Chicago Regional, and Philadelphia instances, for which $T_{i} = 21600$s, and the nonsmooth Goldcoast, Chicago Regional, and Philadelphia instances (each with more than a million OD records), for which $T_{i} = 86400$s.
As a time-dependent aggregate performance metric over set~$\mathcal{I}$ of instances, we compute the fraction of instances for which method~$m$ finds a sufficiently good solution by time~$t$, given by
\begin{equation*}
P_{m}(t) \ \ \defi \ \ \frac{1}{|\mathcal{I}|} \sum_{i \in \mathcal{I}} \mathds{1}\!\left\{\gamma_{i,m}(\min\{t,T_{i}\}) \le \varepsilon_{i}\right\}.
\label{eqn:exp_performance_function}
\end{equation*}
Other aggregate performance metrics for each method~$m$ are
\begin{equation*}
\begin{aligned}
g_{m} & \ \ = \ \ \left[\prod_{i \in \mathcal{I}}(t_{i,m} + 1)\right]^{1/|\mathcal{I}|} - 1,
\qquad
r_{m} \ \ = \ \ \frac{g_{m}}{\min_{m'} g_{m'}}, \\
h_{m} & \ \ = \ \ \frac{1}{|\mathcal{I}|} \sum_{i \in \mathcal{I}} \mathds{1}\left\{\gamma_{i,m}(T_{i}) \le \varepsilon_{i}, \ t_{i,m} = \min_{m'} t_{i,m'}\right\},
\end{aligned}
\label{eqn:exp_sgm}
\end{equation*}
the shifted geometric mean capped completion time, its normalized value, and the fraction of all instances for which method~$m$ reaches the tolerance $\varepsilon_{i}$ by time~$T_{i}$ and is the fastest.

\subsection{Convergence Performance on Benchmark Networks}
\label{sec:exp_convergence_perf}

\begin{figure}[t]
\centering
\begin{tikzpicture}
\begin{axis}[
  width=0.66\textwidth,
  height=0.30\textwidth,
  xmode=log,
  xmin=0.00005,
  xmax=21600,
  ymin=-0.02,
  ymax=1.02,
  xtick={0.0001,0.01,1,100,21600},
  xticklabels={$10^{-4}$,$10^{-2}$,$1$,$10^2$,$21600$},
  ytick={0,0.2,0.4,0.6,0.8,1},
  xlabel={elapsed time $t$ (seconds)},
  ylabel={$P_{s}(t;10^{-5})$},
  ylabel style={font=\small},
  xlabel style={font=\small},
  tick label style={font=\scriptsize},
  label style={font=\small},
  axis line style={line width=0.4pt},
  tick style={line width=0.4pt},
  grid=both,
  major grid style={gray!28,line width=0.3pt},
  minor grid style={gray!12,line width=0.2pt},
  legend style={draw=none,font=\scriptsize,at={(0.02,0.98)},anchor=north west,legend columns=2},
  every axis plot/.append style={mark=none}
]
\addplot+[dag,const plot] coordinates {
  (0.000001,0) (0.472978514,0.0588235294118) (21600,0.0588235294118)
};
\addlegendentry{Daganzo}
\addplot+[dial,const plot] coordinates {
  (0.000001,0) (0.006782194,0.0588235294118) (0.131250667,0.117647058824) (2.134686944,0.176470588235) (2.857909346,0.235294117647)
  (4.010682268,0.294117647059) (13.744234805,0.352941176471) (26.851457809,0.411764705882) (111.513322345,0.470588235294) (610.800880473,0.529411764706)
  (1034.9819414,0.588235294118) (1750.34744454,0.647058823529) (2283.96784453,0.705882352941) (21600,0.705882352941)
};
\addlegendentry{Dial}
\addplot+[mz,const plot] coordinates {
  (0.000001,0) (0.00784215,0.0588235294118) (0.059687523,0.117647058824) (0.086407653,0.176470588235) (0.214067572,0.235294117647)
  (0.456762939,0.294117647059) (0.586413118,0.352941176471) (1.90857612,0.411764705882) (2.353759416,0.470588235294) (161.333954314,0.529411764706)
  (288.802883684,0.588235294118) (593.193098392,0.647058823529) (820.423669641,0.705882352941) (823.444751881,0.764705882353) (21600,0.764705882353)
};
\addlegendentry{Marcotte--Zhu}
\addplot+[sba,const plot] coordinates {
  (0.000001,0) (0.002136541,0.0588235294118) (0.122015333,0.117647058824) (0.143232021,0.176470588235) (0.258668088,0.235294117647)
  (0.470355531,0.294117647059) (1.155380708,0.352941176471) (2.782356946,0.411764705882) (6.4708289,0.470588235294) (38.276883269,0.529411764706)
  (88.691323119,0.588235294118) (112.735489188,0.647058823529) (496.222754574,0.705882352941) (2128.49549319,0.764705882353) (21600,0.764705882353)
};
\addlegendentry{SBA}
\addplot+[bta,const plot] coordinates {
  (0.000001,0) (0.000103703700006,0.0588235294118) (0.0315230400302,0.117647058824) (0.0402438919991,0.176470588235) (0.0445494980086,0.235294117647)
  (0.0703399860067,0.294117647059) (0.168326881016,0.352941176471) (0.321042817959,0.411764705882) (1.037966631,0.470588235294) (1.27543046896,0.529411764706)
  (4.67275369802,0.588235294118) (5.987209243,0.647058823529) (9.935192907,0.705882352941) (14.656834039,0.764705882353) (25.351896361,0.823529411765)
  (21600,0.823529411765)
};
\addlegendentry{BTA}
\end{axis}
\end{tikzpicture}
\caption{Plot of $P_{m}(t)$ for BTA and the methods of \citet{daga:83}, \citet{dial:96}, \citet{marc:96b}, and \citet{xie:26} over all $17$ smooth instances.
BTA reaches $0.82$, missing only Goldcoast, Chicago Regional, and Philadelphia, for which no method reaches the target tolerance; SBA reaches $0.76$, Marcotte-Zhu $0.76$, Dial $0.71$, and Daganzo $0.06$.}
\label{fig:exp_common_gap_profile}
\end{figure}
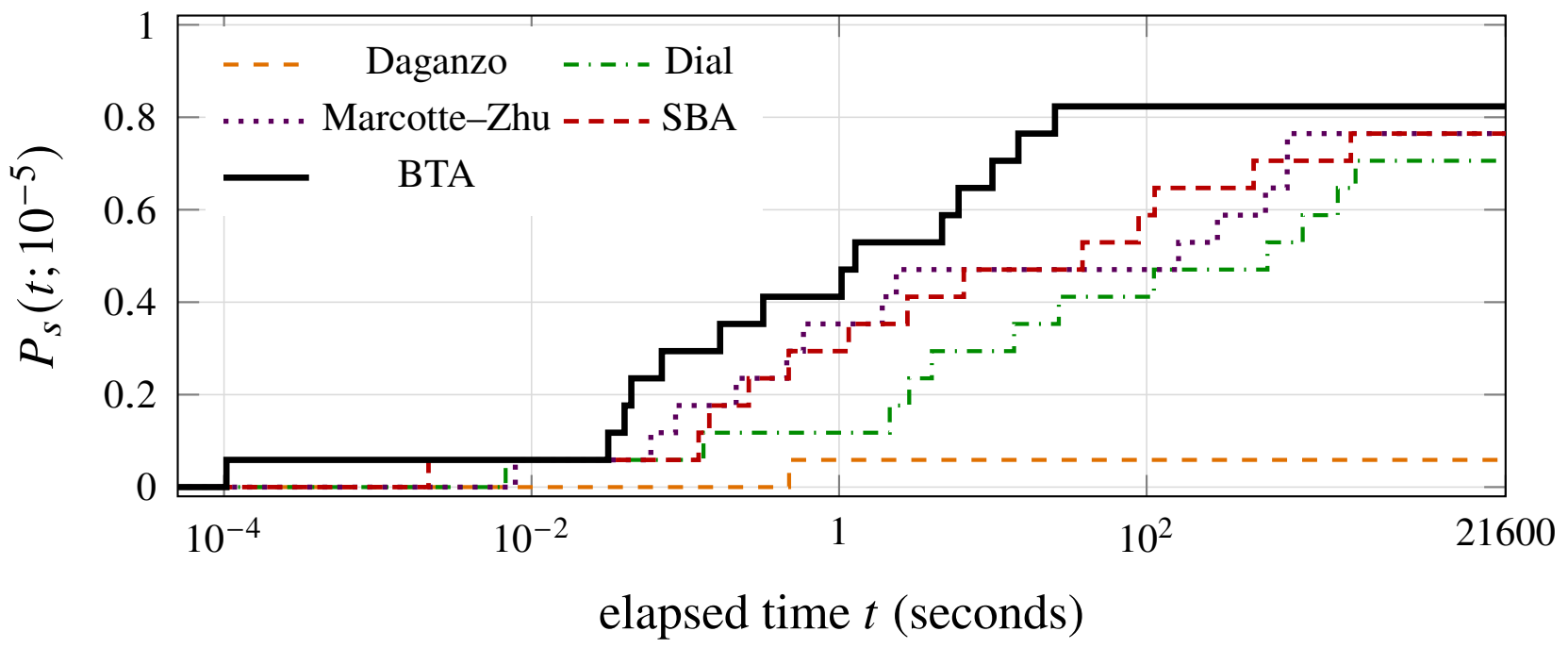

Table~\ref{tab:exp_metrics_summary} and Figure~\ref{fig:exp_common_gap_profile} compare BTA with the methods of \citet{daga:83}, \citet{dial:96}, \citet{marc:96b} and \citet{xie:26} for smooth instances.
No method reached the target tolerance $\varepsilon_{i}$ by time~$T_{i}$ for the Goldcoast, Chicago Regional, or Philadelphia instances, and BTA was the fastest on all $14$ remaining instances.

\begin{table}[t]
\centering
\small
\caption{Comparison of BTA and the methods of \citet{daga:83}, \citet{dial:96}, \citet{marc:96b}, and \citet{xie:26} for smooth instances, in terms of the shifted geometric mean completion time $g_{m}$, its normalized value $r_{m}$, the fastest share $h_{m}$, and the number of instances~$i$ solved in time~$T_{i}$.
}
\label{tab:exp_metrics_summary}
\begin{tabular}{lrrrr}
\toprule
Method &  $g_{m}$ &  $r_{m}$ &  $h_{m}$ & Instances~$i$ solved in time~$T_{i}$ \\
\midrule
Daganzo & $3120.898$ & $250.622$ & $0.000$ & $1/17$ \\
Dial & $196.485$ & $15.779$ & $0.000$ & $12/17$ \\
Marcotte--Zhu & $68.737$ & $5.520$ & $0.000$ & $13/17$ \\
SBA & $59.983$ & $4.817$ & $0.000$ & $13/17$ \\
BTA & $12.453$ & $1.000$ & $0.824$ & $14/17$ \\
\bottomrule
\end{tabular}
\end{table}

Table~\ref{tab:exp_metrics_nonsmooth} compares BTA with the method of \citet{dafe:72} for nonsmooth instances for which $w$ has finite support.
Dafermos reached the target tolerance $\varepsilon_{i}$ by time~$T_{i}$ for all 17 instances, while BTA reached $\varepsilon_{i}$ by time~$T_{i}$ for 15 of the 17 instances (excluding Goldcoast and Philadelphia).
Tables~\ref{tab:exp_smooth_runtimes} and~\ref{tab:exp_dafermos_runtimes} in Appendix~\ref{app:numerical-experiments} report the per-instance computation times for both comparisons.

\begin{table}[t]
\centering
\small
\caption{Comparison of BTA and the method of \citet{dafe:72} for instances with finite support for $w$, in terms of the shifted geometric mean completion time $g_{m}$, its normalized value $r_{m}$, the fastest share $h_{m}$, and the number of instances~$i$ solved in time~$T_{i}$.
}
\label{tab:exp_metrics_nonsmooth}
\begin{tabular}{lrrrr}
\toprule
Method &  $g_{m}$ &  $r_{m}$ &  $h_{m}$ & Instances~$i$ solved in time~$T_{i}$ \\
\midrule
Dafermos & $25.641$ & $1.000$ & $0.882$ & $17/17$ \\
BTA  & $204.737$ & $7.985$ & $0.118$ & $15/17$ \\
\bottomrule
\end{tabular}
\end{table}

Table~\ref{tab:exp_nonsmooth_staircase} reports the solution times $t_{i,\mathrm{BTA}}$ and final gaps $\gamma_{i,\mathrm{BTA}}(T_{i})$ for BTA.
If $t_{i,\mathrm{BTA}} < T_{i}$ then $\gamma_{i,\mathrm{BTA}}(T_{i}) \le \varepsilon_{i}$, and if $t_{i,\mathrm{BTA}} = T_{i}$ then $\gamma_{i,\mathrm{BTA}}(T_{i}) > \varepsilon_{i}$.

\begin{table}[t]
\centering
\scriptsize
\setlength{\tabcolsep}{4pt}
\caption{Solution times $t_{i,\mathrm{BTA}}$ and final gaps $\gamma_{i,\mathrm{BTA}}(T_{i})$ for BTA.}
\label{tab:exp_nonsmooth_staircase}
\begin{tabular}{p{0.235\textwidth}rr@{\hskip 10pt}rr}
\toprule
& \multicolumn{2}{c}{Smooth, $w \sim \mathcal{U}[6,30]$} & \multicolumn{2}{c}{Nonsmooth, $w \sim$ staircase} \\
\cmidrule(lr){2-3}\cmidrule(lr){4-5}
Network & Time & Gap & Time & Gap \\
\midrule
Braess & $<0.001$ & $\le 10^{-5}$ & $0.016$ & $\le 10^{-3}$ \\
Sioux Falls & $0.032$ & $\le 10^{-5}$ & $0.191$ & $\le 10^{-3}$ \\
Berlin Tiergarten & $0.040$ & $\le 10^{-5}$ & $7.915$ & $\le 10^{-3}$ \\
Friedrichshain & $0.045$ & $\le 10^{-5}$ & $15.671$ & $\le 10^{-3}$ \\
Anaheim & $0.321$ & $\le 10^{-5}$ & $140.964$ & $\le 10^{-3}$ \\
Berlin Prenzlauerberg & $0.168$ & $\le 10^{-5}$ & $39.797$ & $\le 10^{-3}$ \\
Berlin Mitte & $0.070$ & $\le 10^{-5}$ & $33.597$ & $\le 10^{-3}$ \\
Winnipeg-Asym & $1.275$ & $\le 10^{-5}$ & $42.728$ & $\le 10^{-3}$ \\
Berlin M.--P.--F. & $1.038$ & $\le 10^{-5}$ & $275.643$ & $\le 10^{-3}$ \\
Terrassa-Asym & $9.935$ & $\le 10^{-5}$ & $280.418$ & $\le 10^{-3}$ \\
Hessen-Asym & $14.657$ & $\le 10^{-5}$ & $314.519$ & $\le 10^{-3}$ \\
Winnipeg & $4.673$ & $\le 10^{-5}$ & $84.343$ & $\le 10^{-3}$ \\
Barcelona & $5.987$ & $\le 10^{-5}$ & $2093.585$ & $\le 10^{-3}$ \\
Berlin Center & $25.352$ & $\le 10^{-5}$ & $2014.743$ & $\le 10^{-3}$ \\
Goldcoast & $21600$ & $1.09\times10^{-2}$ & $86400$ & $4.75\times10^{-3}$ \\
Philadelphia & $21600$ & $5.06\times10^{-4}$ & $86400$ & $1.62\times10^{-2}$ \\
Chicago Regional & $21600$ & $2.03\times10^{-1}$ & $20352.302$ & $\le 10^{-3}$ \\
\bottomrule
\end{tabular}
\end{table}

\section{Conclusion}
\label{sec:conclusion}

The method for the bicriteria traffic assignment problem proposed in this paper has a number of advantages.
First, it can be used with any distribution for the weight that trades off the two criteria, whereas previously proposed methods were either restricted to a distribution with finite support or a distribution with a positive density on an interval.
Second, with the proposed method an equilibrium can be computed by solving a finite dimensional convex optimization problem.
This advantage contributes to the proposed algorithm having better empirical performance than previously proposed algorithms, with the exception that the algorithm of \citet{dafe:72}, designed for instances with a finite number of user classes (the distribution has finite support), often performs better on such instances.
A closer look at the relative performance of the algorithms on these instances reveals a shortcoming of the proposed algorithm --- it spends too much time solving the restricted problem with a small subset of paths to the accuracy required of the final solution, instead of stopping the iteration earlier to generate new paths.
A related enhancement requires a different algorithm, and is the topic of further research.


\bibliographystyle{informs2014}
\bibliography{references}

\ECSwitch


\ECHead{Electronic Companion}

\section{Proofs for Results in Section~\ref{sec:preliminary results}}
\label{app:proofs for preliminary results}

\begin{manuallemma}{\ref{lem:ordering_invariance}}
For any $\pi \in \Pi$, consider two total orders $\prec$ and $\lessdot$ on $\mathcal{R}^{\pi}$ such that for any $r,r' \in \mathcal{R}^{\pi}$ it holds that if $r \prec r'$ or $r \lessdot r'$ then $S^{\pi}_{r} \ge S^{\pi}_{r'}$.
Consider any $\hat{\mathcal{R}}^{\pi} \subset \mathcal{R}^{\pi}$ and any $q^{\pi} \in \Delta_{\hat{\mathcal{R}}^{\pi}}$.
For each $r \in \hat{\mathcal{R}}^{\pi}$, let
\[
p^{\pi}_{r} \ \ \defi \ \ \sum_{\{r' \in \hat{\mathcal{R}}^{\pi} \, : \, r' \preceq r\}} q^{\pi}_{r'}
\qquad \mbox{and} \qquad
\dot{p}^{\pi}_{r} \ \ \defi \ \ \sum_{\{r' \in \hat{\mathcal{R}}^{\pi} \, : \, r' \ledot r\}} q^{\pi}_{r'}
\]
Then
\begin{equation}
\label{eqn:order invariance}
\sum_{r \in \hat{\mathcal{R}}^{\pi}} S^{\pi}_{r} \int_{\mathbb{R}} y \, \diff G^{\pi}_{p^{\pi}_{r-1}, p^{\pi}_{r}}(y)
\ \ = \ \ \sum_{r \in \hat{\mathcal{R}}^{\pi}} S^{\pi}_{r} \int_{\mathbb{R}} y \, \diff G^{\pi}_{\dot{p}^{\pi}_{r-1}, \dot{p}^{\pi}_{r}}(y)
\end{equation}
and thus $\phi(q)$ is invariant with respect to the chosen order.
\end{manuallemma}

\proof{Proof.}
Consider any $r \in \hat{\mathcal{R}}^{\pi}$.
Let $r' \in \hat{\mathcal{R}}^{\pi}$ satisfy $S^{\pi}_{r'} = S^{\pi}_{r}$ and $S^{\pi}_{r'''} > S^{\pi}_{r}$ for all $r''' \in \hat{\mathcal{R}}^{\pi}$ such that $r''' \prec r'$, and let $r'' \in \hat{\mathcal{R}}^{\pi}$ satisfy $S^{\pi}_{r''} = S^{\pi}_{r}$ and $S^{\pi}_{r'''} < S^{\pi}_{r}$ for all $r''' \in \hat{\mathcal{R}}^{\pi}$ such that $r'' \prec r'''$.
Similarly, let $\dot{r}' \in \hat{\mathcal{R}}^{\pi}$ satisfy $S^{\pi}_{\dot{r}'} = S^{\pi}_{r}$ and $S^{\pi}_{r'''} > S^{\pi}_{r}$ for all $r''' \in \hat{\mathcal{R}}^{\pi}$ such that $r''' \lessdot \dot{r}'$, and let $\dot{r}'' \in \hat{\mathcal{R}}^{\pi}$ satisfy $S^{\pi}_{\dot{r}''} = S^{\pi}_{r}$ and $S^{\pi}_{r'''} < S^{\pi}_{r}$ for all $r''' \in \hat{\mathcal{R}}^{\pi}$ such that $\dot{r}'' \lessdot r'''$.
Then $\{r''' \in \hat{\mathcal{R}}^{\pi} \, : \, S^{\pi}_{r'''} > S^{\pi}_{r}\} = \{r''' \in \hat{\mathcal{R}}^{\pi} \, : \, r''' \prec r'\} = \{r''' \in \hat{\mathcal{R}}^{\pi} \, : \, r''' \lessdot \dot{r}'\}$, $\{r''' \in \hat{\mathcal{R}}^{\pi} \, : \, S^{\pi}_{r'''} = S^{\pi}_{r}\} = \{r''' \in \hat{\mathcal{R}}^{\pi} \, : \, r' \preceq r''' \preceq r''\} = \{r''' \in \hat{\mathcal{R}}^{\pi} \, : \, \dot{r}' \ledot r''' \ledot \dot{r}''\}$, and $\{r''' \in \hat{\mathcal{R}}^{\pi} \, : \, S^{\pi}_{r'''} < S^{\pi}_{r}\} = \{r''' \in \hat{\mathcal{R}}^{\pi} \, : \, r'' \prec r'''\} = \{r''' \in \hat{\mathcal{R}}^{\pi} \, : \, \dot{r}'' \lessdot r'''\}$.
Thus,
\[
p^{\pi}_{r'-1} \ \ = \ \ \sum_{\{r''' \in \hat{\mathcal{R}}^{\pi} \, : \, r''' \prec r'\}} q^{\pi}_{r'''}
\ \ = \ \ \sum_{\{r''' \in \hat{\mathcal{R}}^{\pi} \, : \, r''' \lessdot \dot{r}'\}} q^{\pi}_{r'''}
\ \ = \ \ \dot{p}^{\pi}_{\dot{r}'-1}
\]
and
\[
p^{\pi}_{r''} \ \ = \ \ p^{\pi}_{r'-1} + \sum_{\{r''' \in \hat{\mathcal{R}}^{\pi} \, : \, r' \preceq r''' \preceq r''\}} q^{\pi}_{r'''}
\ \ = \ \ \dot{p}^{\pi}_{\dot{r}'-1} + \sum_{\{r''' \in \hat{\mathcal{R}}^{\pi} \, : \, \dot{r}' \ledot r''' \ledot \dot{r}''\}} q^{\pi}_{r'''}
\ \ = \ \ \dot{p}^{\pi}_{\dot{r}''}
\]
Hence,
\begin{align*}
\sum_{\{r''' \in \hat{\mathcal{R}}^{\pi} \, : \, S^{\pi}_{r'''} = S^{\pi}_{r}\}} S^{\pi}_{r'''} \int_{\mathbb{R}} y \, \diff G^{\pi}_{p^{\pi}_{r'''-1}, p^{\pi}_{r'''}}(y)
\ \ & = \ \ S^{\pi}_{r'} \sum_{\{r''' \in \hat{\mathcal{R}}^{\pi} \, : \, r' \preceq r''' \preceq r''\}} \int_{\mathbb{R}} y \, \diff G^{\pi}_{p^{\pi}_{r'''-1}, p^{\pi}_{r'''}}(y) \\
& \overset{\eqref{eqn:additivity}}{=} \ \ S^{\pi}_{r} \int_{\mathbb{R}} y \, \diff G^{\pi}_{p^{\pi}_{r'-1}, p^{\pi}_{r''}}(y) \\
& = \ \ S^{\pi}_{r} \int_{\mathbb{R}} y \, \diff G^{\pi}_{\dot{p}^{\pi}_{\dot{r}'-1}, \dot{p}^{\pi}_{\dot{r}''}}(y) \\
& \overset{\eqref{eqn:additivity}}{=} \ \ S^{\pi}_{r} \sum_{\{r''' \in \hat{\mathcal{R}}^{\pi} \, : \, \dot{r}' \ledot r''' \ledot \dot{r}''\}} \int_{\mathbb{R}} y \, \diff G^{\pi}_{\dot{p}^{\pi}_{r'''-1}, \dot{p}^{\pi}_{r'''}}(y) \\
& = \ \ \sum_{\{r''' \in \hat{\mathcal{R}}^{\pi} \, : \, S^{\pi}_{r'''} = S^{\pi}_{r}\}} S^{\pi}_{r'''} \int_{\mathbb{R}} y \, \diff G^{\pi}_{\dot{p}^{\pi}_{r'''-1}, \dot{p}^{\pi}_{r'''}}(y)
\end{align*}
The equality holds for every value $S^{\pi}_{r}$ and every $\pi \in \Pi$, and therefore~\eqref{eqn:order invariance} holds.
\ignore{
Any two compatible orderings can be connected by a sequence of adjacent swaps involving only paths with equal tolls $S_r$. It is therefore enough to consider one adjacent swap of consecutive paths $r_i$ and $r_{i+1}$ with $S_{r_i}=S_{r_{i+1}}$. Let $p_{i-1}$ and $p_{i+1}$ denote the cumulative sums before and after these two paths in the original ordering, let $p_i=p_{i-1}+q_{r_i}$, and let $p_i'=p_{i-1}+q_{r_{i+1}}$ denote the corresponding intermediate cumulative sum after the swap.

The only terms in \eqref{eq:toll-cost} affected by the swap are the two consecutive contributions associated with $r_i$ and $r_{i+1}$. In the original ordering they equal
\[
S_{r_i}\int_{\mathbb{R}} y\,\diff G^{\pi}_{p_{i-1},p_i}(y)
+
S_{r_{i+1}}\int_{\mathbb{R}} y\,\diff G^{\pi}_{p_i,p_{i+1}}(y),
\]
which, because $S_{r_i}=S_{r_{i+1}}$ and by the additivity property \eqref{eq:additivity}, is
\[
S_{r_i}\int_{\mathbb{R}} y\,\diff G^{\pi}_{p_{i-1},p_{i+1}}(y).
\]
After the swap, the same two terms become
\[
S_{r_{i+1}}\int_{\mathbb{R}} y\,\diff G^{\pi}_{p_{i-1},p_i'}(y)
+
S_{r_i}\int_{\mathbb{R}} y\,\diff G^{\pi}_{p_i',p_{i+1}}(y),
\]
which again equals $S_{r_i}\int_{\mathbb{R}} y\,\diff G^{\pi}_{p_{i-1},p_{i+1}}(y)$ by \eqref{eq:additivity}. All other terms are unchanged, so one adjacent equal-$S$ swap leaves $F^{\pi}(q^{\pi})$ invariant. Repeating this argument along the swap sequence yields the claim.
}
\Halmos
\endproof

\begin{manuallemma}{\ref{lem:varphi technical properties}}
The function $\varphi : [0,1] \to \mathbb{R}$ is convex and Lipschitz continuous, and $\partial \varphi(p) = G^{-1}(p)$ for every $p \in [0,1]$.
\end{manuallemma}

\proof{Proof.}
Recall that $G^{-1}(p) \ne \varnothing$ for all $p \in [0,1]$.
We first show that for every $\hat{p} \in [0,1]$ and every $\hat{w} \in G^{-1}(\hat{p})$, $\hat{w}$ is a subgradient of $\varphi$ at $\hat{p}$, which establishes that $\varphi$ is convex and that $G^{-1}(\hat{p}) \subset \partial \varphi(\hat{p})$.
For every $p \in [\hat{p},1]$ it holds that
\begin{align*}
\varphi(p) - \varphi(\hat{p})
\ \ = \ \ \int_{\mathbb{R}} w \, \diff G_{0,p}(w) - \int_{\mathbb{R}} w \, \diff G_{0,\hat{p}}(w)
\ \ = \ \ \int_{\mathbb{R}} w \, \diff G_{\hat{p},p}(w)
\ \ \ge \ \ \int_{\mathbb{R}} \hat{w} \, \diff G_{\hat{p},p}(w)
\ \ = \ \ \hat{w} (p - \hat{p}).
\end{align*}
Also, for every $p \in [0,\hat{p}]$ it holds that
\begin{align*}
\varphi(p) - \varphi(\hat{p})
\ \ = \ \ \int_{\mathbb{R}} w \, \diff G_{0,p}(w) - \int_{\mathbb{R}} w \, \diff G_{0,\hat{p}}(w)
\ \ = \ \ - \int_{\mathbb{R}} w \, \diff G_{p,\hat{p}}(w)
\ \ \ge \ \ - \int_{\mathbb{R}} \hat{w} \, \diff G_{p,\hat{p}}(w)
\ \ = \ \ \hat{w} (p - \hat{p}).
\end{align*}
Thus $\varphi(p) - \varphi(\hat{p}) \ge \hat{w} (p - \hat{p})$ for all $p \in [0,1]$.
It follows that $\varphi$ is convex and that $G^{-1}(\hat{p}) \subset \partial \varphi(\hat{p})$ for all $\hat{p} \in [0,1]$.
Also, for any $0 \le \hat{p} \le p \le 1$ it holds that
\begin{align*}
\varphi(p) - \varphi(\hat{p})
\ \ = \ \ \int_{\underline{w}}^{\overline{w}} w \, \diff G_{\hat{p},p}(w)
\ \ \le \ \ \overline{w} (p - \hat{p})
\end{align*}
and thus $\varphi$ is Lipschitz continuous with Lipschitz constant $\overline{w}$.

Next we show that $\partial \varphi(\hat{p}) \subset G^{-1}(\hat{p})$ for all $\hat{p} \in [0,1]$.
Consider any $\hat{w} \notin G^{-1}(\hat{p})$.
Then, either $\hat{w} < \min\{w \in G^{-1}(\hat{p})\}$ or $\hat{w} > \max\{w \in G^{-1}(\hat{p})\}$.

First, consider the case in which $\hat{w} < \min\{w \in G^{-1}(\hat{p})\}$, that is, there is a $\check{w} > \hat{w}$ such that $\check{p} \defi G(\check{w}) < \hat{p}$.
Then
\begin{align*}
\varphi(\check{p}) - \varphi(\hat{p})
\ \ = \ \ - \int_{\mathbb{R}} w \, \diff G_{\check{p},\hat{p}}(w)
\ \ \le \ \ - \int_{\mathbb{R}} \check{w} \, \diff G_{\check{p},\hat{p}}(w)
\ \ = \ \ - \check{w} (\hat{p} - \check{p})
\ \ < \ \ \hat{w} (\check{p} - \hat{p})
\end{align*}
and thus $\hat{w} \notin \partial \varphi(\hat{p})$.
Next, consider the case in which $\hat{w} > \max\{w \in G^{-1}(\hat{p})\}$, that is, there is a $\bar{w} < \hat{w}$ such that $\bar{p} \defi G(\bar{w}) > \hat{p}$.
Then
\begin{align*}
\varphi(\bar{p}) - \varphi(\hat{p})
\ \ = \ \ \int_{\mathbb{R}} w \, \diff G_{\hat{p},\bar{p}}(w)
\ \ \le \ \ \int_{\mathbb{R}} \bar{w} \, \diff G_{\hat{p},\bar{p}}(w)
\ \ = \ \ \bar{w} (\bar{p} - \hat{p})
\ \ < \ \ \hat{w} (\bar{p} - \hat{p})
\end{align*}
and thus $\hat{w} \notin \partial \varphi(\hat{p})$.
Hence $\partial \varphi(\hat{p}) \subset G^{-1}(\hat{p})$ for all $\hat{p} \in [0,1]$.

Therefore it has been established that $\partial \varphi(p) = G^{-1}(p)$ for every $p \in [0,1]$.
\Halmos
\endproof

\begin{manuallemma}{\ref{lem:phi technical properties}}
The function $\phi : \prod_{\pi \in \Pi} \Delta_{\hat{\mathcal{R}}^{\pi}} \to \mathbb{R}$ is convex and continuous, and
\[
\partial \phi(q) \ \ = \ \ \prod_{\pi \in \Pi} \left(D^{\pi} \sum_{\{r \in \hat{\mathcal{R}}^{\pi} \, : \, r+1 \in \hat{\mathcal{R}}^{\pi}\}} \left(S^{\pi}_{r} - S^{\pi}_{r+1}\right) \left(\mathds{1}\left\{r' \preceq r\right\}, \; r' \in \hat{\mathcal{R}}^{\pi}\right) \left(G^{\pi}\right)^{-1}\big(L^{\pi}_{r}(q^{\pi})\big) + \mathcal{N}_{\Delta_{\hat{\mathcal{R}}^{\pi}}}(q^{\pi})\right)
\]
\end{manuallemma}

\proof{Proof.}
\ignore{
Consider any $\hat{q} \in \square_{\hat{\mathcal{R}}^{\pi}}$ and any $\hat{w} \in \left(G^{\pi}\right)^{-1}\big(L^{\pi}_{r}(\hat{q})\big)$.
It follows from Lemma~\ref{lem:varphi technical properties} that for any $q \in \square_{\hat{\mathcal{R}}^{\pi}}$ it holds that
\begin{align*}
\varphi^{\pi}\big(L^{\pi}_{r}(q)\big) - \varphi^{\pi}\big(L^{\pi}_{r}(\hat{q})\big)
\ \ \ge \ \ \hat{w} \left[L^{\pi}_{r}(q) - L^{\pi}_{r}(\hat{q})\right]
\ \ = \ \ \hat{w} \sum_{\{r' \in \hat{\mathcal{R}}^{\pi} \, : \, r' \preceq r\}} (q_{r'} - \hat{q}_{r'})
\ \ = \ \ \hat{w} \left(\mathds{1}\left\{r_{1} \preceq r\right\}, \ldots, \mathds{1}\left\{r_{k} \preceq r\right\}\right)^{\top} (q - \hat{q})
\end{align*}
and thus $\left(\mathds{1}\left\{r_{1} \preceq r\right\}, \ldots, \mathds{1}\left\{r_{k} \preceq r\right\}\right) \left(G^{\pi}\right)^{-1}\big(L^{\pi}_{r}(\hat{q})\big) \subset \partial \phi^{\pi}_{r}(\hat{q})$.
}
It follows from $S^{\pi}_{r} \ge S^{\pi}_{r+1}$ and Lemma~\ref{lem:varphi technical properties} that $\phi$ is convex and continuous.
For each $\pi \in \Pi$, let $\phi^{\pi} : \Delta_{\hat{\mathcal{R}}^{\pi}} \to \mathbb{R}$ be given by
\[
\phi^{\pi}(q^{\pi}) \ \ \defi \ \ D^{\pi} \sum_{r \in \hat{\mathcal{R}}^{\pi}} \left[\mathds{1}\{r+1 \in \hat{\mathcal{R}}^{\pi}\} \left(S^{\pi}_{r} - S^{\pi}_{r+1}\right) \phi^{\pi}_{r}(q^{\pi}) + \mathds{1}\{r+1 \notin \hat{\mathcal{R}}^{\pi}\} S^{\pi}_{\min} \varphi^{\pi}(1)\right]
\]
For each $\pi \in \Pi$, it holds that
\[
\bigcap_{r \in \hat{\mathcal{R}}^{\pi}} \text{ri}\big(\text{dom}(\phi^{\pi}_{r})\big) \ \ = \ \ \Delta_{\hat{\mathcal{R}}^{\pi}}^{\circ} \ \ \neq \ \ \varnothing
\]
and thus it follows from the Moreau-Rockafellar Theorem \citep[Theorem 23.8]{rockafellar1997convex} that
\begin{align*}
\partial \phi^{\pi}(q^{\pi}) \ \ & = \ \ D^{\pi} \sum_{\{r \in \hat{\mathcal{R}}^{\pi} \, : \, r+1 \in \hat{\mathcal{R}}^{\pi}\}} \left(S^{\pi}_{r} - S^{\pi}_{r+1}\right) \partial \phi^{\pi}_{r}(q^{\pi}) + \mathcal{N}_{\Delta_{\hat{\mathcal{R}}^{\pi}}}(q^{\pi}) \\
& = \ \ D^{\pi} \sum_{\{r \in \hat{\mathcal{R}}^{\pi} \, : \, r+1 \in \hat{\mathcal{R}}^{\pi}\}} \left(S^{\pi}_{r} - S^{\pi}_{r+1}\right) \left(\mathds{1}\left\{r' \preceq r\right\}, \; r' \in \hat{\mathcal{R}}^{\pi}\right) \left(G^{\pi}\right)^{-1}\big(L^{\pi}_{r}(q^{\pi})\big) + \mathcal{N}_{\Delta_{\hat{\mathcal{R}}^{\pi}}}(q^{\pi})
\end{align*}
Then it follows from $\phi(q) = \sum_{\pi \in \Pi} \phi^{\pi}(q^{\pi})$ being separable that
\begin{align*}
\partial \phi(q) \ \ & = \ \ \prod_{\pi \in \Pi} \partial \phi^{\pi}(q^{\pi}) \\
& = \ \ \prod_{\pi \in \Pi} \left(D^{\pi} \sum_{\{r \in \hat{\mathcal{R}}^{\pi} \, : \, r+1 \in \hat{\mathcal{R}}^{\pi}\}} \left(S^{\pi}_{r} - S^{\pi}_{r+1}\right) \left(\mathds{1}\left\{r' \preceq r\right\}, \; r' \in \hat{\mathcal{R}}^{\pi}\right) \left(G^{\pi}\right)^{-1}\big(L^{\pi}_{r}(q^{\pi})\big) + \mathcal{N}_{\Delta_{\hat{\mathcal{R}}^{\pi}}}(q^{\pi})\right)
\end{align*}
\Halmos
\endproof

\begin{manuallemma}{\ref{lem:psi technical properties}}
The function $\psi$ is convex and continuous on $\prod_{\pi \in \Pi} \Delta_{\hat{\mathcal{R}}^{\pi}}$, and
\[
\partial \psi(q) \ \ = \ \ \left\{g = \left(g^{\pi}_{r}, \pi \in \Pi, r \in \hat{\mathcal{R}}^{\pi}\right) \; : \; g^{\pi}_{r} = D^{\pi} \sum_{a \in \mathcal{A}^{\pi}_{r}} \tau_{a}, \; \tau_{a} \in \left[t_{a}^{-}(x_{a}(q)), t_{a}^{+}(x_{a}(q))\right] \ \forall \; a \in \mathcal{A}\right\} + \prod_{\pi \in \Pi} \mathcal{N}_{\Delta_{\hat{\mathcal{R}}^{\pi}}}(q^{\pi})
\]
with $\left[t_{a}^{-}(0), t_{a}^{+}(0)\right] \defi \left(-\infty, t_{a}^{+}(0)\right]$.
\end{manuallemma}

\proof{Proof.}
It is convenient to extend the domain of $\psi$ before establishing its subdifferential:
Let
\[
Q \ \ \defi \ \ \left\{q = \left(q^{\pi}_{r}, \pi \in \Pi, r \in \hat{\mathcal{R}}^{\pi}\right) \; : \; x_{a}(q) \ge 0 \; \forall \; a \in \mathcal{A}\right\}
\]
Note that $\prod_{\pi \in \Pi} \Delta_{\hat{\mathcal{R}}^{\pi}} \subset \prod_{\pi \in \Pi} \mathbb{R}_{+}^{\hat{\mathcal{R}}^{\pi}} \subset Q \subset \prod_{\pi \in \Pi} \mathbb{R}^{\hat{\mathcal{R}}^{\pi}}$.
Let $\psi : \prod_{\pi \in \Pi} \mathbb{R}^{\hat{\mathcal{R}}^{\pi}} \mapsto \mathbb{R} \cup \{\infty\}$ be given by~\eqref{eqn:travel time objective} if $q \in \prod_{\pi \in \Pi} \Delta_{\hat{\mathcal{R}}^{\pi}}$, and $\psi(q) = \infty$ if $q \notin \prod_{\pi \in \Pi} \Delta_{\hat{\mathcal{R}}^{\pi}}$, that is, $\psi(q) = \hat{\psi}(q) + \sum_{\pi \in \Pi} \delta^{\pi}_{\Delta_{\hat{\mathcal{R}}^{\pi}}}(q^{\pi})$ for all $q \in \prod_{\pi \in \Pi} \mathbb{R}^{\hat{\mathcal{R}}^{\pi}}$, where
\begin{equation*}
\hat{\psi}(q) \ \ \defi \ \ \begin{cases}
\sum_{a \in \mathcal{A}} \int_{0}^{x_{a}(q)} t_{a}(y) \, \diff y & \mbox{ if } q \in Q \\
\infty & \mbox{ if } q \notin Q
\end{cases}
\end{equation*}
and, for each $\pi \in \Pi$, $\delta^{\pi}_{\Delta_{\hat{\mathcal{R}}^{\pi}}} : \mathbb{R}^{\hat{\mathcal{R}}^{\pi}} \mapsto \mathbb{R} \cup \{\infty\}$ is the indicator function of $\Delta_{\hat{\mathcal{R}}^{\pi}}$ given by
\begin{equation*}
\delta^{\pi}_{\Delta_{\hat{\mathcal{R}}^{\pi}}}(q^{\pi}) \ \ \defi \ \ \begin{cases}
0 & \mbox{ if } q^{\pi} \in \Delta_{\hat{\mathcal{R}}^{\pi}} \\
\infty & \mbox{ if } q^{\pi} \notin \Delta_{\hat{\mathcal{R}}^{\pi}}
\end{cases}
\end{equation*}
Note that
\[
\text{ri}\left(\text{dom}\big(\hat{\psi}\big)\right) \, \cap \, \text{ri}\left(\text{dom}\bigg(\sum_{\pi \in \Pi} \delta^{\pi}_{\Delta_{\hat{\mathcal{R}}^{\pi}}}\bigg)\right) \ \ = \ \ \prod_{\pi \in \Pi} \Delta_{\hat{\mathcal{R}}^{\pi}}^{\circ} \ \ \neq \ \ \varnothing
\]
and thus it follows from the Moreau-Rockafellar Theorem \citep[Theorem 23.8]{rockafellar1997convex} and $\psi(q) = \hat{\psi}(q) + \sum_{\pi \in \Pi} \delta^{\pi}_{\Delta_{\hat{\mathcal{R}}^{\pi}}}(q^{\pi})$ that $\partial \psi(q) = \partial \hat{\psi}(q) + \partial \sum_{\pi \in \Pi} \delta^{\pi}_{\Delta_{\hat{\mathcal{R}}^{\pi}}}(q^{\pi})$.
Also, it follows from $\sum_{\pi \in \Pi} \delta^{\pi}_{\Delta_{\hat{\mathcal{R}}^{\pi}}}(q^{\pi})$ being separable that $\partial \sum_{\pi \in \Pi} \delta^{\pi}_{\Delta_{\hat{\mathcal{R}}^{\pi}}}(q^{\pi}) = \prod_{\pi \in \Pi} \partial \delta^{\pi}_{\Delta_{\hat{\mathcal{R}}^{\pi}}}(q^{\pi}) = \prod_{\pi \in \Pi} \mathcal{N}_{\Delta_{\hat{\mathcal{R}}^{\pi}}}(q^{\pi})$.

Next we determine $\partial \hat{\psi}(q)$.
Consider any $\hat{q} \in Q$ and any $\hat{g}$ such that $\hat{g}^{\pi}_{r} = D^{\pi} \sum_{a \in \mathcal{A}^{\pi}_{r}} \tau_{a}$, $\tau_{a} \in \left[t_{a}^{-}(x_{a}(\hat{q})), t_{a}^{+}(x_{a}(\hat{q}))\right]$ $\; \forall \, a \in \mathcal{A}$.
For any $q \in Q$ and $a \in \mathcal{A}$, first consider the case in which $x_{a}(q) \ge x_{a}(\hat{q})$.
Then
\begin{align*}
\int_{0}^{x_{a}(q)} t_{a}(y) \, \diff y - \int_{0}^{x_{a}(\hat{q})} t_{a}(y) \, \diff y
\ \ & = \ \ \int_{x_{a}(\hat{q})}^{x_{a}(q)} t_{a}(y) \, \diff y \\
& \ge \ \ t_{a}^{+}(x_{a}(\hat{q})) \left[x_{a}(q) - x_{a}(\hat{q})\right]
\ \ \ge \ \ \tau_{a} \left[x_{a}(q) - x_{a}(\hat{q})\right]
\end{align*}
Next consider the case in which $x_{a}(q) < x_{a}(\hat{q})$.
Then
\begin{align*}
\int_{0}^{x_{a}(q)} t_{a}(y) \, \diff y - \int_{0}^{x_{a}(\hat{q})} t_{a}(y) \, \diff y
\ \ & = \ \ - \int_{x_{a}(q)}^{x_{a}(\hat{q})} t_{a}(y) \, \diff y \\
& \ge \ \ - t_{a}^{-}(x_{a}(\hat{q})) \left[x_{a}(\hat{q}) - x_{a}(q)\right]
\ \ \ge \ \ \tau_{a} \left[x_{a}(q) - x_{a}(\hat{q})\right]
\end{align*}
Thus,
\begin{align*}
\hat{\psi}(q) - \hat{\psi}(\hat{q}) \ \ & = \ \ \sum_{a \in \mathcal{A}} \int_{0}^{x_{a}(q)} t_{a}(y) \, \diff y - \sum_{a \in \mathcal{A}} \int_{0}^{x_{a}(\hat{q})} t_{a}(y) \, \diff y \\
& \ge \ \ \sum_{a \in \mathcal{A}} \tau_{a} \left[x_{a}(q) - x_{a}(\hat{q})\right] \\
& = \ \ \sum_{a \in \mathcal{A}} \tau_{a} \left[\sum_{\pi \in \Pi} D^{\pi} \sum_{r \in \hat{\mathcal{R}}^{\pi}_{a}} q^{\pi}_{r} - \sum_{\pi \in \Pi} D^{\pi} \sum_{r \in \hat{\mathcal{R}}^{\pi}_{a}} \hat{q}^{\pi}_{r}\right] \\
& = \ \ \sum_{\pi \in \Pi} \sum_{r \in \hat{\mathcal{R}}^{\pi}} D^{\pi} \sum_{a \in \mathcal{A}^{\pi}_{r}} \tau_{a} \left(q^{\pi}_{r} - \hat{q}^{\pi}_{r}\right)
\ \ = \ \ \hat{g}^{\top} \left(q - \hat{q}\right)
\end{align*}
and hence $\hat{g} \in \partial \hat{\psi}(\hat{q})$.
Therefore,
\[
\partial \hat{\psi}(q) \ \ \supset \ \ \left\{g = \left(g^{\pi}_{r}, \pi \in \Pi, r \in \hat{\mathcal{R}}^{\pi}\right) \; : \; g^{\pi}_{r} = D^{\pi} \sum_{a \in \mathcal{A}^{\pi}_{r}} \tau_{a}, \; \tau_{a} \in \left[t_{a}^{-}(x_{a}(q)), t_{a}^{+}(x_{a}(q))\right] \ \forall \; a \in \mathcal{A}\right\}
\]

Next, consider any $\hat{g} \notin \left\{g = \left(g^{\pi}_{r}, \pi \in \Pi, r \in \hat{\mathcal{R}}^{\pi}\right) \; : \; g^{\pi}_{r} = D^{\pi} \sum_{a \in \mathcal{A}^{\pi}_{r}} \tau_{a}, \; \tau_{a} \in \left[t_{a}^{-}(x_{a}(\hat{q})), t_{a}^{+}(x_{a}(\hat{q}))\right] \ \forall \; a \in \mathcal{A}\right\}$.
Let $\mathcal{A}_{0} \defi \{a \in \mathcal{A} \, : \, x_{a}(\hat{q}) = 0\}$, $\mathcal{A}_{=} \defi \{a \in \mathcal{A} \, : \, x_{a}(\hat{q}) > 0, \, t_{a}^{-}(x_{a}(\hat{q})) = t_{a}^{+}(x_{a}(\hat{q}))\}$, and $\mathcal{A}_{+} \defi \{a \in \mathcal{A} \, : \, x_{a}(\hat{q}) > 0, \, t_{a}^{-}(x_{a}(\hat{q})) < t_{a}^{+}(x_{a}(\hat{q}))\}$.
Then
\begin{align}
\label{eqn:auxiliary objective}
\hat{z} \ \ \defi \ \ \min_{\theta,\tau} \ \ & \sum_{\pi \in \Pi} \sum_{r \in \hat{\mathcal{R}}^{\pi}} \theta^{\pi}_{r} \\
\mbox{s.t.} \ \ & \theta^{\pi}_{r} \ \ \ge \ \ \hat{g}^{\pi}_{r} - D^{\pi} \sum_{a \in \mathcal{A}^{\pi}_{r}} \tau_{a} \qquad \forall \ \pi \in \Pi, \ \forall \ r \in \hat{\mathcal{R}}^{\pi}
\label{eqn:auxiliary positive deviation} \\
& \theta^{\pi}_{r} \ \ \ge \ \ D^{\pi} \sum_{a \in \mathcal{A}^{\pi}_{r}} \tau_{a} - \hat{g}^{\pi}_{r} \qquad \forall \ \pi \in \Pi, \ \forall \ r \in \hat{\mathcal{R}}^{\pi}
\label{eqn:auxiliary negative deviation} \\
& \tau_{a} \ \ \ge \ \ t_{a}^{-}(x_{a}(\hat{q})) \qquad \forall \ a \in \mathcal{A}_{+}
\label{eqn:auxiliary lower derivative} \\
& \tau_{a} \ \ \le \ \ t_{a}^{+}(x_{a}(\hat{q})) \qquad \forall \ a \in \mathcal{A}_{0} \cup \mathcal{A}_{+}
\label{eqn:auxiliary upper derivative} \\
& \tau_{a} \ \ = \ \ t_{a}(x_{a}(\hat{q})) \qquad \forall \ a \in \mathcal{A}_{=}
\label{eqn:auxiliary equal derivative}
\end{align}
satisfies $\hat{z} > 0$.
Let $\eta^{\pi}_{r} \ge 0$ denote the dual variables associated with~\eqref{eqn:auxiliary positive deviation}, let $\zeta^{\pi}_{r} \ge 0$ denote the dual variables associated with~\eqref{eqn:auxiliary negative deviation}, let $\lambda^{-}_{a} \ge 0$ denote the dual variables associated with~\eqref{eqn:auxiliary lower derivative}, let $\lambda^{+}_{a} \ge 0$ denote the dual variables associated with~\eqref{eqn:auxiliary upper derivative}, and let $\lambda^{=}_{a}$ denote the dual variables associated with~\eqref{eqn:auxiliary equal derivative}.
Then
\begin{align}
\label{eqn:auxiliary dual objective}
\hat{z} \ \ = \ \ \max_{\eta,\zeta,\lambda^{-},\lambda^{+},\lambda^{=}} \ \ & \left\{\sum_{\pi \in \Pi} \sum_{r \in \hat{\mathcal{R}}^{\pi}} \hat{g}^{\pi}_{r} \left(\eta^{\pi}_{r} - \zeta^{\pi}_{r}\right) + \sum_{a \in \mathcal{A}_{+}} t_{a}^{-}(x_{a}(\hat{q})) \lambda^{-}_{a} - \sum_{a \in \mathcal{A}_{0} \cup \mathcal{A}_{+}} t_{a}^{+}(x_{a}(\hat{q})) \lambda^{+}_{a} + \sum_{a \in \mathcal{A}_{=}} t_{a}(x_{a}(\hat{q})) \lambda^{=}_{a}\right\} \\
\mbox{s.t.} \ \ & \eta^{\pi}_{r} + \zeta^{\pi}_{r} \ \ = \ \ 1 \qquad \forall \ \pi \in \Pi, \ \forall \ r \in \hat{\mathcal{R}}^{\pi}
\label{eqn:auxiliary total deviation} \\
& - \lambda^{+}_{a} + \sum_{\pi \in \Pi} \sum_{r \in \hat{\mathcal{R}}^{\pi}_{a}} D^{\pi} \left(\eta^{\pi}_{r} - \zeta^{\pi}_{r}\right) \ \ = \ \ 0 \qquad \forall \ a \in \mathcal{A}_{0}
\label{eqn:auxiliary upper flow adjustment} \\
& \left(\lambda^{-}_{a} - \lambda^{+}_{a}\right) + \sum_{\pi \in \Pi} \sum_{r \in \hat{\mathcal{R}}^{\pi}_{a}} D^{\pi} \left(\eta^{\pi}_{r} - \zeta^{\pi}_{r}\right) \ \ = \ \ 0 \qquad \forall \ a \in \mathcal{A}_{+}
\label{eqn:auxiliary total flow adjustment} \\
& \lambda^{=}_{a} + \sum_{\pi \in \Pi} \sum_{r \in \hat{\mathcal{R}}^{\pi}_{a}} D^{\pi} \left(\eta^{\pi}_{r} - \zeta^{\pi}_{r}\right) \ \ = \ \ 0 \qquad \forall \ a \in \mathcal{A}_{=}
\label{eqn:auxiliary equal flow adjustment} \\
& \eta,\zeta,\lambda^{-},\lambda^{+} \ \ \ge \ \ 0
\label{eqn:auxiliary sign restrictions}
\end{align}
Let $(\hat{\eta},\hat{\zeta},\hat{\lambda}^{-},\hat{\lambda}^{+},\hat{\lambda}^{=})$ be an optimal solution of \eqref{eqn:auxiliary dual objective}--\eqref{eqn:auxiliary sign restrictions}.
Let $\check{q} : \mathbb{R}_{+} \mapsto \prod_{\pi \in \Pi} \mathbb{R}^{\hat{\mathcal{R}}^{\pi}}$ be given by $\check{q}(\alpha) \defi \hat{q} + \alpha \left(\hat{\eta} - \hat{\zeta}\right)$.
Consider any $\varepsilon \in (0, \hat{z} / \sum_{a \in \mathcal{A}} \sum_{\pi \in \Pi} \sum_{r \in \hat{\mathcal{R}}^{\pi}_{a}} D^{\pi})$.

For any $a \in \mathcal{A}_{0}$, it holds that
\begin{equation}
0 \ \ \overset{\eqref{eqn:auxiliary sign restrictions}}{\le} \ \ \alpha \hat{\lambda}^{+}_{a} \ \ \overset{\eqref{eqn:auxiliary upper flow adjustment}}{=} \ \ \alpha \sum_{\pi \in \Pi} \sum_{r \in \hat{\mathcal{R}}^{\pi}_{a}} D^{\pi} \left(\hat{\eta}^{\pi}_{r} - \hat{\zeta}^{\pi}_{r}\right)
\ \ = \ \ \sum_{\pi \in \Pi} \sum_{r \in \hat{\mathcal{R}}^{\pi}_{a}} D^{\pi} \left(\check{q}^{\pi}_{r}(\alpha) - \hat{q}^{\pi}_{r}\right)
\ \ = \ \ x_{a}(\check{q}(\alpha)) - x_{a}(\hat{q})
\label{eqn:flow adjustment a in A0}
\end{equation}
Thus $x_{a}(\check{q}(\alpha)) \ge x_{a}(\hat{q}) = 0$ for all $\alpha \ge 0$.
Also note that $x_{a}(\check{q}(\cdot))$ is continuous, and thus $\lim_{\alpha \downarrow 0} x_{a}(\check{q}(\alpha)) = x_{a}(\hat{q})$ and $\lim_{\alpha \downarrow 0} t_{a}\big(x_{a}(\check{q}(\alpha))\big) \le t_{a}^{+}\big(x_{a}(\hat{q})\big)$.
(It is possible that $\hat{\lambda}^{+}_{a} = 0$ and $x_{a}(\check{q}(\alpha)) = x_{a}(\hat{q})$ for all $\alpha$, in which case it may hold that $t_{a}\big(x_{a}(\check{q}(\alpha))\big) < t_{a}^{+}\big(x_{a}(\hat{q})\big)$ for all $\alpha$.)
Let $A_{a} > 0$ be such that $t_{a}\big(x_{a}(\check{q}(\alpha))\big) \le t_{a}^{+}\big(x_{a}(\hat{q})\big) + \varepsilon$ for all $\alpha \in [0,A_{a}]$.
Then, for all $\alpha \in [0,A_{a}]$ it holds that
\begin{align}
\hat{\psi}_{a}\big(\check{q}(\alpha)\big) - \hat{\psi}_{a}\big(\hat{q}\big)
\ \ & \defi \ \ \int_{0}^{x_{a}(\check{q}(\alpha))} t_{a}(y) \, \diff y - \int_{0}^{x_{a}(\hat{q})} t_{a}(y) \, \diff y
\ \ = \ \ \int_{x_{a}(\hat{q})}^{x_{a}(\check{q}(\alpha))} t_{a}(y) \, \diff y \nonumber \\
& \le \ \ \left[t_{a}^{+}\big(x_{a}(\hat{q})\big) + \varepsilon\right] \left[x_{a}(\check{q}(\alpha)) - x_{a}(\hat{q})\right]
\ \ \overset{\eqref{eqn:flow adjustment a in A0}}{=} \ \ \left[t_{a}^{+}\big(x_{a}(\hat{q})\big) + \varepsilon\right] \alpha \hat{\lambda}^{+}_{a}
\label{eqn:change bound a in A0}
\end{align}
It follows from~\eqref{eqn:auxiliary total deviation} and~\eqref{eqn:auxiliary sign restrictions} that $\hat{\eta}^{\pi}_{r} - \hat{\zeta}^{\pi}_{r} \le 1$ and thus
\begin{equation}
\varepsilon \hat{\lambda}^{+}_{a} \ \ \overset{\eqref{eqn:auxiliary upper flow adjustment}}{=} \ \ \varepsilon \sum_{\pi \in \Pi} \sum_{r \in \hat{\mathcal{R}}^{\pi}_{a}} D^{\pi} \left(\hat{\eta}^{\pi}_{r} - \hat{\zeta}^{\pi}_{r}\right)
\ \ < \ \ \hat{z} \frac{\sum_{\pi \in \Pi} \sum_{r \in \hat{\mathcal{R}}^{\pi}_{a}} D^{\pi}}{\sum_{a' \in \mathcal{A}_{0} \cup \mathcal{A}_{+}} \sum_{\pi \in \Pi} \sum_{r \in \hat{\mathcal{R}}^{\pi}_{a'}} D^{\pi}}
\label{eqn:error bound a in A0}
\end{equation}

Consider any $a \in \mathcal{A}_{+}$.
It follows from $t_{a}^{-}\big(x_{a}(\hat{q})\big) < t_{a}^{+}\big(x_{a}(\hat{q})\big)$ and complementary slackness that $\hat{\lambda}^{-}_{a} \hat{\lambda}^{+}_{a} = 0$.
First consider $a \in \mathcal{A}_{++} \defi \left\{a \in \mathcal{A}_{+} \, : \, \sum_{\pi \in \Pi} \sum_{r \in \hat{\mathcal{R}}^{\pi}_{a}} D^{\pi} \left(\hat{\eta}^{\pi}_{r} - \hat{\zeta}^{\pi}_{r}\right) \ge 0\right\}$.
Then it holds that
\begin{equation}
\alpha \left(\hat{\lambda}^{+}_{a} - \hat{\lambda}^{-}_{a}\right) \ \ \overset{\eqref{eqn:auxiliary total flow adjustment}}{=} \ \ \alpha \sum_{\pi \in \Pi} \sum_{r \in \hat{\mathcal{R}}^{\pi}_{a}} D^{\pi} \left(\hat{\eta}^{\pi}_{r} - \hat{\zeta}^{\pi}_{r}\right)
\ \ = \ \ \sum_{\pi \in \Pi} \sum_{r \in \hat{\mathcal{R}}^{\pi}_{a}} D^{\pi} \left(\check{q}^{\pi}_{r}(\alpha) - \hat{q}^{\pi}_{r}\right)
\ \ = \ \ x_{a}(\check{q}(\alpha)) - x_{a}(\hat{q})
\ \ \ge \ \ 0
\label{eqn:flow adjustment a in A++}
\end{equation}
Thus $\hat{\lambda}^{-}_{a} = 0$ and $x_{a}(\check{q}(\alpha)) \ge x_{a}(\hat{q}) \ge 0$ for all $\alpha \ge 0$.
Let $A_{a} > 0$ be such that $t_{a}\big(x_{a}(\check{q}(\alpha))\big) \le t_{a}^{+}\big(x_{a}(\hat{q})\big) + \varepsilon$ for all $\alpha \in [0,A_{a}]$.
Then, for all $\alpha \in [0,A_{a}]$ it holds that
\begin{align}
\hat{\psi}_{a}\big(\check{q}(\alpha)\big) - \hat{\psi}_{a}\big(\hat{q}\big)
\ \ & \defi \ \ \int_{0}^{x_{a}(\check{q}(\alpha))} t_{a}(y) \, \diff y - \int_{0}^{x_{a}(\hat{q})} t_{a}(y) \, \diff y
\ \ = \ \ \int_{x_{a}(\hat{q})}^{x_{a}(\check{q}(\alpha))} t_{a}(y) \, \diff y \nonumber \\
& \le \ \ \left[t_{a}^{+}\big(x_{a}(\hat{q})\big) + \varepsilon\right] \left[x_{a}(\check{q}(\alpha)) - x_{a}(\hat{q})\right]
\ \ \overset{\eqref{eqn:flow adjustment a in A++}}{=} \ \ \left[t_{a}^{+}\big(x_{a}(\hat{q})\big) + \varepsilon\right] \alpha \hat{\lambda}^{+}_{a}
\label{eqn:change bound a in A++}
\end{align}
It follows from~\eqref{eqn:auxiliary total deviation} and~\eqref{eqn:auxiliary sign restrictions} that $\hat{\eta}^{\pi}_{r} - \hat{\zeta}^{\pi}_{r} \le 1$ and thus
\begin{equation}
\varepsilon \hat{\lambda}^{+}_{a} \ \ = \ \ \varepsilon \left(\hat{\lambda}^{+}_{a} - \hat{\lambda}^{-}_{a}\right)
\ \ \overset{\eqref{eqn:auxiliary total flow adjustment}}{=} \ \ \varepsilon \sum_{\pi \in \Pi} \sum_{r \in \hat{\mathcal{R}}^{\pi}_{a}} D^{\pi} \left(\hat{\eta}^{\pi}_{r} - \hat{\zeta}^{\pi}_{r}\right)
\ \ < \ \ \hat{z} \frac{\sum_{\pi \in \Pi} \sum_{r \in \hat{\mathcal{R}}^{\pi}_{a}} D^{\pi}}{\sum_{a' \in \mathcal{A}_{0} \cup \mathcal{A}_{+}} \sum_{\pi \in \Pi} \sum_{r \in \hat{\mathcal{R}}^{\pi}_{a'}} D^{\pi}}
\label{eqn:error bound a in A++}
\end{equation}
Next, consider $a \in \mathcal{A}_{+-} \defi \left\{a \in \mathcal{A}_{+} \, : \, \sum_{\pi \in \Pi} \sum_{r \in \hat{\mathcal{R}}^{\pi}_{a}} D^{\pi} \left(\hat{\eta}^{\pi}_{r} - \hat{\zeta}^{\pi}_{r}\right) < 0\right\}$.
Then it holds that
\begin{equation}
\alpha \left(\hat{\lambda}^{+}_{a} - \hat{\lambda}^{-}_{a}\right) \ \ \overset{\eqref{eqn:auxiliary total flow adjustment}}{=} \ \ \alpha \sum_{\pi \in \Pi} \sum_{r \in \hat{\mathcal{R}}^{\pi}_{a}} D^{\pi} \left(\hat{\eta}^{\pi}_{r} - \hat{\zeta}^{\pi}_{r}\right)
\ \ = \ \ \sum_{\pi \in \Pi} \sum_{r \in \hat{\mathcal{R}}^{\pi}_{a}} D^{\pi} \left(\check{q}^{\pi}_{r}(\alpha) - \hat{q}^{\pi}_{r}\right)
\ \ = \ \ x_{a}(\check{q}(\alpha)) - x_{a}(\hat{q})
\ \ \le \ \ 0
\label{eqn:flow adjustment a in A+-}
\end{equation}
Thus $\hat{\lambda}^{+}_{a} = 0$.
Note that $x_{a}(\hat{q}) > 0$.
Also note that $\lim_{\alpha \downarrow 0} x_{a}(\check{q}(\alpha)) = x_{a}(\hat{q})$ and $\lim_{\alpha \downarrow 0} t_{a}\big(x_{a}(\check{q}(\alpha))\big) \ge t_{a}^{-}\big(x_{a}(\hat{q})\big)$.
Let $A_{a} > 0$ be such that $t_{a}\big(x_{a}(\check{q}(\alpha))\big) \ge t_{a}^{-}\big(x_{a}(\hat{q})\big) - \varepsilon$ and $x_{a}(\check{q}(\alpha)) \ge 0$ for all $\alpha \in [0,A_{a}]$.
Then, for all $\alpha \in [0,A_{a}]$ it holds that
\begin{align}
\hat{\psi}_{a}\big(\check{q}(\alpha)\big) - \hat{\psi}_{a}\big(\hat{q}\big)
\ \ & \defi \ \ \int_{0}^{x_{a}(\check{q}(\alpha))} t_{a}(y) \, \diff y - \int_{0}^{x_{a}(\hat{q})} t_{a}(y) \, \diff y
\ \ = \ \ - \int_{x_{a}(\check{q}(\alpha))}^{x_{a}(\hat{q})} t_{a}(y) \, \diff y \nonumber \\
& \le \ \ \left[t_{a}^{-}\big(x_{a}(\hat{q})\big) - \varepsilon\right] \left[x_{a}(\check{q}(\alpha)) - x_{a}(\hat{q})\right]
\ \ \overset{\eqref{eqn:flow adjustment a in A+-}}{=} \ \ - \left[t_{a}^{-}\big(x_{a}(\hat{q})\big) - \varepsilon\right] \alpha \hat{\lambda}^{-}_{a}
\label{eqn:change bound a in A+-}
\end{align}
It follows from~\eqref{eqn:auxiliary total deviation} and~\eqref{eqn:auxiliary sign restrictions} that $\hat{\zeta}^{\pi}_{r} - \hat{\eta}^{\pi}_{r} \le 1$ and thus
\begin{equation}
\varepsilon \hat{\lambda}^{-}_{a} \ \ = \ \ \varepsilon \left(\hat{\lambda}^{-}_{a} - \hat{\lambda}^{+}_{a}\right)
\ \ \overset{\eqref{eqn:auxiliary total flow adjustment}}{=} \ \ \varepsilon \sum_{\pi \in \Pi} \sum_{r \in \hat{\mathcal{R}}^{\pi}_{a}} D^{\pi} \left(\hat{\zeta}^{\pi}_{r} - \hat{\eta}^{\pi}_{r}\right)
\ \ < \ \ \hat{z} \frac{\sum_{\pi \in \Pi} \sum_{r \in \hat{\mathcal{R}}^{\pi}_{a}} D^{\pi}}{\sum_{a' \in \mathcal{A}_{0} \cup \mathcal{A}_{+}} \sum_{\pi \in \Pi} \sum_{r \in \hat{\mathcal{R}}^{\pi}_{a'}} D^{\pi}}
\label{eqn:error bound a in A+-}
\end{equation}

Consider any $a \in \mathcal{A}_{=}$.
First consider $a \in \mathcal{A}_{=+} \defi \left\{a \in \mathcal{A}_{=} \, : \, \sum_{\pi \in \Pi} \sum_{r \in \hat{\mathcal{R}}^{\pi}_{a}} D^{\pi} \left(\hat{\eta}^{\pi}_{r} - \hat{\zeta}^{\pi}_{r}\right) \ge 0\right\}$.
Then it holds that
\begin{equation}
- \alpha \hat{\lambda}^{=}_{a} \ \ \overset{\eqref{eqn:auxiliary equal flow adjustment}}{=} \ \ \alpha \sum_{\pi \in \Pi} \sum_{r \in \hat{\mathcal{R}}^{\pi}_{a}} D^{\pi} \left(\hat{\eta}^{\pi}_{r} - \hat{\zeta}^{\pi}_{r}\right)
\ \ = \ \ \sum_{\pi \in \Pi} \sum_{r \in \hat{\mathcal{R}}^{\pi}_{a}} D^{\pi} \left(\check{q}^{\pi}_{r}(\alpha) - \hat{q}^{\pi}_{r}\right)
\ \ = \ \ x_{a}(\check{q}(\alpha)) - x_{a}(\hat{q})
\ \ \ge \ \ 0
\label{eqn:flow adjustment a in A=+}
\end{equation}
Thus $\hat{\lambda}^{=}_{a} \le 0$ and $x_{a}(\check{q}(\alpha)) \ge x_{a}(\hat{q}) \ge 0$ for all $\alpha \ge 0$.
Let $A_{a} > 0$ be such that $t_{a}\big(x_{a}(\check{q}(\alpha))\big) \le t_{a}\big(x_{a}(\hat{q})\big) + \varepsilon$ for all $\alpha \in [0,A_{a}]$.
Then, for all $\alpha \in [0,A_{a}]$ it holds that
\begin{align}
\hat{\psi}_{a}\big(\check{q}(\alpha)\big) - \hat{\psi}_{a}\big(\hat{q}\big)
\ \ & \defi \ \ \int_{0}^{x_{a}(\check{q}(\alpha))} t_{a}(y) \, \diff y - \int_{0}^{x_{a}(\hat{q})} t_{a}(y) \, \diff y
\ \ = \ \ \int_{x_{a}(\hat{q})}^{x_{a}(\check{q}(\alpha))} t_{a}(y) \, \diff y \nonumber \\
& \le \ \ \left[t_{a}\big(x_{a}(\hat{q})\big) + \varepsilon\right] \left[x_{a}(\check{q}(\alpha)) - x_{a}(\hat{q})\right]
\ \ \overset{\eqref{eqn:flow adjustment a in A=+}}{=} \ \ - \left[t_{a}\big(x_{a}(\hat{q})\big) + \varepsilon\right] \alpha \hat{\lambda}^{=}_{a}
\label{eqn:change bound a in A=+}
\end{align}
It follows from~\eqref{eqn:auxiliary total deviation} and~\eqref{eqn:auxiliary sign restrictions} that $\hat{\eta}^{\pi}_{r} - \hat{\zeta}^{\pi}_{r} \le 1$ and thus
\begin{equation}
- \varepsilon \hat{\lambda}^{=}_{a} \ \ \overset{\eqref{eqn:auxiliary equal flow adjustment}}{=} \ \ \varepsilon \sum_{\pi \in \Pi} \sum_{r \in \hat{\mathcal{R}}^{\pi}_{a}} D^{\pi} \left(\hat{\eta}^{\pi}_{r} - \hat{\zeta}^{\pi}_{r}\right)
\ \ < \ \ \hat{z} \frac{\sum_{\pi \in \Pi} \sum_{r \in \hat{\mathcal{R}}^{\pi}_{a}} D^{\pi}}{\sum_{a' \in \mathcal{A}_{0} \cup \mathcal{A}_{+}} \sum_{\pi \in \Pi} \sum_{r \in \hat{\mathcal{R}}^{\pi}_{a'}} D^{\pi}}
\label{eqn:error bound a in A=+}
\end{equation}
Next consider $a \in \mathcal{A}_{=-} \defi \left\{a \in \mathcal{A}_{=} \, : \, \sum_{\pi \in \Pi} \sum_{r \in \hat{\mathcal{R}}^{\pi}_{a}} D^{\pi} \left(\hat{\eta}^{\pi}_{r} - \hat{\zeta}^{\pi}_{r}\right) < 0\right\}$.
Then it holds that
\begin{equation}
- \alpha \hat{\lambda}^{=}_{a} \ \ \overset{\eqref{eqn:auxiliary equal flow adjustment}}{=} \ \ \alpha \sum_{\pi \in \Pi} \sum_{r \in \hat{\mathcal{R}}^{\pi}_{a}} D^{\pi} \left(\hat{\eta}^{\pi}_{r} - \hat{\zeta}^{\pi}_{r}\right)
\ \ = \ \ \sum_{\pi \in \Pi} \sum_{r \in \hat{\mathcal{R}}^{\pi}_{a}} D^{\pi} \left(\check{q}^{\pi}_{r}(\alpha) - \hat{q}^{\pi}_{r}\right)
\ \ = \ \ x_{a}(\check{q}(\alpha)) - x_{a}(\hat{q})
\ \ \le \ \ 0
\label{eqn:flow adjustment a in A=-}
\end{equation}
Thus $\hat{\lambda}^{=}_{a} \ge 0$.
Note that $x_{a}(\hat{q}) > 0$.
Let $A_{a} > 0$ be such that $t_{a}\big(x_{a}(\check{q}(\alpha))\big) \ge t_{a}\big(x_{a}(\hat{q})\big) - \varepsilon$ for all $\alpha \in [0,A_{a}]$.
Then, for all $\alpha \in [0,A_{a}]$ it holds that
\begin{align}
\hat{\psi}_{a}\big(\check{q}(\alpha)\big) - \hat{\psi}_{a}\big(\hat{q}\big)
\ \ & \defi \ \ \int_{0}^{x_{a}(\check{q}(\alpha))} t_{a}(y) \, \diff y - \int_{0}^{x_{a}(\hat{q})} t_{a}(y) \, \diff y
\ \ = \ \ - \int_{x_{a}(\check{q}(\alpha))}^{x_{a}(\hat{q})} t_{a}(y) \, \diff y \nonumber \\
& \le \ \ \left[t_{a}\big(x_{a}(\hat{q})\big) - \varepsilon\right] \left[x_{a}(\check{q}(\alpha)) - x_{a}(\hat{q})\right]
\ \ \overset{\eqref{eqn:flow adjustment a in A=-}}{=} \ \ - \left[t_{a}\big(x_{a}(\hat{q})\big) - \varepsilon\right] \alpha \hat{\lambda}^{=}_{a}
\label{eqn:change bound a in A=-}
\end{align}
It follows from~\eqref{eqn:auxiliary total deviation} and~\eqref{eqn:auxiliary sign restrictions} that $\hat{\zeta}^{\pi}_{r} - \hat{\eta}^{\pi}_{r} \le 1$ and thus
\begin{equation}
\varepsilon \hat{\lambda}^{=}_{a} \ \ \overset{\eqref{eqn:auxiliary equal flow adjustment}}{=} \ \ \varepsilon \sum_{\pi \in \Pi} \sum_{r \in \hat{\mathcal{R}}^{\pi}_{a}} D^{\pi} \left(\hat{\zeta}^{\pi}_{r} - \hat{\eta}^{\pi}_{r}\right)
\ \ < \ \ \hat{z} \frac{\sum_{\pi \in \Pi} \sum_{r \in \hat{\mathcal{R}}^{\pi}_{a}} D^{\pi}}{\sum_{a' \in \mathcal{A}_{0} \cup \mathcal{A}_{+}} \sum_{\pi \in \Pi} \sum_{r \in \hat{\mathcal{R}}^{\pi}_{a'}} D^{\pi}}
\label{eqn:error bound a in A=-}
\end{equation}

Let $A \defi \min\left\{A_{a} \, : \, a \in \mathcal{A}\right\} > 0$.
Note that $x_{a}(\check{q}(\alpha)) \ge 0$ for all $\alpha \in [0,A]$ and all $a \in \mathcal{A}$, and thus $\check{q}(\alpha) \in Q$ for all $\alpha \in [0,A]$.
Then, for all $\alpha \in (0,A]$ it holds that
\begin{align*}
\hat{\psi}\big(\check{q}(\alpha)\big) - \hat{\psi}\big(\hat{q}\big)
\ \ & = \ \ \sum_{a \in \mathcal{A}} \left[\hat{\psi}_{a}\big(\check{q}(\alpha)\big) - \hat{\psi}_{a}\big(\hat{q}\big)\right] \\
& = \ \ \sum_{a \in \mathcal{A}_{0}} \left[\hat{\psi}_{a}\big(\check{q}(\alpha)\big) - \hat{\psi}_{a}\big(\hat{q}\big)\right] + \sum_{a \in \mathcal{A}_{++}} \left[\hat{\psi}_{a}\big(\check{q}(\alpha)\big) - \hat{\psi}_{a}\big(\hat{q}\big)\right] + \sum_{a \in \mathcal{A}_{+-}} \left[\hat{\psi}_{a}\big(\check{q}(\alpha)\big) - \hat{\psi}_{a}\big(\hat{q}\big)\right] \\
& \hspace{10mm} {} + \sum_{a \in \mathcal{A}_{=+}} \left[\hat{\psi}_{a}\big(\check{q}(\alpha)\big) - \hat{\psi}_{a}\big(\hat{q}\big)\right] + \sum_{a \in \mathcal{A}_{=-}} \left[\hat{\psi}_{a}\big(\check{q}(\alpha)\big) - \hat{\psi}_{a}\big(\hat{q}\big)\right] \\
& \overset{\eqref{eqn:change bound a in A0}\eqref{eqn:change bound a in A++}\eqref{eqn:change bound a in A+-}\eqref{eqn:change bound a in A=+}\eqref{eqn:change bound a in A=-}}{\le} \ \ \sum_{a \in \mathcal{A}_{0}} \left[t_{a}^{+}\big(x_{a}(\hat{q})\big) + \varepsilon\right] \alpha \hat{\lambda}^{+}_{a} \\
& \hspace{40mm} {} + \sum_{a \in \mathcal{A}_{++}} \left[t_{a}^{+}\big(x_{a}(\hat{q})\big) + \varepsilon\right] \alpha \hat{\lambda}^{+}_{a} - \sum_{a \in \mathcal{A}_{+-}} \left[t_{a}^{-}\big(x_{a}(\hat{q})\big) - \varepsilon\right] \alpha \hat{\lambda}^{-}_{a} \\
& \hspace{40mm} {} - \sum_{a \in \mathcal{A}_{=+}} \left[t_{a}\big(x_{a}(\hat{q})\big) + \varepsilon\right] \alpha \hat{\lambda}^{=}_{a} - \sum_{a \in \mathcal{A}_{=-}} \left[t_{a}\big(x_{a}(\hat{q})\big) - \varepsilon\right] \alpha \hat{\lambda}^{=}_{a} \\
& \overset{\eqref{eqn:error bound a in A0}\eqref{eqn:error bound a in A++}\eqref{eqn:error bound a in A+-}\eqref{eqn:error bound a in A=+}\eqref{eqn:error bound a in A=-}}{<} \ \ \alpha \sum_{a \in \mathcal{A}_{0}} t_{a}^{+}\big(x_{a}(\hat{q})\big) \hat{\lambda}^{+}_{a} \\
& \hspace{40mm} {} + \alpha \sum_{a \in \mathcal{A}_{++}} t_{a}^{+}\big(x_{a}(\hat{q})\big) \hat{\lambda}^{+}_{a} - \alpha \sum_{a \in \mathcal{A}_{+-}} t_{a}^{-}\big(x_{a}(\hat{q})\big) \hat{\lambda}^{-}_{a} \\
& \hspace{40mm} {} - \alpha \sum_{a \in \mathcal{A}_{=}} t_{a}\big(x_{a}(\hat{q})\big) \hat{\lambda}^{=}_{a} + \alpha \sum_{a \in \mathcal{A}} \hat{z} \frac{\sum_{\pi \in \Pi} \sum_{r \in \hat{\mathcal{R}}^{\pi}_{a}} D^{\pi}}{\sum_{a' \in \mathcal{A}} \sum_{\pi \in \Pi} \sum_{r \in \hat{\mathcal{R}}^{\pi}_{a'}} D^{\pi}} \\
& = \ \ \alpha \sum_{a \in \mathcal{A}_{0}} t_{a}^{+}\big(x_{a}(\hat{q})\big) \hat{\lambda}^{+}_{a} + \alpha \sum_{a \in \mathcal{A}_{+}} t_{a}^{+}\big(x_{a}(\hat{q})\big) \hat{\lambda}^{+}_{a} - \alpha \sum_{a \in \mathcal{A}_{+}} t_{a}^{-}\big(x_{a}(\hat{q})\big) \hat{\lambda}^{-}_{a} - \alpha \sum_{a \in \mathcal{A}_{=}} t_{a}\big(x_{a}(\hat{q})\big) \hat{\lambda}^{=}_{a} + \alpha \hat{z} \\
& = \ \ \alpha \sum_{a \in \mathcal{A}_{0} \cup \mathcal{A}_{+}} t_{a}^{+}\big(x_{a}(\hat{q})\big) \hat{\lambda}^{+}_{a} - \alpha \sum_{a \in \mathcal{A}_{+}} t_{a}^{-}\big(x_{a}(\hat{q})\big) \hat{\lambda}^{-}_{a} - \alpha \sum_{a \in \mathcal{A}_{=}} t_{a}\big(x_{a}(\hat{q})\big) \hat{\lambda}^{=}_{a} \\
& \hspace{5mm} {} + \alpha \left\{\sum_{\pi \in \Pi} \sum_{r \in \hat{\mathcal{R}}^{\pi}} \hat{g}^{\pi}_{r} \left(\hat{\eta}^{\pi}_{r} - \hat{\zeta}^{\pi}_{r}\right) + \sum_{a \in \mathcal{A}_{+}} t_{a}^{-}(x_{a}(\hat{q})) \hat{\lambda}^{-}_{a} - \sum_{a \in \mathcal{A}_{0} \cup \mathcal{A}_{+}} t_{a}^{+}(x_{a}(\hat{q})) \hat{\lambda}^{+}_{a} + \sum_{a \in \mathcal{A}_{=}} t_{a}(x_{a}(\hat{q})) \hat{\lambda}^{=}_{a}\right\} \\
& = \ \ \alpha \; \hat{g}^{\top} \!\left(\hat{\eta} - \hat{\zeta}\right)
\ \ = \ \ \hat{g}^{\top} \big(\check{q}(\alpha) - \hat{q}\big)
\end{align*}
and thus $\hat{g} \notin \partial \hat{\psi}(\hat{q})$.
Therefore, it has been established that
\[
\partial \hat{\psi}(q) \ \ = \ \ \left\{g = \left(g^{\pi}_{r}, \pi \in \Pi, r \in \hat{\mathcal{R}}^{\pi}\right) \; : \; g^{\pi}_{r} = D^{\pi} \sum_{a \in \mathcal{A}^{\pi}_{r}} \tau_{a}, \; \tau_{a} \in \left[t_{a}^{-}(x_{a}(q)), t_{a}^{+}(x_{a}(q))\right] \ \forall \; a \in \mathcal{A}\right\}
\]
and the result of the lemma follows.
\Halmos
\endproof

\section{Proof for Result in Section~\ref{sec:optimality conditions}}
\label{app:proofs for optimality conditions}
\label{app:convex-reformulation}

\begin{manuallemma}{\ref{lem:KKT equivalence}}
A vector $q^{*} \in \prod_{\pi \in \Pi} \Delta_{\hat{\mathcal{R}}^{\pi}}$ is optimal for problem~\eqref{eqn:q reformulated BiTAP objective} if and only if there exist $\tau^* = (\tau^*_{a}, a \in \mathcal{A})$ with $\tau^*_{a} \in \left[t_{a}^{-}(x_{a}(q^*)), t_{a}^{+}(x_{a}(q^*))\right]$, $w^{*} = \left(w^{\pi*}_{r}, \pi \in \Pi, r \in \hat{\mathcal{R}}^{\pi}\right)$ with $w^{\pi*}_{r} \in \left(G^{\pi}\right)^{-1}\left(L^{\pi}_{r}(q^{\pi*})\right)$, and $\lambda^{*} = (\lambda^{\pi*}, \pi \in \Pi)$, such that the following holds:
\begin{align*}
& q^{\pi*}_{r} \left\{D^{\pi} \left[\sum_{a \in \mathcal{A}^{\pi}_{r}} \tau^*_{a} + \sum_{\{r' \in \hat{\mathcal{R}}^{\pi} \, : \, r \preceq r', r'+1 \in \hat{\mathcal{R}}^{\pi}\}} \left(S^{\pi}_{r'} - S^{\pi}_{r'+1}\right) w^{\pi*}_{r'}\right] + \lambda^{\pi*}\right\} \ \ = \ \ 0 \qquad \forall \ \pi \in \Pi, \ \forall \ r \in \hat{\mathcal{R}}^{\pi} \\
& D^{\pi} \left[\sum_{a \in \mathcal{A}^{\pi}_{r}} \tau^*_{a} + \sum_{\{r' \in \hat{\mathcal{R}}^{\pi} \, : \, r \preceq r', r'+1 \in \hat{\mathcal{R}}^{\pi}\}} \left(S^{\pi}_{r'} - S^{\pi}_{r'+1}\right) w^{\pi*}_{r'}\right] + \lambda^{\pi*} \ \ \ge \ \ 0 \qquad \forall \ \pi \in \Pi, \ \forall \ r \in \hat{\mathcal{R}}^{\pi}
\end{align*}
\end{manuallemma}

\proof{Proof.}
Convexity of $\Upsilon : \prod_{\pi \in \Pi} \Delta_{\hat{\mathcal{R}}^{\pi}} \mapsto \mathbb{R}$ follows from Lemmas~\ref{lem:phi technical properties} and~\ref{lem:psi technical properties}.
It follows from convexity of $\Upsilon$ and $\prod_{\pi \in \Pi} \Delta_{\hat{\mathcal{R}}^{\pi}}$ that $q^{*} \in \prod_{\pi \in \Pi} \Delta_{\hat{\mathcal{R}}^{\pi}}$ is optimal for problem~\eqref{eqn:q reformulated BiTAP objective} if and only if $0 \in \partial \Upsilon(q^*)$.
Note that $\Upsilon(q) = \psi(q) + \phi(q) + \sum_{\pi \in \Pi} \delta^{\pi}_{\Delta_{\hat{\mathcal{R}}^{\pi}}}(q^{\pi})$ and
\[
\text{ri}\left(\text{dom}\big(\psi\big)\right) \, \cap \, \text{ri}\left(\text{dom}\big(\phi\big)\right) \, \cap \, \text{ri}\left(\text{dom}\bigg(\sum_{\pi \in \Pi} \delta^{\pi}_{\Delta_{\hat{\mathcal{R}}^{\pi}}}\bigg)\right) \ \ = \ \ \prod_{\pi \in \Pi} \Delta_{\hat{\mathcal{R}}^{\pi}}^{\circ} \ \ \neq \ \ \varnothing
\]
and thus it follows from the Moreau-Rockafellar Theorem \citep[Theorem 23.8]{rockafellar1997convex} that $\partial \Upsilon(q) = \partial \psi(q) + \partial \phi(q) + \partial \sum_{\pi \in \Pi} \delta^{\pi}_{\Delta_{\hat{\mathcal{R}}^{\pi}}}(q^{\pi})$.
Hence, it follows from Lemmas~\ref{lem:phi technical properties} and~\ref{lem:psi technical properties} that $0 \in \partial \Upsilon(q^*)$ if and only if there exist $\tau^* = (\tau^*_{a}, a \in \mathcal{A})$ with $\tau^*_{a} \in \left[t_{a}^{-}(x_{a}(q^*)), t_{a}^{+}(x_{a}(q^*))\right]$, $w^{*} = \left(w^{\pi*}_{r}, \pi \in \Pi, r \in \hat{\mathcal{R}}^{\pi}\right)$ with $w^{\pi*}_{r} \in \left(G^{\pi}\right)^{-1}\left(L^{\pi}_{r}(q^{\pi*})\right)$, and $\nu^{\pi*} = (\nu^{\pi*}_{r}, r \in \hat{\mathcal{R}}^{\pi}) \in \mathcal{N}_{\Delta_{\hat{\mathcal{R}}^{\pi}}}(q^{\pi*})$ for all $\pi \in \Pi$,  such that
\begin{equation}
0 \ \ = \ \ D^{\pi} \left[\sum_{a \in \mathcal{A}^{\pi}_{r}} \tau^*_{a} + \sum_{\{r' \in \hat{\mathcal{R}}^{\pi} \, : \, r \preceq r', r'+1 \in \hat{\mathcal{R}}^{\pi}\}} \left(S^{\pi}_{r'} - S^{\pi}_{r'+1}\right) w^{\pi*}_{r'}\right] + \nu^{\pi*}_{r} \qquad \forall \ \pi \in \Pi, \ \forall \ r \in \hat{\mathcal{R}}^{\pi}
\label{eqn:optimality condition 1}
\end{equation}
Note that $\nu^{\pi*} \in \mathcal{N}_{\Delta_{\hat{\mathcal{R}}^{\pi}}}(q^{\pi*})$ if and only if there are $\lambda^{\pi*}$ and $\mu^{\pi*}_{r} \ge 0$ for all $r \in \hat{\mathcal{R}}^{\pi}$, such that $\nu^{\pi*}_{r} = \lambda^{\pi*} - \mu^{\pi*}_{r}$ and $\mu^{\pi*}_{r} q^{\pi*}_{r} = 0$ for all $r \in \hat{\mathcal{R}}^{\pi}$.
Therefore, in both the cases $q^{\pi*}_{r} = 0$ and $q^{\pi*}_{r} > 0$, \eqref{eqn:optimality condition 1} is equivalent to \eqref{eqn:q reformulated BiTAP optimality1}--\eqref{eqn:q reformulated BiTAP optimality3}.
\Halmos
\endproof

\section{Proofs for Results in Section~\ref{sec:optimality to equilibrium}}
\label{app:proofs for equilibrium conditions}
\label{app:threshold-proofs}

\begin{manuallemma}{\ref{lem:optimal thresholds}}
Assume that $t_{a} : \mathbb{R}_{+} \mapsto \mathbb{R}_{+}$ is continuous for all $a \in \mathcal{A}$.
Consider any optimal solution $q$ for problem~\eqref{eqn:q reformulated BiTAP objective}, and the resulting output $w$ and $\mathcal{R}^{\pi}_{\min}$, $\pi \in \Pi$, of Algorithm~\ref{algo:threshold}.
Then the following hold:
\begin{enumerate}[label=(\arabic*), ref=(\arabic*)]
\item
The thresholds $w$ satisfy $\underline{w}^{\pi} = w^{\pi}_{0} \le w^{\pi}_{1} \le \cdots \le w^{\pi}_{\hat{R}^{\pi}} = \overline{w}^{\pi}$ for all $\pi \in \Pi$.
\item
Every path with positive flow is an efficient path, that is, $\left\{r \in \hat{\mathcal{R}}^{\pi} \, : \, q^{\pi}_{r} > 0\right\} \subset \mathcal{R}^{\pi}_{\min}$.
\item
The thresholds $w$ satisfy $w^{\pi}_{r} \in \left(G^{\pi}\right)^{-1}\left(L^{\pi}_{r}(q^{\pi})\right)$ for all $\pi \in \Pi$, $r \in \hat{\mathcal{R}}^{\pi}$.
\item
For every $\pi \in \Pi$ and every $r_{\ell} \in \mathcal{R}^{\pi}_{\min}$ it holds that
\begin{equation*}
T^{\pi}_{r_{\ell}}(x(q)) + w S^{\pi}_{r_{\ell}} \ \ \le \ \ T^{\pi}_{r}(x(q)) + w S^{\pi}_{r}
\qquad \forall \ r \in \hat{\mathcal{R}}^{\pi}, \ \forall \ w \in [w^{\pi}_{r_{\ell-1}},w^{\pi}_{r_{\ell}}]
\end{equation*}
where $w^{\pi}_{r_{0}} = w^{\pi}_{0} = \underline{w}^{\pi}$.
\end{enumerate}
\end{manuallemma}

\proof{Proof.}
Consider any $\pi \in \Pi$.
Let $\hat{\mathcal{R}}^{\pi}$ be indexed with $r = 1,\ldots,\hat{R}^{\pi}$, let $\mathcal{R}^{\pi}_{\min} = \{r_{1},\ldots,r_{L}\}$, and let $p^{\pi}_{r} = L^{\pi}_{r}(q^{\pi})$.
All the properties will be established together by induction.

\paragraph{Initialization.}$\mbox{}$

\begin{enumerate}[label=(\arabic*), ref=(\arabic*), leftmargin=7mm]
\item
Note that $\underline{w}^{\pi} = w^{\pi}_{0} = \cdots = w^{\pi}_{r_{1}-1}$, and $w^{\pi}_{r_{\ell}} = w^{\pi}_{r_{\ell}+1} = \cdots = w^{\pi}_{r_{\ell+1}-1}$ for $\ell = 1,\ldots,L-1$, and $w^{\pi}_{r_{L}} = w^{\pi}_{r_{L}+1} = \cdots = w^{\pi}_{\hat{R}^{\pi}} = \overline{w}^{\pi}$.
Thus, to establish \ref{item:thresholds ordered}, it remains to show that $\underline{w}^{\pi} \le w^{\pi}_{r_{\ell}} \le w^{\pi}_{r_{\ell+1}} \le \overline{w}^{\pi}$ for $\ell = 1,\ldots,L-1$.
Next we show that $\underline{w}^{\pi} \le w^{\pi}_{r_{1}}$.
In Algorithm~\ref{algo:threshold}, $w^{\pi}_{r_{1}}$ is determined by~\eqref{eqn:threshold >} or~\eqref{eqn:threshold ==} or~\eqref{eqn:threshold =>} or~\eqref{eqn:threshold =0}.
First consider the case in which $w^{\pi}_{r_{1}}$ is determined by~\eqref{eqn:threshold >} or~\eqref{eqn:threshold =>}.
If $w^{\min}_{1} > \overline{w}^{\pi}$, then $L = 1$ and $\underline{w}^{\pi} \le w^{\pi}_{r_{1}} = \overline{w}^{\pi}$.
Otherwise, $w^{\pi}_{r_{1}} = w^{\min}_{1} \le \overline{w}^{\pi}$, $S^{\pi}_{r_{1}} > S^{\pi}_{r_{2}}$, and $T^{\pi}_{r_{1}}(x(q)) + w^{\pi}_{r_{1}} S^{\pi}_{r_{1}} = T^{\pi}_{r_{2}}(x(q)) + w^{\pi}_{r_{1}} S^{\pi}_{r_{2}}$.
Suppose that $\underline{w}^{\pi} > w^{\pi}_{r_{1}}$.
Then
\[
\left(S^{\pi}_{r_{1}} - S^{\pi}_{r_{2}}\right) \underline{w}^{\pi} \ \ > \ \ \left(S^{\pi}_{r_{1}} - S^{\pi}_{r_{2}}\right) w^{\pi}_{r_{1}} \ \ = \ \ T^{\pi}_{r_{2}}(x(q)) - T^{\pi}_{r_{1}}(x(q))
\quad \Rightarrow \quad
T^{\pi}_{r_{1}}(x(q)) + \underline{w}^{\pi} S^{\pi}_{r_{1}} \ \ > \ \ T^{\pi}_{r_{2}}(x(q)) + \underline{w}^{\pi} S^{\pi}_{r_{2}}
\]
which contradicts~\eqref{eqn:efficient path 1}.
Thus $\underline{w}^{\pi} \le w^{\pi}_{r_{1}} \le \overline{w}^{\pi}$.
If $w^{\pi}_{r_{1}}$ is determined by~\eqref{eqn:threshold ==} or~\eqref{eqn:threshold =0}, then it follows immediately that $\underline{w}^{\pi} \le w^{\pi}_{r_{1}} \le \overline{w}^{\pi}$.

\item
We will show that for any $r \in \hat{\mathcal{R}}^{\pi} \setminus \mathcal{R}^{\pi}_{\min}$ it holds that $q^{\pi}_{r} = 0$.
First consider any $r < r_{1}$.
Then $S^{\pi}_{r} \ge S^{\pi}_{r_{1}}$ and it follows from~\eqref{eqn:efficient path 1} that $T^{\pi}_{r}(x(q)) + \underline{w}^{\pi} S^{\pi}_{r} > T^{\pi}_{r_{1}}(x(q)) + \underline{w}^{\pi} S^{\pi}_{r_{1}}$.
Suppose that $q^{\pi}_{r} > 0$.
Then for all $r' \ge r$ it holds that $p^{\pi}_{r'} \ge q^{\pi}_{r} > 0$, and thus any $\hat{w} \in \left(G^{\pi}\right)^{-1}\left(p^{\pi}_{r'}\right)$ satisfies $\hat{w} \ge \underline{w}^{\pi}$.
Consider the directional derivative of $\Upsilon$ at $q$ in the feasible direction $e^{\pi}_{r_{1}} - e^{\pi}_{r}$:
\begin{align*}
\hspace{-5mm}
& \max\Bigg\{\left(D^{\pi} \left[T^{\pi}_{r_{1}}(x(q)) + \sum_{r'=r_{1}}^{\hat{R}^{\pi}-1} \left(S^{\pi}_{r'} - S^{\pi}_{r'+1}\right) \hat{w}^{\pi}_{r'}\right] + \lambda^{\pi} - \mu^{\pi}_{r_{1}}\right) \\
\hspace{-5mm}
& \hspace{10mm} - \left(D^{\pi} \left[T^{\pi}_{r}(x(q)) + \sum_{r'=r}^{\hat{R}^{\pi}-1} \left(S^{\pi}_{r'} - S^{\pi}_{r'+1}\right) \hat{w}^{\pi}_{r'}\right] + \lambda^{\pi}\right) \; : \; \hat{w}^{\pi}_{1} \in \left(G^{\pi}\right)^{-1}\left(p^{\pi}_{1}\right), \ldots, \hat{w}^{\pi}_{\hat{R}^{\pi}-1} \in \left(G^{\pi}\right)^{-1}\left(p^{\pi}_{\hat{R}^{\pi}-1}\right), \\
\hspace{-5mm}
& \hspace{145mm} \mu^{\pi}_{r_{1}} \ge 0, \; \mu^{\pi}_{r_{1}} q^{\pi}_{r_{1}} = 0\Bigg\} \\
\hspace{-5mm}
& = \ \ \max\Bigg\{D^{\pi} \left[T^{\pi}_{r_{1}}(x(q)) - T^{\pi}_{r}(x(q)) - \sum_{r'=r}^{r_{1}-1} \left(S^{\pi}_{r'} - S^{\pi}_{r'+1}\right) \hat{w}^{\pi}_{r'}\right] \; : \; \hat{w}^{\pi}_{r} \in \left(G^{\pi}\right)^{-1}\left(p^{\pi}_{r}\right), \ldots, \hat{w}^{\pi}_{r_{1}-1} \in \left(G^{\pi}\right)^{-1}\left(p^{\pi}_{r_{1}-1}\right)\Bigg\} \\
\hspace{-5mm}
& \le \ \ D^{\pi} \left[T^{\pi}_{r_{1}}(x(q)) - T^{\pi}_{r}(x(q)) - \sum_{r'=r}^{r_{1}-1} \left(S^{\pi}_{r'} - S^{\pi}_{r'+1}\right) \underline{w}^{\pi}\right] \\
\hspace{-5mm}
& = \ \ D^{\pi} \left[T^{\pi}_{r_{1}}(x(q)) - T^{\pi}_{r}(x(q)) + \left(S^{\pi}_{r_{1}} - S^{\pi}_{r}\right) \underline{w}^{\pi}\right]
\ \ < \ \ 0
\end{align*}
which contradicts the optimality of $q$.
Thus for all $r < r_{1}$, it holds that $q^{\pi}_{r} = 0$.

\item
Note that $p^{\pi}_{0} = 0$ and $w^{\pi}_{0} = \underline{w}^{\pi} \in \left(G^{\pi}\right)^{-1}\left(0\right)$, and thus $w^{\pi}_{0} \in \left(G^{\pi}\right)^{-1}\left(p^{\pi}_{0}\right)$.
By Algorithm~\ref{algo:threshold}, $w^{\pi}_{r} = \underline{w}^{\pi}$ for all $r < r_{1}$.
It was established in \ref{item:positive flow paths are efficient} above that $q^{\pi}_{r} = 0$ for all $r < r_{1}$, and thus $p^{\pi}_{r} = 0$ and $w^{\pi}_{r} \in \left(G^{\pi}\right)^{-1}\left(p^{\pi}_{r}\right)$ for all $r < r_{1}$.
Also, $p^{\pi}_{\hat{R}^{\pi}} = 1$ and $w^{\pi}_{\hat{R}^{\pi}} = \overline{w}^{\pi} \in \left(G^{\pi}\right)^{-1}\left(1\right)$, and thus $w^{\pi}_{\hat{R}^{\pi}} \in \left(G^{\pi}\right)^{-1}\left(p^{\pi}_{\hat{R}^{\pi}}\right)$.

Next we show by contradiction that $w^{\pi}_{r_{1}} \in \left(G^{\pi}\right)^{-1}\left(p^{\pi}_{r_{1}}\right)$.
For each $r \in \{0\} \cup \hat{\mathcal{R}}^{\pi}$, if $p^{\pi}_{r} \in (0,1)$ then let $\left(G^{\pi}\right)^{-1}\left(p^{\pi}_{r}\right) = [\underline{w}^{\pi}_{r},\overline{w}^{\pi}_{r}]$, if $p^{\pi}_{r} = 0$ then let $\left(G^{\pi}\right)^{-1}\left(p^{\pi}_{r}\right) = (-\infty,\overline{w}^{\pi}_{r}]$, and if $p^{\pi}_{r} = 1$ then let $\left(G^{\pi}\right)^{-1}\left(p^{\pi}_{r}\right) = [\underline{w}^{\pi}_{r},\infty)$.

\noindent
\textbf{Case~1}: Suppose that $w^{\pi}_{r_{1}} < \underline{w}^{\pi}_{r_{1}}$. \\
Note that $w^{\pi}_{r_{1}} < \underline{w}^{\pi}_{r_{1}}$ implies that $G^{\pi}(w^{\pi}_{r_{1}}) < p^{\pi}_{r_{1}}$.
Also, it follows from $q^{\pi}_{r} = 0$ for all $r < r_{1}$ and $p^{\pi}_{r_{1}} > 0$ that $q^{\pi}_{r_{1}} > 0$.
In Algorithm~\ref{algo:threshold}, $w^{\pi}_{r_{1}}$ is determined by~\eqref{eqn:threshold >} or~\eqref{eqn:threshold ==} or~\eqref{eqn:threshold =>} or~\eqref{eqn:threshold =0}.
It follows from $w^{\pi}_{r_{1}} < \underline{w}^{\pi}_{r_{1}}$ that $w^{\pi}_{r_{1}}$ is not determined by~\eqref{eqn:threshold ==}.
Also, it follows from $w^{\pi}_{r_{1}} < \underline{w}^{\pi}_{r_{1}} \le \overline{w}^{\pi}$ that $w^{\pi}_{r_{1}}$ is not determined by~\eqref{eqn:threshold =0}.
Thus $w^{\pi}_{r_{1}}$ is determined by~\eqref{eqn:threshold >} or~\eqref{eqn:threshold =>}.
In both cases, $S^{\pi}_{r_{1}} > S^{\pi}_{r_{2}}$, and $w^{\pi}_{r} = w^{\pi}_{r_{1}}$ for all $r \in \{r_{1},\ldots,r_{2}-1\}$.
It follows from $w^{\pi}_{r_{1}} < \underline{w}^{\pi}_{r_{1}} \le \overline{w}^{\pi}$ that
\begin{align}
& w^{\pi}_{r_{1}} \ \ = \ \ w^{\min}_{1} \ \ = \ \ \min\left\{\frac{T^{\pi}_{r}(x(q)) - T^{\pi}_{r_{1}}(x(q))}{S^{\pi}_{r_{1}} - S^{\pi}_{r}} \; : \; S^{\pi}_{r_{1}} > S^{\pi}_{r}, \; r \in \hat{\mathcal{R}}^{\pi}\right\} \ \ = \ \ \frac{T^{\pi}_{r_{2}}(x(q)) - T^{\pi}_{r_{1}}(x(q))}{S^{\pi}_{r_{1}} - S^{\pi}_{r_{2}}} \nonumber \\
\Rightarrow \quad & T^{\pi}_{r_{1}}(x(q)) + w^{\pi}_{r_{1}} S^{\pi}_{r_{1}} \ \ = \ \ T^{\pi}_{r_{2}}(x(q)) + w^{\pi}_{r_{1}} S^{\pi}_{r_{2}} \nonumber \\
\Rightarrow \quad & T^{\pi}_{r_{1}}(x(q)) + w S^{\pi}_{r_{1}} \ \ > \ \ T^{\pi}_{r_{2}}(x(q)) + w S^{\pi}_{r_{2}} \quad \forall \ w > w^{\pi}_{r_{1}}
\label{eqn:least cost interval}
\end{align}
It follows from $q^{\pi}_{r_{1}} > 0$ that direction $e^{\pi}_{r_{2}} - e^{\pi}_{r_{1}}$ is feasible.
Next, consider the directional derivative of $\Upsilon$ at $q$ in the direction $e^{\pi}_{r_{2}} - e^{\pi}_{r_{1}}$:
\begin{align*}
\hspace{-5mm}
& \max\Bigg\{\left(D^{\pi} \left[T^{\pi}_{r_{2}}(x(q)) + \sum_{r=r_{2}}^{\hat{R}^{\pi}-1} \left(S^{\pi}_{r} - S^{\pi}_{r+1}\right) \hat{w}^{\pi}_{r}\right] + \lambda^{\pi} - \mu^{\pi}_{r_{2}}\right) \\
\hspace{-5mm}
& \hspace{10mm} - \left(D^{\pi} \left[T^{\pi}_{r_{1}}(x(q)) + \sum_{r=r_{1}}^{\hat{R}^{\pi}-1} \left(S^{\pi}_{r} - S^{\pi}_{r+1}\right) \hat{w}^{\pi}_{r}\right] + \lambda^{\pi}\right) \; : \; \hat{w}^{\pi}_{1} \in \left(G^{\pi}\right)^{-1}\left(p^{\pi}_{1}\right), \ldots, \hat{w}^{\pi}_{\hat{R}^{\pi}-1} \in \left(G^{\pi}\right)^{-1}\left(p^{\pi}_{\hat{R}^{\pi}-1}\right), \\
\hspace{-5mm}
& \hspace{145mm} \mu^{\pi}_{r_{2}} \ge 0, \; \mu^{\pi}_{r_{2}} q^{\pi}_{r_{2}} = 0\Bigg\} \\
\hspace{-5mm}
& = \ \ \max\Bigg\{D^{\pi} \left[T^{\pi}_{r_{2}}(x(q)) - T^{\pi}_{r_{1}}(x(q)) - \sum_{r=r_{1}}^{r_{2}-1} \left(S^{\pi}_{r} - S^{\pi}_{r+1}\right) \hat{w}^{\pi}_{r}\right] \; : \; \hat{w}^{\pi}_{r_{1}} \in \left(G^{\pi}\right)^{-1}\left(p^{\pi}_{r_{1}}\right), \ldots, \hat{w}^{\pi}_{r_{2}-1} \in \left(G^{\pi}\right)^{-1}\left(p^{\pi}_{r_{2}-1}\right)\Bigg\} \\
\hspace{-5mm}
& \le \ \ D^{\pi} \left[T^{\pi}_{r_{2}}(x(q)) - T^{\pi}_{r_{1}}(x(q)) - \sum_{r=r_{1}}^{r_{2}-1} \left(S^{\pi}_{r} - S^{\pi}_{r+1}\right) \underline{w}^{\pi}_{r_{1}}\right] \\
\hspace{-5mm}
& = \ \ D^{\pi} \left[T^{\pi}_{r_{2}}(x(q)) - T^{\pi}_{r_{1}}(x(q)) + \left(S^{\pi}_{r_{2}} - S^{\pi}_{r_{1}}\right) \underline{w}^{\pi}_{r_{1}}\right]
\ \ \overset{\eqref{eqn:least cost interval}}{<} \ \ 0
\end{align*}
which contradicts the optimality of $q$.

\noindent
\textbf{Case~2}: Suppose that $w^{\pi}_{r_{1}} > \overline{w}^{\pi}_{r_{1}}$. \\
Note that $w^{\pi}_{r_{1}} > \overline{w}^{\pi}_{r_{1}}$ implies that $1 \ge G^{\pi}(w^{\pi}_{r_{1}}) > p^{\pi}_{r_{1}}$.
It follows from $p^{\pi}_{r_{1}} < 1$ that $\big\{r \in \{r_{1}+1,\ldots,\hat{R}^{\pi}\} \, : \, q^{\pi}_{r} > 0\big\} \neq \varnothing$.
Let $\check{r} \defi \min\big\{r \in \{r_{1}+1,\ldots,\hat{R}^{\pi}\} \, : \, q^{\pi}_{r} > 0\big\}$.
It follows that $p^{\pi}_{r} = p^{\pi}_{r_{1}}$, and $\left(G^{\pi}\right)^{-1}\left(p^{\pi}_{r}\right) = \left(G^{\pi}\right)^{-1}\left(p^{\pi}_{r_{1}}\right) = [\underline{w}^{\pi}_{r_{1}},\overline{w}^{\pi}_{r_{1}}]$ for all $r \in \{r_{1},\ldots,\check{r}-1\}$.

\noindent
\textbf{Case~2.1}: $S^{\pi}_{r_{1}} = S^{\pi}_{\check{r}}$. \\
Then $w^{\pi}_{r_{1}}$ is determined by~\eqref{eqn:threshold ==} or~\eqref{eqn:threshold =>} or~\eqref{eqn:threshold =0}.
It follows from $w^{\pi}_{r_{1}-1} = \underline{w}^{\pi} \le \underline{w}^{\pi}_{r_{1}} \le \overline{w}^{\pi}_{r_{1}} < w^{\pi}_{r_{1}}$ that $w^{\pi}_{r_{1}}$ is not determined by~\eqref{eqn:threshold ==}.
If $w^{\pi}_{r_{1}}$ is determined by~\eqref{eqn:threshold =>} or~\eqref{eqn:threshold =0}, then $T^{\pi}_{\check{r}}(x(q)) \ge T^{\min}_{1} > T^{\pi}_{r_{1}}(x(q))$, and thus $T^{\pi}_{\check{r}}(x(q)) + w S^{\pi}_{\check{r}} > T^{\pi}_{r_{1}}(x(q)) + w S^{\pi}_{r_{1}}$ for all~$w$, which contradicts the optimality of $q$ with $q^{\pi}_{\check{r}} > 0$.

\noindent
\textbf{Case~2.2}: $S^{\pi}_{r_{1}} > S^{\pi}_{\check{r}}$. \\
Then $w^{\pi}_{r_{1}}$ is determined by~\eqref{eqn:threshold >} or~\eqref{eqn:threshold ==} or~\eqref{eqn:threshold =>}.
It follows from $w^{\pi}_{r_{1}-1} = \underline{w}^{\pi} \le \underline{w}^{\pi}_{r_{1}} \le \overline{w}^{\pi}_{r_{1}} < w^{\pi}_{r_{1}}$ that $w^{\pi}_{r_{1}}$ is not determined by~\eqref{eqn:threshold ==}.
If $w^{\pi}_{r_{1}}$ is determined by~\eqref{eqn:threshold >} or~\eqref{eqn:threshold =>}, then
\begin{align}
& w^{\pi}_{r_{1}} \ \ \le \ \ w^{\min}_{1} \ \ = \ \ \min\left\{\frac{T^{\pi}_{r}(x(q)) - T^{\pi}_{r_{1}}(x(q))}{S^{\pi}_{r_{1}} - S^{\pi}_{r}} \; : \; S^{\pi}_{r_{1}} > S^{\pi}_{r}, \; r \in \hat{\mathcal{R}}^{\pi}\right\}
\ \ \le \ \ \frac{T^{\pi}_{\check{r}}(x(q)) - T^{\pi}_{r_{1}}(x(q))}{S^{\pi}_{r_{1}} - S^{\pi}_{\check{r}}} \nonumber \\
\Rightarrow \quad & T^{\pi}_{r_{1}}(x(q)) + w^{\pi}_{r_{1}} S^{\pi}_{r_{1}} \ \ \le \ \ T^{\pi}_{\check{r}}(x(q)) + w^{\pi}_{r_{1}} S^{\pi}_{\check{r}} \nonumber \\
\Rightarrow \quad & T^{\pi}_{r_{1}}(x(q)) + w S^{\pi}_{r_{1}} \ \ < \ \ T^{\pi}_{\check{r}}(x(q)) + w S^{\pi}_{\check{r}} \quad \forall \ w < w^{\pi}_{r_{1}}
\label{eqn:least cost interval 2}
\end{align}
It follows from $q^{\pi}_{\check{r}} > 0$ that direction $e^{\pi}_{r_{1}} - e^{\pi}_{\check{r}}$ is feasible.
Next, consider the directional derivative of $\Upsilon$ at $q$ in the direction $e^{\pi}_{r_{1}} - e^{\pi}_{\check{r}}$:
\begin{align*}
& \max\Bigg\{\left(D^{\pi} \left[T^{\pi}_{r_{1}}(x(q)) + \sum_{r=r_{1}}^{\hat{R}^{\pi}-1} \left(S^{\pi}_{r} - S^{\pi}_{r+1}\right) \hat{w}^{\pi}_{r}\right] + \lambda^{\pi} - \mu^{\pi}_{r_{1}}\right) \\
& \hspace{10mm} - \left(D^{\pi} \left[T^{\pi}_{\check{r}}(x(q)) + \sum_{r=\check{r}}^{\hat{R}^{\pi}-1} \left(S^{\pi}_{r} - S^{\pi}_{r+1}\right) \hat{w}^{\pi}_{r}\right] + \lambda^{\pi}\right) \; : \; \hat{w}^{\pi}_{1} \in \left(G^{\pi}\right)^{-1}\left(p^{\pi}_{1}\right), \ldots, \hat{w}^{\pi}_{\hat{R}^{\pi}-1} \in \left(G^{\pi}\right)^{-1}\left(p^{\pi}_{\hat{R}^{\pi}-1}\right), \\
& \hspace{145mm} \mu^{\pi}_{r_{1}} \ge 0, \; \mu^{\pi}_{r_{1}} q^{\pi}_{r_{1}} = 0\Bigg\} \\
& = \ \ \max\Bigg\{D^{\pi} \left[T^{\pi}_{r_{1}}(x(q)) - T^{\pi}_{\check{r}}(x(q)) + \sum_{r=r_{1}}^{\check{r}-1} \left(S^{\pi}_{r} - S^{\pi}_{r+1}\right) \hat{w}^{\pi}_{r}\right] \; : \; \hat{w}^{\pi}_{r_{1}} \in \left(G^{\pi}\right)^{-1}\left(p^{\pi}_{r_{1}}\right), \ldots, \hat{w}^{\pi}_{\check{r}-1} \in \left(G^{\pi}\right)^{-1}\left(p^{\pi}_{\check{r}-1}\right)\Bigg\} \\
& = \ \ D^{\pi} \left[T^{\pi}_{r_{1}}(x(q)) - T^{\pi}_{\check{r}}(x(q)) + \sum_{r=r_{1}}^{\check{r}-1} \left(S^{\pi}_{r} - S^{\pi}_{r+1}\right) \overline{w}^{\pi}_{r_{1}}\right] \\
& = \ \ D^{\pi} \left[T^{\pi}_{r_{1}}(x(q)) - T^{\pi}_{\check{r}}(x(q)) + \left(S^{\pi}_{r_{1}} - S^{\pi}_{\check{r}}\right) \overline{w}^{\pi}_{r_{1}}\right]
\ \ \overset{\eqref{eqn:least cost interval 2}}{<} \ \ 0
\end{align*}
which contradicts the optimality of $q$.

\item
It follows from~\eqref{eqn:efficient path 1} that
\[
T^{\pi}_{r_{1}}(x(q)) + \underline{w}^{\pi} S^{\pi}_{r_{1}} \ \ \le \ \ T^{\pi}_{r}(x(q)) + \underline{w}^{\pi} S^{\pi}_{r} \qquad \forall \ r \in \hat{\mathcal{R}}^{\pi}
\]
Thus, for any $w \ge \underline{w}^{\pi}$ and any $r < r_{1}$ it holds that $S^{\pi}_{r} \ge S^{\pi}_{r_{1}}$, and hence
\begin{equation}
\left(S^{\pi}_{r} - S^{\pi}_{r_{1}}\right) w \ \ \ge \ \ \left(S^{\pi}_{r} - S^{\pi}_{r_{1}}\right) \underline{w}^{\pi}
\qquad \Rightarrow \quad T^{\pi}_{r}(x(q)) + w S^{\pi}_{r} \ \ \overset{\eqref{eqn:efficient path 1}}{\ge} \ \ T^{\pi}_{r_{1}}(x(q)) + w S^{\pi}_{r_{1}}
\label{eqn:first path lower cost than previous paths}
\end{equation}
In Algorithm~\ref{algo:threshold}, $w^{\pi}_{r_{1}}$ is determined by~\eqref{eqn:threshold >} or~\eqref{eqn:threshold ==} or~\eqref{eqn:threshold =>} or~\eqref{eqn:threshold =0}.

\noindent
\textbf{Case~1}: $w^{\pi}_{r_{1}}$ is determined by~\eqref{eqn:threshold >} or~\eqref{eqn:threshold =>}: \\
Then, for any $w \le w^{\min}_{1}$ and any $r > r_{1}$ it holds that either $S^{\pi}_{r} = S^{\pi}_{r_{1}}$ and $T^{\pi}_{r}(x(q)) > T^{\pi}_{r_{1}}(x(q))$, or $S^{\pi}_{r_{1}} > S^{\pi}_{r}$ and
\begin{equation}
w \ \ \le \ \ w^{\min}_{1} \ \ \le \ \ \frac{T^{\pi}_{r}(x(q)) - T^{\pi}_{r_{1}}(x(q))}{S^{\pi}_{r_{1}} - S^{\pi}_{r}}
\qquad \Rightarrow \quad T^{\pi}_{r_{1}}(x(q)) + w S^{\pi}_{r_{1}} \ \ \le \ \ T^{\pi}_{r}(x(q)) + w S^{\pi}_{r}
\label{eqn:first path lower cost than following paths}
\end{equation}
It follows from~\eqref{eqn:first path lower cost than previous paths} and~\eqref{eqn:first path lower cost than following paths} that
\[
T^{\pi}_{r_{1}}(x(q)) + w S^{\pi}_{r_{1}} \ \ \le \ \ T^{\pi}_{r}(x(q)) + w S^{\pi}_{r} \qquad \forall \ r \in \hat{\mathcal{R}}^{\pi}, \ \forall \ w \in [w^{\pi}_{r_{0}},w^{\pi}_{r_{1}}] \subset [\underline{w}^{\pi},w^{\min}_{1}]
\]

\noindent
\textbf{Case~2}: $w^{\pi}_{r_{1}}$ is determined by~\eqref{eqn:threshold ==}: \\
\textbf{Case~2.1}: $q^{\pi}_{r_{1}} = 0$:
It was established above that $q^{\pi}_{r} = 0$ for all $r < r_{1}$, and thus $p^{\pi}_{r_{1}} = 0$ and $w^{\pi}_{r_{1}} = \underline{w}^{\pi} = \max\left\{w^{\pi}_{r_{0}}, \min\left\{w \in \left(G^{\pi}\right)^{-1}\left(p^{\pi}_{r_{1}}\right)\right\}\right\}$.
It follows from~\eqref{eqn:efficient path 1} that
\[
T^{\pi}_{r_{1}}(x(q)) + w S^{\pi}_{r_{1}} \ \ \le \ \ T^{\pi}_{r}(x(q)) + w S^{\pi}_{r} \qquad \forall \ r \in \hat{\mathcal{R}}^{\pi}, \ \forall \ w \in [w^{\pi}_{r_{0}},w^{\pi}_{r_{1}}] = [\underline{w}^{\pi},\underline{w}^{\pi}]
\]
\textbf{Case~2.2}: $q^{\pi}_{r_{1}} > 0$:
For any $r \in \hat{\mathcal{R}}^{\pi}$, if $S^{\pi}_{r} = S^{\pi}_{r_{1}}$, then $T^{\pi}_{r}(x(q)) \ge T^{\pi}_{r_{1}}(x(q))$ and thus $T^{\pi}_{r}(x(q)) + w S^{\pi}_{r} \ge T^{\pi}_{r_{1}}(x(q)) + w S^{\pi}_{r_{1}}$ for all~$w$.
Suppose that $w^{\pi}_{r_{1}} = \min\left\{w \in \left(G^{\pi}\right)^{-1}\left(p^{\pi}_{r_{1}}\right)\right\} > \underline{w}^{\pi}$ and there is $\hat{r} > r_{1}$ such that $S^{\pi}_{\hat{r}} < S^{\pi}_{r_{1}}$ and $T^{\pi}_{r_{1}}(x(q)) + w^{\pi}_{r_{1}} S^{\pi}_{r_{1}} > T^{\pi}_{\hat{r}}(x(q)) + w^{\pi}_{r_{1}} S^{\pi}_{\hat{r}}$.
Then consider the directional derivative of $\Upsilon$ at $q$ in the feasible direction $e^{\pi}_{\hat{r}} - e^{\pi}_{r_{1}}$:
\begin{align*}
& \max\Bigg\{\left(D^{\pi} \left[T^{\pi}_{\hat{r}}(x(q)) + \sum_{r=\hat{r}}^{\hat{R}^{\pi}-1} \left(S^{\pi}_{r} - S^{\pi}_{r+1}\right) \hat{w}^{\pi}_{r}\right] + \lambda^{\pi} - \mu^{\pi}_{\hat{r}}\right) \\
& \hspace{10mm} - \left(D^{\pi} \left[T^{\pi}_{r_{1}}(x(q)) + \sum_{r=r_{1}}^{\hat{R}^{\pi}-1} \left(S^{\pi}_{r} - S^{\pi}_{r+1}\right) \hat{w}^{\pi}_{r}\right] + \lambda^{\pi}\right) \; : \; \hat{w}^{\pi}_{1} \in \left(G^{\pi}\right)^{-1}\left(p^{\pi}_{1}\right), \ldots, \hat{w}^{\pi}_{\hat{R}^{\pi}-1} \in \left(G^{\pi}\right)^{-1}\left(p^{\pi}_{\hat{R}^{\pi}-1}\right), \\
& \hspace{150mm} \mu^{\pi}_{\hat{r}} \ge 0, \; \mu^{\pi}_{\hat{r}} q^{\pi}_{\hat{r}} = 0\Bigg\} \\
& = \ \ \max\Bigg\{D^{\pi} \left[T^{\pi}_{\hat{r}}(x(q)) - T^{\pi}_{r_{1}}(x(q)) - \sum_{r=r_{1}}^{\hat{r}-1} \left(S^{\pi}_{r} - S^{\pi}_{r+1}\right) \hat{w}^{\pi}_{r}\right] \; : \; \hat{w}^{\pi}_{r_{1}} \in \left(G^{\pi}\right)^{-1}\left(p^{\pi}_{r_{1}}\right), \ldots, \hat{w}^{\pi}_{\hat{r}-1} \in \left(G^{\pi}\right)^{-1}\left(p^{\pi}_{\hat{r}-1}\right)\Bigg\} \\
& \le \ \ D^{\pi} \left[T^{\pi}_{\hat{r}}(x(q)) - T^{\pi}_{r_{1}}(x(q)) - \sum_{r=r_{1}}^{\hat{r}-1} \left(S^{\pi}_{r} - S^{\pi}_{r+1}\right) w^{\pi}_{r_{1}}\right] \\
& = \ \ D^{\pi} \left[T^{\pi}_{\hat{r}}(x(q)) - T^{\pi}_{r_{1}}(x(q)) + w^{\pi}_{r_{1}} \left(S^{\pi}_{\hat{r}} - S^{\pi}_{r_{1}}\right)\right]
\ \ < \ \ 0
\end{align*}
which contradicts the optimality of $q$.
(The first inequality follows from $w^{\pi}_{r_{1}} = \min\left\{w \in \left(G^{\pi}\right)^{-1}\left(p^{\pi}_{r_{1}}\right)\right\}$.)
Thus, for all $r \in \hat{\mathcal{R}}^{\pi}$ such that $r > r_{1}$ it holds that
\begin{align}
& T^{\pi}_{r_{1}}(x(q)) + w^{\pi}_{r_{1}} S^{\pi}_{r_{1}} \ \ \le \ \ T^{\pi}_{r}(x(q)) + w^{\pi}_{r_{1}} S^{\pi}_{r} \nonumber \\
\Rightarrow \qquad & T^{\pi}_{r_{1}}(x(q)) + w S^{\pi}_{r_{1}} \ \ \le \ \ T^{\pi}_{r}(x(q)) + w S^{\pi}_{r} \qquad \forall \ w \le w^{\pi}_{r_{1}}
\label{eqn:first path lower cost than following paths2}
\end{align}
It follows from~\eqref{eqn:first path lower cost than previous paths} and~\eqref{eqn:first path lower cost than following paths2} that
\[
T^{\pi}_{r_{1}}(x(q)) + w S^{\pi}_{r_{1}} \ \ \le \ \ T^{\pi}_{r}(x(q)) + w S^{\pi}_{r} \qquad \forall \ r \in \hat{\mathcal{R}}^{\pi}, \ \forall \ w \in [w^{\pi}_{r_{0}},w^{\pi}_{r_{1}}]
\]

\noindent
\textbf{Case~3}: $w^{\pi}_{r_{1}}$ is determined by~\eqref{eqn:threshold =0}. \\
Then $S^{\pi}_{r} \ge S^{\pi}_{r_{1}}$ for all $r \in \hat{\mathcal{R}}^{\pi}$, and thus it follows from~\eqref{eqn:first path lower cost than previous paths} that
\[
T^{\pi}_{r_{1}}(x(q)) + w S^{\pi}_{r_{1}} \ \ \le \ \ T^{\pi}_{r}(x(q)) + w S^{\pi}_{r} \qquad \forall \ r \in \hat{\mathcal{R}}^{\pi}, \ \forall \ w \in [w^{\pi}_{r_{0}},w^{\pi}_{r_{1}}] \subset [\underline{w}^{\pi},\infty)
\]
\end{enumerate}

\paragraph{Induction.}$\mbox{}$

Suppose that the following induction hypothesis holds for some $\ell \in \{1,\ldots,L-1\}$:
\begin{enumerate}[label=(\arabic*), ref=(\arabic*), leftmargin=7mm]
\item
The thresholds $w$ satisfy $\underline{w}^{\pi} = w^{\pi}_{0} \le \cdots \le w^{\pi}_{r_{\ell-1}} \le w^{\pi}_{r_{\ell}} \le \overline{w}^{\pi}$.
\item
For every $r \in \hat{\mathcal{R}}^{\pi} \setminus \mathcal{R}^{\pi}_{\min}$ such that $r < r_{\ell}$, it holds that $q^{\pi}_{r} = 0$.
\item
For every $r \in \hat{\mathcal{R}}^{\pi}$ such that $r \le r_{\ell}$, it holds that $w^{\pi}_{r} \in \left(G^{\pi}\right)^{-1}\left(p^{\pi}_{r}\right)$.
\item
For every $\ell' \in \{1,\ldots,\ell\}$, it holds that $T^{\pi}_{r_{\ell'}}(x(q)) + w S^{\pi}_{r_{\ell'}} \le T^{\pi}_{r}(x(q)) + w S^{\pi}_{r}$ for all $r \in \hat{\mathcal{R}}^{\pi}$ and all $w \in [w^{\pi}_{r_{\ell'-1}},w^{\pi}_{r_{\ell'}}]$.
\end{enumerate}
The initialization step above established that the hypothesis holds for $\ell = 1$.
Next we show how the hypothesis extends to $\ell+1$.

\begin{enumerate}[label=(\arabic*), ref=(\arabic*), leftmargin=7mm]
\item
If $w^{\pi}_{r_{\ell}}$ is determined by~\eqref{eqn:threshold =0} in Algorithm~\ref{algo:threshold}, then $\ell = L$.
If $w^{\pi}_{r_{\ell}}$ is determined by~\eqref{eqn:threshold >} or~\eqref{eqn:threshold =>}, and $w^{\min}_{\ell} > \overline{w}^{\pi}$, then $\ell = L$.
Otherwise, $w^{\pi}_{r_{\ell}}$ is determined by~\eqref{eqn:threshold >} or~\eqref{eqn:threshold =>}, and $w^{\pi}_{r_{\ell}} = w^{\min}_{\ell} \le \overline{w}^{\pi}$, or $w^{\pi}_{r_{\ell}}$ is determined by~\eqref{eqn:threshold ==}.
In all of these cases, it holds that
\begin{equation}
\label{eqn:cost ell = cost ell+1}
T^{\pi}_{r_{\ell}}(x(q)) + w^{\pi}_{r_{\ell}} S^{\pi}_{r_{\ell}} \ \ = \ \ T^{\pi}_{r_{\ell+1}}(x(q)) + w^{\pi}_{r_{\ell}} S^{\pi}_{r_{\ell+1}}
\end{equation}
Also, $w^{\pi}_{r_{\ell+1}}$ is determined by~\eqref{eqn:threshold >} or~\eqref{eqn:threshold ==} or~\eqref{eqn:threshold =>} or~\eqref{eqn:threshold =0} or~\eqref{eqn:last threshold}.

\noindent
\textbf{Case~1}: $w^{\pi}_{r_{\ell+1}}$ is determined by~\eqref{eqn:threshold >} or~\eqref{eqn:threshold =>}: \\
If $w^{\min}_{\ell+1} > \overline{w}^{\pi}$, then $L = \ell+1$ and $w^{\pi}_{r_{\ell}} \le w^{\pi}_{r_{\ell+1}} = \overline{w}^{\pi}$.
Otherwise, $w^{\pi}_{r_{\ell+1}} = w^{\min}_{\ell+1} \le \overline{w}^{\pi}$, $S^{\pi}_{r_{\ell+1}} > S^{\pi}_{r_{\ell+2}}$, and $T^{\pi}_{r_{\ell+1}}(x(q)) + w^{\pi}_{r_{\ell+1}} S^{\pi}_{r_{\ell+1}} = T^{\pi}_{r_{\ell+2}}(x(q)) + w^{\pi}_{r_{\ell+1}} S^{\pi}_{r_{\ell+2}}$.
Suppose that $w^{\pi}_{r_{\ell}} > w^{\pi}_{r_{\ell+1}}$.
Then
\begin{align*}
& \left(S^{\pi}_{r_{\ell+1}} - S^{\pi}_{r_{\ell+2}}\right) w^{\pi}_{r_{\ell}} \ \ > \ \ \left(S^{\pi}_{r_{\ell+1}} - S^{\pi}_{r_{\ell+2}}\right) w^{\pi}_{r_{\ell+1}} \ \ = \ \ T^{\pi}_{r_{\ell+2}}(x(q)) - T^{\pi}_{r_{\ell+1}}(x(q)) \\
\Rightarrow \quad
& T^{\pi}_{r_{\ell}}(x(q)) + w^{\pi}_{r_{\ell}} S^{\pi}_{r_{\ell}} \ \ \overset{\eqref{eqn:cost ell = cost ell+1}}{=} \ \ T^{\pi}_{r_{\ell+1}}(x(q)) + w^{\pi}_{r_{\ell}} S^{\pi}_{r_{\ell+1}} \ \ > \ \ T^{\pi}_{r_{\ell+2}}(x(q)) + w^{\pi}_{r_{\ell}} S^{\pi}_{r_{\ell+2}}
\end{align*}
which contradicts the induction hypothesis~\ref{item:thresholds optimal}.
Thus $w^{\pi}_{r_{\ell}} \le w^{\pi}_{r_{\ell+1}} \le \overline{w}^{\pi}$.

\noindent
\textbf{Case~2}: $w^{\pi}_{r_{\ell+1}}$ is determined by~\eqref{eqn:threshold ==}: \\
Then it follows immediately that $w^{\pi}_{r_{\ell}} \le w^{\pi}_{r_{\ell+1}} \le \overline{w}^{\pi}$.

\noindent
\textbf{Case~3}: $w^{\pi}_{r_{\ell+1}}$ is determined by~\eqref{eqn:threshold =0} or~\eqref{eqn:last threshold}: \\
Then $w^{\pi}_{r_{\ell+1}} = \overline{w}^{\pi}$, and thus $w^{\pi}_{r_{\ell}} \le w^{\pi}_{r_{\ell+1}} = w^{\pi}_{r_{L}} = \overline{w}^{\pi}$.

\item
Suppose that there is a $r \in \hat{\mathcal{R}}^{\pi} \setminus \mathcal{R}^{\pi}_{\min}$ such that $r_{\ell} < r < r_{\ell+1}$ and $q^{\pi}_{r} > 0$.
In Algorithm~\ref{algo:threshold}, $w^{\pi}_{r_{\ell}}$ is determined by~\eqref{eqn:threshold >} or~\eqref{eqn:threshold ==} or~\eqref{eqn:threshold =>} or~\eqref{eqn:threshold =0}.
If $w^{\pi}_{r_{\ell}}$ is determined by~\eqref{eqn:threshold =0}, then $\ell = L$.
The remaining cases are considered next.

\noindent
\textbf{Case~1}: $w^{\pi}_{r_{\ell}}$ is determined by~\eqref{eqn:threshold >}, or by~\eqref{eqn:threshold =>} and $S^{\pi}_{r_{\ell}} > S^{\pi}_{r}$. \\
If $w^{\min}_{\ell} > \overline{w}^{\pi}$, then $\ell = L$.
Otherwise, $w^{\pi}_{r_{\ell}} = w^{\min}_{\ell} \le \overline{w}^{\pi}$, and
\begin{align*}
& w^{\pi}_{r_{\ell}} \ \ = \ \ w^{\min}_{\ell} \ \ = \ \ \frac{T^{\pi}_{r_{\ell+1}}(x(q)) - T^{\pi}_{r_{\ell}}(x(q))}{S^{\pi}_{r_{\ell}} - S^{\pi}_{r_{\ell+1}}}
\ \ < \ \ \frac{T^{\pi}_{r}(x(q)) - T^{\pi}_{r_{\ell}}(x(q))}{S^{\pi}_{r_{\ell}} - S^{\pi}_{r}} \\
\Rightarrow \quad & T^{\pi}_{r_{\ell}}(x(q)) + w^{\pi}_{r_{\ell}} S^{\pi}_{r_{\ell}} \ \ = \ \ T^{\pi}_{r_{\ell+1}}(x(q)) + w^{\pi}_{r_{\ell}} S^{\pi}_{r_{\ell+1}} \ \ < \ \ T^{\pi}_{r}(x(q)) + w^{\pi}_{r_{\ell}} S^{\pi}_{r}
\end{align*}
It follows from the induction hypothesis that $w^{\pi}_{r_{\ell}} \in \left(G^{\pi}\right)^{-1}\left(p^{\pi}_{r_{\ell}}\right)$.
Also, note that if $q^{\pi}_{r} > 0$ then $p^{\pi}_{r_{\ell}} < p^{\pi}_{r} \le p^{\pi}_{r'}$ for all $r' \ge r$.
Thus $w^{\pi}_{r_{\ell}} \le \hat{w}^{\pi}_{r'}$ for all $\hat{w}^{\pi}_{r'} \in \left(G^{\pi}\right)^{-1}\left(p^{\pi}_{r'}\right)$, and all $r' \ge r$.
Consider the directional derivative of $\Upsilon$ at $q$ in the feasible direction $e^{\pi}_{r_{\ell+1}} - e^{\pi}_{r}$:
\begin{align*}
& \max\Bigg\{\left(D^{\pi} \left[T^{\pi}_{r_{\ell+1}}(x(q)) + \sum_{r'=r_{\ell+1}}^{\hat{R}^{\pi}-1} \left(S^{\pi}_{r'} - S^{\pi}_{r'+1}\right) \hat{w}^{\pi}_{r'}\right] + \lambda^{\pi} - \mu^{\pi}_{r_{\ell+1}}\right) \\
& \hspace{10mm} - \left(D^{\pi} \left[T^{\pi}_{r}(x(q)) + \sum_{r'=r}^{\hat{R}^{\pi}-1} \left(S^{\pi}_{r'} - S^{\pi}_{r'+1}\right) \hat{w}^{\pi}_{r'}\right] + \lambda^{\pi}\right) \; : \; \hat{w}^{\pi}_{1} \in \left(G^{\pi}\right)^{-1}\left(p^{\pi}_{1}\right), \ldots, \hat{w}^{\pi}_{\hat{R}^{\pi}-1} \in \left(G^{\pi}\right)^{-1}\left(p^{\pi}_{\hat{R}^{\pi}-1}\right), \\
& \hspace{140mm} \mu^{\pi}_{r_{\ell+1}} \ge 0, \; \mu^{\pi}_{r_{\ell+1}} q^{\pi}_{r_{\ell+1}} = 0\Bigg\} \\
& = \ \ \max\Bigg\{D^{\pi} \left[T^{\pi}_{r_{\ell+1}}(x(q)) - T^{\pi}_{r}(x(q)) - \sum_{r'=r}^{r_{\ell+1}-1} \left(S^{\pi}_{r'} - S^{\pi}_{r'+1}\right) \hat{w}^{\pi}_{r'}\right] \; : \; \hat{w}^{\pi}_{r} \in \left(G^{\pi}\right)^{-1}\left(p^{\pi}_{r}\right), \ldots \\
& \hspace{125mm} \ldots, \hat{w}^{\pi}_{r_{\ell+1}-1} \in \left(G^{\pi}\right)^{-1}\left(p^{\pi}_{r_{\ell+1}-1}\right)\Bigg\} \\
& \le \ \ D^{\pi} \left[T^{\pi}_{r_{\ell+1}}(x(q)) - T^{\pi}_{r}(x(q)) - \sum_{r'=r}^{r_{\ell+1}-1} \left(S^{\pi}_{r'} - S^{\pi}_{r'+1}\right) w^{\pi}_{r_{\ell}}\right] \\
& = \ \ D^{\pi} \left[T^{\pi}_{r_{\ell+1}}(x(q)) - T^{\pi}_{r}(x(q)) + \left(S^{\pi}_{r_{\ell+1}} - S^{\pi}_{r}\right) w^{\pi}_{r_{\ell}}\right]
\ \ < \ \ 0
\end{align*}
which contradicts the optimality of $q$.

\noindent
\textbf{Case~2}: $w^{\pi}_{r_{\ell}}$ is determined by~\eqref{eqn:threshold ==}, or by~\eqref{eqn:threshold =>} and $S^{\pi}_{r_{\ell}} = S^{\pi}_{r}$. \\
Then $T^{\pi}_{r_{\ell}}(x(q)) < T^{\pi}_{r}(x(q))$ and $S^{\pi}_{r_{\ell}} = S^{\pi}_{r'} = S^{\pi}_{r}$ for all $r' = r_{\ell},\ldots,r$.
Consider the directional derivative of $\Upsilon$ at $q$ in the feasible direction $e^{\pi}_{r_{\ell}} - e^{\pi}_{r}$:
\begin{align*}
\hspace{-5mm}
& \max\Bigg\{\left(D^{\pi} \left[T^{\pi}_{r_{\ell}}(x(q)) + \sum_{r'=r_{\ell}}^{\hat{R}^{\pi}-1} \left(S^{\pi}_{r'} - S^{\pi}_{r'+1}\right) \hat{w}^{\pi}_{r'}\right] + \lambda^{\pi} - \mu^{\pi}_{r_{\ell}}\right) \\
\hspace{-5mm}
& \hspace{10mm} - \left(D^{\pi} \left[T^{\pi}_{r}(x(q)) + \sum_{r'=r}^{\hat{R}^{\pi}-1} \left(S^{\pi}_{r'} - S^{\pi}_{r'+1}\right) \hat{w}^{\pi}_{r'}\right] + \lambda^{\pi}\right) \; : \; \hat{w}^{\pi}_{1} \in \left(G^{\pi}\right)^{-1}\left(p^{\pi}_{1}\right), \ldots, \hat{w}^{\pi}_{\hat{R}^{\pi}-1} \in \left(G^{\pi}\right)^{-1}\left(p^{\pi}_{\hat{R}^{\pi}-1}\right), \\
\hspace{-5mm}
& \hspace{145mm} \mu^{\pi}_{r_{\ell}} \ge 0, \; \mu^{\pi}_{r_{\ell}} q^{\pi}_{r_{\ell}} = 0\Bigg\} \\
\hspace{-5mm}
& = \ \ \max\Bigg\{D^{\pi} \left[T^{\pi}_{r_{\ell}}(x(q)) - T^{\pi}_{r}(x(q)) + \sum_{r'=r_{\ell}}^{r-1} \left(S^{\pi}_{r'} - S^{\pi}_{r'+1}\right) \hat{w}^{\pi}_{r'}\right] \; : \; \hat{w}^{\pi}_{r_{\ell}} \in \left(G^{\pi}\right)^{-1}\left(p^{\pi}_{r_{\ell}}\right), \ldots, \hat{w}^{\pi}_{r-1} \in \left(G^{\pi}\right)^{-1}\left(p^{\pi}_{r-1}\right)\Bigg\} \\
\hspace{-5mm}
& = \ \ D^{\pi} \left[T^{\pi}_{r_{\ell}}(x(q)) - T^{\pi}_{r}(x(q))\right]
\ \ < \ \ 0
\end{align*}
which contradicts the optimality of $q$.

\item
By the induction hypothesis, $w^{\pi}_{r_{\ell}} \in \left(G^{\pi}\right)^{-1}\left(p^{\pi}_{r_{\ell}}\right)$.
By Algorithm~\ref{algo:threshold}, $w^{\pi}_{r} = w^{\pi}_{r_{\ell}}$ for all $r = r_{\ell},\ldots,r_{\ell+1}-1$.
It was established in \ref{item:positive flow paths are efficient} above that $q^{\pi}_{r} = 0$ for all $r \in \{r_{\ell}+1,\ldots,r_{\ell+1}-1\}$, and thus $p^{\pi}_{r} = p^{\pi}_{r_{\ell}}$ and $w^{\pi}_{r} \in \left(G^{\pi}\right)^{-1}\left(p^{\pi}_{r}\right)$ for all $r \in \{r_{\ell}+1,\ldots,r_{\ell+1}-1\}$.
Next we show by contradiction that $w^{\pi}_{r_{\ell+1}} \in \left(G^{\pi}\right)^{-1}\left(p^{\pi}_{r_{\ell+1}}\right)$.

\noindent
\textbf{Case~1}: Suppose that $w^{\pi}_{r_{\ell+1}} < \underline{w}^{\pi}_{r_{\ell+1}}$. \\
Note that $w^{\pi}_{r_{\ell+1}} < \underline{w}^{\pi}_{r_{\ell+1}}$ implies that $G^{\pi}(w^{\pi}_{r_{\ell+1}}) < p^{\pi}_{r_{\ell+1}}$.
It was established in \ref{item:thresholds ordered} above that $w^{\pi}_{r_{\ell}} = \cdots = w^{\pi}_{r_{\ell+1}-1} \le w^{\pi}_{r_{\ell+1}}$.
It follows from $\underline{w}^{\pi}_{r_{\ell+1}-1} \le w^{\pi}_{r_{\ell+1}-1} \le w^{\pi}_{r_{\ell+1}} < \underline{w}^{\pi}_{r_{\ell+1}}$ that $p^{\pi}_{r_{\ell+1}-1} < p^{\pi}_{r_{\ell+1}}$, that is, $q^{\pi}_{r_{\ell+1}} > 0$.
In Algorithm~\ref{algo:threshold}, $w^{\pi}_{r_{\ell+1}}$ is determined by~\eqref{eqn:threshold >} or~\eqref{eqn:threshold ==} or~\eqref{eqn:threshold =>} or~\eqref{eqn:threshold =0} or~\eqref{eqn:last threshold}.
It follows from $w^{\pi}_{r_{\ell+1}} < \underline{w}^{\pi}_{r_{\ell+1}}$ that $w^{\pi}_{r_{\ell+1}}$ is not determined by~\eqref{eqn:threshold ==}.
Also, it follows from $w^{\pi}_{r_{\ell+1}} < \underline{w}^{\pi}_{r_{\ell+1}} \le \overline{w}^{\pi}$ that $w^{\pi}_{r_{\ell+1}}$ is not determined by~\eqref{eqn:threshold =0} or~\eqref{eqn:last threshold}.
Thus $w^{\pi}_{r_{\ell+1}}$ is determined by~\eqref{eqn:threshold >} or~\eqref{eqn:threshold =>}.
In both cases, $S^{\pi}_{r_{\ell+1}} > S^{\pi}_{r_{\ell+2}}$, and $w^{\pi}_{r} = w^{\pi}_{r_{\ell+1}}$ for all $r \in \{r_{\ell+1},\ldots,r_{\ell+2}-1\}$.
It follows from $w^{\pi}_{r_{\ell+1}} < \underline{w}^{\pi}_{r_{\ell+1}} \le \overline{w}^{\pi}$ that
\begin{align}
& w^{\pi}_{r_{\ell+1}} \ \ = \ \ w^{\min}_{\ell+1} \ \ = \ \ \min\left\{\frac{T^{\pi}_{r}(x(q)) - T^{\pi}_{r_{\ell+1}}(x(q))}{S^{\pi}_{r_{\ell+1}} - S^{\pi}_{r}} \; : \; S^{\pi}_{r_{\ell+1}} > S^{\pi}_{r}, \; r \in \hat{\mathcal{R}}^{\pi}\right\} \ \ = \ \ \frac{T^{\pi}_{r_{\ell+2}}(x(q)) - T^{\pi}_{r_{\ell+1}}(x(q))}{S^{\pi}_{r_{\ell+1}} - S^{\pi}_{r_{\ell+2}}} \nonumber \\
\Rightarrow \quad & T^{\pi}_{r_{\ell+1}}(x(q)) + w^{\pi}_{r_{\ell+1}} S^{\pi}_{r_{\ell+1}} \ \ = \ \ T^{\pi}_{r_{\ell+2}}(x(q)) + w^{\pi}_{r_{\ell+1}} S^{\pi}_{r_{\ell+2}} \nonumber \\
\Rightarrow \quad & T^{\pi}_{r_{\ell+1}}(x(q)) + w S^{\pi}_{r_{\ell+1}} \ \ > \ \ T^{\pi}_{r_{\ell+2}}(x(q)) + w S^{\pi}_{r_{\ell+2}} \quad \forall \ w > w^{\pi}_{r_{\ell+1}}
\label{eqn:least cost interval 3}
\end{align}
It follows from $q^{\pi}_{r_{\ell+1}} > 0$ that direction $e^{\pi}_{r_{\ell+2}} - e^{\pi}_{r_{\ell+1}}$ is feasible.
Next, consider the directional derivative of $\Upsilon$ at $q$ in the direction $e^{\pi}_{r_{\ell+2}} - e^{\pi}_{r_{\ell+1}}$:
\begin{align*}
& \max\Bigg\{\left(D^{\pi} \left[T^{\pi}_{r_{\ell+2}}(x(q)) + \sum_{r=r_{\ell+2}}^{\hat{R}^{\pi}-1} \left(S^{\pi}_{r} - S^{\pi}_{r+1}\right) \hat{w}^{\pi}_{r}\right] + \lambda^{\pi} - \mu^{\pi}_{r_{\ell+2}}\right) \\
& \hspace{10mm} - \left(D^{\pi} \left[T^{\pi}_{r_{\ell+1}}(x(q)) + \sum_{r=r_{\ell+1}}^{\hat{R}^{\pi}-1} \left(S^{\pi}_{r} - S^{\pi}_{r+1}\right) \hat{w}^{\pi}_{r}\right] + \lambda^{\pi}\right) \; : \; \hat{w}^{\pi}_{1} \in \left(G^{\pi}\right)^{-1}\left(p^{\pi}_{1}\right), \ldots, \hat{w}^{\pi}_{\hat{R}^{\pi}-1} \in \left(G^{\pi}\right)^{-1}\left(p^{\pi}_{\hat{R}^{\pi}-1}\right), \\
& \hspace{140mm} \mu^{\pi}_{r_{\ell+2}} \ge 0, \; \mu^{\pi}_{r_{\ell+2}} q^{\pi}_{r_{\ell+2}} = 0\Bigg\} \\
& = \ \ \max\Bigg\{D^{\pi} \left[T^{\pi}_{r_{\ell+2}}(x(q)) - T^{\pi}_{r_{\ell+1}}(x(q)) - \sum_{r=r_{\ell+1}}^{r_{\ell+2}-1} \left(S^{\pi}_{r} - S^{\pi}_{r+1}\right) \hat{w}^{\pi}_{r}\right] \; : \; \hat{w}^{\pi}_{r_{\ell+1}} \in \left(G^{\pi}\right)^{-1}\left(p^{\pi}_{r_{\ell+1}}\right), \ldots \\
& \hspace{125mm} \ldots, \hat{w}^{\pi}_{r_{\ell+2}-1} \in \left(G^{\pi}\right)^{-1}\left(p^{\pi}_{r_{\ell+2}-1}\right)\Bigg\} \\
& \le \ \ D^{\pi} \left[T^{\pi}_{r_{\ell+2}}(x(q)) - T^{\pi}_{r_{\ell+1}}(x(q)) - \sum_{r=r_{\ell+1}}^{r_{\ell+2}-1} \left(S^{\pi}_{r} - S^{\pi}_{r+1}\right) \underline{w}^{\pi}_{r_{\ell+1}}\right] \\
& = \ \ D^{\pi} \left[T^{\pi}_{r_{\ell+2}}(x(q)) - T^{\pi}_{r_{\ell+1}}(x(q)) + \left(S^{\pi}_{r_{\ell+2}} - S^{\pi}_{r_{\ell+1}}\right) \underline{w}^{\pi}_{r_{\ell+1}}\right]
\ \ \overset{\eqref{eqn:least cost interval 3}}{<} \ \ 0
\end{align*}
which contradicts the optimality of $q$.

\noindent
\textbf{Case~2}: Suppose that $w^{\pi}_{r_{\ell+1}} > \overline{w}^{\pi}_{r_{\ell+1}}$. \\
Note that $w^{\pi}_{r_{\ell+1}} > \overline{w}^{\pi}_{r_{\ell+1}}$ implies that $1 \ge G^{\pi}(w^{\pi}_{r_{\ell+1}}) > p^{\pi}_{r_{\ell+1}}$.
It follows from $p^{\pi}_{r_{\ell+1}} < 1$ that $\big\{r \in \{r_{\ell+1}+1,\ldots,\hat{R}^{\pi}\} \, : \, q^{\pi}_{r} > 0\big\} \neq \varnothing$.
Let $\check{r} \defi \min\big\{r \in \{r_{\ell+1}+1,\ldots,\hat{R}^{\pi}\} \, : \, q^{\pi}_{r} > 0\big\}$.
It follows that $p^{\pi}_{r} = p^{\pi}_{r_{\ell+1}}$, and $\left(G^{\pi}\right)^{-1}\left(p^{\pi}_{r}\right) = \left(G^{\pi}\right)^{-1}\left(p^{\pi}_{r_{\ell+1}}\right) = [\underline{w}^{\pi}_{r_{\ell+1}},\overline{w}^{\pi}_{r_{\ell+1}}]$ for all $r \in \{r_{\ell+1},\ldots,\check{r}-1\}$.

\noindent
\textbf{Case~2.1}: $S^{\pi}_{r_{\ell+1}} = S^{\pi}_{\check{r}}$. \\
Then $w^{\pi}_{r_{\ell+1}}$ is determined by~\eqref{eqn:threshold ==} or~\eqref{eqn:threshold =>} or~\eqref{eqn:threshold =0}.
It follows from $w^{\pi}_{r_{\ell+1}-1} \le \overline{w}^{\pi}_{r_{\ell+1}-1} \le \overline{w}^{\pi}_{r_{\ell+1}} < w^{\pi}_{r_{\ell+1}}$ that $w^{\pi}_{r_{\ell+1}}$ is not determined by~\eqref{eqn:threshold ==}.
If $w^{\pi}_{r_{\ell+1}}$ is determined by~\eqref{eqn:threshold =>} or~\eqref{eqn:threshold =0}, then $T^{\pi}_{\check{r}}(x(q)) \ge T^{\min}_{\ell+1} > T^{\pi}_{r_{\ell+1}}(x(q))$, and thus $T^{\pi}_{\check{r}}(x(q)) + w S^{\pi}_{\check{r}} > T^{\pi}_{r_{\ell+1}}(x(q)) + w S^{\pi}_{r_{\ell+1}}$ for all~$w$, which contradicts the optimality of $q$ with $q^{\pi}_{\check{r}} > 0$.

\noindent
\textbf{Case~2.2}: $S^{\pi}_{r_{\ell+1}} > S^{\pi}_{\check{r}}$. \\
Then $w^{\pi}_{r_{\ell+1}}$ is determined by~\eqref{eqn:threshold >} or~\eqref{eqn:threshold ==} or~\eqref{eqn:threshold =>}.
It follows from $w^{\pi}_{r_{\ell+1}-1} \le \overline{w}^{\pi}_{r_{\ell+1}-1} \le \overline{w}^{\pi}_{r_{\ell+1}} < w^{\pi}_{r_{\ell+1}}$ that $w^{\pi}_{r_{\ell+1}}$ is not determined by~\eqref{eqn:threshold ==}.
If $w^{\pi}_{r_{\ell+1}}$ is determined by~\eqref{eqn:threshold >} or~\eqref{eqn:threshold =>}, then
\begin{align}
& w^{\pi}_{r_{\ell+1}} \ \ \le \ \ w^{\min}_{\ell+1} \ \ = \ \ \min\left\{\frac{T^{\pi}_{r}(x(q)) - T^{\pi}_{r_{\ell+1}}(x(q))}{S^{\pi}_{r_{\ell+1}} - S^{\pi}_{r}} \; : \; S^{\pi}_{r_{\ell+1}} > S^{\pi}_{r}, \; r \in \hat{\mathcal{R}}^{\pi}\right\}
\ \ \le \ \ \frac{T^{\pi}_{\check{r}}(x(q)) - T^{\pi}_{r_{\ell+1}}(x(q))}{S^{\pi}_{r_{\ell+1}} - S^{\pi}_{\check{r}}} \nonumber \\
\Rightarrow \quad & T^{\pi}_{r_{\ell+1}}(x(q)) + w^{\pi}_{r_{\ell+1}} S^{\pi}_{r_{\ell+1}} \ \ \le \ \ T^{\pi}_{\check{r}}(x(q)) + w^{\pi}_{r_{\ell+1}} S^{\pi}_{\check{r}} \nonumber \\
\Rightarrow \quad & T^{\pi}_{r_{\ell+1}}(x(q)) + w S^{\pi}_{r_{\ell+1}} \ \ < \ \ T^{\pi}_{\check{r}}(x(q)) + w S^{\pi}_{\check{r}} \quad \forall \ w < w^{\pi}_{r_{\ell+1}}
\label{eqn:least cost interval 4}
\end{align}
It follows from $q^{\pi}_{\check{r}} > 0$ that direction $e^{\pi}_{r_{\ell+1}} - e^{\pi}_{\check{r}}$ is feasible.
Next, consider the directional derivative of $\Upsilon$ at $q$ in the direction $e^{\pi}_{r_{\ell+1}} - e^{\pi}_{\check{r}}$:
\begin{align*}
& \max\Bigg\{\left(D^{\pi} \left[T^{\pi}_{r_{\ell+1}}(x(q)) + \sum_{r=r_{\ell+1}}^{\hat{R}^{\pi}-1} \left(S^{\pi}_{r} - S^{\pi}_{r+1}\right) \hat{w}^{\pi}_{r}\right] + \lambda^{\pi} - \mu^{\pi}_{r_{\ell+1}}\right) \\
& \hspace{10mm} - \left(D^{\pi} \left[T^{\pi}_{\check{r}}(x(q)) + \sum_{r=\check{r}}^{\hat{R}^{\pi}-1} \left(S^{\pi}_{r} - S^{\pi}_{r+1}\right) \hat{w}^{\pi}_{r}\right] + \lambda^{\pi}\right) \; : \; \hat{w}^{\pi}_{1} \in \left(G^{\pi}\right)^{-1}\left(p^{\pi}_{1}\right), \ldots, \hat{w}^{\pi}_{\hat{R}^{\pi}-1} \in \left(G^{\pi}\right)^{-1}\left(p^{\pi}_{\hat{R}^{\pi}-1}\right), \\
& \hspace{140mm} \mu^{\pi}_{r_{\ell+1}} \ge 0, \; \mu^{\pi}_{r_{\ell+1}} q^{\pi}_{r_{\ell+1}} = 0\Bigg\} \\
& = \ \ \max\Bigg\{D^{\pi} \left[T^{\pi}_{r_{\ell+1}}(x(q)) - T^{\pi}_{\check{r}}(x(q)) + \sum_{r=r_{\ell+1}}^{\check{r}-1} \left(S^{\pi}_{r} - S^{\pi}_{r+1}\right) \hat{w}^{\pi}_{r}\right] \; : \; \hat{w}^{\pi}_{r_{\ell+1}} \in \left(G^{\pi}\right)^{-1}\left(p^{\pi}_{r_{\ell+1}}\right), \ldots \\
& \hspace{130mm} \ldots, \hat{w}^{\pi}_{\check{r}-1} \in \left(G^{\pi}\right)^{-1}\left(p^{\pi}_{\check{r}-1}\right)\Bigg\} \\
& = \ \ D^{\pi} \left[T^{\pi}_{r_{\ell+1}}(x(q)) - T^{\pi}_{\check{r}}(x(q)) + \sum_{r=r_{\ell+1}}^{\check{r}-1} \left(S^{\pi}_{r} - S^{\pi}_{r+1}\right) \overline{w}^{\pi}_{r_{\ell+1}}\right] \\
& = \ \ D^{\pi} \left[T^{\pi}_{r_{\ell+1}}(x(q)) - T^{\pi}_{\check{r}}(x(q)) + \left(S^{\pi}_{r_{\ell+1}} - S^{\pi}_{\check{r}}\right) \overline{w}^{\pi}_{r_{\ell+1}}\right]
\ \ \overset{\eqref{eqn:least cost interval 4}}{<} \ \ 0
\end{align*}
which contradicts the optimality of $q$.

\item
We show that the following inequalities hold, which implies the induction hypothesis for $\ell+1$:
\begin{align}
\label{eqn:optimal intervals <=}
T^{\pi}_{r_{\ell+1}}(x(q)) + w S^{\pi}_{r_{\ell+1}} \ \ \le \ \ T^{\pi}_{r}(x(q)) + w S^{\pi}_{r}
\qquad \forall \ r \in \hat{\mathcal{R}}^{\pi} \mbox{ s.t. } S^{\pi}_{r} \ge S^{\pi}_{r_{\ell+1}}, \ \forall \ w \ge w^{\pi}_{r_{\ell}} \\
\label{eqn:optimal intervals >}
T^{\pi}_{r_{\ell+1}}(x(q)) + w S^{\pi}_{r_{\ell+1}} \ \ \le \ \ T^{\pi}_{r}(x(q)) + w S^{\pi}_{r}
\qquad \forall \ r \in \hat{\mathcal{R}}^{\pi} \mbox{ s.t. } S^{\pi}_{r} < S^{\pi}_{r_{\ell+1}}, \ \forall \ w \le w^{\pi}_{r_{\ell+1}}
\end{align}
In Algorithm~\ref{algo:threshold}, if $w^{\pi}_{r_{\ell}}$ is determined by~\eqref{eqn:threshold >} or~\eqref{eqn:threshold =>}, and $w^{\min}_{\ell} > \overline{w}^{\pi}$, or $w^{\pi}_{r_{\ell}}$ is determined by~\eqref{eqn:threshold =0}, then $\ell = L$.
Otherwise, $w^{\pi}_{r_{\ell}}$ is determined by~\eqref{eqn:threshold >} or~\eqref{eqn:threshold =>}, and $w^{\min}_{\ell} \le \overline{w}^{\pi}$, or $w^{\pi}_{r_{\ell}}$ is determined by~\eqref{eqn:threshold ==}.
In all these cases, $w^{\pi}_{r_{\ell}} = w^{\min}_{\ell} \le \overline{w}^{\pi}$ and
\[
T^{\pi}_{r_{\ell}}(x(q)) + w^{\pi}_{r_{\ell}} S^{\pi}_{r_{\ell}} \ \ = \ \ T^{\pi}_{r_{\ell+1}}(x(q)) + w^{\pi}_{r_{\ell}} S^{\pi}_{r_{\ell+1}}
\]
It follows from the hypothesis that
\[
T^{\pi}_{r_{\ell}}(x(q)) + w^{\pi}_{r_{\ell}} S^{\pi}_{r_{\ell}} \ \ \le \ \ T^{\pi}_{r}(x(q)) + w^{\pi}_{r_{\ell}} S^{\pi}_{r} \qquad \forall \ r \in \hat{\mathcal{R}}^{\pi}
\]
Thus
\begin{align*}
T^{\pi}_{r_{\ell+1}}(x(q)) + w^{\pi}_{r_{\ell}} S^{\pi}_{r_{\ell+1}} \ \ & \le \ \ T^{\pi}_{r}(x(q)) + w^{\pi}_{r_{\ell}} S^{\pi}_{r} \qquad \forall \ r \in \hat{\mathcal{R}}^{\pi} \\
\Rightarrow \quad T^{\pi}_{r_{\ell+1}}(x(q)) + w S^{\pi}_{r_{\ell+1}} \ \ & \le \ \ T^{\pi}_{r}(x(q)) + w S^{\pi}_{r} \qquad \forall \ r \in \hat{\mathcal{R}}^{\pi} \mbox{ s.t. } S^{\pi}_{r} \ge S^{\pi}_{r_{\ell+1}}, \ \forall \ w \ge w^{\pi}_{r_{\ell}}
\end{align*}
Thus, \eqref{eqn:optimal intervals <=} holds for $\ell+1$.
Next we show that~\eqref{eqn:optimal intervals >} holds for $\ell+1$.

\noindent
\textbf{Case~1}: $w^{\pi}_{r_{\ell+1}}$ is determined by~\eqref{eqn:threshold >} or~\eqref{eqn:threshold =>}. \\
Then
\begin{align*}
& w^{\pi}_{r_{\ell+1}} \ \ \le \ \ w^{\min}_{\ell+1} \ \ = \ \ \min\left\{\frac{T^{\pi}_{r}(x(q)) - T^{\pi}_{r_{\ell+1}}(x(q))}{S^{\pi}_{r_{\ell+1}} - S^{\pi}_{r}} \; : \; S^{\pi}_{r_{\ell+1}} > S^{\pi}_{r}, \; r \in \hat{\mathcal{R}}^{\pi}\right\} \\
\Rightarrow \qquad & T^{\pi}_{r_{\ell+1}}(x(q)) + w^{\pi}_{r_{\ell+1}} S^{\pi}_{r_{\ell+1}} \ \ \le \ \ T^{\pi}_{r}(x(q)) + w^{\pi}_{r_{\ell+1}} S^{\pi}_{r} \qquad \forall \ r \in \hat{\mathcal{R}}^{\pi} \mbox{ s.t. } S^{\pi}_{r_{\ell+1}} > S^{\pi}_{r} \\
\Rightarrow \qquad & T^{\pi}_{r_{\ell+1}}(x(q)) + w S^{\pi}_{r_{\ell+1}} \ \ \le \ \ T^{\pi}_{r}(x(q)) + w S^{\pi}_{r} \qquad \forall \ r \in \hat{\mathcal{R}}^{\pi} \mbox{ s.t. } S^{\pi}_{r} < S^{\pi}_{r_{\ell+1}}, \ \forall \ w \le w^{\pi}_{r_{\ell+1}}
\end{align*}
and thus~\eqref{eqn:optimal intervals >} holds for $\ell+1$.

\noindent
\textbf{Case~2}: $w^{\pi}_{r_{\ell+1}}$ is determined by~\eqref{eqn:threshold ==}.

\noindent
\textbf{Case~2.1}: $q^{\pi}_{r_{\ell+1}} = 0$. \\
It was established in~\ref{item:positive flow paths are efficient} above that $q^{\pi}_{r} = 0$ for $r_{\ell} < r < r_{\ell+1}$.
Thus $p^{\pi}_{r_{\ell}} = p^{\pi}_{r_{\ell+1}}$.
It follows from the hypothesis that $w^{\pi}_{r_{\ell}} \in \left(G^{\pi}\right)^{-1}\left(p^{\pi}_{r_{\ell}}\right)$, and thus
\[
w^{\pi}_{r_{\ell+1}} \ \ = \ \ \max\left\{w^{\pi}_{r_{\ell}}, \; \min\left\{w \in \left(G^{\pi}\right)^{-1}\left(p^{\pi}_{r_{\ell+1}}\right)\right\}\right\}
\ \ = \ \ \max\left\{w^{\pi}_{r_{\ell}}, \; \min\left\{w \in \left(G^{\pi}\right)^{-1}\left(p^{\pi}_{r_{\ell}}\right)\right\}\right\}
\ \ = \ \ w^{\pi}_{r_{\ell}}
\]
It also follows from the hypothesis that
\begin{align*}
T^{\pi}_{r_{\ell+1}}(x(q)) + w^{\pi}_{r_{\ell}} S^{\pi}_{r_{\ell+1}} \ \ & = \ \ T^{\pi}_{r_{\ell}}(x(q)) + w^{\pi}_{r_{\ell}} S^{\pi}_{r_{\ell}}
\ \ \le \ \ T^{\pi}_{r}(x(q)) + w^{\pi}_{r_{\ell}} S^{\pi}_{r} \qquad \forall \ r \in \hat{\mathcal{R}}^{\pi} \\
\Rightarrow \quad T^{\pi}_{r_{\ell+1}}(x(q)) + w^{\pi}_{r_{\ell+1}} S^{\pi}_{r_{\ell+1}} \ \ & \le \ \ T^{\pi}_{r}(x(q)) + w^{\pi}_{r_{\ell+1}} S^{\pi}_{r} \qquad \forall \ r \in \hat{\mathcal{R}}^{\pi} \\
\Rightarrow \quad T^{\pi}_{r_{\ell+1}}(x(q)) + w S^{\pi}_{r_{\ell+1}} \ \ & \le \ \ T^{\pi}_{r}(x(q)) + w S^{\pi}_{r} \qquad \forall \ r \in \hat{\mathcal{R}}^{\pi} \mbox{ s.t. } S^{\pi}_{r} < S^{\pi}_{r_{\ell+1}}, \ \forall \ w \le w^{\pi}_{r_{\ell+1}}
\end{align*}
Thus, \eqref{eqn:optimal intervals >} holds for $\ell+1$.

\noindent
\textbf{Case~2.2}: $q^{\pi}_{r_{\ell+1}} > 0$. \\
We will show by contradiction that \eqref{eqn:optimal intervals >} holds for $\ell+1$.
First, suppose there is an $\hat{r} \in \hat{\mathcal{R}}^{\pi}$ such that $S^{\pi}_{\hat{r}} < S^{\pi}_{r_{\ell+1}}$ and $T^{\pi}_{r_{\ell+1}}(x(q)) + w^{\pi}_{r_{\ell+1}} S^{\pi}_{r_{\ell+1}} > T^{\pi}_{\hat{r}}(x(q)) + w^{\pi}_{r_{\ell+1}} S^{\pi}_{\hat{r}}$.
Thus $\hat{r} > r_{\ell+1}$.
It follows from the hypothesis that $w^{\pi}_{r_{\ell}} \in \left(G^{\pi}\right)^{-1}\left(p^{\pi}_{r_{\ell}}\right)$.
Also, $q^{\pi}_{r_{\ell+1}} > 0$ implies that $p^{\pi}_{r_{\ell+1}} > p^{\pi}_{r_{\ell}}$, and thus $w^{\pi}_{r_{\ell}} \le \min\left\{w \in \left(G^{\pi}\right)^{-1}\left(p^{\pi}_{r_{\ell+1}}\right)\right\}$.
Hence
\begin{equation}
w^{\pi}_{r_{\ell+1}} \ \ = \ \ \max\left\{w^{\pi}_{r_{\ell}}, \; \min\left\{w \in \left(G^{\pi}\right)^{-1}\left(p^{\pi}_{r_{\ell+1}}\right)\right\}\right\}
\ \ = \ \ \min\left\{w \in \left(G^{\pi}\right)^{-1}\left(p^{\pi}_{r_{\ell+1}}\right)\right\}
\ \ \le \ \ \min\left\{w \in \left(G^{\pi}\right)^{-1}\left(p^{\pi}_{r}\right)\right\}
\label{eqn:smallest quantile}
\end{equation}
for all $r \in \hat{\mathcal{R}}^{\pi}$ such that $r \ge r_{\ell+1}$.
Then consider the directional derivative of $\Upsilon$ at $q$ in the feasible direction $e^{\pi}_{\hat{r}} - e^{\pi}_{r_{\ell+1}}$:
\begin{align*}
& \max\Bigg\{\left(D^{\pi} \left[T^{\pi}_{\hat{r}}(x(q)) + \sum_{r=\hat{r}}^{\hat{R}^{\pi}-1} \left(S^{\pi}_{r} - S^{\pi}_{r+1}\right) \hat{w}^{\pi}_{r}\right] + \lambda^{\pi} - \mu^{\pi}_{\hat{r}}\right) \\
& \hspace{10mm} - \left(D^{\pi} \left[T^{\pi}_{r_{\ell+1}}(x(q)) + \sum_{r=r_{\ell+1}}^{\hat{R}^{\pi}-1} \left(S^{\pi}_{r} - S^{\pi}_{r+1}\right) \hat{w}^{\pi}_{r}\right] + \lambda^{\pi}\right) \; : \; \hat{w}^{\pi}_{1} \in \left(G^{\pi}\right)^{-1}\left(p^{\pi}_{1}\right), \ldots, \hat{w}^{\pi}_{\hat{R}^{\pi}-1} \in \left(G^{\pi}\right)^{-1}\left(p^{\pi}_{\hat{R}^{\pi}-1}\right), \\
& \hspace{145mm} \mu^{\pi}_{\hat{r}} \ge 0, \; \mu^{\pi}_{\hat{r}} q^{\pi}_{\hat{r}} = 0\Bigg\} \\
& = \ \ \max\Bigg\{D^{\pi} \left[T^{\pi}_{\hat{r}}(x(q)) - T^{\pi}_{r_{\ell+1}}(x(q)) - \sum_{r=r_{\ell+1}}^{\hat{r}-1} \left(S^{\pi}_{r} - S^{\pi}_{r+1}\right) \hat{w}^{\pi}_{r}\right] \; : \; \hat{w}^{\pi}_{r_{\ell+1}} \in \left(G^{\pi}\right)^{-1}\left(p^{\pi}_{r_{\ell+1}}\right), \ldots \\
& \hspace{130mm} \ldots, \hat{w}^{\pi}_{\hat{r}-1} \in \left(G^{\pi}\right)^{-1}\left(p^{\pi}_{\hat{r}-1}\right)\Bigg\} \\
& \overset{\eqref{eqn:smallest quantile}}{\le} \ \ D^{\pi} \left[T^{\pi}_{\hat{r}}(x(q)) - T^{\pi}_{r_{\ell+1}}(x(q)) - \sum_{r=r_{\ell+1}}^{\hat{r}-1} \left(S^{\pi}_{r} - S^{\pi}_{r+1}\right) w^{\pi}_{r_{\ell+1}}\right] \\
& = \ \ D^{\pi} \left[T^{\pi}_{\hat{r}}(x(q)) - T^{\pi}_{r_{\ell+1}}(x(q)) + w^{\pi}_{r_{\ell+1}} \left(S^{\pi}_{\hat{r}} - S^{\pi}_{r_{\ell+1}}\right)\right]
\ \ < \ \ 0
\end{align*}
which contradicts the optimality of $q$.
Thus, for all $r \in \hat{\mathcal{R}}^{\pi}$ such that $S^{\pi}_{r} < S^{\pi}_{r_{\ell+1}}$ it holds that
\begin{align*}
T^{\pi}_{r_{\ell+1}}(x(q)) + w^{\pi}_{r_{\ell+1}} S^{\pi}_{r_{\ell+1}} \ \ & \le \ \ T^{\pi}_{r}(x(q)) + w^{\pi}_{r_{\ell+1}} S^{\pi}_{r} \\
\Rightarrow \qquad T^{\pi}_{r_{\ell+1}}(x(q)) + w S^{\pi}_{r_{\ell+1}} \ \ & \le \ \ T^{\pi}_{r}(x(q)) + w S^{\pi}_{r} \qquad \forall \ w \le w^{\pi}_{r_{\ell+1}}
\end{align*}
Thus, \eqref{eqn:optimal intervals >} holds for $\ell+1$.

\noindent
\textbf{Case~3}: $w^{\pi}_{r_{\ell+1}}$ is determined by~\eqref{eqn:threshold =0}. \\
Then $\left\{r \in \hat{\mathcal{R}}^{\pi} \; : \; S^{\pi}_{r_{\ell+1}} > S^{\pi}_{r}\right\} = \varnothing$, and thus \eqref{eqn:optimal intervals >} holds for $\ell+1$.
\end{enumerate}
Thus the hypothesis has been established.
It remains to verify that~\ref{item:positive flow paths are efficient} and~\ref{item:thresholds quantiles} hold for $r > r_{L}$.

\begin{enumerate}[label=(\arabic*), ref=(\arabic*), start=2, leftmargin=7mm]
\item
In Algorithm~\ref{algo:threshold}, when $\ell = L$, $w^{\min}_{L}$ and $w^{\pi}_{r_{L}}$ are determined by~\eqref{eqn:threshold >} or~\eqref{eqn:threshold =>} or~\eqref{eqn:threshold =0}, and $w^{\min}_{L} > \overline{w}^{\pi}$.
Consider any $r > r_{L}$ and suppose that $q^{\pi}_{r} > 0$.
Then for any $r' < r$ it holds that $p^{\pi}_{r'} < 1$, and thus any $\hat{w} \in \left(G^{\pi}\right)^{-1}\left(p^{\pi}_{r'}\right)$ satisfies $\hat{w} \le \overline{w}^{\pi}$.

\noindent
\textbf{Case~1}: $w^{\pi}_{r_{L}}$ is determined by~\eqref{eqn:threshold >} or~\eqref{eqn:threshold =>} and $S^{\pi}_{r_{L}} > S^{\pi}_{r}$. \\
Then
\begin{align}
& \overline{w}^{\pi} \ \ < \ \ w^{\min}_{L} \ \ = \ \ \min\left\{\frac{T^{\pi}_{r'}(x(q)) - T^{\pi}_{r_{L}}(x(q))}{S^{\pi}_{r_{L}} - S^{\pi}_{r'}} \; : \; S^{\pi}_{r_{L}} > S^{\pi}_{r'}, \; r' \in \hat{\mathcal{R}}^{\pi}\right\} \ \ \le \ \ \frac{T^{\pi}_{r}(x(q)) - T^{\pi}_{r_{L}}(x(q))}{S^{\pi}_{r_{L}} - S^{\pi}_{r}} \nonumber \\
\Rightarrow & \quad T^{\pi}_{r_{L}}(x(q)) + \overline{w}^{\pi} S^{\pi}_{r_{L}} \ \ < \ \ T^{\pi}_{r}(x(q)) + \overline{w}^{\pi} S^{\pi}_{r}
\label{eqn:least cost interval 5}
\end{align}
Next, consider the directional derivative of $\Upsilon$ at $q$ in the feasible direction $e^{\pi}_{r_{L}} - e^{\pi}_{r}$:
\begin{align*}
\hspace{-5mm}
& \max\Bigg\{\left(D^{\pi} \left[T^{\pi}_{r_{L}}(x(q)) + \sum_{r'=r_{L}}^{\hat{R}^{\pi}-1} \left(S^{\pi}_{r'} - S^{\pi}_{r'+1}\right) \hat{w}^{\pi}_{r'}\right] + \lambda^{\pi} - \mu^{\pi}_{r_{L}}\right) \\
\hspace{-5mm}
& \hspace{10mm} - \left(D^{\pi} \left[T^{\pi}_{r}(x(q)) + \sum_{r'=r}^{\hat{R}^{\pi}-1} \left(S^{\pi}_{r'} - S^{\pi}_{r'+1}\right) \hat{w}^{\pi}_{r'}\right] + \lambda^{\pi}\right) \; : \; \hat{w}^{\pi}_{1} \in \left(G^{\pi}\right)^{-1}\left(p^{\pi}_{1}\right), \ldots, \hat{w}^{\pi}_{\hat{R}^{\pi}-1} \in \left(G^{\pi}\right)^{-1}\left(p^{\pi}_{\hat{R}^{\pi}-1}\right), \\
\hspace{-5mm}
& \hspace{145mm} \mu^{\pi}_{r_{L}} \ge 0, \; \mu^{\pi}_{r_{L}} q^{\pi}_{r_{L}} = 0\Bigg\} \\
\hspace{-5mm}
& = \ \ \max\Bigg\{D^{\pi} \left[T^{\pi}_{r_{L}}(x(q)) - T^{\pi}_{r}(x(q)) + \sum_{r'=r_{L}}^{r-1} \left(S^{\pi}_{r'} - S^{\pi}_{r'+1}\right) \hat{w}^{\pi}_{r'}\right] \; : \; \hat{w}^{\pi}_{r_{L}} \in \left(G^{\pi}\right)^{-1}\left(p^{\pi}_{r_{L}}\right), \ldots, \hat{w}^{\pi}_{r-1} \in \left(G^{\pi}\right)^{-1}\left(p^{\pi}_{r-1}\right)\Bigg\} \\
\hspace{-5mm}
& \le \ \ D^{\pi} \left[T^{\pi}_{r_{L}}(x(q)) - T^{\pi}_{r}(x(q)) + \sum_{r'=r_{L}}^{r-1} \left(S^{\pi}_{r'} - S^{\pi}_{r'+1}\right) \overline{w}^{\pi}\right] \\
\hspace{-5mm}
& = \ \ D^{\pi} \left[T^{\pi}_{r_{L}}(x(q)) - T^{\pi}_{r}(x(q)) + \left(S^{\pi}_{r_{L}} - S^{\pi}_{r}\right) \overline{w}^{\pi}\right]
\ \ \overset{\eqref{eqn:least cost interval 5}}{<} \ \ 0
\end{align*}
which contradicts the optimality of $q$.

\noindent
\textbf{Case~2}: $w^{\pi}_{r_{L}}$ is determined by~\eqref{eqn:threshold =>} and $S^{\pi}_{r_{L}} = S^{\pi}_{r}$, or by~\eqref{eqn:threshold =0}. \\
Then $T^{\pi}_{r_{L}}(x(q)) < T^{\pi}_{r}(x(q))$ and $S^{\pi}_{r_{L}} = S^{\pi}_{r'} = S^{\pi}_{r}$ for all $r' = r_{L},\ldots,r$.
Next, consider the directional derivative of $\Upsilon$ at $q$ in the feasible direction $e^{\pi}_{r_{L}} - e^{\pi}_{r}$: 
\begin{align*}
\hspace{-5mm}
& \max\Bigg\{\left(D^{\pi} \left[T^{\pi}_{r_{L}}(x(q)) + \sum_{r'=r_{L}}^{\hat{R}^{\pi}-1} \left(S^{\pi}_{r'} - S^{\pi}_{r'+1}\right) \hat{w}^{\pi}_{r'}\right] + \lambda^{\pi} - \mu^{\pi}_{r_{L}}\right) \\
\hspace{-5mm}
& \hspace{10mm} - \left(D^{\pi} \left[T^{\pi}_{r}(x(q)) + \sum_{r'=r}^{\hat{R}^{\pi}-1} \left(S^{\pi}_{r'} - S^{\pi}_{r'+1}\right) \hat{w}^{\pi}_{r'}\right] + \lambda^{\pi}\right) \; : \; \hat{w}^{\pi}_{1} \in \left(G^{\pi}\right)^{-1}\left(p^{\pi}_{1}\right), \ldots, \hat{w}^{\pi}_{\hat{R}^{\pi}-1} \in \left(G^{\pi}\right)^{-1}\left(p^{\pi}_{\hat{R}^{\pi}-1}\right), \\
\hspace{-5mm}
& \hspace{145mm} \mu^{\pi}_{r_{L}} \ge 0, \; \mu^{\pi}_{r_{L}} q^{\pi}_{r_{L}} = 0\Bigg\} \\
\hspace{-5mm}
& = \ \ \max\Bigg\{D^{\pi} \left[T^{\pi}_{r_{L}}(x(q)) - T^{\pi}_{r}(x(q)) + \sum_{r'=r_{L}}^{r-1} \left(S^{\pi}_{r'} - S^{\pi}_{r'+1}\right) \hat{w}^{\pi}_{r'}\right] \; : \; \hat{w}^{\pi}_{r_{L}} \in \left(G^{\pi}\right)^{-1}\left(p^{\pi}_{r_{L}}\right), \ldots, \hat{w}^{\pi}_{r-1} \in \left(G^{\pi}\right)^{-1}\left(p^{\pi}_{r-1}\right)\Bigg\} \\
\hspace{-5mm}
& = \ \ D^{\pi} \left[T^{\pi}_{r_{L}}(x(q)) - T^{\pi}_{r}(x(q))\right]
\ \ < \ \ 0
\end{align*}
which contradicts the optimality of $q$.
Thus if $r > r_{L}$ then $q^{\pi}_{r} = 0$.

\item
It was established by induction that $w^{\pi}_{r_{L}} \in \left(G^{\pi}\right)^{-1}\left(p^{\pi}_{r_{L}}\right)$.
By Algorithm~\ref{algo:threshold}, $w^{\pi}_{r} = w^{\pi}_{r_{L}} = \overline{w}^{\pi}$ for all $r \ge r_{L}$.
It was established in \ref{item:positive flow paths are efficient} above that $q^{\pi}_{r} = 0$ for all $r > r_{L}$, and thus $p^{\pi}_{r} = p^{\pi}_{r_{L}} = 1$ and $w^{\pi}_{r} \in \left(G^{\pi}\right)^{-1}\left(p^{\pi}_{r}\right)$ for all $r > r_{L}$.
\end{enumerate}
\Halmos
\endproof

\begin{manuallemma}{\ref{lem:thresholds optimality conditions}}
Assume that $t_{a} : \mathbb{R}_{+} \mapsto \mathbb{R}_{+}$ is continuous for all $a \in \mathcal{A}$.
Consider any optimal solution $q$ for problem~\eqref{eqn:q reformulated BiTAP objective}, and the resulting output $w$ and $\mathcal{R}^{\pi}_{\min} = \{r_{1},\ldots,r_{L}\}$ of Algorithm~\ref{algo:threshold}.
Let $\hat{\mathcal{R}}^{\pi}$ be indexed with $r = 1,\ldots,\hat{R}^{\pi}$, and let
\begin{equation*}
\lambda^{\pi} \ \ \defi \ \ -D^{\pi} \left[T^{\pi}_{r_{L}}(x(q)) + \overline{w}^{\pi} \left(S^{\pi}_{r_{L}} - S^{\pi}_{\hat{R}^{\pi}}\right)\right]
\qquad \forall \ \pi \in \Pi
\end{equation*}
Then $(q,\lambda,w)$ satisfies the optimality conditions~\eqref{eqn:q continuous BiTAP optimality1}--\eqref{eqn:q continuous BiTAP optimality3}.
\end{manuallemma}

\proof{Proof.}
Consider any $\pi \in \Pi$.
Let $p^{\pi}_{r} = L^{\pi}_{r}(q^{\pi})$.
Consider any successive pair $r_{\ell},r_{\ell+1} \in \mathcal{R}^{\pi}_{\min}$.
If $w^{\pi}_{r_{\ell}}$ is determined by~\eqref{eqn:threshold =0} in Algorithm~\ref{algo:threshold}, then $\ell = L$, and thus this case does not apply here.
Similarly, if $w^{\pi}_{r_{\ell}}$ is determined by~\eqref{eqn:threshold >} or~\eqref{eqn:threshold =>}, and $w^{\min}_{\ell} > \overline{w}^{\pi}$, then $\ell = L$, and thus these cases also do not apply here.
Otherwise, $w^{\pi}_{r_{\ell}}$ is determined by~\eqref{eqn:threshold >} or~\eqref{eqn:threshold =>}, and $w^{\pi}_{r_{\ell}} = w^{\min}_{\ell} \le \overline{w}^{\pi}$, or $w^{\pi}_{r_{\ell}}$ is determined by~\eqref{eqn:threshold ==}.
In all of these cases it holds that
\begin{equation}
\label{eqn:successive paths}
T^{\pi}_{r_{\ell}}(x(q)) + w^{\pi}_{r_{\ell}} S^{\pi}_{r_{\ell}} \ \ = \ \ T^{\pi}_{r_{\ell+1}}(x(q)) + w^{\pi}_{r_{\ell}} S^{\pi}_{r_{\ell+1}}
\end{equation}
that is, $w^{\pi}_{r_{\ell}}$ is a boundary or indifference point between the interval $[w^{\pi}_{r_{\ell}-1},w^{\pi}_{r_{\ell}}]$ where $r_{\ell}$ has least bicriteria cost and the interval $[w^{\pi}_{r_{\ell}},w^{\pi}_{r_{\ell+1}}]$ where $r_{\ell+1}$ has least bicriteria cost.

Also, note that $w^{\pi}_{r} = w^{\pi}_{r_{\ell}}$ for all $r \in \{r_{\ell},\ldots,r_{\ell+1}-1\}$, and thus
\begin{equation}
\sum_{r=r_{\ell}}^{r_{\ell+1}-1} \left(S^{\pi}_{r} - S^{\pi}_{r+1}\right) w^{\pi}_{r}
\ \ = \ \ \sum_{r=r_{\ell}}^{r_{\ell+1}-1} \left(S^{\pi}_{r} - S^{\pi}_{r+1}\right) w^{\pi}_{r_{\ell}}
\ \ = \ \ \left(S^{\pi}_{r_{\ell}} - S^{\pi}_{r_{\ell+1}}\right) w^{\pi}_{r_{\ell}}
\label{eqn:cancellation unused paths}
\end{equation}
Similarly, $w^{\pi}_{r} = w^{\pi}_{r_{L}} = \overline{w}^{\pi}$ for all $r \in \{r_{L},\ldots,\hat{R}^{\pi}\}$, and thus
\begin{equation}
\sum_{r=r_{L}}^{\hat{R}^{\pi}-1} \left(S^{\pi}_{r} - S^{\pi}_{r+1}\right) w^{\pi}_{r}
\ \ = \ \ \sum_{r=r_{L}}^{\hat{R}^{\pi}-1} \left(S^{\pi}_{r} - S^{\pi}_{r+1}\right) \overline{w}^{\pi}
\ \ = \ \ \left(S^{\pi}_{r_{L}} - S^{\pi}_{\hat{R}^{\pi}}\right) \overline{w}^{\pi}
\label{eqn:cancellation unused paths2}
\end{equation}

\noindent
\textbf{Case~1}: $r \in \mathcal{R}^{\pi}_{\min}$. \\
Let $r = r_{\ell}$ for some $\ell \in \{1,\ldots,L\}$.
Then
\begin{align*}
& D^{\pi} \left[T^{\pi}_{r}(x(q)) + \sum_{\{r' \in \hat{\mathcal{R}}^{\pi} \, : \, r \preceq r', r'+1 \in \hat{\mathcal{R}}^{\pi}\}} \left(S^{\pi}_{r'} - S^{\pi}_{r'+1}\right) w^{\pi}_{r'}\right] + \lambda^{\pi} \\
& \overset{\eqref{eqn:choice of dual},\eqref{eqn:cancellation unused paths},\eqref{eqn:cancellation unused paths2}}{=} \ \ D^{\pi} \left[T^{\pi}_{r_{\ell}}(x(q)) + \sum_{\ell' = \ell}^{L-1} \left(S^{\pi}_{r_{\ell'}} - S^{\pi}_{r_{\ell'+1}}\right) w^{\pi}_{r_{\ell'}} + \left(S^{\pi}_{r_{L}} - S^{\pi}_{\hat{R}^{\pi}}\right) \overline{w}^{\pi}\right] - D^{\pi} \left[T^{\pi}_{r_{L}}(x(q)) + \overline{w}^{\pi} \left(S^{\pi}_{r_{L}} - S^{\pi}_{\hat{R}^{\pi}}\right)\right] \\
& = \ \ D^{\pi} \sum_{\ell' = \ell}^{L-1} \left[T^{\pi}_{r_{\ell'}}(x(q)) - T^{\pi}_{r_{\ell'+1}}(x(q)) + \left(S^{\pi}_{r_{\ell'}} - S^{\pi}_{r_{\ell'+1}}\right) w^{\pi}_{r_{\ell'}}\right] \ \ \overset{\eqref{eqn:successive paths}}{=} \ \ 0
\end{align*}
Thus conditions~\eqref{eqn:q continuous BiTAP optimality1}--\eqref{eqn:q continuous BiTAP optimality3} hold for all $r \in \mathcal{R}^{\pi}_{\min}$.

\noindent
\textbf{Case~2}: $r \notin \mathcal{R}^{\pi}_{\min}$. \\
It follows from Lemma~\ref{lem:optimal thresholds}\ref{item:positive flow paths are efficient} that $q^{\pi}_{r} = 0$, and thus condition~\eqref{eqn:q continuous BiTAP optimality1} holds.
Next we show that condition~\eqref{eqn:q continuous BiTAP optimality3} also holds.

\noindent
\textbf{Case~2.1}: $r < r_{L}$. \\
Then either $r < r_{1}$, or there is a successive pair $r_{\ell},r_{\ell+1} \in \mathcal{R}^{\pi}_{\min}$ such that $r_{\ell} < r < r_{\ell+1}$.
If $r < r_{1}$, then $\ell = 0$ and $w^{\pi}_{r'} = w^{\pi}_{r_{0}} = \underline{w}^{\pi}$ for all $r' \in \{r,\ldots,r_{1}-1\}$, and if $r_{\ell} < r < r_{\ell+1}$, then $w^{\pi}_{r'} = w^{\pi}_{r_{\ell}}$ for all $r' \in \{r,\ldots,r_{\ell+1}-1\}$.
Then
\begin{align*}
& D^{\pi} \left[T^{\pi}_{r}(x(q)) + \sum_{\{r' \in \hat{\mathcal{R}}^{\pi} \, : \, r \preceq r', r'+1 \in \hat{\mathcal{R}}^{\pi}\}} \left(S^{\pi}_{r'} - S^{\pi}_{r'+1}\right) w^{\pi}_{r'}\right] + \lambda^{\pi} \\
& \overset{\eqref{eqn:cancellation unused paths},\eqref{eqn:cancellation unused paths2}}{=} \ \ D^{\pi} \left[T^{\pi}_{r}(x(q)) + \sum_{r'=r}^{r_{\ell+1}-1} \left(S^{\pi}_{r'} - S^{\pi}_{r'+1}\right) w^{\pi}_{r_{\ell}} + \sum_{\ell' = \ell+1}^{L-1} \left(S^{\pi}_{r_{\ell'}} - S^{\pi}_{r_{\ell'+1}}\right) w^{\pi}_{r_{\ell'}} + \left(S^{\pi}_{r_{L}} - S^{\pi}_{\hat{R}^{\pi}}\right) \overline{w}^{\pi}\right] \\
& \hspace{15mm} {} - D^{\pi} \left[T^{\pi}_{r_{L}}(x(q)) + \overline{w}^{\pi} \left(S^{\pi}_{r_{L}} - S^{\pi}_{\hat{R}^{\pi}}\right)\right] \\
& = \ \ D^{\pi} \left[T^{\pi}_{r}(x(q)) - T^{\pi}_{r_{\ell+1}}(x(q)) + T^{\pi}_{r_{\ell+1}}(x(q)) - T^{\pi}_{r_{L}}(x(q)) + \left(S^{\pi}_{r} - S^{\pi}_{r_{\ell+1}}\right) w^{\pi}_{r_{\ell}} + \sum_{\ell' = \ell+1}^{L-1} \left(S^{\pi}_{r_{\ell'}} - S^{\pi}_{r_{\ell'+1}}\right) w^{\pi}_{r_{\ell'}}\right] \\
& = \ \ D^{\pi} \left[T^{\pi}_{r}(x(q)) - T^{\pi}_{r_{\ell+1}}(x(q)) + \left(S^{\pi}_{r} - S^{\pi}_{r_{\ell+1}}\right) w^{\pi}_{r_{\ell}}\right] + D^{\pi} \sum_{\ell' = \ell+1}^{L-1} \left[T^{\pi}_{r_{\ell'}}(x(q)) - T^{\pi}_{r_{\ell'+1}}(x(q)) + \left(S^{\pi}_{r_{\ell'}} - S^{\pi}_{r_{\ell'+1}}\right) w^{\pi}_{r_{\ell'}}\right] \\
& \overset{\eqref{eqn:successive paths}}{=} \ \ D^{\pi} \left[T^{\pi}_{r}(x(q)) - T^{\pi}_{r_{\ell+1}}(x(q)) + \left(S^{\pi}_{r} - S^{\pi}_{r_{\ell+1}}\right) w^{\pi}_{r_{\ell}}\right]
\end{align*}

\noindent
\textbf{Case~2.1.0}: $r < r_{1}$. \\
Then $\ell = 0$ and $w^{\pi}_{r_{\ell}} = \underline{w}^{\pi}$.
Thus
\[
T^{\pi}_{r}(x(q)) - T^{\pi}_{r_{\ell+1}}(x(q)) + \left(S^{\pi}_{r} - S^{\pi}_{r_{\ell+1}}\right) w^{\pi}_{r_{\ell}}
\ \ = \ \ T^{\pi}_{r}(x(q)) - T^{\pi}_{r_{1}}(x(q)) + \left(S^{\pi}_{r} - S^{\pi}_{r_{1}}\right) \underline{w}^{\pi}
\ \ \overset{\eqref{eqn:efficient path 1}}{\ge} \ \ 0
\]

\noindent
\textbf{Case~2.1.1}: $r_{1} < r < r_{L}$. \\
Then there is a successive pair $r_{\ell},r_{\ell+1} \in \mathcal{R}^{\pi}_{\min}$ such that $r_{\ell} < r < r_{\ell+1}$.

\noindent
\textbf{Case~2.1.1.1}: $w^{\pi}_{r_{\ell}}$ is determined by~\eqref{eqn:threshold >} or~\eqref{eqn:threshold =>} and $S^{\pi}_{r_{\ell}} > S^{\pi}_{r}$. \\
Then
\begin{align*}
& w^{\pi}_{r_{\ell}} \ \ = \ \ w^{\min}_{\ell} \ \ = \ \ \frac{T^{\pi}_{r_{\ell+1}}(x(q)) - T^{\pi}_{r_{\ell}}(x(q))}{S^{\pi}_{r_{\ell}} - S^{\pi}_{r_{\ell+1}}}
\ \ < \ \ \frac{T^{\pi}_{r}(x(q)) - T^{\pi}_{r_{\ell}}(x(q))}{S^{\pi}_{r_{\ell}} - S^{\pi}_{r}} \\
\Rightarrow \qquad & T^{\pi}_{r_{\ell+1}}(x(q)) + w^{\pi}_{r_{\ell}} S^{\pi}_{r_{\ell+1}} \ \ < \ \ T^{\pi}_{r}(x(q)) + w^{\pi}_{r_{\ell}} S^{\pi}_{r} \\
\Rightarrow \qquad & D^{\pi} \left[T^{\pi}_{r}(x(q)) - T^{\pi}_{r_{\ell+1}}(x(q)) + \left(S^{\pi}_{r} - S^{\pi}_{r_{\ell+1}}\right) w^{\pi}_{r_{\ell}}\right]
\ \ > \ \ 0
\end{align*}
(The strict inequality holds because $r_{\ell+1}$ is chosen in~\eqref{eqn:efficient path >} and~\eqref{eqn:efficient path =>} as the \emph{smallest} minimizer.)

\noindent
\textbf{Case~2.1.1.2}: $w^{\pi}_{r_{\ell}}$ is determined by~\eqref{eqn:threshold =>} and $S^{\pi}_{r_{\ell}} = S^{\pi}_{r}$. \\
Then $T^{\pi}_{r_{\ell}}(x(q)) < T^{\pi}_{r}(x(q))$.
Thus
\begin{align*}
& D^{\pi} \left[T^{\pi}_{r}(x(q)) - T^{\pi}_{r_{\ell+1}}(x(q)) + \left(S^{\pi}_{r} - S^{\pi}_{r_{\ell+1}}\right) w^{\pi}_{r_{\ell}}\right] \\
= \ \ & D^{\pi} \left[T^{\pi}_{r}(x(q)) - T^{\pi}_{r_{\ell}}(x(q)) + T^{\pi}_{r_{\ell}}(x(q)) - T^{\pi}_{r_{\ell+1}}(x(q)) + \left(S^{\pi}_{r_{\ell}} - S^{\pi}_{r_{\ell+1}}\right) w^{\pi}_{r_{\ell}}\right] \\
\overset{\eqref{eqn:successive paths}}{=} \ \ & D^{\pi} \left[T^{\pi}_{r}(x(q)) - T^{\pi}_{r_{\ell}}(x(q))\right]
\ \ > \ \ 0
\end{align*}

\noindent
\textbf{Case~2.1.1.3}: $w^{\pi}_{r_{\ell}}$ is determined by~\eqref{eqn:threshold ==}. \\
Then $S^{\pi}_{r_{\ell}} = S^{\pi}_{r} = S^{\pi}_{r_{\ell+1}}$ and $T^{\pi}_{r_{\ell}}(x(q)) = T^{\pi}_{r_{\ell+1}}(x(q)) < T^{\pi}_{r}(x(q))$, and thus
\[
D^{\pi} \left[T^{\pi}_{r}(x(q)) - T^{\pi}_{r_{\ell+1}}(x(q)) + \left(S^{\pi}_{r} - S^{\pi}_{r_{\ell+1}}\right) w^{\pi}_{r_{\ell}}\right]
\ \ = \ \ D^{\pi} \left[T^{\pi}_{r}(x(q)) - T^{\pi}_{r_{\ell+1}}(x(q))\right]
\ \ > \ \ 0
\]
Therefore condition~\eqref{eqn:q continuous BiTAP optimality3} holds for all $r < r_{L}$.

\noindent
\textbf{Case~2.2}: $r > r_{L}$. \\
Then
\begin{align*}
& D^{\pi} \left[T^{\pi}_{r}(x(q)) + \sum_{\{r' \in \hat{\mathcal{R}}^{\pi} \, : \, r \preceq r', r'+1 \in \hat{\mathcal{R}}^{\pi}\}} \left(S^{\pi}_{r'} - S^{\pi}_{r'+1}\right) w^{\pi}_{r'}\right] + \lambda^{\pi} \\
& = \ \ D^{\pi} \left[T^{\pi}_{r}(x(q)) + \sum_{r'=r}^{\hat{R}^{\pi}-1} \left(S^{\pi}_{r'} - S^{\pi}_{r'+1}\right) \overline{w}^{\pi}\right] - D^{\pi} \left[T^{\pi}_{r_{L}}(x(q)) + \overline{w}^{\pi} \left(S^{\pi}_{r_{L}} - S^{\pi}_{\hat{R}^{\pi}}\right)\right] \\
& = \ \ D^{\pi} \left[T^{\pi}_{r}(x(q)) - T^{\pi}_{r_{L}}(x(q)) + \left(S^{\pi}_{r} - S^{\pi}_{r_{L}}\right) \overline{w}^{\pi}\right]
\end{align*}

\noindent
\textbf{Case~2.2.1}: $w^{\pi}_{r_{L}}$ is determined by~\eqref{eqn:threshold >} or~\eqref{eqn:threshold =>} and $S^{\pi}_{r_{L}} > S^{\pi}_{r}$. \\
Then
\begin{align*}
& w^{\pi}_{r_{L}} \ \ = \ \ \overline{w}^{\pi} \ \ < \ \ w^{\min}_{L} \ \ = \ \ \min\left\{\frac{T^{\pi}_{r'}(x(q)) - T^{\pi}_{r_{L}}(x(q))}{S^{\pi}_{r_{L}} - S^{\pi}_{r'}} \; : \; S^{\pi}_{r_{L}} > S^{\pi}_{r'}\right\}
\ \ \le \ \ \frac{T^{\pi}_{r}(x(q)) - T^{\pi}_{r_{L}}(x(q))}{S^{\pi}_{r_{L}} - S^{\pi}_{r}} \\
\Rightarrow \qquad & T^{\pi}_{r_{L}}(x(q)) + \overline{w}^{\pi} S^{\pi}_{r_{L}} \ \ < \ \ T^{\pi}_{r}(x(q)) + \overline{w}^{\pi} S^{\pi}_{r} \\
\Rightarrow \qquad & D^{\pi} \left[T^{\pi}_{r}(x(q)) - T^{\pi}_{r_{L}}(x(q)) + \left(S^{\pi}_{r} - S^{\pi}_{r_{L}}\right) \overline{w}^{\pi}\right]
\ \ > \ \ 0
\end{align*}

\noindent
\textbf{Case~2.2.2}: $w^{\pi}_{r_{L}}$ is determined by~\eqref{eqn:threshold =>} and $S^{\pi}_{r_{L}} = S^{\pi}_{r}$. \\
Then $T^{\pi}_{r_{L}}(x(q)) < T^{\pi}_{r}(x(q))$.
Thus
\begin{align*}
D^{\pi} \left[T^{\pi}_{r}(x(q)) - T^{\pi}_{r_{L}}(x(q)) + \left(S^{\pi}_{r} - S^{\pi}_{r_{L}}\right) \overline{w}^{\pi}\right]
\ \ = \ \ D^{\pi} \left[T^{\pi}_{r}(x(q)) - T^{\pi}_{r_{L}}(x(q))\right]
\ \ > \ \ 0
\end{align*}

\noindent
\textbf{Case~2.2.3}: $w^{\pi}_{r_{L}}$ is determined by~\eqref{eqn:threshold ==}. \\
Then $r_{L+1}$ set in~\eqref{eqn:efficient path ==} would satisfy $r_{L+1} > r_{L}$ and $r_{L+1} \in \mathcal{R}^{\pi}_{\min}$, which cannot hold, and thus $w^{\pi}_{r_{L}}$ is not determined by~\eqref{eqn:threshold ==}.

\noindent
\textbf{Case~2.2.4}: $w^{\pi}_{r_{L}}$ is determined by~\eqref{eqn:threshold =0}. \\
Then $S^{\pi}_{r_{L}} = S^{\pi}_{r} = S^{\pi}_{\hat{R}^{\pi}}$ and $T^{\pi}_{r}(x(q)) > T^{\pi}_{r_{L}}(x(q))$.
Thus
\begin{align*}
D^{\pi} \Bigg[T^{\pi}_{r}(x(q)) - T^{\pi}_{r_{L}}(x(q)) + \left(S^{\pi}_{r} - S^{\pi}_{r_{L}}\right) \overline{w}^{\pi}\Bigg]
\ \ = \ \ & D^{\pi} \left[T^{\pi}_{r}(x(q)) - T^{\pi}_{r_{L}}(x(q))\right]
\ \ > \ \ 0
\end{align*}
Therefore conditions~\eqref{eqn:q continuous BiTAP optimality1}--\eqref{eqn:q continuous BiTAP optimality3} hold for all $r \in \hat{\mathcal{R}}^{\pi}$.
\Halmos
\endproof

\begin{manuallemma}{\ref{lem:optimality conditions}}
Assume that $t_{a} : \mathbb{R}_{+} \mapsto \mathbb{R}_{+}$ is continuous for all $a \in \mathcal{A}$.
For each $\pi \in \Pi$, let $\mathcal{R}^{\pi}$ denote the complete path set and let $\hat{\mathcal{R}}^{\pi} \subset \mathcal{R}^{\pi}$.
Consider any optimal solution $\hat{q}$ for problem~\eqref{eqn:q reformulated BiTAP objective} with path sets $\hat{\mathcal{R}}^{\pi}$, $\pi \in \Pi$.
Let $\hat{w} = \left(\hat{w}^{\pi}_{r}, \pi \in \Pi, r \in \{0\} \cup \hat{\mathcal{R}}^{\pi}\right)$ and $\hat{\mathcal{R}}^{\pi}_{\min} = \{\hat{r}_{1},\ldots,\hat{r}_{\hat{L}}\}$ denote the output of Algorithm~\ref{algo:threshold} with input $\hat{\mathcal{R}}^{\pi}$ and $\hat{q}$.

Let $q^* = \left(q^{\pi*}_{r}, \pi \in \Pi, r \in \mathcal{R}^{\pi}\right)$ be as follows:
For each $r \in \hat{\mathcal{R}}^{\pi}$ let $q^{\pi*}_{r} = \hat{q}^{\pi}_{r}$, and for each $r \in \mathcal{R}^{\pi} \setminus \hat{\mathcal{R}}^{\pi}$ let $q^{\pi*}_{r} = 0$.
Then $x(q^*) = x(\hat{q})$.

Let $w^* = \left(w^{\pi*}_{r}, \pi \in \Pi, r \in \{0\} \cup \mathcal{R}^{\pi}\right)$ and $\mathcal{R}^{\pi*}_{\min} = \{r^*_{1},\ldots,r^*_{L^*}\}$ denote the output of Algorithm~\ref{algo:threshold} with input $\mathcal{R}^{\pi}$ and $q^*$.
Let $\mathcal{R}^{\pi}$ be indexed with $r = 1,\ldots,R^{\pi}$, and let
\begin{equation*}
\lambda^{\pi*} \ \ \defi \ \ -D^{\pi} \left[T^{\pi}_{r^*_{L^*}}(x(q^*)) + \overline{w}^{\pi} \left(S^{\pi}_{r^*_{L^*}} - S^{\pi}_{R^{\pi}}\right)\right]
\qquad \forall \ \pi \in \Pi
\end{equation*}
For each $\pi \in \Pi$, let $\mathcal{R}^{\pi+}$ denote an optimal set of paths for the parametric shortest path problem (over the complete path set $\mathcal{R}^{\pi}$) with arc costs $t_{a}(x_{a}(\hat{q}))$ and $s_{a}$.
Suppose that $\mathcal{R}^{\pi+} \subset \hat{\mathcal{R}}^{\pi}$ for each $\pi \in \Pi$.
Then it holds that $(q^*,\lambda^*,w^*)$ satisfies the optimality conditions for the complete path set:
\begin{align*}
& q^{\pi*}_{r} \left\{D^{\pi} \left[T^{\pi}_{r}(x(q^*)) + \sum_{r' = r}^{R^{\pi}-1} \left(S^{\pi}_{r'} - S^{\pi}_{r'+1}\right) w^{\pi*}_{r'}\right] + \lambda^{\pi*}\right\} \ \ = \ \ 0 \qquad \forall \ \pi \in \Pi, \ \forall \ r \in \mathcal{R}^{\pi} \\
& D^{\pi} \left[T^{\pi}_{r}(x(q^*)) + \sum_{r' = r}^{R^{\pi}-1} \left(S^{\pi}_{r'} - S^{\pi}_{r'+1}\right) w^{\pi*}_{r'}\right] + \lambda^{\pi*} \ \ \ge \ \ 0 \qquad \forall \ \pi \in \Pi, \ \forall \ r \in \mathcal{R}^{\pi}
\end{align*}
\end{manuallemma}

\proof{Proof.}
It is easy to check that $x(q^*) = x(\hat{q})$.

Consider any $\pi \in \Pi$.
It follows from Lemma~\ref{lem:optimal thresholds}\ref{item:thresholds optimal} that for each $\ell \in \{1,\ldots,\hat{L}\}$ it holds that $T^{\pi}_{\hat{r}_{\ell}}(x(\hat{q})) + w S^{\pi}_{\hat{r}_{\ell}} \le T^{\pi}_{r}(x(\hat{q})) + w S^{\pi}_{r}$ for all $r \in \hat{\mathcal{R}}^{\pi} \subset \mathcal{R}^{\pi}$ and all $w \in [\hat{w}^{\pi}_{\hat{r}_{\ell-1}},\hat{w}^{\pi}_{\hat{r}_{\ell}}]$.
Next we show that this inequality holds for the complete set of paths $r \in \mathcal{R}^{\pi}$.
Suppose that for some $\ell \in \{1,\ldots,\hat{L}\}$ and some $w^{+} \in [\hat{w}^{\pi}_{\hat{r}_{\ell-1}},\hat{w}^{\pi}_{\hat{r}_{\ell}}]$ there is an $r \in \mathcal{R}^{\pi} \setminus \hat{\mathcal{R}}^{\pi}$ such that $T^{\pi}_{r}(x(q^*)) + w^{+} S^{\pi}_{r} < T^{\pi}_{\hat{r}_{\ell}}(x(q^*)) + w^{+} S^{\pi}_{\hat{r}_{\ell}}$.
Let $r^{+} \in \mathcal{R}^{\pi+}$ denote the shortest path produced by the parametric shortest path algorithm when the arc costs are $c_{a}(w^{+},x_{a}(q^*)) = t_{a}(x_{a}(q^*)) + w^{+} s_{a}$.
Then $r^{+} \in \mathcal{R}^{\pi+} \setminus \hat{\mathcal{R}}^{\pi}$, which contradicts $\mathcal{R}^{\pi+} \subset \hat{\mathcal{R}}^{\pi}$.
Thus
\begin{equation}
T^{\pi}_{\hat{r}_{\ell}}(x(q^*)) + w S^{\pi}_{\hat{r}_{\ell}} \ \ \le \ \ T^{\pi}_{r}(x(q^*)) + w S^{\pi}_{r} \qquad \forall \ r \in \mathcal{R}^{\pi}, \ \forall \ w \in [\hat{w}^{\pi}_{\hat{r}_{\ell-1}},\hat{w}^{\pi}_{\hat{r}_{\ell}}]
\label{eqn:min cost ell}
\end{equation}

Let $\hat{r}_{0} = r^*_{0} = 0$, and note that $\hat{w}^{\pi}_{0} = w^{\pi*}_{0} = \underline{w}^{\pi}$.
Note that $\mathcal{R}^{\pi+} \subset \hat{\mathcal{R}}^{\pi}$ does not imply that all efficient paths are in $\hat{\mathcal{R}}^{\pi}$, because there may be efficient paths $r^*_{\ell} \in \mathcal{R}^{\pi*}_{\min} \setminus \hat{\mathcal{R}}^{\pi}$ (and thus $q^{\pi*}_{r^*_{\ell}} = 0$).
Next we show that $\hat{\mathcal{R}}^{\pi}_{\min} \subset \mathcal{R}^{\pi*}_{\min}$, and that for each $\hat{r}_{\ell} \in \hat{\mathcal{R}}^{\pi}_{\min}$, Algorithm~\ref{algo:threshold} produces the same interval $[\hat{w}^{\pi}_{\hat{r}_{\ell-1}},\hat{w}^{\pi}_{\hat{r}_{\ell}}]$ whether the input is $\hat{\mathcal{R}}^{\pi}$ and $\hat{q}$ or $\mathcal{R}^{\pi}$ and $q^*$.
Thereafter the proof can proceed similarly to the proof of Lemma~\ref{lem:thresholds optimality conditions}.
Specifically, next we show by induction on $\ell \in \{1,\ldots,\hat{L}\}$ that
\begin{enumerate}[label=(\arabic*), ref=(\arabic*)]
\item
\label{item:efficient paths in complete list}
there are indices $\ell^*(1) < \cdots < \ell^*(\hat{L})$ such that for each $\ell \in \{1,\ldots,\hat{L}\}$ it holds that $\hat{r}_{\ell} = r^*_{\ell^*(\ell)}$, and
\item
\label{item:thresholds in complete list}
for each $\ell \in \{1,\ldots,\hat{L}\}$ it holds that $\hat{w}^{\pi}_{\hat{r}_{\ell-1}} = w^{\pi*}_{r^*_{\ell^*(\ell)-1}}$.
\end{enumerate}

It follows from $\hat{\mathcal{R}}^{\pi} \subset \mathcal{R}^{\pi}$, \eqref{eqn:efficient path 1}, and~\eqref{eqn:min cost ell} that
\begin{equation}
T^{\pi}_{r^*_{1}}(x(q^*)) + \underline{w}^{\pi} S^{\pi}_{r^*_{1}} \ \ = \ \ T^{\pi}_{\hat{r}_{1}}(x(q^*)) + \underline{w}^{\pi} S^{\pi}_{\hat{r}_{1}}
\label{eqn:equal cost ell = 1}
\end{equation}
and thus $r^*_{1} \le \hat{r}_{1}$.
If $r^*_{1} = \hat{r}_{1}$, then claim~\ref{item:efficient paths in complete list} of the hypothesis holds for $\ell = 1$ with $\ell^*(1) = 1$.
Also, $\hat{w}^{\pi}_{\hat{r}_{0}} = \underline{w}^{\pi} = w^{\pi*}_{r^*_{0}}$, and hence claim~\ref{item:thresholds in complete list} of the hypothesis holds for $\ell = 1$.
Otherwise, $r^*_{1} < \hat{r}_{1}$, and thus $S^{\pi}_{r^*_{1}} \ge S^{\pi}_{\hat{r}_{1}}$.
Let $\check{\ell} \defi \max\{\ell \in \{1,\ldots,L^*\} \, : \, r^*_{\ell} < \hat{r}_{1}\}$.
Then it follows from Lemma~\ref{lem:optimal thresholds}\ref{item:positive flow paths are efficient} that $\hat{q}^{\pi}_{r} = 0$ for all $r \in \hat{\mathcal{R}}^{\pi}$ such that $r < \hat{r}_{1}$, and thus $q^{\pi*}_{r} = 0$ for all $r \in \mathcal{R}^{\pi}$ such that $r < \hat{r}_{1}$, and $L^{\pi}_{r^*_{\ell}}(q^{\pi*}) = 0$ for all $\ell \in \{1,\ldots,\check{\ell}\}$.
Also, $S^{\pi}_{r^*_{1}} \ge S^{\pi}_{r^*_{\check{\ell}}} \ge S^{\pi}_{\hat{r}_{1}}$.

\noindent
\textbf{Case~A.1}: $S^{\pi}_{r^*_{1}} = S^{\pi}_{r^*_{\check{\ell}}}$. \\
Then it follows from $r^*_{\ell} \in \mathcal{R}^{\pi*}_{\min}$ and $S^{\pi}_{r^*_{1}} = S^{\pi}_{r^*_{\ell}}$ for all $\ell \in \{1,\ldots,\check{\ell}\}$ that $T^{\pi}_{r^*_{1}}(x(q^*)) = T^{\pi}_{r^*_{\ell}}(x(q^*))$ for all $\ell \in \{1,\ldots,\check{\ell}\}$.
Then, in Algorithm~\ref{algo:threshold}, $w^{\pi*}_{r^*_{\ell}}$ is determined by~\eqref{eqn:threshold ==} for all $\ell \in \{1,\ldots,\check{\ell}-1\}$.
Thus $w^{\pi*}_{r^*_{\ell}} = w^{\pi*}_{0} = \underline{w}^{\pi}$ for all $\ell \in \{1,\ldots,\check{\ell}-1\}$.

\noindent
\textbf{Case~A.1.1}: $S^{\pi}_{r^*_{1}} = S^{\pi}_{r^*_{\check{\ell}}} = S^{\pi}_{\hat{r}_{1}}$. \\
Then it follows from~\eqref{eqn:equal cost ell = 1} that $T^{\pi}_{r^*_{1}}(x(q^*)) = T^{\pi}_{r^*_{\check{\ell}}}(x(q^*)) = T^{\pi}_{\hat{r}_{1}}(x(q^*))$.
Then, in Algorithm~\ref{algo:threshold}, $w^{\pi*}_{r^*_{\check{\ell}}}$ is determined by~\eqref{eqn:threshold ==}, $w^{\pi*}_{r^*_{\check{\ell}}} = w^{\pi*}_{r^*_{\check{\ell}-1}} = \underline{w}^{\pi}$, and $r^*_{\check{\ell}+1} = \hat{r}_{1}$, hence claim~\ref{item:efficient paths in complete list} of the hypothesis holds for $\ell = 1$ with $\ell^*(1) = \check{\ell}+1$.
Also, $\hat{w}^{\pi}_{\hat{r}_{0}} = \underline{w}^{\pi} = w^{\pi*}_{r^*_{\check{\ell}}}$, and hence claim~\ref{item:thresholds in complete list} of the hypothesis holds for $\ell = 1$.

\noindent
\textbf{Case~A.1.2}: $S^{\pi}_{r^*_{1}} = S^{\pi}_{r^*_{\check{\ell}}} > S^{\pi}_{\hat{r}_{1}}$. \\
Then, in Algorithm~\ref{algo:threshold}, $w^{\pi*}_{r^*_{\check{\ell}}}$ is determined by~\eqref{eqn:threshold >} or~\eqref{eqn:threshold =>}.
Then
\begin{align*}
w^{\min}_{\check{\ell}} \ \ & = \ \ \min\left\{\frac{T^{\pi}_{r}(x(q^*)) - T^{\pi}_{r^*_{\check{\ell}}}(x(q^*))}{S^{\pi}_{r^*_{\check{\ell}}} - S^{\pi}_{r}} \; : \; S^{\pi}_{r^*_{\check{\ell}}} > S^{\pi}_{r}\right\} \\
& \le \ \ \frac{T^{\pi}_{\hat{r}_{1}}(x(q^*)) - T^{\pi}_{r^*_{\check{\ell}}}(x(q^*))}{S^{\pi}_{r^*_{\check{\ell}}} - S^{\pi}_{\hat{r}_{1}}}
\ \ = \ \ \frac{T^{\pi}_{\hat{r}_{1}}(x(q^*)) - T^{\pi}_{r^*_{1}}(x(q^*))}{S^{\pi}_{r^*_{1}} - S^{\pi}_{\hat{r}_{1}}}
\ \ \overset{\eqref{eqn:equal cost ell = 1}}{=} \ \ \underline{w}^{\pi}
\end{align*}
Suppose that $w^{\min}_{\check{\ell}} < \underline{w}^{\pi}$.
Then
\begin{align*}
& w^{\min}_{\check{\ell}} \ \ = \ \ \frac{T^{\pi}_{r^*_{\check{\ell}+1}}(x(q^*)) - T^{\pi}_{r^*_{\check{\ell}}}(x(q^*))}{S^{\pi}_{r^*_{\check{\ell}}} - S^{\pi}_{r^*_{\check{\ell}+1}}}
\ \ < \ \ \underline{w}^{\pi} \\
\Rightarrow \qquad & T^{\pi}_{r^*_{\check{\ell}+1}}(x(q^*)) + \underline{w}^{\pi} S^{\pi}_{r^*_{\check{\ell}+1}} \ \ < \ \ T^{\pi}_{r^*_{\check{\ell}}}(x(q^*)) + \underline{w}^{\pi} S^{\pi}_{r^*_{\check{\ell}}}
\ \ = \ \ T^{\pi}_{r^*_{1}}(x(q^*)) + \underline{w}^{\pi} S^{\pi}_{r^*_{1}}
\end{align*}
which contradicts $T^{\pi}_{r^*_{1}}(x(q^*)) + \underline{w}^{\pi} S^{\pi}_{r^*_{1}} \le T^{\pi}_{r}(x(q^*)) + \underline{w}^{\pi} S^{\pi}_{r}$ for all $r \in \mathcal{R}^{\pi}$ from~\eqref{eqn:efficient path 1}.
Thus $w^{\min}_{\check{\ell}} = \underline{w}^{\pi}$, and hence
\[
w^{\pi*}_{r^*_{\check{\ell}}} \ \ = \ \ w^{\min}_{\check{\ell}} \ \ = \ \ \min\left\{\frac{T^{\pi}_{r}(x(q^*)) - T^{\pi}_{r^*_{\check{\ell}}}(x(q^*))}{S^{\pi}_{r^*_{\check{\ell}}} - S^{\pi}_{r}} \; : \; S^{\pi}_{r^*_{\check{\ell}}} > S^{\pi}_{r}\right\}
\ \ = \ \ \frac{T^{\pi}_{\hat{r}_{1}}(x(q^*)) - T^{\pi}_{r^*_{\check{\ell}}}(x(q^*))}{S^{\pi}_{r^*_{\check{\ell}}} - S^{\pi}_{\hat{r}_{1}}}
\]
and $r^*_{\check{\ell}+1} = \hat{r}_{1}$.
Hence claim~\ref{item:efficient paths in complete list} of the hypothesis holds for $\ell = 1$ with $\ell^*(1) = \check{\ell}+1$.
Also, $\hat{w}^{\pi}_{\hat{r}_{0}} = \underline{w}^{\pi} = w^{\min}_{\check{\ell}} = w^{\pi*}_{r^*_{\check{\ell}}}$, and hence claim~\ref{item:thresholds in complete list} of the hypothesis holds for $\ell = 1$.

\noindent
\textbf{Case~A.2}: $S^{\pi}_{r^*_{1}} > S^{\pi}_{r^*_{\check{\ell}}}$. \\
Then $r^*_{1} < r^*_{\check{\ell}} < \hat{r}_{1}$.
Let $\tilde{\ell} \defi \max\left\{\ell \in \{1,\ldots,\check{\ell}-1\} \, : \, S^{\pi}_{r^*_{1}} = S^{\pi}_{r^*_{\ell}}\right\}$.
Thus, $S^{\pi}_{r^*_{1}} = S^{\pi}_{r^*_{\ell}}$ and $T^{\pi}_{r^*_{1}}(x(q^*)) = T^{\pi}_{r^*_{\ell}}(x(q^*))$ for all $\ell \in \{1,\ldots,\tilde{\ell}\}$, and in Algorithm~\ref{algo:threshold}, $w^{\pi*}_{r^*_{\ell}}$ is determined by~\eqref{eqn:threshold ==} for all $\ell = 1,\ldots,\tilde{\ell}-1$, and $\underline{w}^{\pi} = w^{\pi*}_{r^*_{0}} = \cdots = w^{\pi*}_{r^*_{\tilde{\ell}-1}}$.
Also, in Algorithm~\ref{algo:threshold}, $w^{\pi*}_{r^*_{\tilde{\ell}}}$ is determined by~\eqref{eqn:threshold >} or~\eqref{eqn:threshold =>}.
Then
\begin{align*}
w^{\min}_{\tilde{\ell}} \ \ & = \ \ \min\left\{\frac{T^{\pi}_{r}(x(q^*)) - T^{\pi}_{r^*_{\tilde{\ell}}}(x(q^*))}{S^{\pi}_{r^*_{\tilde{\ell}}} - S^{\pi}_{r}} \; : \; S^{\pi}_{r^*_{\tilde{\ell}}} > S^{\pi}_{r}\right\} \\
& \le \ \ \frac{T^{\pi}_{\hat{r}_{1}}(x(q^*)) - T^{\pi}_{r^*_{\tilde{\ell}}}(x(q^*))}{S^{\pi}_{r^*_{\tilde{\ell}}} - S^{\pi}_{\hat{r}_{1}}}
\ \ = \ \ \frac{T^{\pi}_{\hat{r}_{1}}(x(q^*)) - T^{\pi}_{r^*_{1}}(x(q^*))}{S^{\pi}_{r^*_{1}} - S^{\pi}_{\hat{r}_{1}}}
\ \ \overset{\eqref{eqn:equal cost ell = 1}}{=} \ \ \underline{w}^{\pi}
\end{align*}
Suppose that $w^{\min}_{\tilde{\ell}} < \underline{w}^{\pi}$.
Then
\begin{align*}
& w^{\min}_{\tilde{\ell}} \ \ = \ \ \frac{T^{\pi}_{r^*_{\tilde{\ell}+1}}(x(q^*)) - T^{\pi}_{r^*_{\tilde{\ell}}}(x(q^*))}{S^{\pi}_{r^*_{\tilde{\ell}}} - S^{\pi}_{r^*_{\tilde{\ell}+1}}}
\ \ < \ \ \underline{w}^{\pi} \\
\Rightarrow \qquad & T^{\pi}_{r^*_{\tilde{\ell}+1}}(x(q^*)) + \underline{w}^{\pi} S^{\pi}_{r^*_{\tilde{\ell}+1}} \ \ < \ \ T^{\pi}_{r^*_{\tilde{\ell}}}(x(q^*)) + \underline{w}^{\pi} S^{\pi}_{r^*_{\tilde{\ell}}}
\ \ = \ \ T^{\pi}_{r^*_{1}}(x(q^*)) + \underline{w}^{\pi} S^{\pi}_{r^*_{1}}
\end{align*}
which contradicts $T^{\pi}_{r^*_{1}}(x(q^*)) + \underline{w}^{\pi} S^{\pi}_{r^*_{1}} \le T^{\pi}_{r}(x(q^*)) + \underline{w}^{\pi} S^{\pi}_{r}$ for all $r \in \mathcal{R}^{\pi}$ from~\eqref{eqn:efficient path 1}.
Thus $w^{\pi*}_{r^*_{\tilde{\ell}}} = w^{\min}_{\tilde{\ell}} = \underline{w}^{\pi}$ and $T^{\pi}_{r^*_{\tilde{\ell}+1}}(x(q^*)) + \underline{w}^{\pi} S^{\pi}_{r^*_{\tilde{\ell}+1}} = T^{\pi}_{r^*_{\tilde{\ell}}}(x(q^*)) + \underline{w}^{\pi} S^{\pi}_{r^*_{\tilde{\ell}}} = T^{\pi}_{r^*_{1}}(x(q^*)) + \underline{w}^{\pi} S^{\pi}_{r^*_{1}}$.
If $\tilde{\ell}+1 < \check{\ell}$, then let $\tilde{\tilde{\ell}} \defi \max\left\{\ell \in \{\tilde{\ell}+1,\ldots,\check{\ell}\} \, : \, S^{\pi}_{r^*_{\tilde{\ell}+1}} = S^{\pi}_{r^*_{\ell}}\right\}$.
Then the argument is repeated by induction, to obtain that $w^{\pi*}_{r^*_{\check{\ell}-1}} = \underline{w}^{\pi}$ and
\begin{equation}
T^{\pi}_{r^*_{\check{\ell}}}(x(q^*)) + \underline{w}^{\pi} S^{\pi}_{r^*_{\check{\ell}}} \ \ = \ \ T^{\pi}_{r^*_{1}}(x(q^*)) + \underline{w}^{\pi} S^{\pi}_{r^*_{1}}
\label{eqn:equal cost ell = check ell}
\end{equation}

\noindent
\textbf{Case~A.2.1}: $S^{\pi}_{r^*_{1}} > S^{\pi}_{r^*_{\check{\ell}}} = S^{\pi}_{\hat{r}_{1}}$. \\
Then
\begin{align*}
& T^{\pi}_{\hat{r}_{1}}(x(q^*)) + \underline{w}^{\pi} S^{\pi}_{\hat{r}_{1}} \ \ \overset{\eqref{eqn:equal cost ell = 1}}{=} \ \ T^{\pi}_{r^*_{1}}(x(q^*)) + \underline{w}^{\pi} S^{\pi}_{r^*_{1}} \ \ \overset{\eqref{eqn:equal cost ell = check ell}}{=} \ \ T^{\pi}_{r^*_{\check{\ell}}}(x(q^*)) + \underline{w}^{\pi} S^{\pi}_{r^*_{\check{\ell}}} \\
& \Rightarrow \qquad T^{\pi}_{\hat{r}_{1}}(x(q^*)) \ \ = \ \ T^{\pi}_{r^*_{\check{\ell}}}(x(q^*))
\end{align*}
Thus, in Algorithm~\ref{algo:threshold}, $w^{\pi*}_{r^*_{\check{\ell}}}$ is determined by~\eqref{eqn:threshold ==}, and $r^*_{\check{\ell}+1} = \hat{r}_{1}$.
Hence claim~\ref{item:efficient paths in complete list} of the hypothesis holds for $\ell = 1$ with $\ell^*(1) = \check{\ell}+1$.
Also, $w^{\pi*}_{r^*_{\check{\ell}}} = \max\left\{w^{\pi*}_{r^*_{\check{\ell}}-1}, \; \min\left\{w \in \left(G^{\pi}\right)^{-1}\left(L^{\pi}_{r^*_{\check{\ell}}}(q^{\pi*})\right)\right\}\right\} = \max\left\{w^{\pi*}_{r^*_{\check{\ell}-1}}, \; \min\left\{w \in \left(G^{\pi}\right)^{-1}\left(0\right)\right\}\right\} = \underline{w}^{\pi}$.
Thus $\hat{w}^{\pi}_{\hat{r}_{0}} = \underline{w}^{\pi} = w^{\pi*}_{r^*_{\check{\ell}}}$, and hence claim~\ref{item:thresholds in complete list} of the hypothesis holds for $\ell = 1$.

\noindent
\textbf{Case~A.2.2}: $S^{\pi}_{r^*_{1}} > S^{\pi}_{r^*_{\check{\ell}}} > S^{\pi}_{\hat{r}_{1}}$. \\
Then, in Algorithm~\ref{algo:threshold}, $w^{\pi*}_{r^*_{\check{\ell}}}$ is determined by~\eqref{eqn:threshold >} or~\eqref{eqn:threshold =>}.
Then
\[
w^{\min}_{\check{\ell}} \ \ = \ \ \min\left\{\frac{T^{\pi}_{r}(x(q^*)) - T^{\pi}_{r^*_{\check{\ell}}}(x(q^*))}{S^{\pi}_{r^*_{\check{\ell}}} - S^{\pi}_{r}} \; : \; S^{\pi}_{r^*_{\check{\ell}}} > S^{\pi}_{r}\right\} \ \ \le \ \ \frac{T^{\pi}_{\hat{r}_{1}}(x(q^*)) - T^{\pi}_{r^*_{\check{\ell}}}(x(q^*))}{S^{\pi}_{r^*_{\check{\ell}}} - S^{\pi}_{\hat{r}_{1}}}
\ \ \overset{\eqref{eqn:equal cost ell = 1},\eqref{eqn:equal cost ell = check ell}}{=} \ \ \underline{w}^{\pi}
\]
Suppose that $w^{\min}_{\check{\ell}} < \underline{w}^{(u,v)}$.
Then
\begin{align*}
& w^{\min}_{\check{\ell}} \ \ = \ \ \frac{T^{\pi}_{r^*_{\check{\ell}+1}}(x(q^*)) - T^{\pi}_{r^*_{\check{\ell}}}(x(q^*))}{S^{\pi}_{r^*_{\check{\ell}}} - S^{\pi}_{r^*_{\check{\ell}+1}}}
\ \ < \ \ \underline{w}^{\pi} \\
\Rightarrow \qquad & T^{\pi}_{r^*_{\check{\ell}+1}}(x(q^*)) + \underline{w}^{\pi} S^{\pi}_{r^*_{\check{\ell}+1}} \ \ < \ \ T^{\pi}_{r^*_{\check{\ell}}}(x(q^*)) + \underline{w}^{\pi} S^{\pi}_{r^*_{\check{\ell}}}
\ \ \overset{\eqref{eqn:equal cost ell = check ell}}{=} \ \ T^{\pi}_{r^*_{1}}(x(q^*)) + \underline{w}^{\pi} S^{\pi}_{r^*_{1}}
\end{align*}
which contradicts $T^{\pi}_{r^*_{1}}(x(q^*)) + \underline{w}^{\pi} S^{\pi}_{r^*_{1}} \le T^{\pi}_{r}(x(q^*)) + \underline{w}^{\pi} S^{\pi}_{r}$ for all $r \in \mathcal{R}^{\pi}$ from~\eqref{eqn:efficient path 1}.
Thus $w^{\min}_{\check{\ell}} = \underline{w}^{\pi}$, and hence
\[
w^{\pi*}_{r^*_{\check{\ell}}} \ \ = \ \ w^{\min}_{\check{\ell}} \ \ = \ \ \min\left\{\frac{T^{\pi}_{r}(x(q^*)) - T^{\pi}_{r^*_{\check{\ell}}}(x(q^*))}{S^{\pi}_{r^*_{\check{\ell}}} - S^{\pi}_{r}} \; : \; S^{\pi}_{r^*_{\check{\ell}}} > S^{\pi}_{r}\right\}
\ \ = \ \ \frac{T^{\pi}_{\hat{r}_{1}}(x(q^*)) - T^{\pi}_{r^*_{\check{\ell}}}(x(q^*))}{S^{\pi}_{r^*_{\check{\ell}}} - S^{\pi}_{\hat{r}_{1}}}
\]
and $r^*_{\check{\ell}+1} = \hat{r}_{1}$.
Hence claim~\ref{item:efficient paths in complete list} of the hypothesis holds for $\ell = 1$ with $\ell^*(1) = \check{\ell}+1$.
Also, $\hat{w}^{\pi}_{\hat{r}_{0}} = \underline{w}^{\pi} = w^{\min}_{\check{\ell}} = w^{\pi*}_{r^*_{\check{\ell}}}$, and hence claim~\ref{item:thresholds in complete list} of the hypothesis holds for $\ell = 1$.

Suppose that the hypothesis holds for $\ell \in \{1,\ldots,\hat{L}-1\}$.
Then
\begin{equation}
T^{\pi}_{\hat{r}_{\ell}}(x(q^*)) + \hat{w}^{\pi}_{\hat{r}_{\ell-1}} S^{\pi}_{\hat{r}_{\ell}}
\ \ = \ \ T^{\pi}_{r^*_{\ell^*(\ell)}}(x(q^*)) + w^{\pi*}_{r^*_{\ell^*(\ell)-1}} S^{\pi}_{r^*_{\ell^*(\ell)}}
\ \ \overset{\eqref{eqn:min cost ell}}{\le} \ \ T^{\pi}_{r}(x(q^*)) + \hat{w}^{\pi}_{\hat{r}_{\ell-1}} S^{\pi}_{r} \qquad \forall \ r \in \mathcal{R}^{\pi}
\label{eqn:min cost ell*}
\end{equation}
Also, note that $L^{\pi}_{\hat{r}_{\ell}}(\hat{q}^{\pi}) = L^{\pi}_{r^*_{\ell^*(\ell)}}(q^{\pi*})$.
Next we show that $r^*_{\ell^*(\ell)+1} \le \hat{r}_{\ell+1}$.
In Algorithm~\ref{algo:threshold}, $\hat{r}_{\ell+1}$ is determined by~\eqref{eqn:efficient path >} or~\eqref{eqn:efficient path ==} or~\eqref{eqn:efficient path =>}.

\noindent
\textbf{Case~B.1}: $\hat{r}_{\ell+1}$ is determined by~\eqref{eqn:efficient path >} or~\eqref{eqn:efficient path =>}. \\
Suppose that $r^*_{\ell^*(\ell)+1} > \hat{r}_{\ell+1}$.
Then $S^{\pi}_{r^*_{\ell^*(\ell)+1}} \le S^{\pi}_{\hat{r}_{\ell+1}} < S^{\pi}_{\hat{r}_{\ell}} = S^{\pi}_{r^*_{\ell^*(\ell)}}$.
Then
\begin{align*}
& w^{\min}_{\ell^*(\ell)} \ \ = \ \ \frac{T^{\pi}_{r^*_{\ell^*(\ell)+1}}(x(q^*)) - T^{\pi}_{r^*_{\ell^*(\ell)}}(x(q^*))}{S^{\pi}_{r^*_{\ell^*(\ell)}} - S^{\pi}_{r^*_{\ell^*(\ell)+1}}}
\ \ < \ \ \frac{T^{\pi}_{\hat{r}_{\ell+1}}(x(q^*)) - T^{\pi}_{r^*_{\ell^*(\ell)}}(x(q^*))}{S^{\pi}_{r^*_{\ell^*(\ell)}} - S^{\pi}_{\hat{r}_{\ell+1}}} \\
& \hspace{80mm} = \ \ \frac{T^{\pi}_{\hat{r}_{\ell+1}}(x(\hat{q})) - T^{\pi}_{\hat{r}_{\ell}}(x(\hat{q}))}{S^{\pi}_{\hat{r}_{\ell}} - S^{\pi}_{\hat{r}_{\ell+1}}}
\ \ = \ \ \hat{w}^{\pi}_{\hat{r}_{\ell}} \\
\Rightarrow \qquad & T^{\pi}_{r^*_{\ell^*(\ell)+1}}(x(q^*)) + \hat{w}^{\pi}_{\hat{r}_{\ell}} S^{\pi}_{r^*_{\ell^*(\ell)+1}} \ \ < \ \ T^{\pi}_{r^*_{\ell^*(\ell)}}(x(q^*)) + \hat{w}^{\pi}_{\hat{r}_{\ell}} S^{\pi}_{r^*_{\ell^*(\ell)}}
\ \ = \ \ T^{\pi}_{\hat{r}_{\ell}}(x(q^*)) + \hat{w}^{\pi}_{\hat{r}_{\ell}} S^{\pi}_{\hat{r}_{\ell}}
\end{align*}
which contradicts~\eqref{eqn:min cost ell}.

\noindent
\textbf{Case~B.2}: $\hat{r}_{\ell+1}$ is determined by~\eqref{eqn:efficient path ==}. \\
Then $S^{\pi}_{\hat{r}_{\ell+1}} = S^{\pi}_{\hat{r}_{\ell}} = S^{\pi}_{r^*_{\ell^*(\ell)}}$ and $T^{\pi}_{\hat{r}_{\ell+1}}(x(\hat{q})) = T^{\pi}_{\hat{r}_{\ell}}(x(\hat{q})) = T^{\pi}_{r^*_{\ell^*(\ell)}}(x(q^*))$.
Then $\left\{r \in \{r^*_{\ell^*(\ell)}+1,\ldots,R^{\pi}\} \, : \, S^{\pi}_{r} = S^{\pi}_{r^*_{\ell^*(\ell)}}\right\} \neq \varnothing$.
Suppose that there is an $\check{r} \in \{r^*_{\ell^*(\ell)}+1,\ldots,R^{\pi}\}$ such that $S^{\pi}_{\check{r}} = S^{\pi}_{r^*_{\ell^*(\ell)}}$ and $T^{\pi}_{\check{r}}(x(q^*)) < T^{\pi}_{\hat{r}_{\ell+1}}(x(q^*))$.
Then $T^{\pi}_{\check{r}}(x(q^*)) + \hat{w}^{\pi}_{\hat{r}_{\ell}} S^{\pi}_{\check{r}} < T^{\pi}_{\hat{r}_{\ell+1}}(x(q^*)) + \hat{w}^{\pi}_{\hat{r}_{\ell}} S^{\pi}_{\hat{r}_{\ell+1}}$ which contradicts~\eqref{eqn:min cost ell}.
Thus it follows from $\hat{\mathcal{R}}^{\pi} \subset \mathcal{R}^{\pi}$ that $\min\left\{T^{\pi}_{r}(x(q^*)) \, : \, S^{\pi}_{r} = S^{\pi}_{r^*_{\ell^*(\ell)}}, \; r \in \{r^*_{\ell^*(\ell)}+1,\ldots,R^{\pi}\}\right\} = T^{\pi}_{\hat{r}_{\ell+1}}(x(q^*))$, and
\newline $r^*_{\ell^*(\ell)+1} = \min\left\{\arg\min\left\{T^{\pi}_{r}(x(q^*)) \; : \; S^{\pi}_{r} = S^{\pi}_{r^*_{\ell^*(\ell)}}, \; r \in \{r^*_{\ell^*(\ell)}+1,\ldots,R^{\pi}\}\right\}\right\} \le \hat{r}_{\ell+1}$.
Thus $r^*_{\ell^*(\ell)+1} \le \hat{r}_{\ell+1}$.

If $r^*_{\ell^*(\ell)+1} = \hat{r}_{\ell+1}$, then claim~\ref{item:efficient paths in complete list} of the hypothesis holds for $\ell + 1$ with $\ell^*(\ell+1) = \ell^*(\ell) + 1$.
Also, $\hat{w}^{\pi}_{\hat{r}_{\ell}} = w^{\pi*}_{r^*_{\ell^*(\ell)}} = w^{\pi*}_{r^*_{\ell^*(\ell+1)-1}}$, and hence claim~\ref{item:thresholds in complete list} of the hypothesis holds for $\ell + 1$.
Otherwise, $r^*_{\ell^*(\ell)+1} < \hat{r}_{\ell+1}$, and thus $S^{\pi}_{r^*_{\ell^*(\ell)+1}} \ge S^{\pi}_{\hat{r}_{\ell+1}}$.
Let $\check{\ell} \defi \max\{\ell' \in \{\ell^*(\ell)+1,\ldots,L^*\} \, : \, r^*_{\ell'} < \hat{r}_{\ell+1}\}$.
It follows from Lemma~\ref{lem:optimal thresholds}\ref{item:positive flow paths are efficient} that $\hat{q}^{\pi}_{r} = 0$ for all $r \in \hat{\mathcal{R}}^{\pi}$ such that $r \in \{\hat{r}_{\ell}+1,\ldots,\hat{r}_{\ell+1}-1\}$, and thus $L^{\pi}_{r^*_{\ell'}}(q^{\pi*}) = L^{\pi}_{r^*_{\ell^*(\ell)}}(q^{\pi*}) = L^{\pi}_{\hat{r}_{\ell}}(\hat{q}^{\pi})$ for all $\ell' \in \{\ell^*(\ell),\ldots,\check{\ell}\}$.
Also, $S^{\pi}_{\hat{r}_{\ell}} = S^{\pi}_{r^*_{\ell^*(\ell)}} \ge S^{\pi}_{r^*_{\check{\ell}}} \ge S^{\pi}_{\hat{r}_{\ell+1}}$.

\noindent
\textbf{Case~C.1}: $S^{\pi}_{r^*_{\ell^*(\ell)}} = S^{\pi}_{r^*_{\check{\ell}}}$. \\
Then it follows from $r^*_{\ell^*(\ell)} \in \mathcal{R}^{\pi*}_{\min}$ and $S^{\pi}_{r^*_{\ell^*(\ell)}} = S^{\pi}_{r^*_{\ell'}}$ for all $\ell' \in \{\ell^*(\ell),\ldots,\check{\ell}\}$ that $T^{\pi}_{r^*_{\ell^*(\ell)}}(x(q^*)) = T^{\pi}_{r^*_{\ell'}}(x(q^*))$ for all $\ell' \in \{\ell^*(\ell),\ldots,\check{\ell}\}$.
Then, in Algorithm~\ref{algo:threshold}, $w^{\pi*}_{r^*_{\ell'}}$ is determined by~\eqref{eqn:threshold ==} for all $\ell' \in \{\ell^*(\ell),\ldots,\check{\ell}-1\}$.
Also, it follows from $L^{\pi}_{r^*_{\ell'}}(q^{\pi*}) = L^{\pi}_{r^*_{\ell^*(\ell)}}(q^{\pi*}) = L^{\pi}_{\hat{r}_{\ell}}(\hat{q}^{\pi})$ for all $\ell' \in \{\ell^*(\ell),\ldots,\check{\ell}\}$ and induction on $\ell'$ that
\begin{align*}
w^{\pi*}_{r^*_{\ell'}} \ \ & = \ \ \max\left\{w^{\pi*}_{r^*_{\ell'}-1}, \; \min\left\{w \in \left(G^{\pi}\right)^{-1}\left(L^{\pi}_{r^*_{\ell'}}(q^{\pi*})\right)\right\}\right\}
\ \ = \ \ \max\left\{w^{\pi*}_{r^*_{\ell^*(\ell)-1}}, \; \min\left\{w \in \left(G^{\pi}\right)^{-1}\left(L^{\pi}_{r^*_{\ell^*(\ell)}}(q^{\pi*})\right)\right\}\right\} \\
& = \ \ \max\left\{\hat{w}^{\pi}_{\hat{r}_{\ell-1}}, \; \min\left\{w \in \left(G^{\pi}\right)^{-1}\left(L^{\pi}_{\hat{r}_{\ell}}(\hat{q}^{\pi})\right)\right\}\right\}
\qquad \forall \ \ell' \in \{\ell^*(\ell),\ldots,\check{\ell}-1\}
\end{align*}

\noindent
\textbf{Case~C.1.1}: $S^{\pi}_{r^*_{\ell^*(\ell)}} = S^{\pi}_{r^*_{\check{\ell}}} = S^{\pi}_{\hat{r}_{\ell+1}}$. \\
Then $T^{\pi}_{r^*_{\ell^*(\ell)}}(x(q^*)) = T^{\pi}_{r^*_{\check{\ell}}}(x(q^*)) \le T^{\pi}_{\hat{r}_{\ell+1}}(x(q^*))$.
Suppose $T^{\pi}_{r^*_{\check{\ell}}}(x(q^*)) < T^{\pi}_{\hat{r}_{\ell+1}}(x(q^*))$.
Then $T^{\pi}_{r^*_{\check{\ell}}}(x(q^*)) + \hat{w}^{\pi}_{\hat{r}_{\ell}} S^{\pi}_{r^*_{\check{\ell}}} < T^{\pi}_{\hat{r}_{\ell+1}}(x(q^*)) + \hat{w}^{\pi}_{\hat{r}_{\ell}} S^{\pi}_{\hat{r}_{\ell+1}}$, which contradicts~\eqref{eqn:min cost ell}.
Thus $T^{\pi}_{r^*_{\ell^*(\ell)}}(x(q^*)) = T^{\pi}_{r^*_{\check{\ell}}}(x(q^*)) = T^{\pi}_{\hat{r}_{\ell+1}}(x(q^*))$.
Then, in Algorithm~\ref{algo:threshold}, $w^{\pi*}_{r^*_{\check{\ell}}}$ is determined by~\eqref{eqn:threshold ==}, and
\begin{align*}
w^{\pi*}_{r^*_{\check{\ell}}} \ \ & = \ \ \max\left\{w^{\pi*}_{r^*_{\check{\ell}}-1}, \; \min\left\{w \in \left(G^{\pi}\right)^{-1}\left(L^{\pi}_{r^*_{\check{\ell}}}(q^{\pi*})\right)\right\}\right\}
\ \ = \ \ \max\left\{w^{\pi*}_{r^*_{\ell^*(\ell)-1}}, \; \min\left\{w \in \left(G^{\pi}\right)^{-1}\left(L^{\pi}_{r^*_{\ell^*(\ell)}}(q^{\pi*})\right)\right\}\right\} \\
& = \ \ \max\left\{\hat{w}^{\pi}_{\hat{r}_{\ell-1}}, \; \min\left\{w \in \left(G^{\pi}\right)^{-1}\left(L^{\pi}_{\hat{r}_{\ell}}(\hat{q}^{\pi})\right)\right\}\right\}
\ \ = \ \ \hat{w}^{\pi}_{\hat{r}_{\ell}}
\end{align*}
Thus $r^*_{\check{\ell}+1} = \hat{r}_{\ell+1}$, hence claim~\ref{item:efficient paths in complete list} of the hypothesis holds for $\ell + 1$ with $\ell^*(\ell+1) = \check{\ell}+1$.
Also, $\hat{w}^{\pi}_{\hat{r}_{\ell}} = w^{\pi*}_{r^*_{\check{\ell}}} = w^{\pi*}_{r^*_{\ell^*(\ell+1)-1}}$, and hence claim~\ref{item:thresholds in complete list} of the hypothesis holds for $\ell + 1$.

\noindent
\textbf{Case~C.1.2}: $S^{\pi}_{r^*_{\ell^*(\ell)}} = S^{\pi}_{r^*_{\check{\ell}}} > S^{\pi}_{\hat{r}_{\ell+1}}$. \\
Then, in Algorithm~\ref{algo:threshold}, $w^{\pi*}_{r^*_{\check{\ell}}}$ is determined by~\eqref{eqn:threshold >} or~\eqref{eqn:threshold =>}.
Then
\begin{align*}
w^{\min}_{\check{\ell}} \ \ & = \ \ \min\left\{\frac{T^{\pi}_{r}(x(q^*)) - T^{\pi}_{r^*_{\check{\ell}}}(x(q^*))}{S^{\pi}_{r^*_{\check{\ell}}} - S^{\pi}_{r}} \; : \; S^{\pi}_{r^*_{\check{\ell}}} > S^{\pi}_{r}\right\} \\
& \le \ \ \frac{T^{\pi}_{\hat{r}_{\ell+1}}(x(q^*)) - T^{\pi}_{r^*_{\check{\ell}}}(x(q^*))}{S^{\pi}_{r^*_{\check{\ell}}} - S^{\pi}_{\hat{r}_{\ell+1}}}
\ \ = \ \ \frac{T^{\pi}_{\hat{r}_{\ell+1}}(x(q^*)) - T^{\pi}_{r^*_{\ell^*(\ell)}}(x(q^*))}{S^{\pi}_{r^*_{\ell^*(\ell)}} - S^{\pi}_{\hat{r}_{\ell+1}}}
\ \ = \ \ \frac{T^{\pi}_{\hat{r}_{\ell+1}}(x(\hat{q})) - T^{\pi}_{\hat{r}_{\ell}}(x(\hat{q}))}{S^{\pi}_{\hat{r}_{\ell}} - S^{\pi}_{\hat{r}_{\ell+1}}}
\ \ = \ \ \hat{w}^{\pi}_{\hat{r}_{\ell}}
\end{align*}
Suppose that $w^{\min}_{\check{\ell}} < \hat{w}^{\pi}_{\hat{r}_{\ell}}$.
Then
\begin{align*}
w^{\min}_{\check{\ell}} \ \ & = \ \ \frac{T^{\pi}_{r^*_{\check{\ell}+1}}(x(q^*)) - T^{\pi}_{r^*_{\check{\ell}}}(x(q^*))}{S^{\pi}_{r^*_{\check{\ell}}} - S^{\pi}_{r^*_{\check{\ell}+1}}}
\ \ = \ \ \frac{T^{\pi}_{r^*_{\check{\ell}+1}}(x(q^*)) - T^{\pi}_{r^*_{\ell^*(\ell)}}(x(q^*))}{S^{\pi}_{r^*_{\ell^*(\ell)}} - S^{\pi}_{r^*_{\check{\ell}+1}}}
\ \ = \ \ \frac{T^{\pi}_{r^*_{\check{\ell}+1}}(x(q^*)) - T^{\pi}_{\hat{r}_{\ell}}(x(q^*))}{S^{\pi}_{\hat{r}_{\ell}} - S^{\pi}_{r^*_{\check{\ell}+1}}} \\
& < \ \ \frac{T^{\pi}_{\hat{r}_{\ell+1}}(x(\hat{q})) - T^{\pi}_{\hat{r}_{\ell}}(x(\hat{q}))}{S^{\pi}_{\hat{r}_{\ell}} - S^{\pi}_{\hat{r}_{\ell+1}}}
\ \ = \ \ \hat{w}^{\pi}_{\hat{r}_{\ell}} \\
\Rightarrow \qquad & T^{\pi}_{r^*_{\check{\ell}+1}}(x(q^*)) + \hat{w}^{\pi}_{\hat{r}_{\ell}} S^{\pi}_{r^*_{\check{\ell}+1}} \ \ < \ \ T^{\pi}_{\hat{r}_{\ell}}(x(q^*)) + \hat{w}^{\pi}_{\hat{r}_{\ell}} S^{\pi}_{\hat{r}_{\ell}}
\ \ = \ \ T^{\pi}_{\hat{r}_{\ell+1}}(x(q^*)) + \hat{w}^{\pi}_{\hat{r}_{\ell}} S^{\pi}_{\hat{r}_{\ell+1}}
\end{align*}
which contradicts~\eqref{eqn:min cost ell}.
Thus $w^{\min}_{\check{\ell}} = \hat{w}^{\pi}_{\hat{r}_{\ell}}$, and hence
\[
w^{\pi*}_{r^*_{\check{\ell}}} \ \ = \ \ w^{\min}_{\check{\ell}} \ \ = \ \ \min\left\{\frac{T^{\pi}_{r}(x(q^*)) - T^{\pi}_{r^*_{\check{\ell}}}(x(q^*))}{S^{\pi}_{r^*_{\check{\ell}}} - S^{\pi}_{r}} \; : \; S^{\pi}_{r^*_{\check{\ell}}} > S^{\pi}_{r}\right\}
\ \ = \ \ \frac{T^{\pi}_{\hat{r}_{\ell+1}}(x(q^*)) - T^{\pi}_{r^*_{\check{\ell}}}(x(q^*))}{S^{\pi}_{r^*_{\check{\ell}}} - S^{\pi}_{\hat{r}_{\ell+1}}}
\]
and $r^*_{\check{\ell}+1} = \hat{r}_{\ell+1}$.
Hence claim~\ref{item:efficient paths in complete list} of the hypothesis holds for $\ell + 1$ with $\ell^*(\ell+1) = \check{\ell}+1$.
Also, $\hat{w}^{\pi}_{\hat{r}_{\ell}} = w^{\pi*}_{r^*_{\check{\ell}}} = w^{\pi*}_{r^*_{\ell^*(\ell+1)-1}}$, and hence claim~\ref{item:thresholds in complete list} of the hypothesis holds for $\ell + 1$.

\noindent
\textbf{Case~C.2}: $S^{\pi}_{r^*_{\ell^*(\ell)}} > S^{\pi}_{r^*_{\check{\ell}}}$. \\
Then $r^*_{\ell^*(\ell)} < r^*_{\check{\ell}} < \hat{r}_{\ell+1}$.
Let $\tilde{\ell} \defi \max\left\{\ell' \in \{\ell^*(\ell),\ldots,\check{\ell}-1\} \, : \, S^{\pi}_{r^*_{\ell^*(\ell)}} = S^{\pi}_{r^*_{\ell'}}\right\}$.
Thus, $S^{\pi}_{r^*_{\ell^*(\ell)}} = S^{\pi}_{r^*_{\ell'}}$ and $T^{\pi}_{r^*_{\ell^*(\ell)}}(x(q^*)) = T^{\pi}_{r^*_{\ell'}}(x(q^*))$ for all $\ell' \in \{\ell^*(\ell),\ldots,\tilde{\ell}\}$.
Then, in Algorithm~\ref{algo:threshold}, $w^{\pi*}_{r^*_{\ell'}}$ is determined by~\eqref{eqn:threshold ==} for all $\ell' = \ell^*(\ell),\ldots,\tilde{\ell}-1$.
Also, it follows from $L^{\pi}_{r^*_{\ell'}}(q^{\pi*}) = L^{\pi}_{r^*_{\ell^*(\ell)}}(q^{\pi*}) = L^{\pi}_{\hat{r}_{\ell}}(\hat{q}^{\pi})$ for all $\ell' \in \{\ell^*(\ell),\ldots,\tilde{\ell}\}$ and induction on $\ell'$ that
\begin{align*}
w^{\pi*}_{r^*_{\ell'}} \ \ & = \ \ \max\left\{w^{\pi*}_{r^*_{\ell'}-1}, \; \min\left\{w \in \left(G^{\pi}\right)^{-1}\left(L^{\pi}_{r^*_{\ell'}}(q^{\pi*})\right)\right\}\right\}
\ \ = \ \ \max\left\{w^{\pi*}_{r^*_{\ell^*(\ell)-1}}, \; \min\left\{w \in \left(G^{\pi}\right)^{-1}\left(L^{\pi}_{r^*_{\ell^*(\ell)}}(q^{\pi*})\right)\right\}\right\} \\
& = \ \ \max\left\{\hat{w}^{\pi}_{\hat{r}_{\ell-1}}, \; \min\left\{w \in \left(G^{\pi}\right)^{-1}\left(L^{\pi}_{\hat{r}_{\ell}}(\hat{q}^{\pi})\right)\right\}\right\}
\qquad \forall \ \ell' \in \{\ell^*(\ell),\ldots,\tilde{\ell}-1\}
\end{align*}
Also, in Algorithm~\ref{algo:threshold}, $w^{\pi*}_{r^*_{\tilde{\ell}}}$ is determined by~\eqref{eqn:threshold >} or~\eqref{eqn:threshold =>}.
Then
\begin{align*}
w^{\min}_{\tilde{\ell}} \ \ & = \ \ \min\left\{\frac{T^{\pi}_{r}(x(q^*)) - T^{\pi}_{r^*_{\tilde{\ell}}}(x(q^*))}{S^{\pi}_{r^*_{\tilde{\ell}}} - S^{\pi}_{r}} \; : \; S^{\pi}_{r^*_{\tilde{\ell}}} > S^{\pi}_{r}\right\} \\
& \le \ \ \frac{T^{\pi}_{\hat{r}_{\ell+1}}(x(q^*)) - T^{\pi}_{r^*_{\tilde{\ell}}}(x(q^*))}{S^{\pi}_{r^*_{\tilde{\ell}}} - S^{\pi}_{\hat{r}_{\ell+1}}}
\ \ = \ \ \frac{T^{\pi}_{\hat{r}_{\ell+1}}(x(q^*)) - T^{\pi}_{r^*_{\ell^*(\ell)}}(x(q^*))}{S^{\pi}_{r^*_{\ell^*(\ell)}} - S^{\pi}_{\hat{r}_{\ell+1}}}
\ \ = \ \ \frac{T^{\pi}_{\hat{r}_{\ell+1}}(x(\hat{q})) - T^{\pi}_{\hat{r}_{\ell}}(x(\hat{q}))}{S^{\pi}_{\hat{r}_{\ell}} - S^{\pi}_{\hat{r}_{\ell+1}}}
\ \ = \ \ \hat{w}^{\pi}_{\hat{r}_{\ell}}
\end{align*}
Suppose that $w^{\min}_{\tilde{\ell}} < \hat{w}^{\pi}_{\hat{r}_{\ell}}$.
Then
\begin{align*}
& w^{\min}_{\tilde{\ell}} \ \ = \ \ \frac{T^{\pi}_{r^*_{\tilde{\ell}+1}}(x(q^*)) - T^{\pi}_{r^*_{\tilde{\ell}}}(x(q^*))}{S^{\pi}_{r^*_{\tilde{\ell}}} - S^{\pi}_{r^*_{\tilde{\ell}+1}}}
\ \ < \ \ \hat{w}^{\pi}_{\hat{r}_{\ell}} \\
\Rightarrow \qquad & T^{\pi}_{r^*_{\tilde{\ell}+1}}(x(q^*)) + \hat{w}^{\pi}_{\hat{r}_{\ell}} S^{\pi}_{r^*_{\tilde{\ell}+1}} \ \ < \ \ T^{\pi}_{r^*_{\tilde{\ell}}}(x(q^*)) + \hat{w}^{\pi}_{\hat{r}_{\ell}} S^{\pi}_{r^*_{\tilde{\ell}}}
\ \ = \ \ T^{\pi}_{r^*_{\ell^*(\ell)}}(x(q^*)) + \hat{w}^{\pi}_{\hat{r}_{\ell}} S^{\pi}_{r^*_{\ell^*(\ell)}}
\ \ = \ \ T^{\pi}_{\hat{r}_{\ell}}(x(q^*)) + \hat{w}^{\pi}_{\hat{r}_{\ell}} S^{\pi}_{\hat{r}_{\ell}}
\end{align*}
which contradicts~\eqref{eqn:min cost ell}.
Thus $w^{\pi*}_{r^*_{\tilde{\ell}}} = w^{\min}_{\tilde{\ell}} = \hat{w}^{\pi}_{\hat{r}_{\ell}}$ and $T^{\pi}_{r^*_{\tilde{\ell}+1}}(x(q^*)) + \hat{w}^{\pi}_{\hat{r}_{\ell}} S^{\pi}_{r^*_{\tilde{\ell}+1}} = T^{\pi}_{r^*_{\tilde{\ell}}}(x(q^*)) + \hat{w}^{\pi}_{\hat{r}_{\ell}} S^{\pi}_{r^*_{\tilde{\ell}}} = T^{\pi}_{r^*_{\ell^*(\ell)}}(x(q^*)) + \hat{w}^{\pi}_{\hat{r}_{\ell}} S^{\pi}_{r^*_{\ell^*(\ell)}} = T^{\pi}_{\hat{r}_{\ell}}(x(q^*)) + \hat{w}^{\pi}_{\hat{r}_{\ell}} S^{\pi}_{\hat{r}_{\ell}}$.
If $\tilde{\ell}+1 < \check{\ell}$, then let $\tilde{\tilde{\ell}} \defi \max\left\{\ell' \in \{\tilde{\ell}+1,\ldots,\check{\ell}\} \, : \, S^{\pi}_{r^*_{\tilde{\ell}+1}} = S^{\pi}_{r^*_{\ell'}}\right\}$.
Then the argument is repeated by induction, to obtain that $w^{\pi*}_{r^*_{\check{\ell}-1}} = \hat{w}^{\pi}_{\hat{r}_{\ell}}$ and
\begin{equation}
T^{\pi}_{r^*_{\check{\ell}}}(x(q^*)) + \hat{w}^{\pi}_{\hat{r}_{\ell}} S^{\pi}_{r^*_{\check{\ell}}}
\ \ = \ \ T^{\pi}_{r^*_{\ell^*(\ell)}}(x(q^*)) + \hat{w}^{\pi}_{\hat{r}_{\ell}} S^{\pi}_{r^*_{\ell^*(\ell)}}
\ \ = \ \ T^{\pi}_{\hat{r}_{\ell}}(x(q^*)) + \hat{w}^{\pi}_{\hat{r}_{\ell}} S^{\pi}_{\hat{r}_{\ell}}
\ \ \overset{\eqref{eqn:min cost ell}}{=} \ \ T^{\pi}_{\hat{r}_{\ell+1}}(x(q^*)) + \hat{w}^{\pi}_{\hat{r}_{\ell}} S^{\pi}_{\hat{r}_{\ell+1}}
\label{eqn:equal cost ell = check ell 2}
\end{equation}

\noindent
\textbf{Case~C.2.1}: $S^{\pi}_{r^*_{\ell^*(\ell)}} > S^{\pi}_{r^*_{\check{\ell}}} = S^{\pi}_{\hat{r}_{\ell+1}}$. \\
It follows from~\eqref{eqn:equal cost ell = check ell 2} that $T^{\pi}_{r^*_{\check{\ell}}}(x(q^*)) = T^{\pi}_{\hat{r}_{\ell+1}}(x(q^*))$.
Thus, in Algorithm~\ref{algo:threshold}, $w^{\pi*}_{r^*_{\check{\ell}}}$ is determined by~\eqref{eqn:threshold ==}, and $r^*_{\check{\ell}+1} = \hat{r}_{\ell+1}$.
Hence claim~\ref{item:efficient paths in complete list} of the hypothesis holds for $\ell + 1$ with $\ell^*(\ell+1) = \check{\ell}+1$.
Also, $w^{\pi*}_{r^*_{\check{\ell}}} = \max\left\{w^{\pi*}_{r^*_{\check{\ell}}-1}, \; \min\left\{w \in \left(G^{\pi}\right)^{-1}\left(L^{\pi}_{r^*_{\check{\ell}}}(q^{\pi*})\right)\right\}\right\} = \max\left\{w^{\pi*}_{r^*_{\check{\ell}-1}}, \; \min\left\{w \in \left(G^{\pi}\right)^{-1}\left(L^{\pi}_{r^*_{\check{\ell}}}(q^{\pi*})\right)\right\}\right\} = \max\left\{\hat{w}^{\pi}_{\hat{r}_{\ell}}, \; \min\left\{w \in \left(G^{\pi}\right)^{-1}\left(L^{\pi}_{\hat{r}_{\ell}}(\hat{q}^{\pi})\right)\right\}\right\} = \hat{w}^{\pi}_{\hat{r}_{\ell}}$, where the last equality follows from Lemma~\ref{lem:optimal thresholds}\ref{item:thresholds quantiles}.
Thus $\hat{w}^{\pi}_{\hat{r}_{\ell}} = w^{\pi*}_{r^*_{\check{\ell}}} = w^{\pi*}_{r^*_{\ell^*(\ell+1)-1}}$, and hence claim~\ref{item:thresholds in complete list} of the hypothesis holds for $\ell + 1$.

\noindent
\textbf{Case~C.2.2}: $S^{\pi}_{r^*_{\ell^*(\ell)}} > S^{\pi}_{r^*_{\check{\ell}}} > S^{\pi}_{\hat{r}_{\ell+1}}$. \\
Then, in Algorithm~\ref{algo:threshold}, $w^{\pi*}_{r^*_{\check{\ell}}}$ is determined by~\eqref{eqn:threshold >} or~\eqref{eqn:threshold =>}.
Then
\[
w^{\min}_{\check{\ell}} \ \ = \ \ \min\left\{\frac{T^{\pi}_{r}(x(q^*)) - T^{\pi}_{r^*_{\check{\ell}}}(x(q^*))}{S^{\pi}_{r^*_{\check{\ell}}} - S^{\pi}_{r}} \; : \; S^{\pi}_{r^*_{\check{\ell}}} > S^{\pi}_{r}\right\} \ \ \le \ \ \frac{T^{\pi}_{\hat{r}_{\ell+1}}(x(q^*)) - T^{\pi}_{r^*_{\check{\ell}}}(x(q^*))}{S^{\pi}_{r^*_{\check{\ell}}} - S^{\pi}_{\hat{r}_{\ell+1}}}
\ \ \overset{\eqref{eqn:equal cost ell = check ell 2}}{=} \ \ \hat{w}^{\pi}_{\hat{r}_{\ell}}
\]
Suppose that $w^{\min}_{\check{\ell}} < \hat{w}^{\pi}_{\hat{r}_{\ell}}$.
Then
\begin{align*}
& w^{\min}_{\check{\ell}} \ \ = \ \ \frac{T^{\pi}_{r^*_{\check{\ell}+1}}(x(q^*)) - T^{\pi}_{r^*_{\check{\ell}}}(x(q^*))}{S^{\pi}_{r^*_{\check{\ell}}} - S^{\pi}_{r^*_{\check{\ell}+1}}}
\ \ < \ \ \hat{w}^{\pi}_{\hat{r}_{\ell}} \\
\Rightarrow \qquad & T^{\pi}_{r^*_{\check{\ell}+1}}(x(q^*)) + \hat{w}^{\pi}_{\hat{r}_{\ell}} S^{\pi}_{r^*_{\check{\ell}+1}} \ \ < \ \ T^{\pi}_{r^*_{\check{\ell}}}(x(q^*)) + \hat{w}^{\pi}_{\hat{r}_{\ell}} S^{\pi}_{r^*_{\check{\ell}}}
\ \ \overset{\eqref{eqn:equal cost ell = check ell 2}}{=} \ \ T^{\pi}_{\hat{r}_{\ell}}(x(q^*)) + \hat{w}^{\pi}_{\hat{r}_{\ell}} S^{\pi}_{\hat{r}_{\ell}} \\
& \hspace{80mm} \overset{\eqref{eqn:min cost ell}}{=} \ \ T^{\pi}_{\hat{r}_{\ell+1}}(x(q^*)) + \hat{w}^{\pi}_{\hat{r}_{\ell}} S^{\pi}_{\hat{r}_{\ell+1}}
\end{align*}
which contradicts~\eqref{eqn:min cost ell}.
Thus $w^{\min}_{\check{\ell}} = \hat{w}^{\pi}_{\hat{r}_{\ell}}$, and hence
\[
w^{\pi*}_{r^*_{\check{\ell}}} \ \ = \ \ w^{\min}_{\check{\ell}} \ \ = \ \ \min\left\{\frac{T^{\pi}_{r}(x(q^*)) - T^{\pi}_{r^*_{\check{\ell}}}(x(q^*))}{S^{\pi}_{r^*_{\check{\ell}}} - S^{\pi}_{r}} \; : \; S^{\pi}_{r^*_{\check{\ell}}} > S^{\pi}_{r}\right\}
\ \ = \ \ \frac{T^{\pi}_{\hat{r}_{\ell+1}}(x(q^*)) - T^{\pi}_{r^*_{\check{\ell}}}(x(q^*))}{S^{\pi}_{r^*_{\check{\ell}}} - S^{\pi}_{\hat{r}_{\ell+1}}}
\]
and $r^*_{\check{\ell}+1} = \hat{r}_{\ell+1}$.
Hence claim~\ref{item:efficient paths in complete list} of the hypothesis holds for $\ell + 1$ with $\ell^*(\ell+1) = \check{\ell}+1$.
Also, $\hat{w}^{\pi}_{\hat{r}_{\ell}} = w^{\pi*}_{r^*_{\check{\ell}}} = w^{\pi*}_{r^*_{\ell^*(\ell+1)-1}}$, and hence claim~\ref{item:thresholds in complete list} of the hypothesis holds for $\ell + 1$.
Therefore the hypothesis has been established.
In particular, it follows that $\hat{\mathcal{R}}^{\pi}_{\min} \subset \mathcal{R}^{\pi*}_{\min}$.

Next we show that $(q^*,\lambda^*,w^*)$ satisfies the optimality conditions~\eqref{eqn:q complete BiTAP optimality1}--\eqref{eqn:q complete BiTAP optimality3}.
Consider any $\pi \in \Pi$ and any successive pair $r^*_{\ell},r^*_{\ell+1} \in \mathcal{R}^{\pi*}_{\min}$.
If $w^{\pi*}_{r^*_{\ell}}$ is determined by~\eqref{eqn:threshold >} or~\eqref{eqn:threshold =>}, then it holds that $S^{\pi}_{r^*_{\ell}} > S^{\pi}_{r^*_{\ell+1}}$ and
\[
w^{\pi*}_{r^*_{\ell}} \ \ = \ \ \frac{T^{\pi}_{r^*_{\ell+1}}(x(q^*)) - T^{\pi}_{r^*_{\ell}}(x(q^*))}{S^{\pi}_{r^*_{\ell}} - S^{\pi}_{r^*_{\ell+1}}}
\]
(Note that if $w^{\min}_{\ell} > \overline{w}^{\pi}$ or $w^{\pi*}_{r^*_{\ell}}$ is determined by~\eqref{eqn:threshold =0}, then there is no $r^*_{\ell+1} \in \mathcal{R}^{\pi*}_{\min}$.)
If $w^{\pi*}_{r^*_{\ell}}$ is determined by~\eqref{eqn:threshold ==}, then $S^{\pi}_{r^*_{\ell}} = S^{\pi}_{r^*_{\ell+1}}$ and $T^{\pi}_{r^*_{\ell}}(x(q^*)) = T^{\pi}_{r^*_{\ell+1}}(x(q^*))$, and thus $T^{\pi}_{r^*_{\ell}}(x(q^*)) + w S^{\pi}_{r^*_{\ell}} = T^{\pi}_{r^*_{\ell+1}}(x(q^*)) + w S^{\pi}_{r^*_{\ell+1}}$ for all~$w$.
Thus, in all cases it holds that
\begin{equation}
\label{eqn:successive paths 2}
T^{\pi}_{r^*_{\ell}}(x(q^*)) + w^{\pi*}_{r^*_{\ell}} S^{\pi}_{r^*_{\ell}} \ \ = \ \ T^{\pi}_{r^*_{\ell+1}}(x(q^*)) + w^{\pi*}_{r^*_{\ell}} S^{\pi}_{r^*_{\ell+1}}
\end{equation}
Also, note that $w^{\pi*}_{r} = w^{\pi*}_{r^*_{\ell}}$ for all $r \in \{r^*_{\ell},\ldots,r^*_{\ell+1}-1\}$, and thus
\begin{equation}
\sum_{r=r^*_{\ell}}^{r^*_{\ell+1}-1} \left(S^{\pi}_{r} - S^{\pi}_{r+1}\right) w^{\pi*}_{r}
\ \ = \ \ \sum_{r=r^*_{\ell}}^{r^*_{\ell+1}-1} \left(S^{\pi}_{r} - S^{\pi}_{r+1}\right) w^{\pi*}_{r^*_{\ell}}
\ \ = \ \ \left(S^{\pi}_{r^*_{\ell}} - S^{\pi}_{r^*_{\ell+1}}\right) w^{\pi*}_{r^*_{\ell}}
\label{eqn:cancellation unused optimal paths}
\end{equation}
Similarly, $w^{\pi*}_{r} = w^{\pi*}_{r^*_{L^*}} = \overline{w}^{\pi}$ for all $r \in \{r^*_{L^*},\ldots,R^{\pi}\}$, and thus
\begin{equation}
\sum_{r=r^*_{L^*}}^{R^{\pi}-1} \left(S^{\pi}_{r} - S^{\pi}_{r+1}\right) w^{\pi*}_{r}
\ \ = \ \ \sum_{r=r^*_{L^*}}^{R^{\pi}-1} \left(S^{\pi}_{r} - S^{\pi}_{r+1}\right) \overline{w}^{\pi}
\ \ = \ \ \left(S^{\pi}_{r^*_{L^*}} - S^{\pi}_{R^{\pi}}\right) \overline{w}^{\pi}
\label{eqn:cancellation unused optimal paths2}
\end{equation}

\noindent
\textbf{Case~D.1}: $r \in \mathcal{R}^{\pi*}_{\min}$. \\
Let $r = r^*_{\ell}$ for some $\ell \in \{1,\ldots,L^*\}$.
Then
\begin{align*}
& D^{\pi} \left[T^{\pi}_{r}(x(q^*)) + \sum_{r' = r}^{R^{\pi}-1} \left(S^{\pi}_{r'} - S^{\pi}_{r'+1}\right) w^{\pi*}_{r'}\right] + \lambda^{\pi*} \\
& \overset{\eqref{eqn:choice of dual 2},\eqref{eqn:cancellation unused optimal paths},\eqref{eqn:cancellation unused optimal paths2}}{=} \ \ D^{\pi} \left[T^{\pi}_{r^*_{\ell}}(x(q^*)) + \sum_{\ell' = \ell}^{L^*-1} \left(S^{\pi}_{r^*_{\ell'}} - S^{\pi}_{r^*_{\ell'+1}}\right) w^{\pi*}_{r^*_{\ell'}} + \left(S^{\pi}_{r^*_{L^*}} - S^{\pi}_{R^{\pi}}\right) \overline{w}^{\pi}\right] \\
& \hspace{30mm} {} - D^{\pi} \left[T^{\pi}_{r^*_{L^*}}(x(q^*)) + \overline{w}^{\pi} \left(S^{\pi}_{r^*_{L^*}} - S^{\pi}_{R^{\pi}}\right)\right] \\
& = \ \ D^{\pi} \sum_{\ell' = \ell}^{L^*-1} \left[T^{\pi}_{r^*_{\ell'}}(x(q^*)) - T^{\pi}_{r^*_{\ell'+1}}(x(q^*)) + \left(S^{\pi}_{r^*_{\ell'}} - S^{\pi}_{r^*_{\ell'+1}}\right) w^{\pi*}_{r^*_{\ell'}}\right] \ \ \overset{\eqref{eqn:successive paths 2}}{=} \ \ 0
\end{align*}
Thus conditions~\eqref{eqn:q complete BiTAP optimality1}--\eqref{eqn:q complete BiTAP optimality3} hold for all $r \in \mathcal{R}^{\pi*}_{\min}$.

\noindent
\textbf{Case~D.2}: $r \notin \mathcal{R}^{\pi*}_{\min}$. \\
It follows from $\hat{\mathcal{R}}^{\pi}_{\min} \subset \mathcal{R}^{\pi*}_{\min}$ that $r \notin \hat{\mathcal{R}}^{\pi}_{\min}$.
If $r \in \hat{\mathcal{R}}^{\pi}$, then it follows from Lemma~\ref{lem:optimal thresholds}\ref{item:positive flow paths are efficient} that $\hat{q}^{\pi}_{r} = 0$, and thus $q^{\pi*}_{r} = 0$.
If $r \in \mathcal{R}^{\pi} \setminus \hat{\mathcal{R}}^{\pi}$ then $q^{\pi*}_{r} = 0$.
Thus condition~\eqref{eqn:q complete BiTAP optimality1} holds.
Next we show that condition~\eqref{eqn:q complete BiTAP optimality3} also holds.

\noindent
\textbf{Case~D.2.1}: $r < r^*_{L^*}$. \\
Then either $r < r^*_{1}$, or there is a successive pair $r^*_{\ell},r^*_{\ell+1} \in \mathcal{R}^{\pi*}_{\min}$ such that $r^*_{\ell} < r < r^*_{\ell+1}$.
If $r < r^*_{1}$, then $\ell = 0$ and $w^{\pi*}_{r'} = w^{\pi*}_{r^*_{0}} = \underline{w}^{\pi}$ for all $r' \in \{r,\ldots,r^*_{1}-1\}$, and if $r^*_{\ell} < r < r^*_{\ell+1}$, then $w^{\pi*}_{r'} = w^{\pi*}_{r^*_{\ell}}$ for all $r' \in \{r,\ldots,r^*_{\ell+1}-1\}$.
Then
\begin{align*}
& D^{\pi} \left[T^{\pi}_{r}(x(q^*)) + \sum_{r' = r}^{R^{\pi}-1} \left(S^{\pi}_{r'} - S^{\pi}_{r'+1}\right) w^{\pi*}_{r'}\right] + \lambda^{\pi*} \\
& \overset{\eqref{eqn:choice of dual 2},\eqref{eqn:cancellation unused optimal paths},\eqref{eqn:cancellation unused optimal paths2}}{=} \ \ D^{\pi} \left[T^{\pi}_{r}(x(q^*)) + \sum_{r'=r}^{r^*_{\ell+1}-1} \left(S^{\pi}_{r'} - S^{\pi}_{r'+1}\right) w^{\pi*}_{r^*_{\ell}} + \sum_{\ell' = \ell+1}^{L^*-1} \left(S^{\pi}_{r^*_{\ell'}} - S^{\pi}_{r^*_{\ell'+1}}\right) w^{\pi*}_{r^*_{\ell'}} + \left(S^{\pi}_{r^*_{L^*}} - S^{\pi}_{R^{\pi}}\right) \overline{w}^{\pi}\right] \\
& \hspace{20mm} {} - D^{\pi} \left[T^{\pi}_{r^*_{L^*}}(x(q^*)) + \overline{w}^{\pi} \left(S^{\pi}_{r^*_{L^*}} - S^{\pi}_{R^{\pi}}\right)\right] \\
& = \ \ D^{\pi} \left[T^{\pi}_{r}(x(q^*)) - T^{\pi}_{r^*_{\ell+1}}(x(q^*)) + T^{\pi}_{r^*_{\ell+1}}(x(q^*)) - T^{\pi}_{r^*_{L^*}}(x(q^*)) + \left(S^{\pi}_{r} - S^{\pi}_{r^*_{\ell+1}}\right) w^{\pi*}_{r^*_{\ell}} + \sum_{\ell' = \ell+1}^{L^*-1} \left(S^{\pi}_{r^*_{\ell'}} - S^{\pi}_{r^*_{\ell'+1}}\right) w^{\pi*}_{r^*_{\ell'}}\right] \\
& = \ \ D^{\pi} \left[T^{\pi}_{r}(x(q^*)) - T^{\pi}_{r^*_{\ell+1}}(x(q^*)) + \left(S^{\pi}_{r} - S^{\pi}_{r^*_{\ell+1}}\right) w^{\pi*}_{r^*_{\ell}}\right] + D^{\pi} \sum_{\ell' = \ell+1}^{L^*-1} \left[T^{\pi}_{r^*_{\ell'}}(x(q^*)) - T^{\pi}_{r^*_{\ell'+1}}(x(q^*)) + \left(S^{\pi}_{r^*_{\ell'}} - S^{\pi}_{r^*_{\ell'+1}}\right) w^{\pi*}_{r^*_{\ell'}}\right] \\
& \overset{\eqref{eqn:successive paths 2}}{=} \ \ D^{\pi} \left[T^{\pi}_{r}(x(q^*)) - T^{\pi}_{r^*_{\ell+1}}(x(q^*)) + \left(S^{\pi}_{r} - S^{\pi}_{r^*_{\ell+1}}\right) w^{\pi*}_{r^*_{\ell}}\right]
\end{align*}
If $r < r^*_{1}$, then $w^{\pi*}_{r^*_{\ell}} = w^{\pi*}_{r^*_{0}} = \underline{w}^{\pi}$ and $r^*_{\ell+1} = r^*_{1}$.
Then it follows from~\eqref{eqn:efficient path 1} that
\[
T^{\pi}_{r}(x(q^*)) - T^{\pi}_{r^*_{\ell+1}}(x(q^*)) + \left(S^{\pi}_{r} - S^{\pi}_{r^*_{\ell+1}}\right) w^{\pi*}_{r^*_{\ell}} \ \ = \ \ \left[T^{\pi}_{r}(x(q^*)) + \underline{w}^{\pi} S^{\pi}_{r}\right] - \left[T^{\pi}_{r^*_{1}}(x(q^*)) + \underline{w}^{\pi} S^{\pi}_{r^*_{1}}\right] \ \ > \ \ 0
\]
Otherwise, if $r > r^*_{1}$, then consider the following cases.

\noindent
\textbf{Case~D.2.1.1}: $w^{\pi*}_{r^*_{\ell}}$ is determined by~\eqref{eqn:threshold >} or~\eqref{eqn:threshold =>} and $S^{\pi}_{r^*_{\ell}} > S^{\pi}_{r}$. \\
Then
\begin{align*}
& w^{\pi*}_{r^*_{\ell}} \ \ = \ \ w^{\min}_{\ell} \ \ = \ \ \frac{T^{\pi}_{r^*_{\ell+1}}(x(q^*)) - T^{\pi}_{r^*_{\ell}}(x(q^*))}{S^{\pi}_{r^*_{\ell}} - S^{\pi}_{r^*_{\ell+1}}}
\ \ < \ \ \frac{T^{\pi}_{r}(x(q^*)) - T^{\pi}_{r^*_{\ell}}(x(q^*))}{S^{\pi}_{r^*_{\ell}} - S^{\pi}_{r}} \\
\Rightarrow \qquad & T^{\pi}_{r^*_{\ell+1}}(x(q^*)) + w^{\pi*}_{r^*_{\ell}} S^{\pi}_{r^*_{\ell+1}} \ \ < \ \ T^{\pi}_{r}(x(q^*)) + w^{\pi*}_{r^*_{\ell}} S^{\pi}_{r} \\
\Rightarrow \qquad & D^{\pi} \left[T^{\pi}_{r}(x(q^*)) - T^{\pi}_{r^*_{\ell+1}}(x(q^*)) + \left(S^{\pi}_{r} - S^{\pi}_{r^*_{\ell+1}}\right) w^{\pi*}_{r^*_{\ell}}\right]
\ \ > \ \ 0
\end{align*}

\noindent
\textbf{Case~D.2.1.2}: $w^{\pi*}_{r^*_{\ell}}$ is determined by~\eqref{eqn:threshold =>} and $S^{\pi}_{r^*_{\ell}} = S^{\pi}_{r}$. \\
Then $T^{\pi}_{r^*_{\ell}}(x(q^*)) < T^{\pi}_{r}(x(q^*))$.
Thus
\begin{align*}
& D^{\pi} \left[T^{\pi}_{r}(x(q^*)) - T^{\pi}_{r^*_{\ell+1}}(x(q^*)) + \left(S^{\pi}_{r} - S^{\pi}_{r^*_{\ell+1}}\right) w^{\pi*}_{r^*_{\ell}}\right] \\
= \ \ & D^{\pi} \left[T^{\pi}_{r}(x(q^*)) - T^{\pi}_{r^*_{\ell}}(x(q^*)) + T^{\pi}_{r^*_{\ell}}(x(q^*)) - T^{\pi}_{r^*_{\ell+1}}(x(q^*)) + \left(S^{\pi}_{r^*_{\ell}} - S^{\pi}_{r^*_{\ell+1}}\right) w^{\pi*}_{r^*_{\ell}}\right] \\
\overset{\eqref{eqn:successive paths 2}}{=} \ \ & D^{\pi} \left[T^{\pi}_{r}(x(q^*)) - T^{\pi}_{r^*_{\ell}}(x(q^*))\right]
\ \ > \ \ 0
\end{align*}

\noindent
\textbf{Case~D.2.1.3}: $w^{\pi*}_{r^*_{\ell}}$ is determined by~\eqref{eqn:threshold ==}. \\
Then $S^{\pi}_{r^*_{\ell}} = S^{\pi}_{r} = S^{\pi}_{r^*_{\ell+1}}$ and $T^{\pi}_{r^*_{\ell}}(x(q^*)) = T^{\pi}_{r^*_{\ell+1}}(x(q^*)) < T^{\pi}_{r}(x(q^*))$, and thus
\[
D^{\pi} \left[T^{\pi}_{r}(x(q^*)) - T^{\pi}_{r^*_{\ell+1}}(x(q^*)) + \left(S^{\pi}_{r} - S^{\pi}_{r^*_{\ell+1}}\right) w^{\pi*}_{r^*_{\ell}}\right]
\ \ = \ \ D^{\pi} \left[T^{\pi}_{r}(x(q^*)) - T^{\pi}_{r^*_{\ell+1}}(x(q^*))\right]
\ \ > \ \ 0
\]
Therefore condition~\eqref{eqn:q complete BiTAP optimality3} holds for all $r < r^*_{L^*}$.

\noindent
\textbf{Case~D.2.2}: $r > r^*_{L^*}$. \\
Then
\begin{align*}
& D^{\pi} \left[T^{\pi}_{r}(x(q^*)) + \sum_{r' = r}^{R^{\pi}-1} \left(S^{\pi}_{r'} - S^{\pi}_{r'+1}\right) w^{\pi*}_{r'}\right] + \lambda^{\pi*} \\
& \overset{\eqref{eqn:choice of dual 2}}{=} \ \ D^{\pi} \left[T^{\pi}_{r}(x(q^*)) + \sum_{r'=r}^{R^{\pi}-1} \left(S^{\pi}_{r'} - S^{\pi}_{r'+1}\right) \overline{w}^{\pi}\right] - D^{\pi} \left[T^{\pi}_{r^*_{L^*}}(x(q^*)) + \overline{w}^{\pi} \left(S^{\pi}_{r^*_{L^*}} - S^{\pi}_{R^{\pi}}\right)\right] \\
& = \ \ D^{\pi} \left[T^{\pi}_{r}(x(q^*)) - T^{\pi}_{r^*_{L^*}}(x(q^*)) + \left(S^{\pi}_{r} - S^{\pi}_{r^*_{L^*}}\right) \overline{w}^{\pi}\right]
\end{align*}

\noindent
\textbf{Case~D.2.2.1}: $w^{\pi*}_{r^*_{L^*}}$ is determined by~\eqref{eqn:threshold >} or~\eqref{eqn:threshold =>} and $S^{\pi}_{r^*_{L^*}} > S^{\pi}_{r}$. \\
Then
\begin{align*}
& w^{\pi*}_{r^*_{L^*}} \ \ = \ \ \overline{w}^{\pi} \ \ < \ \ w^{\min}_{L^*} \ \ = \ \ \min\left\{\frac{T^{\pi}_{r'}(x(q^*)) - T^{\pi}_{r^*_{L^*}}(x(q^*))}{S^{\pi}_{r^*_{L^*}} - S^{\pi}_{r'}} \; : \; S^{\pi}_{r^*_{L^*}} > S^{\pi}_{r'}\right\}
\ \ \le \ \ \frac{T^{\pi}_{r}(x(q^*)) - T^{\pi}_{r^*_{L^*}}(x(q^*))}{S^{\pi}_{r^*_{L^*}} - S^{\pi}_{r}} \\
\Rightarrow \qquad & T^{\pi}_{r^*_{L^*}}(x(q^*)) + \overline{w}^{\pi} S^{\pi}_{r^*_{L^*}} \ \ < \ \ T^{\pi}_{r}(x(q^*)) + \overline{w}^{\pi} S^{\pi}_{r} \\
\Rightarrow \qquad & D^{\pi} \left[T^{\pi}_{r}(x(q^*)) - T^{\pi}_{r^*_{L^*}}(x(q^*)) + \left(S^{\pi}_{r} - S^{\pi}_{r^*_{L^*}}\right) \overline{w}^{\pi}\right]
\ \ > \ \ 0
\end{align*}

\noindent
\textbf{Case~D.2.2.2}: $w^{\pi*}_{r^*_{L^*}}$ is determined by~\eqref{eqn:threshold =>} and $S^{\pi}_{r^*_{L^*}} = S^{\pi}_{r}$. \\
Then $T^{\pi}_{r^*_{L^*}}(x(q^*)) < T^{\pi}_{r}(x(q^*))$.
Thus
\begin{align*}
D^{\pi} \left[T^{\pi}_{r}(x(q^*)) - T^{\pi}_{r^*_{L^*}}(x(q^*)) + \left(S^{\pi}_{r} - S^{\pi}_{r^*_{L^*}}\right) \overline{w}^{\pi}\right]
\ \ = \ \ D^{\pi} \left[T^{\pi}_{r}(x(q^*)) - T^{\pi}_{r^*_{L^*}}(x(q^*))\right]
\ \ > \ \ 0
\end{align*}

\noindent
\textbf{Case~D.2.2.3}: $w^{\pi*}_{r^*_{L^*}}$ is determined by~\eqref{eqn:threshold ==}. \\
Then $r^*_{L^*+1}$ set in~\eqref{eqn:efficient path ==} would satisfy $r^*_{L^*+1} > r^*_{L^*}$ and $r^*_{L^*+1} \in \mathcal{R}^{\pi*}_{\min}$, which cannot hold, and thus $w^{\pi}_{r^*_{L^*}}$ is not determined by~\eqref{eqn:threshold ==}.

\noindent
\textbf{Case~D.2.2.4}: $w^{\pi*}_{r^*_{L^*}}$ is determined by~\eqref{eqn:threshold =0}. \\
Then $S^{\pi}_{r^*_{L^*}} = S^{\pi}_{r} = S^{\pi}_{R^{\pi}}$ and $T^{\pi}_{r}(x(q^*)) > T^{\pi}_{r^*_{L^*}}(x(q^*))$.
Thus
\begin{align*}
D^{\pi} \Bigg[T^{\pi}_{r}(x(q^*)) - T^{\pi}_{r^*_{L^*}}(x(q^*)) + \left(S^{\pi}_{r} - S^{\pi}_{r^*_{L^*}}\right) \overline{w}^{\pi}\Bigg]
\ \ = \ \ & D^{\pi} \left[T^{\pi}_{r}(x(q^*)) - T^{\pi}_{r^*_{L^*}}(x(q^*))\right]
\ \ > \ \ 0
\end{align*}
Therefore conditions~\eqref{eqn:q complete BiTAP optimality1}--\eqref{eqn:q complete BiTAP optimality3} hold for all $r \in \mathcal{R}^{\pi}$.
\Halmos
\endproof

\section{Supplementary Numerical Data and Results}
\label{app:numerical-experiments}
\renewcommand{\arraystretch}{1.18}

In this section of the appendix, we provide additional data and numerical results to evaluate the performance of BTA.

\subsection{Test Networks}

Table~\ref{tab:exp_networks} gives summary network data for the benchmark networks used in the numerical tests.

\begin{table}[t]
\centering
\small
\caption{Summary network data for the benchmark networks.}
\label{tab:exp_networks}
\begin{tabular}{p{0.39\textwidth}rrrrrrr}
\toprule
Network & Nodes & Links & Zones & OD records & Total demand & $\alpha$ & $\beta$ \\
\midrule
Braess & 4 & 5 & 1 & 2 & 6 & $0.15$ & $4$ \\
Sioux Falls & 24 & 76 & 24 & 576 & 360,600 & $0.15$ & $4$ \\
Friedrichshain & 224 & 523 & 23 & 506 & 11,205 & $0.15$ & $4$ \\
Berlin Prenzlauerberg  & 352 & 749 & 38 & 1,406 & 16,660 & $0.15$ & $4$ \\
Berlin Tiergarten  & 359 & 766 & 26 & 644 & 10,755 & $0.15$ & $4$ \\
Berlin Mitte & 397 & 871 & 36 & 1,260 & 11,482 & $0.15$ & $4$ \\
Anaheim & 416 & 914 & 38 & 1,406 & 104,694 & $0.15$ & $4$ \\
Barcelona & 930 & 2,522 & 97 & 7,922 & 184,680 & $0.15$ & $4$ \\
Winnipeg-Asym & 948 & 2,535 & 135 & 4,345 & 1,361,475 & $0.1$ & $1.5$ \\
Berlin M.--P.--F. & 974 & 2,184 & 98 & 9,505 & 23,648 & $0.15$ & $4$ \\
Winnipeg & 1,040 & 2,836 & 135 & 4,345 & 64,784 & $0.15$ & $4$ \\
Terrassa-Asym & 1,603 & 3,264 & 55 & 2,215 & 25,225,747 & $0.1$ & $1.5$ \\
Hessen-Asym & 4,660 & 6,674 & 195 & 17,213 & 71,250,600 & $0.1$ & $1.5$ \\
Goldcoast & 4,783 & 11,140 & 1,068 & 1,140,624 & 139,256 & $0.15$ & $4$ \\
Chicago Regional & 12,979 & 39,018 & 1,771 & 3,136,441 & 1,360,428 & $0.15$ & $4$ \\
Berlin Center & 12,981 & 28,376 & 865 & 49,688 & 168,222 & $0.15$ & $4$ \\
Philadelphia & 13,389 & 40,003 & 1,489 & 1,151,166 & 18,503,872 & $0.15$ & $4$ \\
\bottomrule
\end{tabular}
\end{table}

\subsection{BTA Subroutine Parameters}

BTA was used with two solvers for solving problem~\eqref{eqn:q reformulated BiTAP objective} with subsets of paths~$\hat{\mathcal{R}}^{\pi}$ in Line~10 of Algorithm~\ref{algo:BTA}.
The FISTA method of \citet{sujanani2024efficient} was used for the smooth instances, and the two-cut adaptive proximal bundle method (Ad-GPB) of \citet{monteiro2024parameter} was used for the nonsmooth instances.
Table~\ref{tab:exp_master_params} gives the parameter settings chosen for these solvers, including the network-specific initializations for Ad-GPB.
FISTA estimates the Lipschitz constant from an initial Lipschitz estimate $M_{0}$ with backtracking, using backtracking factor $\beta$.
We use the unknown-optimum variant of Ad-GPB with a two-cut bundle update; see \citet{monteiro2024parameter} for details.


\begin{table}[t]
\centering
\small
\caption{Parameter settings for the two solvers.}
\label{tab:exp_master_params}
\begin{tabular}{llr}
\toprule
Solver & Parameter & Value \\
\midrule
FISTA  & Restricted-master tolerance $\widehat\varepsilon$ & $10^{-6}$ \\
       & Backtracking factor $\beta$        & $1.25$ \\
       & Initial Lipschitz estimate $M_{0}$ & $50$ \\
       & Metric                             & fraction ($q$) coordinates \\
       & Lipschitz estimate                 & allowed to decrease \\
       & Backtracking acceptance            & scaled to iterate magnitude \\
\addlinespace
Ad-GPB & Initial tolerance $\widehat{\varepsilon}$ & $10^{-3}$ \\
       & Initial acceptance parameter $\beta_1$ & $1/4$ \\
       & $(\tau,\lambda_1)$         & $(0.5,10^{-3})$ \\
       & Maximum number of affine cuts in the model & $2$ \\
\addlinespace
  Both & Path pool                          & monotone retention \\
       & Outer stopping rule                & PSPP solution adds no new path \\
\bottomrule
\end{tabular}
\end{table}

\subsection{Numerical Results Comparing Different Methods}

Table~\ref{tab:exp_smooth_runtimes} compares the computational performance of Algorithm~BTA and the methods of \citet{daga:83}, \citet{dial:96}, \citet{marc:96b}, and \citet{xie:26} on smooth instances~$i$.
If method~$m$ achieved relative gap $\gamma_{i,m}(t) \le \varepsilon_{i} = 10^{-5}$ for instance~$i$ at some time $t \le T_{i}$, where $T_{i} = 3600$ seconds for all instances except Goldcoast, Chicago Regional, and Philadelphia, which used $T_{i} = 21600$ seconds, then the table shows the solve time $t_{i,m}$ (in seconds), otherwise the table shows the final relative gap $\gamma_{i,m}(T_{i})$ achieved by method~$m$ in parentheses.
Dial saved no feasible solution within the $T_{i} = 21600$ second cutoff for Chicago Regional and Philadelphia, because it did not complete its initial minimum-path assignment.
A summary of these results is given in Table~\ref{tab:exp_metrics_summary} and Figure~\ref{fig:exp_common_gap_profile}.

\begin{table}[h]
\centering
\small
\setlength{\tabcolsep}{4pt}
\caption{Computational performance on smooth instances~$i$.
The table shows either the solve time $t_{i,m}$ (in seconds) for instance~$i$ if method~$m$ achieved relative gap $\gamma_{i,m}(t) \le \varepsilon_{i} = 10^{-5}$ at some time $t \le T_{i}$ ($T_{i} = 6$h for Goldcoast, Chicago Regional, and Philadelphia; $T_{i} = 1$h otherwise), or otherwise the final relative gap $\gamma_{i,m}(T_{i})$ achieved by method~$m$ in parentheses, or ``---'' meaning that the method saved no feasible solution within $T_{i}$.}
\label{tab:exp_smooth_runtimes}
\begin{tabular}{lrrrrr}
\toprule
Network & Daganzo & Dial & Marcotte--Zhu & SBA & BTA \\
\midrule
Braess & $(2.38\times10^{-1})$ & $0.00678$ & $0.00784$ & $0.00214$ & $<0.001$ \\
Sioux Falls & $0.473$ & $0.131$ & $0.0597$ & $0.122$ & $0.032$ \\
Friedrichshain & $(2.09\times10^{-2})$ & $2.86$ & $0.457$ & $0.259$ & $0.045$ \\
Berlin Prenzlauerberg & $(1.79\times10^{-1})$ & $13.7$ & $0.586$ & $2.78$ & $0.168$ \\
Berlin Tiergarten & $(4.41\times10^{-1})$ & $2.13$ & $0.0864$ & $0.143$ & $0.040$ \\
Berlin Mitte & $(2.00\times10^{-2})$ & $4.01$ & $0.214$ & $0.470$ & $0.070$ \\
Anaheim & $(2.61\times10^{-1})$ & $26.9$ & $2.35$ & $1.16$ & $0.321$ \\
Barcelona & $(24.1)$ & $1750$ & $289$ & $38.3$ & $5.987$ \\
Winnipeg-Asym & $(6.22\times10^{-2})$ & $611$ & $161$ & $113$ & $1.275$ \\
Berlin M.-P.-F. & $(3.72\times10^{-2})$ & $112$ & $1.91$ & $6.47$ & $1.038$ \\
Winnipeg & $(33.5)$ & $2280$ & $823$ & $88.7$ & $4.673$ \\
Terrassa-Asym & $(3.59\times10^{-1})$ & $1030$ & $593$ & $496$ & $9.935$ \\
Hessen-Asym & $(3.59\times10^{-2})$ & $(1.94\times10^{-5})$ & $820$ & $2130$ & $14.657$ \\
Goldcoast & $(112)$ & $(4.71)$ & $(1.37\times10^{-2})$ & $(1.01\times10^{-2})$ & $(1.09\times10^{-2})$ \\
Chicago Regional & $(3.28\times10^{3})$ & --- & $(8.67\times10^{-1})$ & $(1.85\times10^{4})$ & $(2.03\times10^{-1})$ \\
Berlin Center & $(3.44\times10^{4})$ & $(1.29\times10^{-2})$ & $(1.00)$ & $(2.60\times10^{-3})$ & $25.352$ \\
Philadelphia & $(2.21)$ & --- & $(1.00)$ & $(24.7)$ & $(5.06\times10^{-4})$ \\
\bottomrule
\end{tabular}

\end{table}


For each smooth instance~$i$, Figures~\ref{fig:exp_common_objective_i}--\ref{fig:exp_common_objective_iii} show the relative gaps $\gamma_{i,m}(t)$ of each of the five methods~$m$ as a function of computation time~$t$, with both axes using a logarithmic scale.

\begin{figure}[h]
\centering

\caption{For each smooth instances~$i$, the curves show the relative gaps $\gamma_{i,m}(t)$ of the five methods~$m$ as a function of computation time~$t$.
The dashed line shows the tolerance $\varepsilon_{i} = 10^{-5}$.
Each curve stops at its first recorded point with $\gamma_{i,m}(t) \le \varepsilon_{i}$, or at $T_{i}$.}
\label{fig:exp_common_objective_i}
\end{figure}

\begin{figure}[p]
\centering

\caption{For each smooth instances~$i$, the curves show the relative gaps $\gamma_{i,m}(t)$ of the five methods~$m$ as a function of computation time~$t$.
The dashed line shows the tolerance $\varepsilon_{i} = 10^{-5}$.
Each curve stops at its first recorded point with $\gamma_{i,m}(t) \le \varepsilon_{i}$, or at $T_{i}$.}
\label{fig:exp_common_objective_ii}
\end{figure}

\begin{figure}[p]
\centering

\caption{For each smooth instances~$i$, the curves show the relative gaps $\gamma_{i,m}(t)$ of the five methods~$m$ as a function of computation time~$t$.
The dashed line shows the tolerance $\varepsilon_{i} = 10^{-5}$.
Each curve stops at its first recorded point with $\gamma_{i,m}(t) \le \varepsilon_{i}$, or at $T_{i}$.}
\label{fig:exp_common_objective_iii}
\end{figure}

\FloatBarrier

Table~\ref{tab:exp_dafermos_runtimes} compares the computational performance of Algorithm~BTA and the method of \citet{dafe:72} on nonsmooth instances~$i$.
If method~$m$ achieved relative gap $\gamma_{i,m}(t) \le \varepsilon_{i} = 10^{-3}$ for instance~$i$ at some time $t \le T_{i}$, where $T_{i} = 3600$ seconds for all instances except Goldcoast, Chicago Regional, and Philadelphia, which used $T_{i} = 86400$ seconds, then the table shows the solve time $t_{i,m}$ (in seconds), otherwise the table shows the final relative gap $\gamma_{i,m}(T_{i})$ achieved by method~$m$ in parentheses.
A summary of these results is given in Table~\ref{tab:exp_metrics_nonsmooth}.

\begin{table}[h]
\centering
\small
\caption{Computational performance on nonsmooth instances~$i$.
The table shows either the solve time $t_{i,m}$ (in seconds) for instance~$i$ if method~$m$ achieved relative gap $\gamma_{i,m}(t) \le \varepsilon_{i} = 10^{-3}$ at some time $t \le T_{i}$ ($T_{i} = 24$h for Goldcoast, Chicago Regional, and Philadelphia; $T_{i} = 1$h otherwise), or otherwise the final relative gap $\gamma_{i,m}(T_{i})$ achieved by method~$m$ in parentheses.
}
\label{tab:exp_dafermos_runtimes}
\begin{tabular}{lrr}
\toprule
Network & Dafermos & BTA \\
\midrule
Braess & $0.0545$ & $0.016$ \\
Sioux Falls & $0.031$ & $0.191$ \\
Friedrichshain & $0.060$ & $15.671$ \\
Berlin Prenzlauerberg & $0.315$ & $39.797$ \\
Berlin Tiergarten & $0.130$ & $7.915$ \\
Berlin Mitte & $0.220$ & $33.597$ \\
Anaheim & $0.670$ & $140.964$ \\
Barcelona & $26.931$ & $2093.585$ \\
Winnipeg-Asym & $26.739$ & $42.728$ \\
Berlin M.-P.-F. & $2.500$ & $275.643$ \\
Winnipeg & $21.447$ & $84.343$ \\
Terrassa-Asym & $25.665$ & $280.418$ \\
Hessen-Asym & $493.420$ & $314.519$ \\
Goldcoast & $3216.493$ & $(4.75\times10^{-3})$ \\
Chicago Regional & $12385.983$ & $20352.302$ \\
Berlin Center & $420.221$ & $2014.743$ \\
Philadelphia & $36485.144$ & $(1.62\times10^{-2})$ \\
\bottomrule
\end{tabular}
\end{table}

\begin{figure}[p]
\centering
\begin{tikzpicture}
\begin{groupplot}[
  objectiveaxis,
  group style={group size=2 by 3,horizontal sep=1.5cm,vertical sep=1.35cm}
]
\nextgroupplot[
  title={Sioux Falls},
  xmin=0.0001,
  xmax=3600,
  xtick={0.0001,0.01,1,100,3600},
  xticklabels={$10^{-4}$,$10^{-2}$,$1$,$10^{2}$,$3600$},
  ymin=0.0001,
  ymax=10,
  ytick={0.0001,0.01,1},
  yticklabels={$10^{-4}$,$10^{-2}$,$1$},
  legend columns=2,
  legend style={draw=none,font=\scriptsize,at={(1.10,1.30)},anchor=south,/tikz/every even column/.append style={column sep=7pt}}
]
\addlegendimage{daf}\addlegendentry{Dafermos}
\addlegendimage{bta}\addlegendentry{BTA}
\addplot+[daf,const plot,forget plot] coordinates {
  (0.003981933,1.78309811497) (0.009009273,1.44642603186) (0.013552881,0.0227121288805) (0.020348983,0.00138658551805)
  (0.030592107,0.000178472870335)
};
\addplot+[bta,const plot,forget plot] coordinates {
  (0.000420285039581,1) (0.00562017300399,1) (0.00911156705115,0.392525530611) (0.00948574801441,0.273569808106)
  (0.00980769802118,0.173180196032) (0.0100676290458,0.0939641549306) (0.0101574980072,0.0539012183629) (0.01023746602,0.047083867761)
  (0.0103639030131,0.0440700642075) (0.0108051260468,0.0349676684962) (0.0188891610014,0.0309626286484) (0.0253263600171,0.03008486197)
  (0.0351176430122,0.0298504204129) (0.0446946480079,0.0293398113985) (0.0598535830504,0.0291408336263) (0.0799625010113,0.0287819753532)
  (0.106656870048,0.0286076713228) (0.149037708004,0.0284478240675) (0.175615087035,0.0219354709812) (0.175649432,0.0165593522071)
  (0.175875732035,0.0130599768125) (0.176374504052,0.00808645846644) (0.176611950039,0.00631403711597) (0.178138411022,0.00489103130244)
  (0.178755073051,0.00341060357672) (0.180609395029,0.00242340543028) (0.182150760025,0.00178728934322) (0.187340892036,0.0013431835196)
  (0.190932997037,0.000928404864481)
};
\addplot[gray!70,densely dashed,line width=0.45pt,forget plot]
coordinates {(1e-9,0.001) (3600,0.001)};
\node[anchor=south west,font=\tiny,text=gray!70]
at ({rel axis cs:0,0}|-{axis cs:1,0.001}) {$10^{-3}$};

\nextgroupplot[
  title={Anaheim},
  xmin=0.0001,
  xmax=3600,
  xtick={0.0001,0.01,1,100,3600},
  xticklabels={$10^{-4}$,$10^{-2}$,$1$,$10^{2}$,$3600$},
  ymin=0.0001,
  ymax=10,
  ytick={0.0001,0.01,1},
  yticklabels={$10^{-4}$,$10^{-2}$,$1$}
]
\addplot+[daf,const plot,forget plot] coordinates {
  (0.044482224,0.299483501226) (0.085914217,0.26673299544) (0.129492838,0.0630934732859) (0.195233214,0.0489349413824)
  (0.294978103,0.0123105674189) (0.4459272,0.00121647814601) (0.67047849,0.000264381323336)
};
\addplot+[bta,const plot,forget plot] coordinates {
  (0.000861518026795,1) (0.00862648698967,1) (0.0143244320061,0.557837599143) (0.0149754000013,0.369187642818)
  (0.0154538570205,0.346219443988) (0.015901923005,0.258228605635) (0.0173324039788,0.156617323166) (0.0201028350275,0.126666051767)
  (0.0267668470042,0.102290026589) (0.0272450459888,0.0996977618015) (0.0435102880001,0.0711159098042) (0.047608740977,0.0656767276738)
  (0.0817039040267,0.0513226320847) (0.0850607170141,0.0495328809129) (0.155648102984,0.0451761893201) (0.265439092007,0.0411475164173)
  (0.51199408801,0.0379131087826) (0.623216268024,0.0369685886169) (0.848372810986,0.0364943090883) (1.49215978401,0.0351560406089)
  (2.64503950998,0.0344030949109) (4.85982738802,0.0337422974512) (8.28608878999,0.0334483383609) (14.767689183,0.0332224624452)
  (27.133868411,0.0329459943887) (27.21320509,0.0245149074393) (27.220988879,0.0189206940571) (27.275471097,0.0134028368507)
  (27.372029353,0.00997711027042) (27.652356621,0.00705399633199) (28.176540579,0.00510094229344) (29.555225046,0.00370210154425)
  (32.150117404,0.00267445256194) (43.378086468,0.00192805080494) (45.087735383,0.00178441589205) (79.592261492,0.00140240668028)
  (84.693451102,0.00138379533493) (140.3077209,0.00100533160088) (140.963924732,0.00099875073188)
};
\addplot[gray!70,densely dashed,line width=0.45pt,forget plot]
coordinates {(1e-9,0.001) (3600,0.001)};
\node[anchor=south west,font=\tiny,text=gray!70]
at ({rel axis cs:0,0}|-{axis cs:1,0.001}) {$10^{-3}$};

\nextgroupplot[
  title={Barcelona},
  xmin=0.001,
  xmax=3600,
  xtick={0.01,1,100,3600},
  xticklabels={$10^{-2}$,$1$,$10^{2}$,$3600$},
  ymin=0.0001,
  ymax=100,
  ytick={0.0001,0.01,1,100},
  yticklabels={$10^{-4}$,$10^{-2}$,$1$,$10^{2}$}
]
\addplot+[daf,const plot,forget plot] coordinates {
  (0.0708763,9.80150611375) (0.180237665,9.58098962546) (0.280424193,9.50131138845) (0.443607787,9.42820266631)
  (0.681119911,9.410167389) (1.02411219,7.33498319024) (1.557255841,5.99122120794) (2.348690143,3.86960386908)
  (3.529077374,0.125305958361) (5.30786978,0.0497925006241) (7.963462072,0.0178579731607) (11.960431123,0.00809140392022)
  (17.944208943,0.00220208243691) (26.930844268,0.000481011410145)
};
\addplot+[bta,const plot,forget plot] coordinates {
  (0.0056397630251,1) (0.0539088270161,1) (0.072934122989,0.702106839388) (0.0823945840239,0.491546712566)
  (0.0971804970177,0.366578496683) (0.120681481028,0.329698673031) (0.175502652011,0.264728105468) (0.225933348993,0.247506878317)
  (0.419194052985,0.192395531213) (0.429755938996,0.190770320802) (0.752652076015,0.165901137843) (1.37606641103,0.154024705653)
  (2.54083946801,0.147891976056) (4.78499189799,0.142128716376) (8.67773600802,0.139118055328) (11.830992006,0.137897477156)
  (15.950113919,0.136850083997) (29.579053859,0.135450329171) (54.18900139,0.134480974762) (98.834264941,0.133741045979)
  (182.242676988,0.133302022671) (335.778220244,0.133001944038) (629.139318406,0.132815394107) (936.062780678,0.0908874482468)
  (936.1202236,0.0699578774488) (936.275850886,0.0499461123573) (936.796112702,0.0300020029645) (937.2941719,0.0241724090147)
  (938.278938424,0.0187850824285) (940.544223154,0.0136159056657) (951.33286652,0.00978869776082) (996.293005054,0.00705966253624)
  (1142.16525208,0.00599134813897) (1943.74905137,0.0049545822083) (1946.46728513,0.00363992551936) (1953.6222407,0.00260918952627)
  (1969.48969189,0.00187541129116) (2014.75056817,0.00134904760379) (2091.94413332,0.00101619537893) (2093.58535858,0.000974324781893)
};
\addplot[gray!70,densely dashed,line width=0.45pt,forget plot]
coordinates {(1e-9,0.001) (3600,0.001)};
\node[anchor=south west,font=\tiny,text=gray!70]
at ({rel axis cs:0,0}|-{axis cs:1,0.001}) {$10^{-3}$};

\nextgroupplot[
  title={Winnipeg},
  xmin=0.001,
  xmax=3600,
  xtick={0.01,1,100,3600},
  xticklabels={$10^{-2}$,$1$,$10^{2}$,$3600$},
  ymin=0.0001,
  ymax=100,
  ytick={0.0001,0.01,1,100},
  yticklabels={$10^{-4}$,$10^{-2}$,$1$,$10^{2}$}
]
\addplot+[daf,const plot,forget plot] coordinates {
  (0.088907555,58.4700942825) (0.141823172,58.4691689367) (0.2218696,50.4226134335) (0.347248853,19.9194966567)
  (0.533846269,9.07815291867) (0.820173081,4.24857062908) (1.239981616,0.0723730581607) (1.874172139,0.0435540479366)
  (2.813487346,0.0259418584533) (4.221221591,0.0162266628287) (6.339119226,0.00915772628957) (9.530970586,0.00433437633119)
  (14.297755137,0.00127151660117) (21.447397372,0.000471297222858)
};
\addplot+[bta,const plot,forget plot] coordinates {
  (0.00513286096975,1) (0.0405876479927,1) (0.0520967370248,0.339911909317) (0.0583586490247,0.203433507471)
  (0.0610821670271,0.152314505679) (0.0645274629933,0.134467716832) (0.0888531350065,0.110120438595) (0.106831546,0.098732734366)
  (0.137159297999,0.0895426757052) (0.218804650009,0.07648881778) (0.273571307014,0.0718933472789) (0.33980041102,0.068562941295)
  (0.529978529026,0.0636053012481) (0.840022897988,0.0606372199164) (1.340165539,0.0587905330903) (2.12840858201,0.0579614229846)
  (3.32520959101,0.0574671605931) (5.35273172997,0.0571552877972) (8.98977005098,0.0569009606409) (13.506190374,0.056712584754)
  (21.319088087,0.0566334682916) (21.384333957,0.0475730918285) (21.436495572,0.0359514981138) (21.469705472,0.0265354898391)
  (21.590583709,0.0185589114076) (21.789386244,0.0136299962419) (22.054948684,0.00994395055634) (22.604347333,0.00718315701938)
  (23.836043393,0.00516602822706) (28.204468041,0.00371511718909) (33.990251553,0.0032353714485) (54.029861844,0.00295405293631)
  (70.417956242,0.00265708816416) (72.043459371,0.00191241819381) (76.315673989,0.00138883504985) (84.291824392,0.00100180855496)
  (84.343232506,0.000999682864331)
};
\addplot[gray!70,densely dashed,line width=0.45pt,forget plot]
coordinates {(1e-9,0.001) (3600,0.001)};
\node[anchor=south west,font=\tiny,text=gray!70]
at ({rel axis cs:0,0}|-{axis cs:1,0.001}) {$10^{-3}$};

\nextgroupplot[
  title={Chicago Regional},
  xmin=10,
  xmax=100000,
  xtick={10,100,1000,86400},
  xticklabels={$10$,$10^{2}$,$10^{3}$,$86400$},
  ymin=0.0001,
  ymax=100,
  ytick={0.0001,0.01,1,100},
  yticklabels={$10^{-4}$,$10^{-2}$,$1$,$10^{2}$}
]
\addplot+[daf,const plot,forget plot] coordinates {
  (34.235553897,9.85931197295) (133.992668317,9.70370544919) (205.676740829,9.2601563878) (313.562022454,8.98637097309)
  (475.657074619,8.02142145664) (718.226429075,7.40917284215) (1080.922696,6.84217673871) (1623.77386319,3.7885233493)
  (2439.57188352,3.19115926759) (3664.46914153,1.93736705086) (5503.7580709,1.23787904623) (8256.98049815,0.171010343977)
  (12385.9826711,0.00086938089531)
};
\addplot+[bta,const plot,forget plot] coordinates {
  (11.544162257,1) (81.625693213,1) (109.172388957,0.00467040235767) (159.40886264,0.00467040235767)
  (229.769942717,0.0042711280567) (244.38793063,0.00299553897551) (287.838233585,0.00292006734578) (310.069442223,0.0026554627595)
  (406.837440089,0.00204786448136) (511.588704314,0.00184444403038) (586.099036643,0.00181461841331) (830.83740451,0.00171173010176)
  (1215.58784796,0.00162129734133) (1692.32750202,0.00154394280778) (2438.30975692,0.00143112518564) (2925.39098757,0.00138287779334)
  (3644.21611328,0.00130917064835) (5031.32543987,0.00122935590798) (7378.35787804,0.00114654165364) (10081.2698836,0.00108913988399)
  (14322.8023817,0.00104241788762) (20232.4865217,0.00100011837939) (20352.3017289,0.000999403395175)
};
\addplot[gray!70,densely dashed,line width=0.45pt,forget plot]
coordinates {(1,0.001) (86400,0.001)};
\node[anchor=south west,font=\tiny,text=gray!70]
at ({rel axis cs:0,0}|-{axis cs:1,0.001}) {$10^{-3}$};

\nextgroupplot[hide axis]
\end{groupplot}
\node[font=\small] at ($(group c1r3.south)!0.5!(group c2r3.south)+(0,-0.62cm)$)
  {elapsed time $t$ (seconds)};
\node[font=\small,rotate=90] at ($(group c1r1.west)!0.5!(group c1r3.west)+(-0.95cm,0)$)
  {relative gap $\gamma_{i,m}(t)$};
\end{tikzpicture}
\caption{For nonsmooth instances~$i$, the curves show the relative gaps $\gamma_{i,m}(t)$ of Dafermos and BTA as a function of computation time~$t$.
The dashed line shows the tolerance $\varepsilon_{i} = 10^{-3}$.
Each curve stops at its first recorded point with $\gamma_{i,m}(t) \le \varepsilon_{i}$, or at $T_{i}$.
}
\label{fig:exp_nonsmooth_common_objective_i}
\end{figure}
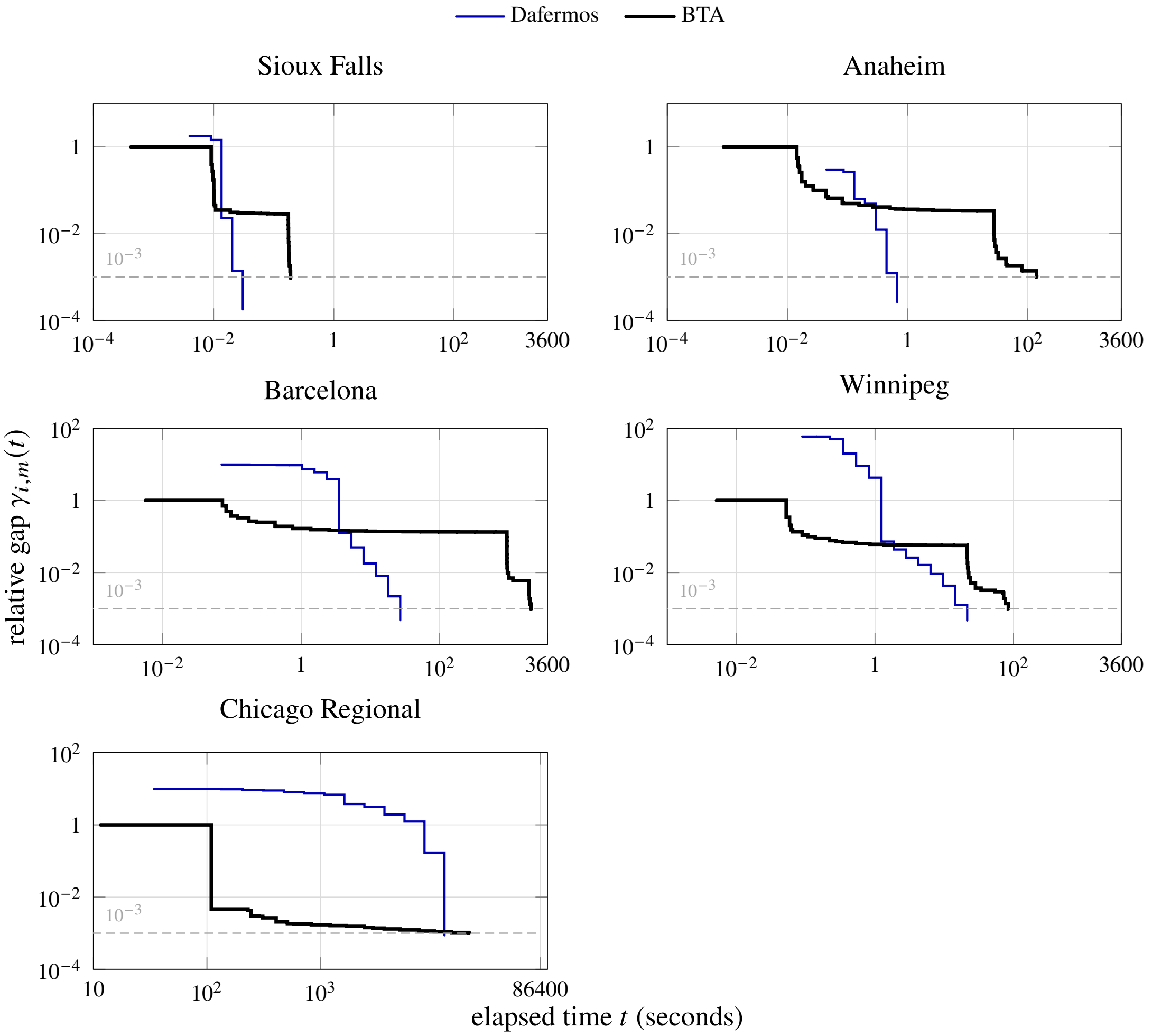

\begin{figure}[p]
\centering
\begin{tikzpicture}
\begin{groupplot}[
  objectiveaxis,
  group style={group size=2 by 3,horizontal sep=1.5cm,vertical sep=1.35cm}
]
\nextgroupplot[
  title={Braess},
  xmin=1e-05,
  xmax=3600,
  xtick={0.0001,0.01,1,100,3600},
  xticklabels={$10^{-4}$,$10^{-2}$,$1$,$10^{2}$,$3600$},
  ymin=1e-10,
  ymax=10,
  ytick={1e-09,1e-06,0.001,1},
  yticklabels={$10^{-9}$,$10^{-6}$,$10^{-3}$,$1$},
  legend columns=2,
  legend style={draw=none,font=\scriptsize,at={(1.10,1.30)},anchor=south,/tikz/every even column/.append style={column sep=7pt}}
]
\addlegendimage{daf}\addlegendentry{Dafermos}
\addlegendimage{bta}\addlegendentry{BTA}
\addplot+[daf,const plot,forget plot] coordinates {
  (0.030239778,0.327568776787) (0.054373708,0.0552965362167) (0.054469735,7.31537083831e-10)
};
\addplot+[bta,const plot,forget plot] coordinates {
  (1.4903023839e-05,1) (0.00521865501651,1) (0.00891023900476,0.417231997011) (0.00891617502202,0.364551294038)
  (0.00892310601193,0.359434336536) (0.00892966199899,0.33107828518) (0.00894039901323,0.328683865738) (0.0089563660149,0.328617265846)
  (0.0089626040135,0.327582758454) (0.00898988201516,0.327582758454) (0.0126770550269,0.327582758454) (0.0163797400019,0.174506139087)
  (0.0163849740056,0.0788784975467) (0.0163895470032,0.0380009894378) (0.0163937180187,0.0191379912875) (0.0163982270169,0.00901097445328)
  (0.0164021470118,0.00300459496388) (0.0164113490027,0.00176271287326) (0.0164199340215,0.00163029879464) (0.0164266840147,0.000243339224615)
};
\addplot[gray!70,densely dashed,line width=0.45pt,forget plot]
coordinates {(1e-06,0.001) (3600,0.001)};
\node[anchor=south west,font=\tiny,text=gray!70]
at ({rel axis cs:0,0}|-{axis cs:1,0.001}) {$10^{-3}$};

\nextgroupplot[
  title={Goldcoast},
  xmin=1,
  xmax=100000,
  xtick={1,10,100,1000,86400},
  xticklabels={$1$,$10$,$10^{2}$,$10^{3}$,$86400$},
  ymin=0.0001,
  ymax=100,
  ytick={0.0001,0.01,1,100},
  yticklabels={$10^{-4}$,$10^{-2}$,$1$,$10^{2}$}
]
\addplot+[daf,const plot,forget plot] coordinates {
  (5.538462757,12.9478710403) (23.835613206,12.8754636614) (36.645182178,2.64821964023) (55.016888762,2.31863355746)
  (82.840245018,2.23484388917) (124.549171614,1.3490585849) (187.24442426,1.06718991656) (281.653895983,0.874386144881)
  (423.133211946,0.7612302221) (635.049061435,0.00871648816946) (952.638747034,0.00518236382645) (1428.98623038,0.00253525760451)
  (2144.10427681,0.00138086773062) (3216.49302542,0.000570485493424)
};
\addplot+[bta,const plot,forget plot] coordinates {
  (1.35576979,1) (10.885752461,1) (14.32517609,0.0244437287308) (20.9754894,0.0244437287308)
  (32.940589085,0.0244142861867) (34.758308491,0.0194186717657) (42.062782123,0.0168247145005) (54.293880347,0.0150563399315)
  (98.080083219,0.0131122966083) (160.914224225,0.0125337774769) (271.535197935,0.0118721152969) (465.343977804,0.0108471249734)
  (754.484639101,0.0103323596745) (788.12077486,0.0101099872685) (1280.22838739,0.00941065036331) (2163.40788166,0.00875110640987)
  (3760.45893409,0.00795721558621) (3788.35298485,0.00788710367549) (6368.40389444,0.00719768730912) (10500.1190076,0.00674116283627)
  (17786.4947706,0.00620095966429) (19180.6481885,0.00612726093574) (30119.1491119,0.00570787830339) (51293.8124494,0.00512044196283)
  (86159.2703211,0.00475179879145) (86346.4119436,0.00475016030378) (86400,0.00475016030378)
};
\addplot[gray!70,densely dashed,line width=0.45pt,forget plot]
coordinates {(1,0.001) (86400,0.001)};
\node[anchor=south west,font=\tiny,text=gray!70]
at ({rel axis cs:0,0}|-{axis cs:1,0.001}) {$10^{-3}$};

\nextgroupplot[
  title={Hessen-Asym},
  xmin=0.01,
  xmax=3600,
  xtick={0.01,1,100,3600},
  xticklabels={$10^{-2}$,$1$,$10^{2}$,$3600$},
  ymin=0.0001,
  ymax=10,
  ytick={0.0001,0.01,1},
  yticklabels={$10^{-4}$,$10^{-2}$,$1$}
]
\addplot+[daf,const plot,forget plot] coordinates {
  (0.545212531,0.579441771844) (0.990106325,0.569007936925) (1.543381309,0.548206883308) (2.402385364,0.509312783042)
  (3.656643916,0.462951936343) (5.575770175,0.394746815905) (8.380412303,0.262136966227) (12.722055203,0.0932451013997)
  (19.138055512,0.0591522804189) (28.813014968,0.0218995729993) (43.231316366,0.0120490578902) (64.89196094,0.00743510527343)
  (97.372490968,0.0051710309496) (146.107546901,0.00339759800897) (219.171741145,0.00215643695178) (328.860478024,0.00129402150912)
  (493.419513753,0.000692358507683)
};
\addplot+[bta,const plot,forget plot] coordinates {
  (0.0171917550033,1) (0.151254378026,1) (0.173749786045,0.289033141359) (0.232116669999,0.285357183789)
  (0.330595510022,0.2402528188) (0.359056767018,0.186643923077) (0.476610321028,0.165828467937) (0.733512743027,0.142979170711)
  (0.990995338012,0.138558593858) (1.15559723304,0.137418642437) (1.84570615104,0.135593206367) (2.93849977804,0.134708012751)
  (4.70536273904,0.13421077458) (7.49256975803,0.13389420577) (11.948830237,0.133680070146) (19.116560857,0.133531030769)
  (30.414886818,0.13342860363) (31.184966377,0.0887295781011) (31.222060463,0.0647072295997) (31.284457374,0.0511174583886)
  (31.366423881,0.033987748215) (31.569573989,0.026440928385) (32.364471196,0.0188253686314) (40.2119144271,0.0138358783414)
  (48.584205666,0.0132349223929) (77.482534108,0.0127244791774) (119.950484773,0.00818507811623) (120.021758442,0.00677717567659)
  (120.441721255,0.00505663299844) (122.587002995,0.00370271194933) (123.718144053,0.00348297595223) (135.861787288,0.0026638688375)
  (197.327793931,0.00217657254759) (311.546916767,0.0016453830422) (313.184429717,0.00137825141374) (314.510039722,0.00100184361934)
  (314.519349688,0.000994062369484)
};
\addplot[gray!70,densely dashed,line width=0.45pt,forget plot]
coordinates {(1e-9,0.001) (3600,0.001)};
\node[anchor=south west,font=\tiny,text=gray!70]
at ({rel axis cs:0,0}|-{axis cs:1,0.001}) {$10^{-3}$};

\nextgroupplot[
  title={Philadelphia},
  xmin=1,
  xmax=100000,
  xtick={1,10,100,1000,86400},
  xticklabels={$1$,$10$,$10^{2}$,$10^{3}$,$86400$},
  ymin=0.0001,
  ymax=100,
  ytick={0.0001,0.01,1,100},
  yticklabels={$10^{-4}$,$10^{-2}$,$1$,$10^{2}$}
]
\addplot+[daf,const plot,forget plot] coordinates {
  (24.134862436,26.6263895279) (77.564251523,26.2637425711) (117.775195925,25.1218084492) (182.934607728,22.4407899784)
  (277.937007793,21.4386959985) (418.164080459,20.3622778812) (627.41934429,18.9236824129) (943.648626542,16.8101163563)
  (1417.50429125,13.5913660566) (2128.51587262,10.1063573278) (3196.14773185,7.29085593544) (4799.68212403,0.081732240814)
  (7201.96359935,0.0165121886786) (10807.1797715,0.00606114119872) (16212.9585444,0.00307245624432) (24321.9687841,0.00132527645746)
  (36485.1439161,0.000499741681468)
};
\addplot+[bta,const plot,forget plot] coordinates {
  (7.147242739,1) (47.360150229,1) (62.353107924,0.124693671727) (92.842709341,0.124693671727)
  (119.716080042,0.122781853952) (151.832452712,0.10135660671) (167.095310263,0.0800905261353) (176.534752573,0.0462943087952)
  (201.951199075,0.0419717244357) (248.835824287,0.0419717244357) (372.716656794,0.03507249352) (393.98708898,0.0337709816294)
  (591.139164718,0.0290093033104) (607.641412893,0.0290093033104) (973.856188266,0.0258555915038) (1336.62834255,0.0239113900714)
  (1489.10296333,0.0234344694247) (2262.93565982,0.0216502737669) (3519.26051591,0.0204306766111) (4747.68960642,0.0197145114301)
  (5474.80377024,0.0194037850527) (8521.44317992,0.0185380173506) (13266.6701978,0.0176928567748) (20696.2039632,0.0171310196761)
  (32196.4920462,0.0167113818359) (50094.596033,0.016408915347) (77892.100143,0.0162017906542) (77954.3803868,0.0162008093521)
  (86400,0.0162008093521)
};
\addplot[gray!70,densely dashed,line width=0.45pt,forget plot]
coordinates {(1,0.001) (86400,0.001)};
\node[anchor=south west,font=\tiny,text=gray!70]
at ({rel axis cs:0,0}|-{axis cs:1,0.001}) {$10^{-3}$};

\nextgroupplot[
  title={Terrassa-Asym},
  xmin=0.001,
  xmax=3600,
  xtick={0.01,1,100,3600},
  xticklabels={$10^{-2}$,$1$,$10^{2}$,$3600$},
  ymin=0.0001,
  ymax=10,
  ytick={0.0001,0.01,1},
  yticklabels={$10^{-4}$,$10^{-2}$,$1$}
]
\addplot+[daf,const plot,forget plot] coordinates {
  (0.053292499,2.07420707514) (0.107565041,1.9190185857) (0.174298162,1.59534199447) (0.277245408,1.25416013582)
  (0.433860015,0.933533768839) (0.652507322,0.532773581961) (0.983505348,0.18062794173) (1.478731092,0.0877005031329)
  (2.233809404,0.0405237062506) (3.371032946,0.0220445409186) (5.056564588,0.0114173935417) (7.597569427,0.00635469618432)
  (11.396885606,0.00312966765033) (17.102340962,0.00159634060043) (25.664853,0.000677916989915)
};
\addplot+[bta,const plot,forget plot] coordinates {
  (0.00365657103248,1) (0.0264993280289,1) (0.0414736199891,0.915644616969) (0.0502635960001,0.691023188087)
  (0.0607349009952,0.65861103029) (0.0669336040155,0.487924551625) (0.082688509021,0.369838085722) (0.0954376300215,0.356925280493)
  (0.155375534028,0.324257057751) (0.266250725021,0.311938757402) (0.453631815035,0.307805765726) (0.774853651004,0.305154813267)
  (1.32499326504,0.303090149943) (2.26131210104,0.301599429935) (3.86327698402,0.300660896058) (6.60149465903,0.300076874521)
  (11.281457013,0.299708496733) (12.876474219,0.263403675096) (12.890737245,0.178700097024) (12.893100898,0.127658806156)
  (12.900477355,0.09682524057) (12.942982274,0.0717887849688) (14.187562015,0.0516635808601) (19.28147915,0.0491327101375)
  (31.166292283,0.0427377194542) (31.173028567,0.0288476652867) (31.20512108,0.0256690101266) (31.479066064,0.0192420368952)
  (32.914330168,0.0147649154106) (34.417911394,0.0138444888147) (56.24787441,0.0123306680493) (78.083758457,0.00993827757897)
  (78.870216183,0.0071520246139) (84.998794237,0.00459719892964) (97.451774341,0.00413893704274) (136.721936722,0.00356296117533)
  (137.979727064,0.00265941120949) (144.211950722,0.00192047024555) (164.150342593,0.00142025150722) (168.241016968,0.00138130709924)
  (280.415595912,0.00101247231493) (280.417824856,0.000994191913807)
};
\addplot[gray!70,densely dashed,line width=0.45pt,forget plot]
coordinates {(1e-9,0.001) (3600,0.001)};
\node[anchor=south west,font=\tiny,text=gray!70]
at ({rel axis cs:0,0}|-{axis cs:1,0.001}) {$10^{-3}$};

\nextgroupplot[
  title={Winnipeg-Asym},
  xmin=0.001,
  xmax=3600,
  xtick={0.01,1,100,3600},
  xticklabels={$10^{-2}$,$1$,$10^{2}$,$3600$},
  ymin=0.0001,
  ymax=10,
  ytick={0.0001,0.01,1},
  yticklabels={$10^{-4}$,$10^{-2}$,$1$}
]
\addplot+[daf,const plot,forget plot] coordinates {
  (0.078209694,2.22978377755) (0.12564086,2.20621515852) (0.194839888,1.98920208365) (0.297743999,1.54577518513)
  (0.455954872,1.33464874474) (0.685237333,1.01547269789) (1.03493379,0.565199126927) (1.552758045,0.240682225143)
  (2.33637206,0.0986172877239) (3.513051103,0.0506251802297) (5.272917513,0.0258391342755) (7.911411563,0.0117277128877)
  (11.878668121,0.00531383952691) (17.826065205,0.00220874001689) (26.739451533,0.00088864315129)
};
\addplot+[bta,const plot,forget plot] coordinates {
  (0.0069171069772,1) (0.0289210499614,1) (0.0388854089542,0.841942236049) (0.0427882449585,0.555824373577)
  (0.0500106809777,0.41551339492) (0.0615947399638,0.385423370772) (0.0649947079946,0.371697518148) (0.0849166839616,0.328323501848)
  (0.127113537979,0.312470401512) (0.193563827954,0.306379047131) (0.291668603953,0.303498388528) (0.44733186299,0.301930118984)
  (0.692025735974,0.301388167359) (1.03200146096,0.301072970638) (1.63236189599,0.30093088769) (2.132263429,0.283240669154)
  (2.13671180495,0.192635061845) (2.16007671197,0.137302736836) (2.192043819,0.0970976724938) (2.30667906196,0.0719072191435)
  (2.33167554101,0.0689723195932) (2.60721660795,0.0516002288753) (3.23203161696,0.037238311372) (3.52927091398,0.034111742258)
  (5.34756049298,0.0307649291264) (8.18795086397,0.0301393568191) (12.724603753,0.029882492427) (12.841127746,0.0239689278803)
  (12.910049492,0.0192156512761) (13.276028055,0.0138337336651) (13.940799178,0.00998296230029) (15.191988563,0.00717826466437)
  (18.649565124,0.00537554380828) (19.85619637,0.0051654965949) (28.204224862,0.00478384944944) (33.044292128,0.00371140731348)
  (33.349551004,0.00267536170253) (34.405422613,0.00192660071793) (36.784155153,0.00138382363085) (42.511516956,0.00100551253994)
  (42.728042129,0.000998016283433)
};
\addplot[gray!70,densely dashed,line width=0.45pt,forget plot]
coordinates {(1e-9,0.001) (3600,0.001)};
\node[anchor=south west,font=\tiny,text=gray!70]
at ({rel axis cs:0,0}|-{axis cs:1,0.001}) {$10^{-3}$};

\end{groupplot}
\node[font=\small] at ($(group c1r3.south)!0.5!(group c2r3.south)+(0,-0.62cm)$)
  {elapsed time $t$ (seconds)};
\node[font=\small,rotate=90] at ($(group c1r1.west)!0.5!(group c1r3.west)+(-0.95cm,0)$)
  {relative gap $\gamma_{i,m}(t)$};
\end{tikzpicture}
\caption{For nonsmooth instances~$i$, the curves show the relative gaps $\gamma_{i,m}(t)$ of Dafermos and BTA as a function of computation time~$t$.
The dashed line shows the tolerance $\varepsilon_{i} = 10^{-3}$.
Each curve stops at its first recorded point with $\gamma_{i,m}(t) \le \varepsilon_{i}$, or at $T_{i}$.}
\label{fig:exp_nonsmooth_common_objective_ii}
\end{figure}

\begin{figure}[p]
\centering
\begin{tikzpicture}
\begin{groupplot}[
  objectiveaxis,
  group style={group size=2 by 3,horizontal sep=1.5cm,vertical sep=1.35cm}
]
\nextgroupplot[
  title={Berlin Center},
  xmin=0.1,
  xmax=3600,
  xtick={0.1,1,10,100,3600},
  xticklabels={$10^{-1}$,$1$,$10$,$10^{2}$,$3600$},
  ymin=0.0001,
  ymax=100,
  ytick={0.0001,0.01,1,100},
  yticklabels={$10^{-4}$,$10^{-2}$,$1$,$10^{2}$},
  legend columns=2,
  legend style={draw=none,font=\scriptsize,at={(1.10,1.30)},anchor=south,/tikz/every even column/.append style={column sep=7pt}}
]
\addlegendimage{daf}\addlegendentry{Dafermos}
\addlegendimage{bta}\addlegendentry{BTA}
\addplot+[daf,const plot,forget plot] coordinates {
  (4.313883154,43.7573204232) (6.527882024,43.6947777707) (10.112043733,43.6083212935) (16.055942878,43.4984134741)
  (24.34424455,10.7856957253) (36.667603747,10.7769069897) (55.105542757,0.0451500670953) (82.819725415,0.0289067330415)
  (124.449497059,0.0115937881818) (186.679503443,0.00541106320375) (280.091490455,0.00183829052155) (420.221425801,0.00039054304668)
};
\addplot+[bta,const plot,forget plot] coordinates {
  (0.281282195006,1) (0.585270637006,1) (0.636744987976,0.426089398655) (0.804270816967,0.403727520451)
  (0.840404638962,0.357528979624) (0.901946244005,0.255643878578) (1.00087140099,0.205663749368) (1.02892888797,0.190895704285)
  (1.25378334901,0.137911808678) (1.53447758895,0.114394640599) (1.788239107,0.0990954506572) (2.35307324299,0.0806638849578)
  (3.02431397,0.070905131735) (3.59165395895,0.0664162963739) (5.457184139,0.0556869985232) (6.87407138298,0.0515028501616)
  (8.49964864296,0.0483859022752) (12.677993366,0.044246141694) (19.315096361,0.0407745313877) (29.715295905,0.038732052965)
  (43.047268808,0.0372114829292) (44.950986572,0.0371298102496) (68.649011386,0.0366186493931) (104.562677515,0.0361183984573)
  (185.906966852,0.0357250335146) (246.9983879,0.0356059879989) (389.55234535,0.0354329013219) (572.749201178,0.0353304110647)
  (783.380179099,0.0236273780986) (783.460210184,0.0183434203302) (783.755801846,0.0138278365414) (786.583107648,0.00995894187979)
  (793.723296068,0.00715495718287) (806.26597075,0.0051473054015) (828.422027934,0.00371594090759) (865.152096507,0.00272641307478)
  (868.687682636,0.00266499474818) (954.707104079,0.00190430634045) (1156.71562577,0.00137564576968) (1328.31160286,0.00119992257153)
  (1987.78838305,0.00100462123671) (2014.74250067,0.000996413262828)
};
\addplot[gray!70,densely dashed,line width=0.45pt,forget plot]
coordinates {(1e-9,0.001) (3600,0.001)};
\node[anchor=south west,font=\tiny,text=gray!70]
at ({rel axis cs:0,0}|-{axis cs:1,0.001}) {$10^{-3}$};

\nextgroupplot[
  title={Berlin Mitte},
  xmin=0.0001,
  xmax=3600,
  xtick={0.0001,0.01,1,100,3600},
  xticklabels={$10^{-4}$,$10^{-2}$,$1$,$10^{2}$,$3600$},
  ymin=0.0001,
  ymax=10,
  ytick={0.0001,0.01,1},
  yticklabels={$10^{-4}$,$10^{-2}$,$1$}
]
\addplot+[daf,const plot,forget plot] coordinates {
  (0.010785764,0.0420397407358) (0.025975591,0.0391304634469) (0.039904432,0.0363698236961) (0.063199121,0.0224130814294)
  (0.095425032,0.00859175315877) (0.145431181,0.00495373869841) (0.219953656,0.00026934593271)
};
\addplot+[bta,const plot,forget plot] coordinates {
  (0.000859484018292,1) (0.00942617200781,1) (0.014603728021,0.748163603933) (0.015250254015,0.487780411846)
  (0.015866731992,0.287725953106) (0.0170682379976,0.244484424447) (0.0175306940218,0.177822904878) (0.0180116660194,0.158566836228)
  (0.0189342969679,0.128090859621) (0.0215109210112,0.0972132587985) (0.0293276779703,0.0779603007207) (0.0315086379997,0.0700403865022)
  (0.0473286620108,0.0511059583273) (0.0484261650126,0.0501733844461) (0.0712610129849,0.0367626652049) (0.080977560021,0.0323959346399)
  (0.116072152974,0.0264960449102) (0.135859176982,0.0242070867749) (0.191510866978,0.0189325233419) (0.22707596398,0.0176401917013)
  (0.367495133018,0.0150266724551) (0.431787804002,0.0137220247377) (0.59866186697,0.01183840342) (1.00935480301,0.0103434803347)
  (1.25038034999,0.00975153707732) (1.65933406702,0.0091623420662) (2.831053601,0.00819139344089) (4.51559913601,0.00785500684462)
  (7.74176184298,0.00743088805001) (11.504297952,0.0071052373554) (12.546298199,0.00709779120769) (21.103810679,0.00698027321899)
  (30.55788024,0.005067149613) (30.580426626,0.00361866190597) (30.689454024,0.0025860895746) (30.966243146,0.00172225455127)
  (31.903900069,0.00134505697607) (33.463418168,0.00102195292858) (33.5965078,0.000982377563464)
};
\addplot[gray!70,densely dashed,line width=0.45pt,forget plot]
coordinates {(1e-9,0.001) (3600,0.001)};
\node[anchor=south west,font=\tiny,text=gray!70]
at ({rel axis cs:0,0}|-{axis cs:1,0.001}) {$10^{-3}$};

\nextgroupplot[
  title={Berlin M.--P.--F.},
  xmin=0.001,
  xmax=3600,
  xtick={0.01,1,100,3600},
  xticklabels={$10^{-2}$,$1$,$10^{2}$,$3600$},
  ymin=0.0001,
  ymax=10,
  ytick={0.0001,0.01,1},
  yticklabels={$10^{-4}$,$10^{-2}$,$1$}
]
\addplot+[daf,const plot,forget plot] coordinates {
  (0.076388877,0.0109311321213) (0.20833075,0.0107183977674) (0.312569572,0.00952644340146) (0.474568688,0.00934926125071)
  (0.728454565,0.00827416641261) (1.100853907,0.00586017506004) (1.662179131,0.00285545966073) (2.50016376,0.000948971134041)
};
\addplot+[bta,const plot,forget plot] coordinates {
  (0.00469831901137,1) (0.0485048590926,1) (0.0597908290802,0.485801528161) (0.061958659091,0.337909923367)
  (0.0654921140522,0.227244379763) (0.0724347210489,0.17427687245) (0.0841420020442,0.133410376051) (0.108574233018,0.102854018447)
  (0.111487861024,0.0968325270579) (0.161269530072,0.0694354278066) (0.183232968091,0.0639406019631) (0.224843438016,0.0516326361704)
  (0.321142384084,0.0381645870857) (0.339251074009,0.0372195240217) (0.52278142306,0.0280332873282) (0.566409614054,0.0266032444816)
  (0.876632214058,0.0203577633107) (0.934592967038,0.0192496185359) (1.48116679711,0.0147359723516) (1.72885821504,0.0137763945008)
  (2.60413982708,0.0110119408322) (3.14131555206,0.00992533049514) (4.23296101904,0.0087960474256) (6.51123580011,0.00710192247858)
  (7.4474088361,0.00696719287459) (12.2019293041,0.00572340300425) (18.2517004401,0.00507004035944) (20.35827065,0.00490459604265)
  (34.7605782311,0.00443889036837) (57.667748051,0.00410284603283) (101.223092941,0.00385003778596) (167.791597894,0.00374805116886)
  (175.323358454,0.00372423096924) (241.467518649,0.00266481753577) (247.019959548,0.00191248873713) (258.218519845,0.00138541784043)
  (275.600515302,0.00100785486672) (275.643324316,0.000999309917385)
};
\addplot[gray!70,densely dashed,line width=0.45pt,forget plot]
coordinates {(1e-9,0.001) (3600,0.001)};
\node[anchor=south west,font=\tiny,text=gray!70]
at ({rel axis cs:0,0}|-{axis cs:1,0.001}) {$10^{-3}$};

\nextgroupplot[
  title={Berlin Prenzlauerberg},
  xmin=0.0001,
  xmax=3600,
  xtick={0.0001,0.01,1,100,3600},
  xticklabels={$10^{-4}$,$10^{-2}$,$1$,$10^{2}$,$3600$},
  ymin=0.0001,
  ymax=10,
  ytick={0.0001,0.01,1},
  yticklabels={$10^{-4}$,$10^{-2}$,$1$}
]
\addplot+[daf,const plot,forget plot] coordinates {
  (0.009610472,1.04922822201) (0.025826628,1.04637371888) (0.039929438,0.237841857914) (0.060453115,0.0788840961531)
  (0.091480923,0.0376543239391) (0.137496572,0.00981342970265) (0.208436099,0.00275517589983) (0.314777316,0.000630389332749)
};
\addplot+[bta,const plot,forget plot] coordinates {
  (0.000701672979631,1) (0.0107039390132,1) (0.0168816960067,0.659290014014) (0.0176036799676,0.51125545287)
  (0.0183459389955,0.370748442692) (0.0186718589975,0.263623069095) (0.0199458900024,0.182085755546) (0.0225988699822,0.137894602484)
  (0.0306886889739,0.11526256399) (0.0365441119648,0.0987040907018) (0.0512688820018,0.0855832359953) (0.081412332016,0.070870023026)
  (0.128801201994,0.0640546787842) (0.218625832989,0.0598741957087) (0.368084339018,0.0566325978605) (0.618366611016,0.0541030625169)
  (1.03687968501,0.0529717926605) (1.749392848,0.0518676558128) (2.802829261,0.0512136555861) (2.99720672501,0.0511431359621)
  (4.98496409302,0.0507469288732) (8.88826709799,0.0503790372846) (14.1140487,0.0501161937326) (23.68290296,0.0500296292146)
  (26.60124149,0.0326702278014) (26.607332778,0.024578714307) (26.621225785,0.0188331244929) (26.660641031,0.0132063262782)
  (26.738034121,0.00974492831804) (26.892804371,0.00692808925541) (27.087931326,0.00504537202055) (27.658840952,0.00362549505538)
  (28.462701445,0.00255025378498) (30.126991253,0.0018818256067) (34.157967493,0.00132014577523) (39.492529718,0.00100637362615)
  (39.79692514,0.000974117922866)
};
\addplot[gray!70,densely dashed,line width=0.45pt,forget plot]
coordinates {(1e-9,0.001) (3600,0.001)};
\node[anchor=south west,font=\tiny,text=gray!70]
at ({rel axis cs:0,0}|-{axis cs:1,0.001}) {$10^{-3}$};

\nextgroupplot[
  title={Berlin Tiergarten},
  xmin=0.0001,
  xmax=3600,
  xtick={0.0001,0.01,1,100,3600},
  xticklabels={$10^{-4}$,$10^{-2}$,$1$,$10^{2}$,$3600$},
  ymin=1e-05,
  ymax=10,
  ytick={0.0001,0.01,1},
  yticklabels={$10^{-4}$,$10^{-2}$,$1$}
]
\addplot+[daf,const plot,forget plot] coordinates {
  (0.007082815,0.143438957634) (0.015597647,0.142605378595) (0.024489446,0.0326208548533) (0.038024438,0.027299901537)
  (0.057120984,0.0204285537207) (0.086141354,0.0042217068469) (0.130088103,6.220449855e-05)
};
\addplot+[bta,const plot,forget plot] coordinates {
  (0.000468757993076,1) (0.00720565399388,1) (0.0159171089763,0.791302586886) (0.0165814280044,0.628274456455)
  (0.0170184340095,0.396666177514) (0.0173023559619,0.286257520053) (0.0188198339893,0.199462800209) (0.0193035149714,0.129321682067)
  (0.0202576699667,0.0926760225389) (0.0221907839878,0.0700457232143) (0.0227640049998,0.0492682405118) (0.0295084069949,0.0366146971963)
  (0.0304931339924,0.0297802131659) (0.0313287089812,0.0254595348779) (0.0496895010001,0.0191878948144) (0.0773329850053,0.0174353017555)
  (0.122369983001,0.0160384476032) (0.201227943006,0.0147722208846) (0.311484035978,0.0138396942685) (0.321257475007,0.0138056633749)
  (0.49668389099,0.012941649754) (0.793048572959,0.0121993000206) (1.26342148497,0.0115463366038) (1.99459909997,0.0111841104528)
  (3.23192475497,0.0110033564185) (5.27092467301,0.0108049291398) (7.70505759201,0.00954033853229) (7.70706787601,0.00699325584728)
  (7.712003714,0.00489994580759) (7.72002280201,0.00352793599546) (7.73110436299,0.00264030371839) (7.76099469699,0.00175935203432)
  (7.80412055901,0.00127589995706) (7.88510847196,0.00105238540312) (7.914832282,0.000990114674164)
};
\addplot[gray!70,densely dashed,line width=0.45pt,forget plot]
coordinates {(1e-9,0.001) (3600,0.001)};
\node[anchor=south west,font=\tiny,text=gray!70]
at ({rel axis cs:0,0}|-{axis cs:1,0.001}) {$10^{-3}$};

\nextgroupplot[
  title={Berlin Friedrichshain},
  xmin=0.0001,
  xmax=3600,
  xtick={0.0001,0.01,1,100,3600},
  xticklabels={$10^{-4}$,$10^{-2}$,$1$,$10^{2}$,$3600$},
  ymin=0.0001,
  ymax=10,
  ytick={0.0001,0.01,1},
  yticklabels={$10^{-4}$,$10^{-2}$,$1$}
]
\addplot+[daf,const plot,forget plot] coordinates {
  (0.004856988,0.112822905826) (0.011410393,0.110157982507) (0.017257792,0.0427479456004) (0.026606003,0.023221708654)
  (0.040204143,0.00309450075116) (0.060475488,0.000471844212896)
};
\addplot+[bta,const plot,forget plot] coordinates {
  (0.000374272000045,1) (0.017801687005,1) (0.0288397810073,0.477635515352) (0.0298668960459,0.277433957898)
  (0.0300180980121,0.222279146063) (0.0303127410007,0.170612702526) (0.030780770001,0.136926424027) (0.0323580140248,0.0941210963346)
  (0.0347543700482,0.0705704215583) (0.0359713929938,0.0670409113853) (0.0417659259983,0.0509138428227) (0.0599274080014,0.0397293941898)
  (0.0634656900074,0.036895724758) (0.0989732680027,0.0305523967841) (0.141497847042,0.0264891042351) (0.164740820008,0.0248378222003)
  (0.274826498993,0.02198877861) (0.45359175303,0.0199891469152) (0.577229192015,0.0191886491613) (0.755558626028,0.0186167201487)
  (1.25046029105,0.0178733073268) (2.12015887804,0.0173073963682) (3.430520088,0.0169650319393) (5.72889834904,0.0167511671599)
  (9.48822951905,0.0166250903765) (11.327699553,0.013635115599) (11.352425173,0.0093080366566) (11.419072425,0.00700234684699)
  (11.514428076,0.00483363728313) (11.719717382,0.00349206932429) (12.090480688,0.0025867647937) (12.770935959,0.00189143181264)
  (13.879832949,0.00136766050939) (15.542015564,0.00102804484901) (15.670814115,0.000990512258337)
};
\addplot[gray!70,densely dashed,line width=0.45pt,forget plot]
coordinates {(1e-9,0.001) (3600,0.001)};
\node[anchor=south west,font=\tiny,text=gray!70]
at ({rel axis cs:0,0}|-{axis cs:1,0.001}) {$10^{-3}$};

\end{groupplot}
\node[font=\small] at ($(group c1r3.south)!0.5!(group c2r3.south)+(0,-0.62cm)$)
  {elapsed time $t$ (seconds)};
\node[font=\small,rotate=90] at ($(group c1r1.west)!0.5!(group c1r3.west)+(-0.95cm,0)$)
  {relative gap $\gamma_{i,m}(t)$};
\end{tikzpicture}
\caption{For nonsmooth instances~$i$, the curves show the relative gaps $\gamma_{i,m}(t)$ of Dafermos and BTA as a function of computation time~$t$.
The dashed line shows the tolerance $\varepsilon_{i} = 10^{-3}$.
Each curve stops at its first recorded point with $\gamma_{i,m}(t) \le \varepsilon_{i}$, or at $T_{i}$.}
\label{fig:exp_nonsmooth_common_objective_iii}
\end{figure}
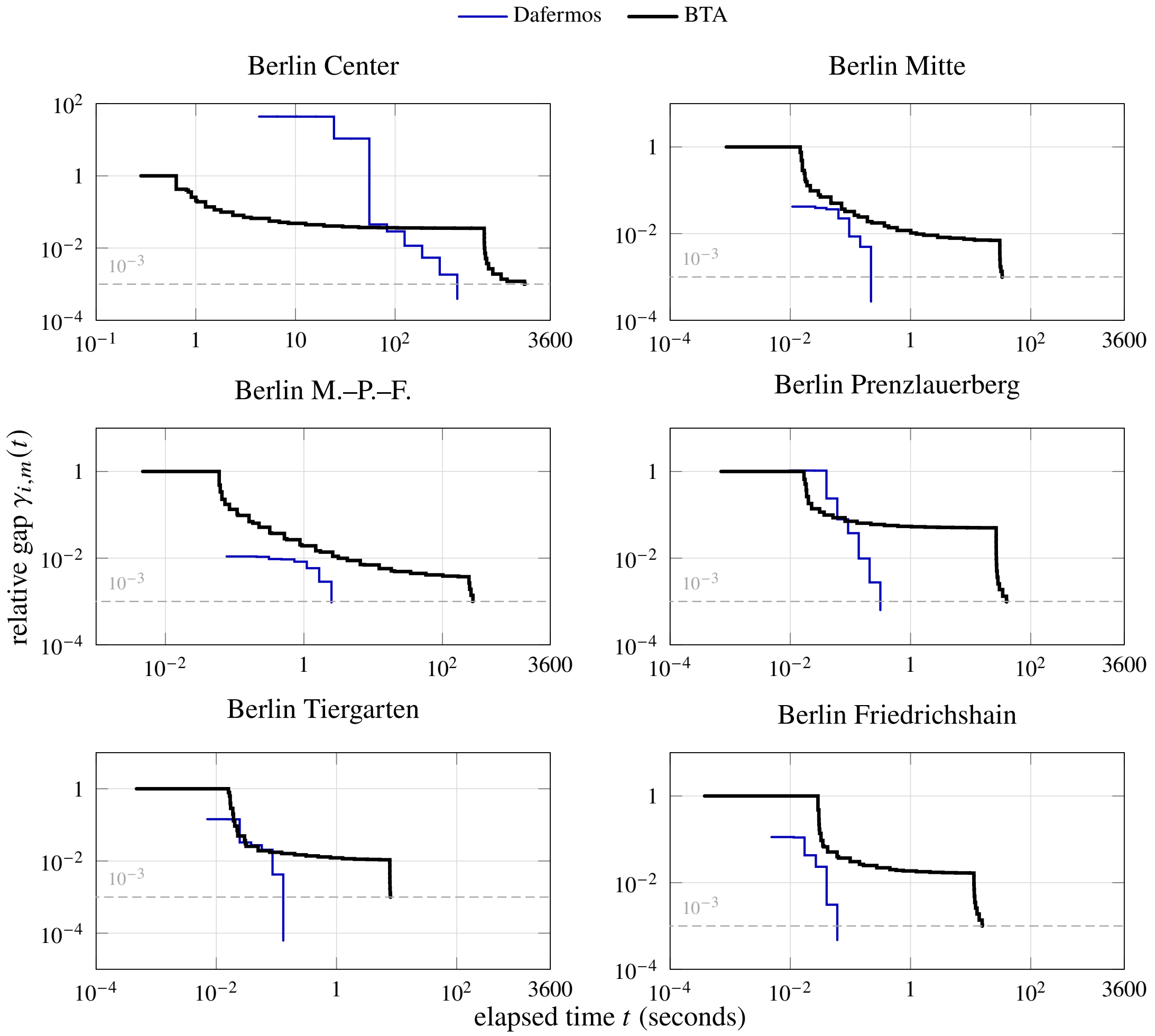


\FloatBarrier

\subsection{Path-support Sparsity}

As mentioned before, in typical networks only a small fraction of the paths from an origin to a destination carry positive flow.
This observation motivated the design of an algorithm that optimizes the objective for (small) subsets $\hat{\mathcal{R}}^{\pi} \subset \mathcal{R}^{\pi}$ of paths.
A related statistic that is of interest is the number of paths generated by the algorithm ($\sum_{\pi \in \Pi} |\hat{\mathcal{R}}^{\pi}|$ when the algorithm stops) that carry positive flow in the final solution.
We define a {\em used path} as a path~$r$ that has flow $q^{\pi*}_{r} > 10^{-6}$ in the final solution.
Table~\ref{tab:bta_smooth_sparsity} reports the number of generated paths, the maximum number of generated paths for any OD pair, the number of used paths, and the average number of used paths per OD pair, when the BTA algorithm stopped for Sioux Falls, Anaheim, Barcelona, and Winnipeg.
For these instances, the fraction of paths generated by the algorithm that carry positive flow in the final solution is quite large, but the fraction tends to decrease as the network size increases.

\begin{table}[h]
\centering
\normalsize
\caption{The number of paths generated by the algorithm ($\sum_{\pi \in \Pi} |\hat{\mathcal{R}}^{\pi}|$ when the algorithm stops), the maximum number of generated paths for an OD pair ($\max_{\pi \in \Pi} |\hat{\mathcal{R}}^{\pi}|$), the number of used paths, and the average number of used paths per OD pair, for BTA on four smooth instances.}
\label{tab:bta_smooth_sparsity}
\begin{tabular}{lrrrr}
\toprule
Network & Generated paths & Maximum generated paths & Used paths & Average used paths \\
 &  & per OD pair &  & per OD pair \\
\midrule
Sioux Falls & $922$ & $5$ & $779$ & $1.48$ \\
Anaheim & $4513$ & $22$ & $3542$ & $2.52$ \\
Barcelona & $60435$ & $59$ & $31268$ & $3.95$ \\
Winnipeg & $57510$ & $85$ & $24642$ & $5.67$ \\
\bottomrule
\end{tabular}
\end{table}

\FloatBarrier

\subsection{Impact of Tolls on Path Diversity}
\label{sec:exp_toll_impact}

It is of interest to determine how much the introduction of a second criterion such as tolls affect the diversity of paths used in equilibrium.
We compare the fraction of OD pairs with positive flow in equilibrium on more than one path in four settings: no tolls, lower marginal-cost tolls with $\alpha_{0} = 1$, higher marginal-cost tolls with $\alpha_{0} = 18$, and selective capacity-UE tolls (specified in Appendix~\ref{app:numerical-experiments}).
The no-toll case is the single-criterion baseline.
For Capacity-UE, we computed a reference flow $\widehat{x}^{\mathrm{UE}}$ for the single criterion (time) equilibrium problem on each network with the same OD demands as the bicriteria instances.
Recall that ``arc capacities'' $c_{a}$ are parameters of the BPR travel-time functions.
Then we assign monetary tolls
\[
s_{a} \ \ = \ \
\begin{cases}
\operatorname{round}_{0.01}\!\left(
t_{a}(\widehat{x}^{\mathrm{UE}}_{a})\right)
& \mbox{ if } \widehat{x}^{\mathrm{UE}}_{a} > c_{a} \\
0 & \mbox{ if } \widehat{x}^{\mathrm{UE}}_{a} \le c_{a}
\end{cases}
\qquad a \in \mathcal{A},
\]
where $\operatorname{round}_{0.01}$ denotes rounding to the nearest hundredth.
Table~\ref{tab:exp_toll_schemes} gives the results under the four tolling schemes, including computation times, final gaps, used-path counts, and fraction of OD pairs with multiple paths with positive flow in equilibrium.
The fraction of OD pairs with positive flow on more than one path is greater with a second criterion (tolls) than with a single criterion (time), and in most, but not all, cases, the fraction tends to increase as the tolls increase.

\ignore{
\begin{table}[t]
\centering
\small
\caption{The fraction of OD pairs with positive flow in equilibrium on more than one path under four tolling schemes.
Table~\ref{tab:exp_toll_schemes} in Appendix~\ref{app:numerical-experiments} also gives elapsed times, final gaps, and used-path counts for the same instances.}
\label{tab:exp_toll_multi}
\begin{tabular}{llrrrr}
\toprule
Network & Instance & No toll & MC $\alpha_{0} = 1$ & MC $\alpha_{0} = 18$ & Capacity-UE \\
\midrule
Sioux Falls & Smooth & $0.265$ & $0.394$ & $0.286$ & $0.345$ \\
 & Nonsmooth & $0.222$ & $0.570$ & $0.697$ & $0.540$ \\
\addlinespace
Anaheim & Smooth & $0.286$ & $0.583$ & $0.762$ & $0.323$ \\
 & Nonsmooth & $0.337$ & $0.764$ & $0.925$ & $0.698$ \\
\addlinespace
Berlin Tiergarten & Smooth & $0.222$ & $0.484$ & $0.630$ & $0.247$ \\
 & Nonsmooth & $0.238$ & $0.618$ & $0.835$ & $0.315$ \\
\addlinespace
Winnipeg-Asym & Smooth & $0.788$ & $0.827$ & $0.722$ & $0.789$ \\
 & Nonsmooth & $0.786$ & $0.882$ & $0.919$ & $0.829$ \\
\bottomrule
\end{tabular}
\end{table}
}

\begin{table}[t]
\centering
\small
\caption{Computation times, final gaps, used-path counts, and fraction of OD pairs with multiple used paths, for the four tolling schemes.
``Used'' paths are paths with equilibrium flow above $10^{-6}$.
}
\label{tab:exp_toll_schemes}
\begin{tabular}{lrrrr@{\hskip 10pt}rrrr}
\toprule
& \multicolumn{4}{c}{Smooth instances} & \multicolumn{4}{c}{Nonsmooth instances} \\
\cmidrule(lr){2-5}\cmidrule(lr){6-9}
Toll scheme & Time & Gap & Used & Multi & Time & Gap & Used & Multi \\
\midrule
\multicolumn{9}{l}{\textit{Sioux Falls}} \\
None & $0.024$ & $5.52\times10^{-6}$ & $734$ & $0.265$ & $0.475$ & $7.11\times10^{-4}$ & $671$ & $0.222$ \\
MC $\alpha_{0}=1$ & $0.021$ & $9.13\times10^{-6}$ & $817$ & $0.394$ & $0.204$ & $9.93\times10^{-4}$ & $975$ & $0.570$ \\
MC $\alpha_{0}=18$ & $0.011$ & $8.19\times10^{-6}$ & $692$ & $0.286$ & $0.286$ & $9.62\times10^{-4}$ & $1176$ & $0.697$ \\
Capacity-UE & $0.013$ & $9.55\times10^{-6}$ & $762$ & $0.345$ & $2.170$ & $8.57\times10^{-4}$ & $904$ & $0.540$ \\
\midrule
\multicolumn{9}{l}{\textit{Anaheim}} \\
None & $0.074$ & $9.45\times10^{-6}$ & $1949$ & $0.286$ & $0.076$ & $5.09\times10^{-4}$ & $2034$ & $0.337$ \\
MC $\alpha_{0}=1$ & $0.270$ & $9.86\times10^{-6}$ & $2757$ & $0.583$ & $0.532$ & $9.73\times10^{-4}$ & $3631$ & $0.764$ \\
MC $\alpha_{0}=18$ & $1.046$ & $9.89\times10^{-6}$ & $3965$ & $0.762$ & $17.722$ & $1.43\times10^{-3}$ & $7250$ & $0.925$ \\
Capacity-UE & $0.288$ & $9.99\times10^{-6}$ & $2172$ & $0.323$ & $11.198$ & $2.83\times10^{-3}$ & $3552$ & $0.698$ \\
\midrule
\multicolumn{9}{l}{\textit{Berlin Tiergarten}} \\
None & $0.028$ & $8.51\times10^{-6}$ & $863$ & $0.222$ & $0.180$ & $6.49\times10^{-4}$ & $809$ & $0.238$ \\
MC $\alpha_{0}=1$ & $0.044$ & $9.45\times10^{-6}$ & $1144$ & $0.484$ & $0.259$ & $9.93\times10^{-4}$ & $1565$ & $0.618$ \\
MC $\alpha_{0}=18$ & $0.153$ & $9.61\times10^{-6}$ & $1383$ & $0.630$ & $4.118$ & $2.04\times10^{-2}$ & $2534$ & $0.835$ \\
Capacity-UE & $0.023$ & $9.20\times10^{-6}$ & $810$ & $0.247$ & $1.740$ & $9.40\times10^{-4}$ & $975$ & $0.315$ \\
\midrule
\multicolumn{9}{l}{\textit{Winnipeg-Asym}} \\
None & $10.783$ & $9.45\times10^{-6}$ & $19017$ & $0.788$ & $69.067$ & $8.95\times10^{-4}$ & $14306$ & $0.786$ \\
MC $\alpha_{0}=1$ & $5.919$ & $9.84\times10^{-6}$ & $20061$ & $0.827$ & $16.850$ & $9.72\times10^{-4}$ & $24380$ & $0.882$ \\
MC $\alpha_{0}=18$ & $1.187$ & $9.80\times10^{-6}$ & $14898$ & $0.722$ & $2.034$ & $9.99\times10^{-4}$ & $19598$ & $0.919$ \\
Capacity-UE & $6.098$ & $9.86\times10^{-6}$ & $19423$ & $0.789$ & $5.959$ & $9.90\times10^{-4}$ & $14293$ & $0.829$ \\
\bottomrule
\end{tabular}
\end{table}

\end{document}